\documentclass[10pt,twoside,openright,a4paper]{report} 
\usepackage[utf8]{inputenc} 
\usepackage[T1]{fontenc}
\usepackage{lmodern}

\usepackage{amsmath}
\usepackage{amsfonts}
\usepackage{stixgras}
\usepackage{mathrsfs}
\usepackage{mathtools}
\usepackage{stmaryrd}
\usepackage[all]{xy}
\usepackage{makeidx}
\usepackage[french]{babel}
\usepackage{csquotes}

\usepackage{amssymb}
\usepackage{xspace}
\usepackage{amsthm}
\makeatletter\def\th@plain{\slshape}\makeatother
\makeatletter\patchcmd{\th@remark}{\itshape}{\slshape}{}{}\makeatother

\usepackage{xr}

\usepackage{titlesec}
\titleformat{\chapter}[hang]{\Huge\bfseries}{\thechapter.}{.5em}{\vspace*{.05em}}[\vskip 0em]
\renewcommand\thechapter{\Alph{chapter}}
\titlelabel{\thetitle.\quad}

\usepackage[url=true,isbn=false,doi=true,backend=biber,defernumbers=true,backref=true,safeinputenc,maxnames=4,sorting=anyt]{biblatex}

\newcounter{bidon}
\newcommand{\rdb}{\refstepcounter{bidon}}

\usepackage[french,nohints]{minitoc}

\theoremstyle{plain}
\newtheorem{theorem}{Théorème}[section]

\newtheorem{pstc}[theorem]{Positivstellensatz concret}
\newtheorem{pstf}[theorem]{Positivstellensatz formel}

\newtheorem{lemma}[theorem]{Lemme}%[section]
\newtheorem{corollary}[theorem]{Corolaire}%[section]
\newtheorem{proposition}[theorem]{Proposition}%[section]
\newtheorem{prpta}[theorem]{Propriétés attendues}%[section]
\newtheorem{propdef}[theorem]{Proposition et définition}
\newtheorem{fact}[theorem]{Fait}%[section]

\newtheorem{plcc}[theorem]{Principe local-global concret}

\newtheorem{theoremc}[theorem]{Th\'{e}or\`{e}me\etoz}

\newtheorem{corollaryc}[theorem]{Corolaire\etoz}

\newtheorem{propositionc}[theorem]{Proposition\etoz}

\theoremstyle{definition}
\newtheorem{rstra}[theorem]{Règles structurelles admissibles}%[section]
\newtheorem{definition}[theorem]{Définition}%[section]
\newtheorem{dfni}[theorem]{Définition informelle}%[section]
\newtheorem{example}[theorem]{Exemple}%[section]
\newtheorem{examples}[theorem]{Exemples}%[section]
\newtheorem{question}[theorem]{Question}
\newtheorem{questions}[theorem]{Questions}
\newtheorem{definota}[theorem]{Définition et notation} 
\newtheorem{definotas}[theorem]{Définitions et notations} 

\newtheorem{definitionc}[theorem]{Définition\etoz}%[section]

\theoremstyle{remark}
\newtheorem{notE}[theorem]{Note}

\newtheorem{remark}[theorem]{Remarque}%[section]
\newtheorem{remarks}[theorem]{Remarques}%[section]

\newcommand \hum[1] {\sibrouillon{\noindent {\sf hum: #1}}}

\newcommand\Subsection[1]{%\goodbreak
\rdb\addcontentsline{toc}{subsection}{#1} \subsection*{#1}}

\newcommand\Subsectio[2]{%\goodbreak
\rdb\addcontentsline{toc}{subsection}{#2} 
\subsection*{#1}}

\newcommand\Subsubsection[1]{%\goodbreak
\rdb\addcontentsline{toc}{subsubsection}{#1} \subsubsection*{#1}}

\newcommand{\Tp}{{\sa T\,'}}

\renewcommand\paragraph[1]{

\rdb\addcontentsline{toc}{subsubsection}{#1} \medskip \noindent $\bullet$ \textbf{#1}}

\newcommand{\vou}{\MA{\tsbf{ ou }}}
\newcommand{\Vou}{\MA{\tsbf{OU}}}
\newcommand \EXists[1] {\tsbf{Introduire }{#1\,}\tsbf{ tel que }\,}
\newcommand \vet {\tsbf{,}\;}
\newcommand \Atcl {\mathrm{Atcl}}
\newcommand \Tcl {\mathrm{Tcl}}
\newcommand \Atclv {\mathrm{Atclv}}

\newcounter{MF}
\newcommand\stMF{\stepcounter{MF}}

\newcommand{\lec}{\stMF\ifodd\value{MF}lecteur\xspace \else 
lectrice\xspace \fi}

\newcommand{\lecs}{\stMF\ifodd\value{MF}lecteurs\xspace \else 
lectrices\xspace \fi}

\newcommand{\alec}{\stMF\ifodd\value{MF}au lecteur\xspace \else%
à la lectrice\xspace \fi}

\newcommand{\dlec}{\stMF\ifodd\value{MF}du lecteur\xspace \else%
de la lectrice\xspace \fi}

\newcommand{\llec}{\stMF\ifodd\value{MF}le lecteur\xspace \else la lectrice\xspace \fi}

\newcommand{\Llec}{\stMF\ifodd\value{MF}Le lecteur\xspace \else La lectrice\xspace \fi}

\newcommand{\lui}{\ifodd\value{MF}lui\xspace \else
elle\xspace \fi}

\newcommand{\du}{\ifodd\value{MF}du\xspace \else
de la\xspace \fi}

\newcommand{\celui}{\ifodd\value{MF}celui\xspace \else
celle\xspace \fi}

\newcommand{\ceux}{\ifodd\value{MF}ceux\xspace \else
celles\xspace \fi}

\newcommand{\er}{\ifodd\value{MF}er\xspace \else
ère\xspace \fi}

\newcommand{\eux}{\ifodd\value{MF}eux\xspace \else
elles\xspace \fi}

\newcommand{\eUx}{\ifodd\value{MF}eux\xspace \else
euse\xspace \fi}

\newcommand{\leux}{\ifodd\value{MF}leux \else
leuse \fi}

\newcommand{\il}{\ifodd\value{MF}il\xspace \else
elle\xspace \fi}

\newcommand{\ien}{\ifodd\value{MF}ien\xspace \else
ienne\xspace \fi}

\newcommand{\iens}{\ifodd\value{MF}iens\xspace\else
iennes\xspace\fi}

\newcommand{\e}{\ifodd\value{MF}\xspace \else e\xspace \fi}

\newcommand{\n}{\ifodd\value{MF}n \else nne \fi}
\newcommand{\nz}{\ifodd\value{MF}n\else nne\fi}

\makeatletter
\newcommand{\la}{\@ifstar{\ifodd\value{MF}le\xspace\else
la\xspace\fi}{\stMF\ifodd\value{MF}le\xspace \else la\xspace \fi}}
\makeatother

\newcommand \Note{\rdb
\noi{\sl Note. }}

\newcommand \rem{\rdb
\noi{\sl Remarque. }}

\newcommand \rems{\rdb
\noi{\sl Remarques. }}

\newcommand \thref[1] {théorème~\ref{#1}}
\newcommand \paref[1] {page~\pageref{#1}}
\newcommand \pstfref[1] {Posi\-tiv\-stel\-lensatz formel~\ref{#1}}
\newcommand \pstref[1] {Posi\-tiv\-stel\-lensatz~\ref{#1}}

\newcommand\oge{\leavevmode\raise.3ex\hbox{$\scriptscriptstyle\langle\!\langle\,$}}
\newcommand\feg{\leavevmode\raise.3ex\hbox{$\scriptscriptstyle\,\rangle\!\rangle$}}

\newcommand\gui[1]{\oge{#1}\feg}

\newcommand\comm{\rdb
\noi{\sl Commentaire. }}

\newcommand \Cad {C'est-à-dire\xspace}
\newcommand \recu {récur\-rence\xspace}
\newcommand \hdr {hypo\-thèse de \recu}
\newcommand \ie {i.e.\xspace}
\newcommand \cad {c'est-à-dire\xspace}
\newcommand \cade {c'est-à-dire en\-co\-re\xspace}
\newcommand \ssi {si, et seu\-lement si, }

\newcommand \spdg {sans per\-te de géné\-ra\-lité\xspace}

\newcommand \Propeq {Les pro\-pri\-é\-tés sui\-van\-tes sont 
équi\-va\-len\-tes.}

\newcommand \Kev {$\gK$-\evc}

\newcommand \QQlg {$\QQ$-\alg}
\newcommand \QQlgs {$\QQ$-\algs}

\newcommand \Rlgs {$\gR$-\algs}

\newcommand \RRlg {$\RR$-\alg}
\newcommand \RRlgs {$\RR$-\algs}

\newcommand \Amos {$\gA$-mo\-du\-les\xspace}

\newcommand \ac{algé\-bri\-quement clos\xspace}

\newcommand \acos{anneaux commutatifs\xspace}

\newcommand \adv {anneau de valuation\xspace}
\newcommand \advs {anneaux de valuation\xspace}

\newcommand \Afrs {Anneaux \frls}
\newcommand \afr {anneau \frl}
\newcommand \aFr {\hyperref[theorieAfr]{anneau \frl}\xspace}
\newcommand \afrs {anneaux \frls}

\newcommand \afrdc {\afr $2$-clos\xspace}

\newcommand \afrrs {\afrs réduits\xspace}

\newcommand \afrvr {\afr avec \ravs}

\newcommand \afrvrs {\afrs avec \ravs}

\newcommand \aftr {anneau \ftm réel\xspace}
\newcommand \aftrs {anneaux \ftm réels\xspace}

\newcommand \arftrs {anneaux réticulés \ftm réels\xspace}

\newcommand \afrb {anneau de fonctions réelles bornées\xspace}
\newcommand \afrbs {anneaux de fonctions réelles bornées\xspace}

\newcommand \agB {\alg de Boole\xspace}
\newcommand \agBs {\algs de Boole\xspace}

\newcommand \agq{algé\-bri\-que\xspace}
\newcommand \agqs{algé\-bri\-ques\xspace}

\newcommand \alg {algè\-bre\xspace}
\newcommand \algs {algè\-bres\xspace}

\newcommand \algo{algo\-rithme\xspace}
\newcommand \algos{algo\-rithmes\xspace}

\newcommand \algq{al\-go\-rith\-mi\-que\xspace}

\newcommand \alo {an\-neau lo\-cal\xspace}
\newcommand \alos {an\-neaux lo\-caux\xspace}

\newcommand \alrd {\alo \dcd}
\newcommand \alrds {\alos \dcds}

\newcommand \arc {anneau réel clos\xspace}
\newcommand \aRc {\hyperref[theorieArc]{\arc}\xspace}
\newcommand \arcs {anneaux réels clos\xspace}

\newcommand \ari{arith\-mé\-tique\xspace}

\newcommand \asr {anneau \str}
\newcommand \asrs {anneaux \strs}

\newcommand \asrvr {\asr avec \ravs}

\newcommand \bif {borne infé\-rieure\xspace} %

\newcommand \bsp {borne supé\-rieure\xspace} %
\newcommand \carns{carac\-té\-ri\-sations\xspace}

\newcommand \cdacs{\cdis \ac}  

\newcommand \cdi{corps discret\xspace}
\newcommand \cdis{corps discrets\xspace}
  
\newcommand \cdv{changement de variables\xspace}

\newcommand \cli {clô\-ture inté\-grale\xspace}

\newcommand \codi {corps ordonné discret\xspace}
\newcommand \codis {corps ordonnés discrets\xspace}

\newcommand \coe {coef\-fi\-cient\xspace}
\newcommand \coes {coef\-fi\-cients\xspace}

\newcommand \com {co\-ma\-xi\-maux\xspace}

\newcommand \coo {coor\-donnée\xspace}
\newcommand \coos {coor\-données\xspace}

\newcommand \cop {complé\-men\-taire\xspace}
\newcommand \cops {complé\-men\-taires\xspace}

\newcommand \corl {coro\-laire\xspace}
\newcommand \corls {coro\-laires\xspace}

\newcommand \cosv {conser\-vative\xspace}

\newcommand \covr {corps ordonné avec \ravs}
\newcommand \covrs {corps ordonnés avec \ravs}

\newcommand \crc {corps réel clos\xspace}
\newcommand \crcs {corps réels clos\xspace}

\newcommand \crcd {corps réel clos discret\xspace}
\newcommand \crcds {corps réels clos discrets\xspace}

\newcommand \dcd {rési\-duel\-lement dis\-cret\xspace}
\newcommand \dcds {rési\-duel\-lement dis\-crets\xspace}

\newcommand \ddk {dimension de~Krull\xspace}

\newcommand \demo{démon\-stra\-tion\xspace}     
\newcommand \demos{démon\-stra\-tions\xspace}

\newcommand \Dfn{Défi\-nition\xspace}  
  
\newcommand \dfn{défi\-nition\xspace}  
\newcommand \dfns{défi\-nitions\xspace}

\newcommand \dij{disjonc\-tive\xspace}  
\newcommand \dijs{disjonc\-tives\xspace}

\newcommand \dvz {di\-viseur de zéro\xspace}

\newcommand \dve {divi\-si\-bi\-lité\xspace}

\newcommand \eco {\elts \com}

\newcommand \egmt {éga\-lement\xspace}

\newcommand \egt {éga\-li\-té\xspace}
\newcommand \egts {éga\-li\-tés\xspace}

\newcommand \elr{élé\-men\-taire\xspace}  
\newcommand \elrs{élé\-men\-taires\xspace}

\newcommand \elt{élé\-ment\xspace}  
\newcommand \elts{élé\-ments\xspace}

\newcommand \entrel {rela\-tion impli\-ca\-tive\xspace}
\newcommand \entrels {rela\-tions impli\-ca\-tives\xspace}

\newcommand\evc{es\-pa\-ce vec\-to\-riel\xspace}

\newcommand \eqv {équivalent\xspace}  
\newcommand \eqve {équivalente\xspace}  
\newcommand \eqvs {équivalents\xspace}  
\newcommand \eqves {équivalentes\xspace}  

\newcommand \eqvc {équivalence\xspace}  
\newcommand \eqvcs {équivalences\xspace}  

\newcommand \esid {essentiellement identique\xspace}  
\newcommand \esids {essentiellement identiques\xspace}

\newcommand \eseq {essentiellement \eqve}  
\newcommand \eseqs {essentiellement \eqves}

\newcommand \fsa {fermé \sagq}
\newcommand \fsas {fermés \sagqs}

\newcommand \fsagc {fonction \sagc}
\newcommand \fsagcs {fonctions \sagcs}

\newcommand \fsagces {\fsagcs entières\xspace}

\newcommand \frl {for\-tement réticulé\xspace}
\newcommand \frle {for\-tement réticulée\xspace}
\newcommand \frls {for\-tement réticulés\xspace}

\newcommand \ftm {fortement\xspace}

\newcommand\gmq{géomé\-trique\xspace}  
\newcommand\gmqs{géomé\-triques\xspace}

\newcommand\gnl{géné\-ral\xspace}  
\newcommand\gnle{géné\-rale\xspace}  
  
\newcommand\gnles{géné\-rales\xspace}  

\newcommand\gnlt{géné\-ra\-lement\xspace}  

\newcommand\gnn{géné\-ra\-li\-sa\-tion\xspace}

\newcommand\gnqs {géné\-riques\xspace}

\newcommand \grl{groupe \rtl}
\newcommand \grls{groupes \rtls}

\newcommand\gtrs{géné\-ra\-teurs\xspace}  

\newcommand \homos {ho\-mo\-mor\-phismes\xspace}

\newcommand \icl {inté\-gra\-lement clos\xspace}

\newcommand \icr {intervalle compact réel\xspace}

\newcommand \icrc {intervalle compact réel clos\xspace}
\newcommand \icrcs {intervalles compacts réels clos\xspace}

\newcommand \id {idéal\xspace}
\newcommand \ids {idéaux\xspace}

\newcommand \idas {\idts \agqs}

\newcommand \idemas {idéaux maxi\-maux\xspace}

\newcommand \idep {idéal pre\-mier\xspace}
\newcommand \ideps {idéaux pre\-miers\xspace}

\newcommand \idp {idéal prin\-ci\-pal\xspace}
\newcommand \idps {idé\-aux prin\-ci\-paux\xspace}

\newcommand \idt {iden\-ti\-té\xspace}
\newcommand \idts {iden\-ti\-tés\xspace}

\newcommand \idtr {in\-dé\-ter\-mi\-née\xspace}
\newcommand \idtrs {in\-dé\-ter\-mi\-nées\xspace}

\newcommand \imd {immé\-diat\xspace}

\newcommand \imdt {immé\-dia\-te\-ment\xspace}

\newcommand \inteq {intui\-ti\-vement \eqve}
\newcommand \inteqs {intui\-ti\-vement \eqves}

\newcommand \irds {irré\-duc\-tibles\xspace}

\newcommand \iso {iso\-mor\-phisme\xspace}

\newcommand \itf {idéal \tf}
\newcommand \itfs {idé\-aux \tf}

\newcommand \iv {inversible\xspace}

\newcommand \lgb {local-global\xspace}

\newcommand \lgbs {local-globals\xspace}

\newcommand \lot {loca\-lement\xspace}

\newcommand \lon {loca\-li\-sation\xspace}

\newcommand \lsdz {\lot \sdz}

\newcommand \mo {mo\-no\"{\i}de\xspace}
\newcommand \mos {mo\-no\"{\i}des\xspace}

\newcommand \moco {\mos \com}

\newcommand \ncr{néces\-saire\xspace}       
\newcommand \ncrs{néces\-saires\xspace}       

\newcommand \ncrt{néces\-sai\-rement\xspace}

\newcommand \nds {\textsl{non} discret\xspace}
\newcommand \ndss {\textsl{non} discrets\xspace}

\newcommand \ndsof {corps ordonné \nds}
\newcommand \ndsofs {corps ordonnés \ndss}

\newcommand \ndz {régu\-lier\xspace}

\newcommand \noe {noethé\-rien\xspace}

\newcommand \nst {Null\-stellen\-satz\xspace}

\newcommand \oqc {ouvert \qc}

\newcommand \pb{pro\-blè\-me\xspace}  
\newcommand \pbs{pro\-blè\-mes\xspace}  

\newcommand \peq {purement équa\-tion\-nelle\xspace}
\newcommand \peqs {purement équa\-tion\-nelles\xspace}

\newcommand \pf {de \pn finie\xspace}

\newcommand \plg {prin\-cipe \lgb}
\newcommand \plgs {prin\-cipes \lgbs}

\newcommand \plgc {\plg con\-cret\xspace}
\newcommand \plgcs {\plgs con\-crets\xspace}

\newcommand \pn {présen\-ta\-tion\xspace}
\newcommand \pns {présen\-ta\-tions\xspace}

\newcommand \pol {poly\-nôme\xspace}
\newcommand \pols {poly\-nômes\xspace}

\newcommand \polle{poly\-no\-miale\xspace}  
\newcommand \polles{poly\-no\-miales\xspace}

\newcommand \prmt {préci\-sé\-ment\xspace}
\newcommand \Prmt {Préci\-sé\-ment\xspace}

\newcommand \prt {pro\-pri\-été\xspace}
\newcommand \prts {pro\-pri\-étés\xspace}

\newcommand \Pst {Positiv\-stellen\-satz\xspace}
\newcommand \Psts {Positiv\-stellen\-sätze\xspace}

\newcommand \qc {quasi-compact\xspace}

\newcommand \qtfs {quanti\-fi\-cateurs\xspace}

\newcommand \ralg {règle \agq}
\newcommand \ralgs {règles \agqs}

\newcommand \rav {racine virtuelle\xspace}
\newcommand \ravs {racines virtuelles\xspace}

\newcommand \rdij {règle \dij}
\newcommand \rdijs {règles \dijs}

\newcommand \rdy {règle dyna\-mique\xspace}
\newcommand \rdys {règles dyna\-miques\xspace}

\newcommand \red {règle directe\xspace}
\newcommand \reds {règles directes\xspace}

\newcommand \rex {règle exis\-ten\-tielle simple\xspace}

\newcommand \rexris {règles exis\-ten\-tielles rigides\xspace}

\newcommand \rsim {règle de simplification\xspace}
\newcommand \rsims {règles de simplification\xspace}

\newcommand \rtl {réti\-culé\xspace}
\newcommand \rtls {réti\-culés\xspace}

\newcommand \sad {\salg dynamique\xspace}
\newcommand \sads {\salgs dynamiques\xspace}

\newcommand \sagq {semialgébrique\xspace}
\newcommand \sagqs {semialgébriques\xspace}

\newcommand \sagc {\sagq continue\xspace}
\newcommand \sagcs {\sagqs continues\xspace}

\newcommand \salg {structure \agq}
\newcommand \salgs {structures \agqs}

\newcommand \sdz {sans \dvz}

\newcommand \spo {semipolynôme\xspace}
\newcommand \spos {semipolynômes\xspace}

\newcommand \stm {strictement\xspace}

\newcommand \str {\stm réticulé\xspace}

\newcommand \strs {\stm réticulés\xspace}

\newcommand \sul {supplé\-men\-taire\xspace}
\newcommand \suls {supplé\-men\-taires\xspace}

\newcommand \sys {sys\-tème\xspace}
\newcommand \syss {sys\-tèmes\xspace}

\newcommand \talg {théorie \agq}
\newcommand \talgs {théories \agqs}

\newcommand \tco {théorie cohé\-rente\xspace}

\newcommand \tdij {théorie \dij}

\newcommand \tdy {théorie dyna\-mique\xspace}
\newcommand \tdys {théories dyna\-miques\xspace}

\newcommand \tel {théorie exis\-ten\-tielle\xspace}

\newcommand \telri {théorie exis\-ten\-tiel\-lement rigide\xspace}

\newcommand \tf {de type fini\xspace}

\newcommand \tfo {théorie formelle\xspace}

\newcommand \tgm {théorie \gmq}
\newcommand \tgms {théories \gmqs}

\newcommand \Tho {Théo\-rème\xspace}

\newcommand \tho {théo\-rème\xspace}
\newcommand \thos {théo\-rèmes\xspace}

\newcommand \tpe {théorie \peq}
\newcommand \tpes {théories \peqs}

\newcommand \trdi {treil\-lis dis\-tri\-bu\-tif\xspace}
\newcommand \trdis {treil\-lis dis\-tri\-bu\-tifs\xspace}

\newcommand \unt {uni\-taire\xspace}

\newcommand \zed {z\'{e}ro-di\-men\-sionnel\xspace}

\newcommand \zedr {\zed réduit\xspace}

\newcommand \cof {cons\-truc\-tif\xspace}
\newcommand \cofs {cons\-truc\-tifs\xspace}

\newcommand \cov {cons\-truc\-tive\xspace}
\newcommand \covs {cons\-truc\-tives\xspace}

\newcommand \coma {\maths\covs}
\newcommand \clama {\maths clas\-siques\xspace}

\renewcommand \cot {cons\-truc\-ti\-vement\xspace}

\newcommand \maths {mathé\-ma\-tiques\xspace}
\newcommand \mathe {mathé\-ma\-tique\xspace}

\newcommand \prco {démons\-tration \cov}
\newcommand \prcos {démons\-trations \covs}

\input{MathMacrosRgeom.tex}

\makeindex

\usepackage[bookmarksopen=false,pdftex=true,breaklinks=true,%
 plainpages=false,%
 hyperindex=true,pdfstartview=FitH,%
 pdfpagelabels=true,
 linkcolor=blue,%
 citecolor=blue,urlcolor=red,
 colorlinks=true%
 ]%
 {hyperref}

\newcommand \sibrouillon[1]{}
\newcommand \Today {\sibrouillon{\vspace{-18em}\hspace{10cm}\fbox{\today}\vspace{15.6em}}}

\RequirePackage{etoolbox}

\patchcmd{\sectionmark}{\MakeUppercase}{}{}{}
\patchcmd{\chaptermark}{\MakeUppercase}{}{}{}

\begin{document} 

\title{Théories géométriques pour l'\alg des nombres réels sans test de signe ni axiome de choix dépendant}
\author{Henri Lombardi et Assia Mahboubi}

\thispagestyle{empty}
% center
~
\vspace{1em}
\begin{center} 
{\bf \LARGE
Théories géométriques 
\\[.3em]
 pour l'\alg des nombres réels
\\[.5em]
sans test de signe ni axiome de choix dépendant}
\\[.8em]{\large
Henri Lombardi et Assia Mahboubi
\\[.8em]\normalsize
Version préliminaire
\\[.8em]
la dernière version en cours peut être trouvée en\\[.3em] \url{http://hlombardi.free.fr/Reels-Geometriques.pdf}
}
\vspace{2em}

\today
\end{center}

\vspace{2em}

\centerline{\bf Résumé}

\medskip 
%: quotation
\begin{quotation} \label{quotation}

Dans ce mémoire, nous cherchons à construire une théorie dynamique aussi complète que possible pour décrire les propriétés algébriques du corps des réels en mathématiques constructives sans axiome du choix dépendant.

\medskip Dans une première partie, nous donnons quelques généralités sur les \tgms et leur version dynamique, les \tdys.

\smallskip La deuxième partie est consacrée à l'étude d'une \tgm finitaire qui a pour ambition de décrire de manière exhaustive les \prts \agqs du corps des réels, et plus \gnlt d'un corps réel clos non discret, du moins celles exprimables dans un langage restreint, proche du langage des anneaux ordonnés. 
On obtient une théorie qui s'avère être, en mathématiques classiques, la théorie des anneaux locaux réels clos. La théorie des anneaux réels clos est présentée ici sous une forme constructive comme une théorie purement équationnelle naturelle, basée sur les fonctions racines virtuelles introduites dans des travaux antérieurs.
Tout ceci constitue un développement, avec quelques modifications terminologiques mineures, des idées données dans l'article \cite{LM2017}.
Enfin, on pose la question de savoir si un axiome d'archimédianité infinitaire permettrait de mieux comprendre la théorie finitaire proposée 

\smallskip Dans la troisième partie, on introduit une théorie plus ambitieuse dans laquelle les fonctions semi-algébriques continues ont droit de cité à part entière: des sortes sont créées pour elles. Cela permet de parler \und{à l'intérieur de la théorie} du module de continuité uniforme d'une \fsagc sur un fermé borné de $\RR^n$. 
Cette nouvelle théorie est résolument infinitaire. On obtient alors une meilleure description des \prts \agqs de $\RR$, mais aussi une première ébauche pour une théorie constructive de certaines structures o-minimales.

\end{quotation}

\rdb
%\addcontentsline{toc}{section}{Table des matières}
\setcounter{tocdepth}{1}
\setcounter{minitocdepth}{3}
\dominitoc

\newpage
\setcounter{page}{0}

\thispagestyle{empty}
~
\newpage

\small
\tableofcontents
\normalsize

\thispagestyle{empty}
%~
%
%
%\section*{\huge Préface}
%\addcontentsline{toc}{chapter}{Préface}
%\markboth{}{Préface}
%\vskip 1cm

\chapter*{Avant-propos}
\addstarredchapter{Avant-propos}
\markboth{Avant-propos}{Avant-propos}

Le mémoire présenté ici est un développement, inachevé, de l'article \cite{LM2017}. Par rapport à cet article, nous avons cependant modifié la \dfn des \fsagcs (\dfn \ref{defiFSAGC2}), dans le même esprit où Bishop définit une fonction réelle continue comme une fonction uniformément continue sur tout intervalle borné. 

Malgré son caractère inachevé et les nombreuses questions que nous ne savons pas actuellement résoudre, nous espérons que ce mémoire suscitera de l'intérêt pour son approche originale du sujet. 

\medskip 
Le mémoire est écrit dans le style des \coma à la Bishop,
i.e.\ les mathématiques avec la logique intuitionniste (voir \cite{Bi67,BB85,BR1987,CACM,MRR,CCAPM}).

\medskip Définissons \textsl{l'\alg réelle} comme l'étude les propriétés algébriques des nombres réels,
i.e.,
les propriétés de $\RR$ formulables au premier ordre sur le langage des anneaux strictement ordonnés, défini par la signature 
\Sigt{\Aso}{\cdot=0,\cdot\geq 0,\cdot>0\mathrel{;}\cdot+\cdot, \cdot\times\cdot,-\cdot, 0,1} 

\noindent avec éventuellement comme constantes tout ou partie des réels constructifs. On peut en outre envisager d'introduire de nouveaux symboles de fonctions pour des fonctions $\RR^n\to\RR$ bien définies (d'un point de vue \cof) et dont la description est purement \agq, comme les fonctions $\sup$,
$\inf$ et de nombreuses \fsagcs définies sur $\QQ$.

L'\textsl{\alg réelle \cov} n'est pas bien comprise! 
L'\textsl{analyse \cov} ($\simeq$ les méthodes certifiées en analyse numérique)
est nettement mieux étudiée.

D'un point de vue \cof, l'\alg réelle est \textsl{assez éloignée} 
de la théorie usuelle classique
des corps réels clos à la Artin-Schreyer-Tarski, 
dans laquelle on suppose que l'on a un 
\textsl{test de signe} pour les réels.

La plupart des \algos de l'\alg réelle classique échouent avec les nombres réels,
parce qu'ils requièrent un {test de signe}.

Même en analyse \cov, on pourrait avoir des retombées intéres\-santes
d'une étude plus approfondie de l'algèbre réelle.
Par exemple cela permettrait de mieux comprendre comment éviter le recours
à l'axiome du choix dépendant (fréquent chez Bishop).
 
\smallskip La compréhension de l'\alg réelle \cov peut \egmt être un premier
pas pour une théorie \cov (et donc \algq) des \textsl{structures o-minimales} (cf.~\cite{cos99}, \hbox{\cite{vdD}}).
La droite réelle et les espaces $\RR^n$ étudiés d'un point de vue purement \agq peuvent être vus comme constituant la plus simple des structures o-minimales. La théorie classique (non \algq) des structures o-minimales donne en effet des pseudo-\algos qui nécessitent pour fonctionner correctement d'avoir au moins un test de signe sur les réels
(il faut aussi introduire des sortes pour les parties définissables de $\RR^n$). Et la théorie des structures o-minimales a à priori un champ d'application très important en analyse.

\medskip Ainsi nous cherchons une \tdy aussi complète que possible pour décrire les \prts \agqs du corps des réels en \coma sans axiome du choix dépendant.

Dans l'étude que nous présentons ici, nous évitons aussi l'usage de la négation. 
Fred Richman \cite{Ric2001} montre que les \coma sont plus élégantes lorsque l'on se passe de l'axiome du choix dépendant. Nous pensons qu'elles sont \egmt plus élégantes si l'on se passe de la négation.

\medskip Dans la première partie, constituée des chapitres \ref{sectgmq}
et \ref{subsecgeominfini} nous donnons quelques généralités sur les \tgms et leur version dynamique, les \tdys. Nous renvoyons pour l'essentiel aux sections 1 à 3 de l'article \cite{LM2022} 

\smallskip La deuxième partie est consacrée à l'étude d'une \tgm qui a pour ambition de décrire de manière exhaustive les \prts \agqs du corps des réels, et plus \gnlt d'un corps réel clos non discret, du moins celles exprimables dans un langage restreint, proche du langage des anneaux ordonnés. Cela constitue un développement, avec quelques modifications terminologiques mineures, des idées données dans l'article \cite{LM2017}.

Le chapitre \ref{chapcoo} propose une définition de la structure de corps ordonné en l'absence d'un test de signe.

 Le chapitre \ref{chap-afr} traite les \afrs et quelques structures dérivées.

 Le chapitre \ref{chapreelclos} essaie de définir la structure de corps ordonné réel clos en l'absence d'un test de signe.

 Le chapitre \ref{secGeomReelsArchi} aborde une \tgm infinitaire pour ajouter l'axiome selon lequel le corps des nombres réels est \textsl{archimédien}.

Donc, à la fin de cette deuxième partie, nous proposons pour la \tdy convoitée celle de la structure d'anneau réel clos local archimédien. La théorie des anneaux réels clos est ici présentée sous une forme \elr, \peq, dans le style de \cite{Tre2007}.

\smallskip Dans la troisième partie, on ajoute les sortes correspondant aux \fsagcs sur les sous-ensembles semialgébriques fermés bornés de $\RR^n$. Nous espérons obtenir ainsi une description plus précise de l'\alg réelle et pouvoir esquisser une première théorie \cot satisfaisante pour les structures \hbox{o-minimales}.

\smallskip Dans tout le texte, les \thos ou lemmes de \clama qui n'ont pas de \demo \cov connue, et qui souvent ne peuvent pas en avoir, sont indiqués avec une étoile.

\smallskip Indiquons enfin que l'article \cite{Lom11} contient des réflexions, dans un cadre plus philosophique, analogues à celles proposées ici.

\medskip \noindent {\bf Remerciements.} Nous remercions Michel Coste et Marcus Tressl pour leurs patientes réponses à nos nombreuses questions.

\bigskip 
\begin{flushright}
%\centerline{
Henri Lombardi, Assia Mahboubbi, \today
\end{flushright}

\part{Théories \gmqs}

\chapter*{Introduction}
\addstarredchapter{Introduction}

\rdb

Cette première partie fait l'objet d'un article plus développé
en préparation \cite{Lom-tgac}, que l'on trouve en: \url{http://hlombardi.free.fr/Theories-geometriques.pdf}

Nous donnons les principales \dfns et renvoyons pour l'essentiel aux sections 1 à 3 de l'article \cite{LM2022}

Une \tdy peut être comprise comme une formalisation d'un morceau bien délimité de \maths intuitives. Ces \maths intuitives, pratiquées par la communauté \mathe, sont étudiées sous une forme complètement calculatoire indépendante de tout point de vue philosophique. Mais là où le point de vue classique utilise librement le principe du tiers exclu et l'axiome du choix, les \tdys remplacent ces outils non calculatoires par le point de vue dynamique des structures incomplètement spécifiées, qui est le point de vue de l'évaluation paresseuses en Calcul Formel. 

\smallskip 
Le chapitre \ref{chap-gmqfini} traite les \tdys finitaires. 

Une \tgm finitaire 
correspond à ce que l'on appelle en \clama une \tfo cohérente.
Mais la \tgm que nous considérons est gouvernée par la logique intuitionniste alors que la théorie cohérente est en général gouvernée par la logique classique. 

En outre la \tdy correspondante est une version minimaliste de la \tgm: 
c'est une pure machinerie calculatoire sans logique, un peu dans le même esprit que l'\ari récursive de Goodstein (\cite{Goo1957}). 

Une \tdy finitaire est donc la mise en forme purement calculatoire, sans logique, d'une \tfo du premier ordre cohérente. La logique classique qui est normalement utilisée dans la \tfo est remplacée par un calcul arborescent sans logique. 

On peut \egmt voir la \tdy comme une version partielle de la déduction naturelle dans laquelle les formules examinées sont toutes d'un type très simple, sans le connecteur d'implication (donc sans la négation) et avec une utilisation très limitée des quantificateurs. 

La surprise est que les \tdys sont cependant très expressives (en \clama toute \tfo du premier ordre peut être vue comme une théorie cohérente) et qu'elles effacent la distinction entre logique classique et logique intuitionniste.

Dans le cas fréquent où la signature est un ensemble dénombrable et où les axiomes forment une partie décidable du langage, le monde \mathe extérieur à la théorie, qui est celui où nous nous situons pour étudier une structure donnée, voir fonctionner le système formel qui la décrit, et établir des \thos le concernant, n'a aucune interférence avec la \tdy elle-même. C'est ce qui est confirmé de manière \gnle par le \tho fondamental \ref{thFond}: si l'on force une \tgm finitaire à se comporter de manière classique, les \rdys écrites dans le langage de départ qui sont valides après étaient déjà valides avant.

Cela correspond au fait que les topos cohérents de Grothendieck, qui sont une autre mise en forme des théories cohérentes, ont une logique interne intuitionniste mais qu'ils peuvent cependant être compris de différentes manières selon que l'on se situe dans un monde \mathe extérieur \cof ou non\footnote{\Cad essentiellement selon que l'on accepte le principe du tiers exclu ou non dans ce monde extérieur.}. 

En \coma, seules certaines structures relèvent des \tdys finitaires. Par exemple la structure de corps discret mais pas celle de corps de Heyting\footnote{Un corps discret est un anneau non trivial dans lequel tout \elt est nul ou inversible, un corps de Heyting est un anneau local dans lequel tout \elt non \iv est nul. Les \clama ne connaissent pas la distinction pertinente entre la notion de corps discret et celle de corps de Heyting.}\index{corps!de Heyting}\index{Heyting!corps de ---}.
En cela, le point de vue restreint des \tdys ouvre la voie à une classification pertinente, invisible en \clama, des degrés de complexité des êtres mathématiques inventés par les êtres humains.

\smallskip 
Le chapitre \ref{chap-gmqinfini} traite les \tdys infinitaires, dans lesquelles des disjonctions infinies sont autorisées dans la conclusion d'une \rdy.

Une restriction essentielle est à noter: les variables libres présentes dans une telle disjonction doivent être précisées d'avance et en nombre fini.

Intuitivement, on utilise
de telles règles dans le \sys de preuves des \tdys en \gui{ouvrant les branches de calcul correspondant à la disjonction infinie}. Qu'est-ce que cela signifie \prmt? Cela signifie qu'une conclusion sera déclarée valide si elle
est valide dans chacune des branches. 

 Ces théories sont plus expressives que les théories finitaires et permettent d'axiomatiser un très grand nombre de structures \maths usuelles. 
 
 Contrairement à ce qui se passe pour les \tdys finitaires, le monde \mathe extérieur intervient de manière inévitable pour certifier la validité d'une \rdy.
 
 Prenons un exemple simple, et indiquons ce qui se passe si on a dans les axiomes une règle infinitaire du type

\Regles{\lab {~} $\Vdi{x_1,\dots,x_k} \Vou_{i\in I}\; \Gamma_i$
}

\noindent avec un ensemble infini $I$ et les $\Gamma_i$ des listes de formules atomiques sans autres variables libres que celles mentionnées (i.e. $\xk$). Si pour chaque $i\in I$ on a une règle valide $\,\,\Gamma_i\vd B(\ux)$, alors on déclare valide la règle $\Vdi{x_1,\dots,x_k}B(\ux)$. 

Il intervient donc \ncrt une \demo
intuitive \und{externe} à la \tdy pour certifier que la conclusion souhaitée est valide dans chacune des branches. En effet le \sys de calcul \gui{sans logique} à l'œuvre dans la \tdy ne peut pas prendre en charge
une telle \gui{infinité} de \demos. Un calcul purement mécanique ne saurait ouvrir une infinité de branches! Par exemple, avec $I=\NN$ la \demo
intuitive externe pourra être une preuve par \recu.

Notons par contre que la \demo interne doit démontrer la validité de la conclusion souhaitée selon les règles de \demo \gui{sans logique} de la \tdy.

\medskip \noindent 
{\bf Terminologie.} 
Comme nous nous situons en \coma il apparait inévi\-tablement des \pbs terminologiques, du simple fait par exemple qu'en \gnl un même concept classique donne lieu à plusieurs concepts \cofs intéressants non \eqvs, mais \eqvs en \clama. 

 Nous donnons ci-dessous des petits tableaux comparatifs entre notre terminologie (en \coma) et la terminologie anglaise la plus usuelle (en \clama) pour ce qui concerne les \tgms. Celle que l'on trouve dans
\cite[Chapitre~D1]{Joh2-02}, \cite{Car2017} et dans~\cite{BH2017}.

La comparaison est un peu biaisée par le fait que dans les \tdys on n'utilise pas la logique à proprement parler. Ce sont de pures machines de calcul. Ainsi, bien qu'une \tdy finitaire \gui{engendre} une théorie (formelle du premier ordre) cohérente
et bien que toute théorie cohérente admette une version \gui{\tdy finitaire}, il ne s'agit pas des mêmes objets formels. Témoin le fait qu'une théorie cohérente ne fonctionne pas de la même manière avec la logique classique et avec la logique intuitionniste, alors qu'une \tdy est insensible à cette distinction, car structurellement les \demos dynamiques sont toujours \covs. 

%:\newpage
\newpage

%%%%%%%%%%%%%%%%%%%%%%%%%%%%%%%%%%%%%%%%%%%%%%%%%%%%%%%%%%%%%%%%%%%%
\begin{center}
 \tabcolsep0pt\renewcommand{\arraystretch}{0}%
\begin{tabular}{|c|}
\hline 
\Boite{.7}{12}{Théories finitaires}\\
\hline 
\end{tabular}
\begin{tabular}{|c|c|}
%\hline
\Boite{.7}{7.5} {Théorie}&
\Boite{.7}{4.5}{Theory}\\
\hline
\Boite{.6}{7.5} {purement équationnelle}&
\Boite{.6}{4.5}{algebraic}\\
\hline
\Boite{.6}{7.5} {directe}&
\Boite{.6}{4.5}{}\\
\hline
\Boite{.6}{7.5} {algébrique}&
\Boite{.6}{4.5}{Horn}\\
\hline
\Boite{.6}{7.5} {disjonctive}&
\Boite{.6}{4.5}{}\\
\hline
\Boite{.6}{7.5} {propositionnelle}&
\Boite{.6}{4.5}{propositional}\\
\hline
\Boite{.6}{7.5} {existentielle}&
\Boite{.6}{4.5}{regular}\\
\hline
\Boite{.6}{7.5} {existentiellement rigide}&
\Boite{.6}{4.5}{}\\
\hline
\Boite{1}{7.5} {existentielle existentiellement rigide, ou\\ cartésienne}&
\Boite{1}{4.5}{cartesian}\\
\hline
\Boite{.6}{7.5}{rigide}&
\Boite{.6}{4.5}{disjunctive, \cite{Johnstone79}}\\
\hline
\Boite{.6}{7.5} {dynamique finitaire}&
\Boite{.6}{4.5}{}\\
\hline
\Boite{.6}{7.5} {cohérente intuitionniste}&
\Boite{.6}{4.5}{}\\
\hline
\Boite{.6}{7.5} {cohérente classique}&
\Boite{.6}{4.5}{coherent}\\
\hline
\end{tabular}
\end{center}

%\pagebreak	
%%%%%%%%%%%%%%%%%%%%%%%%%%%%%%%%%%%%%%%%%%%%%%%%%%%%%%%%%%%%%%%%%%%%
\begin{center}
\tabcolsep0pt\renewcommand{\arraystretch}{0}%
\begin{tabular}{|c|}
\hline 
\Boite{.7}{12}{Théories générales (infinitaires)}\\
\hline 
\end{tabular}
\begin{tabular}{|c|c|}
\Boite{.7}{7.5} {Théorie}&
\Boite{.7}{4.5}{Theory}\\
\hline
\Boite{.6}{7.5} {dynamique (infinitaire)}&
\Boite{.6}{4.5}{}\\
\hline
\Boite{.6}{7.5} {\quad géométrique}&
\Boite{.6}{4.5}{geometric intuitionnist}\\
\hline
\Boite{.6}{7.5} {\quad géométrique classique}&
\Boite{.6}{4.5}{geometric}\\
\hline
\end{tabular}
\end{center}

%%%%%%%%%%%%%%%%%%%%%%%%%%%%%%%%%%%%%%%%%%%%%%%%%%%%%%%%%%%%%%%%%%%%
\begin{center}
\tabcolsep0pt\renewcommand{\arraystretch}{0}%
\begin{tabular}{|c|c|}
\hline
\Boite{.7}{7.5}{Théories géométriques}&
\Boite{.7}{4.5}{Geometric theories}\\
\hline
\Boite{.6}{7.5}{identiques (même signature)}&
\Boite{.6}{4.5}{equivalent}\\
\hline
\Boite{.6}{7.5}{essentiellement identiques (mêmes sortes)}&
\Boite{.6}{4.5}{~}\\
\hline
\Boite{.6}{7.5}{classiquement essentiellement identiques}&
\Boite{.6}{4.5}{definitionally equivalent}\\
\hline
\Boite{.6}{7.5}{essentiellement équivalentes}&
\Boite{.6}{4.5}{}\\
\hline
\Boite{.6}{7.5}{classiquement essentiellement équivalentes?}&
\Boite{.6}{4.5}{Morita equivalent}\\
\hline
\end{tabular}
\end{center}

\newpage \thispagestyle{empty}

\chapter{Théories \gmqs 
finitaires}\label{sectgmq}\label{chap-gmqfini}
\index{theorie@théorie!géométrique finitaire}\index{theorie@théorie!géométrique}

\Today

\minitoc

\section{Théories cohérentes et théories dynamiques}\label{subsectdy}
%%%%%%%%%%%%%%%%%%%%%%%%%%%%%%%%%%%%%%%%%

%: Subsection{Théories cohérentes}
\Subsection{Théories cohérentes}\index{theorie@théorie!cohérente}
Une \textsl{théorie cohérente} $\sa{T}=(\cL,\cA)$ est une \tfo 
du premier ordre basée sur le langage~$\cL$ dans laquelle les axiomes (les \elts de $\cA$)
sont tous \gui{\gmqs}, \cad de la forme suivante:
%
%---- equation {eqAgeom} ----
\begin{equation} \label{eqAgeom}
\forall \und x \;\;\big(C\; \Longrightarrow \; \exists\, \und{y^1} \,D_1\;
\vee\;\cdots\;\vee\;\exists\,\und{y^m}\,D_m\big)
\end{equation}
%------end equation----
où $C$ et les $D_j$ sont des \textsl{conjonctions de formules atomiques} du langage 
$\cL$ de la \tfo, les $\und{y^j}$ sont des listes de variables,
et $\und x$ la liste des autres variables présentes (ces listes sont éventuellement vides). Les variables de $C$ sont uniquement dans la liste $\und x$. Les variables de $D_j$ sont uniquement dans les listes disjointes $\und x$ et $\und{y^j}$.
Une disjonction vide au second membre peut être
remplacée par le symbole $\bot$ représentant le $\Faux$.

On dit aussi \textsl{\tgm finitaire} à la place de théorie cohérente lorsqu'on utilise la logique intuitionniste.\index{theory!coherent ---}\index{theorie@théorie!géométrique finitaire}

\Cadre{.9}{\noindent Dans la suite du chapitre \ref{sectgmq}, nous omettons presque toujours le qualificatif \gui{finitaire} après \gui{\tgm} ou \gui{\tdy}.}

%: Subsection{Théories dynamiques finitaires}
\Subsection{Théories dynamiques finitaires}\index{theorie@théorie!dynamique}\label{sectdyfinitaire}

Référence principale \cite{CLR01}. Dans l'article en question sont introduites les notions de \gui{dyanamical theory} et de \gui{dynamical proof}.
Voir \egmt: l'article \cite[Bezem \& Coquand, 2005]{BC2005} qui décrit un certain nombre d'avantages fournis par cette approche, et les articles précurseurs \cite[Prawitz 1971, sections 1.5 et 4.2]{pra1971}, \cite[Matijasevi\v c 1975]{Mat75} et \cite[Lifschitz, 1980]{Lif80}.

Si \sa{T} est une \tgm finitaire, la \textsl{\tdy} finitaire correspondante s'en différencie seulement 
par un usage extrêmement limité des méthodes de \demo:
%-------begin item---
\begin{itemize}
\item Premièrement, on n'utilise jamais d'autres formules que les formules 
atomiques: on n'introduit jamais de nouveau prédicat utilisant des connecteurs 
logiques ou des quantificateurs. Seules sont manipulées des listes de formules 
atomiques du langage~$\cL$.
\item Deuxièmement, et conformément au point précédent, les axiomes ne sont 
pas vus comme des formules vraies, mais comme des \textsl{règles de déduction}: un axiome tel que \pref{eqAgeom} est utilisé en tant que règle 
\pref{eqRgeom}\label{NOTAvou} 
%---- equation {eqRgeom} ----
\begin{equation} \label{eqRgeom}
\Gamma \vd \EXists{\und{y^1}} \Delta_1
\vou \cdots \vou \EXists{\und{y^m}} \Delta_m
\end{equation}
%------end equation----
Ici les conjonctions de formules atomiques $C$, $D_1$, \dots, $D_m$ de \pref{eqAgeom} ont été remplacées par les listes correspondantes
$\Gamma$, $\Delta_1$, \dots, $\Delta_m$.
\item Troisièmement, on ne prouve que des \textsl{\rdys}, \cad des \thos qui sont 
de la forme des règles de déduction ci-dessus.\index{regle@règle!dynamique} 
\item Quatrièmement, la seule
manière de prouver une \rdy est un calcul arborescent
\gui{sans logique}. À la racine de l'arbre se trouvent 
toutes les hypothèses
du \tho que l'on veut valider. L'arbre se développe en appliquant les axiomes 
selon une pure machinerie de calcul algébrique dans la structure.
Voir l'exemple \ref{exasaCd} ci-après. Les \dfns formelles précises sont données dans \cite{CLR01}, on les étend au cas où il y a plusieurs types d'objets comme dans la théorie des modules sur un anneau commutatif avec les objets du type \gui{\elts de l'anneau} et les objets du type \gui{\elts du module}. 
\end{itemize}
%-------end item---

\smallskip Lorsque l'on applique un axiome tel que \pref{eqRgeom}, on substitue aux variables libres $(x_i)$ présentes dans l'hypotèse des termes arbitraires $(t_i)$ du langage. Si les hypothèses, réécrites avec ces termes, sont déjà prouvées, alors on ouvre $m$ branches de calcul dans chacune desquelles on introduit des variables fraiches correspondant aux
variables muettes $\und{y^k}$ (il faut éventuellement changer leurs noms pour éviter un conflit avec les variables libres présentes dans les termes $t_i$) et chaque conclusion~$B_k$ est valide dans sa branche\footnote{$B_k$ est la liste $\Delta_k$ dans laquelle les variables $x_i$ ont été remplacées par les termes $t_i$.}.

Les exemples \ref{exasaCd} très \elrs montrent comment on valide une \rdy dans une \tdy donnée. On développe un arbre de calcul en utilisant les axiomes de la \tdy comme indiqué précédemment et on a gagné lorsqu'à chacune des feuilles de l'arbre, la conclusion est validée.

%: Example{exasaCd}
\begin{example} \label{exasaCd} 
La \tdy \SA{Cd} des \textsl{corps discrets} est basée sur le langage des 
anneaux commutatifs et elle a pour axiomes ceux des anneaux commutatifs non triviaux (théorie \Sa{Ac} dans l'exemple \ref{exaAc}) et 
la \rdy des corps discrets:

\Regles{\Lab{CD} $\vd x=0\vou\EXists {y} xy=1$}%\label{AxCO}

\smallskip \noindent 1) Pour démontrer la \rdy

\Regles{\Lab{ASDZ} $\,\,xy=0\vd x=0 \vou y=0$}

\noindent on ouvre deux branches conformément à l'axiome \tsbf{CD}.
Dans la première on a $x=0$ et la conclusion est prouvée.
Dans la deuxième on introduit un \gui{paramètre} (une variable fraiche) $z$ avec la 
relation~\hbox{$xz=1$}. Les axiomes des anneaux commutatifs permettent alors de démontrer les \egts $y=1\times y=(xz) y=(xy)z=0\times z=0$.\\
La conclusion est donc prouvée à chacune des deux feuilles de l'arbre de calcul. 

\smallskip \noindent 2) Ensuite par exemple, on déduit de la \rdy précédente la \ralg

\Regles{\Lab{Anz} $z^2=0 \vd z=0$}

\smallskip \noindent car cette fois-ci aux deux feuilles de l'arbre on a la même conclusion $z=0$.

\smallskip \noindent 3) La théorie \SA{Al} des anneaux locaux est basée sur le langage des 
anneaux commutatifs et elle a pour axiomes ceux des anneaux commutatifs (théorie \Sa{Ac0} dans l'exemple \ref{exaAc}) et 
la \rdy des anneaux locaux

\Regles{\Lab{AL} $\,\, (x+y)z=1 \vd \EXists u \; xu=1 \;\vou\;\EXists u\;yu=1$}

\noindent Pour démontrer qu'un corps discret vérifie la règle \tsbf{AL},
on ouvre deux branches conformément à l'axiome \tsbf{CD}.
Dans la première on a $x=0$ et la conclusion est prouvée car $(x+y)z=1$ donne $yz=1$.
Dans la deuxième on introduit un \gui{paramètre} (une variable fraiche) $v$ avec la relation~\hbox{$xv=1$}.\\
La conclusion dans la règle \tsbf{AL} est donc prouvée aux deux feuilles de l'arbre de calcul.
\eoe
\end{example}
%--------- fin example ----------------------------------------

Notons aussi que la validité de la règle suivante, que l'on pourrait appeler \gui{L'existence concrète implique l'existence formelle}, est purement tautologique. 
 
On considère une liste $\Gamma(\ux,\uy)$ de formules atomiques dans une \tdy $\sa{T}$. On note $\Gamma(\ux,\underline t)$ la liste de ces formules dans lesquelles on a substitué à chaque variable $y_j$ un terme $t_j$ construit sur les $x_i$ et sur les constantes de la théorie. Alors la règle existentielle suivante est valide. 
\[
\Gamma(\ux,\underline t)\vd \Exists y_1,\dots,y_m \;\Gamma(\ux,\uy).
\]
%

%:paragraph{La logique remplacée par le calcul}
\paragraph{La logique remplacée par le calcul}~

\smallskip En pratique, démontrer une \rdy dans le cadre d'une \tdy suit toujours un raisonnement naturel intuitif et l'on peut voir cette gymnastique comme une version simplifiée de la déduction naturelle de Gentzen.
Le symbole $\vou$ doit être compris comme une abréviation pour 
\gui{{\bf ou}vrir (des branches dans le calcul)}.

\smallskip Les symboles $\vou$ et $\;\EXists{\cdot}$ ont été 
préférés à $\vuu$ et $\exists$, pour bien marquer que leur utilisation
dans les règles de déduction n'est pas l'utilisation de nouvelles formules construites à partir des formules atomiques.
Le symbole $\,\vd\,$ a été préféré à $\vdash$ pour éviter la confusion avec le symbole utilisé pour les relations implicatives dans les \trdis. Notons aussi qu'il n'a pas la même interprétation que le symbole analogue utilisé dans les calculs de séquents à la Gentzen. 

\smallskip Ainsi le langage d'une \tdy ne comporte aucun symbole logique (connecteur ou quantificateur) permettant de construire des formules compliquées à partir des formules atomiques. La \gui{logique} est remplacée par les symboles $\vdi$, $\vou$ et $\;\EXists{\cdot}$ et par le séparateur \gui{$\vet$}, mais ces symboles sont utilisés pour décrire une machinerie de calculs arborescents et non pour former des formules.
La partie non logique d'une \tdy est constituée de symboles pour les variables, et de la \textsl{signature}, qui contient les symboles pour les sortes, les prédicats et les fonctions (ou lois) définies dans la structure.

\Cadre{.9}{\noindent Dans la suite, nous remplaçons \gui{$\EXists{\cdot}$} 
par le symbole moins encombrant \gui{$\Exists$}, plus proche et néanmoins différent du traditionnel \gui{$\exists$}.\label{NOTAExists}}

%:paragraph{Le prédicat d'égalité}
\paragraph{Le prédicat d'égalité} ~

\smallskip 
Dans une \tdy chaque sorte doit être munie d'un prédicat d'\egt $\cdot=\cdot$ et l'on donne les axiomes qui autorisent la substitution d'un terme $t$ par un terme $t'$ lorsque $\vd t=t'$ est valide dans la théorie\footnote{Cela exclut le cas où $t$ contient une variable $x$ sous la dépendance d'un $\Exists \,x$.} en n'importe quelle occurrence d'une formule atomique présente dans une \rdy valide.

On pourrait aussi bien ne pas donner d'axiomes relatifs à cette substitution et considérer qu'elle est tout simplement une procédure de calcul légitime. 
 
%\vspace{-.2em} 
\paragraph{Extension simple d'une \tdy}
%d
%: Definition{defiextsimple}
\begin{definition} \label{defiextsimple}
On dit que \textsl{la \tdy $\Tp=(\cL',\cA')$ est une extension simple de la \tdy $\sa T=(\cL,\cA)$} si $\cL\subseteq \cL'$ et $\cA\subseteq \cA'$.
Dans ce cas les \rdys formulées dans $\cL$ et valides (i.e. démontrables) dans \sa T sont valides dans $\Tp$.\index{extension --- d'une théorie dynamique!simple} 
\end{definition}
%----------- fin definition -------------------------------- 

%r
%: Remark{remdefiextsimple}
\begin{remark} \label{remdefiextsimple} 
 Dans la \dfn précédente, l'expression \gui{extension simple} peut être remise en question.
Si $\cL,\cA,\cL',\cA'$ sont des ensembles finis, ou 
si ce sont des ensembles dénombrables discrets, on peut considérer que tout est clair intuitivement. Il peut arriver cependant que l'on désire 
utiliser des ensembles plus compliqués, par exemple introduire tous les réels comme des constantes dans une théorie dont une sorte vise à décrire les nombres réels. Dans un tel cas, le mot \gui{extension simple} est contestable car il n'y pas de monomorphisme canonique dans la catégorie~$\slbSet$ de Bishop: dans la conception de Bishop, une partie d'un ensemble correspond à la notion catégorique de sous-objet.
La \gui{simplicité} n'est donc pas dans ce cadre une notion objective, ou si l'on préfère, elle n'a pas de \dfn \mathe précise. 
 \eoe
\end{remark}
%----------- fin remark ---------------------------------- 

%:Subsection{Règles structurelles}
\Subsection{Règles structurelles}

Ici on donne des \textsl{règles structurelles admissibles} pour une \tdy. Il s'agit de règles de déduction \textsl{externes} (différentes des \rdys, qui sont internes à la théorie). Elles disent que si certaines \rdys sont valides, alors d'autres \rdys sont automatiquement valides.\index{regle@règle!de déduction externe}\index{regle@règle!admissible}\index{regle@règle!structurelle}

Voici les règles structurelles admissibles qui nous semblent les plus importantes.

%: Regles structurelles {rstr1}
\begin{rstra} \label{rstr1} ~ 
\begin{enumerate}\setcounter{enumi}{-1}
\item \textsl{Variables libres, variables muettes}
\begin{enumerate}
\item \label{6rstr} \textsl{Substitution.} On peut dans une \rdy remplacer toutes les occurrences d'une variable libre par un terme à condition de ne jamais créer de conflit entre variables libres et muettes. 
\item \label{5rstr} \textsl{Renommage.} On peut dans une \rdy renommer les variables libres 
ou les variables muettes (celles présentes dans les $\Exists$) à condition de ne jamais créer de conflit entre variables libres et muettes. 
\end{enumerate}

\item \textsl{Profiter du travail déjà accompli}
\begin{enumerate} 
\item \label{12rstr} \textsl{Raccourcis}. Lorsque l'on a démontré la validité d'une \rdy on peut la rajouter dans les axiomes de la théorie.
\item \label{11rstr} \textsl{Renforcement simultané de l'hypothèse et des conclusions}.
On considère dans une \tdy une règle valide 
\[
\Gamma \vd \Exists{\und{y^1}} \Delta_1 \vou \cdots\vou \Exists{\und{y^m}} \Delta_m
\]
Soit $A$ une formule atomique qui ne fait intervenir aucune des variables
existentielles du second membre. On note $\Gamma'$ la liste $\Gamma$ suivie de $A$ et $\Delta'_i$ la liste $\Delta_i$ suivie de $A$. Alors la règle
suivante est \egmt valide:
\[
\Gamma' \vd \Exists{\und{y^1}}\, \Delta'_1
\vou \cdots\vou \Exists{\und{y^m}}\,\Delta'_m
\]
\end{enumerate}
\item \textsl{Les listes fonctionnent comme des ensembles finis}
\begin{enumerate} 
\item \label{1rstr} \textsl{Permutation des formules atomiques figurant dans une liste.} 
\item \label{2rstr} \textsl{Contraction.} Si deux formules atomiques identiques figurent dans une liste, on peut supprimer l'une des deux. \\
Inversement, on peut dupliquer une formule atomique dans une liste arbitraire.
\item \label{3rstr} \textsl{Monotonie}. On peut ajouter dans la liste à gauche de $\vd$ des formules atomiques comme il nous sied. 
\item \label{4rstr} \textsl{Permutation, contraction et monotonie pour les $\vou $ à droite du $\vd$}. 
\end{enumerate}
\item \label{10rstr} \textsl{Les listes de formules atomiques comme conjonctions}
\begin{enumerate}
\item \label{10rstr1} \textsl{Prouver une liste de formules atomiques, c'est prouver chacune d'entre elles.}\\
On considère dans une théorie~\sa{T} une \rdy 
$
\Gamma\vd (A_1\vet \dots\vet A_n).
$ 
Cette \rdy est valide \ssi sont valides les règles 
$
\,\,\Gamma\vd A_k \quad (k\in\lrbn).
$ 
\item \label{10rstr2}
\textsl{Distributivité du $\vou $ sur le \gui{et} implicite des listes.}\\
On considère dans une théorie~\sa{T} une \rdy 
$$
\Gamma\vd (A_1\vet \dots\vet A_n) \vou \Exists{\und{y^1}}\,\Delta_1\vou \dots \vou \Exists{\und{y^m}}\,\Delta_m.
$$ 
Cette \rdy est valide \ssi sont valides les \rdys 
$$
\Gamma\vd A_k \vou \Exists{\und{y^1}}\,\Delta_1\vou \dots \vou \Exists{\und{y^m}}\,\Delta_m\quad \quad (k\in\lrbn).
$$ 
\end{enumerate}

\item \textsl{Transitivité et variantes}
\begin{enumerate}
\item \label{7rstr} \textsl{Transitivité.}\label{Transitivite} On donne un exemple en laissant le soin \alec de donner la formulation \gnle. Supposons qu'on ait dans une \tdy \sa{T} des \rdys valides 
\[ 
\begin{array}{rcl} 
\Gamma(\ux) & \vd & \Exists y,z\; \Delta_1(\ux,y,z)\vou \Exists u \;\Delta_2(\ux,u), \\[.3em] 
\Gamma(\ux),\Delta_1(\ux,y,z) & \vd & \Exists r,s,t\; \Delta_3(\ux,y,z,r,s,t), \\[.3em] 
\Gamma(\ux),\Delta_2(\ux,u) & \vd & \Exists v\; \Delta_4(\ux,u,v)
\vou \Exists w\; \Delta_5(\ux,u,w). 
 \end{array}
\]
Alors est \egmt valide la règle
$$
\Gamma(\ux)\vd\Exists y,z,r,s,t\; \Delta_3(\ux,y,z,r,s,t)\vou \Exists u,v\; \Delta_4(\ux,u,v)\vou \Exists u,w\; \Delta_5(\ux,u,w).
$$
\item \label{8rstr} \textsl{Coupure.} Considérons des listes de formules atomiques $\Gamma(\ux),\Delta_0(\ux),\Delta_1(\ux),\dots,\Delta_m(\ux)$ $(m\geq 1)$ dans une \tdy $\sa{T}$. Si les deux \rdys $$\Gamma\vd\Delta_0\vou \Delta_1\vou\dots\vou\Delta_m\quad \hbox{ et }\quad \Gamma,\Delta_0\vd\Delta_1\vou\dots\vou\Delta_m $$ sont valides dans $\sa{T}$, alors la règle
$\Gamma\vd\Delta_1\vou \dots\vou \Delta_m$ est \egmt valide.
\item \label{8rstr2} \textsl{Coupure avec existence.} 
Une version plus \gnle est la suivante. 
Considérons des listes de formules atomiques $\Gamma(\ux),\Delta_0(\ux,\und{y^0}),\Delta_1(\ux,\und{y^1}),\dots,\Delta_m(\ux,\und{y^m})$ $(m\geq 1)$ dans une \tdy $\sa{T}$.
 Si les deux \rdys 
\[\hspace{-1.5em}
\Gamma\vd\Exists{\und{y^0}}\,\Delta_0\vou \Exists{\und{y^1}}\,\Delta_1\vou \dots\vou \Exists{\und{y^m}}\,\Delta_m \; \hbox{ et }\; \Gamma,\Delta_0\vd\Exists{\und{y^1}}\,\Delta_1\vou \dots\vou \Exists{\und{y^m}}\,\Delta_m 
\]
sont valides dans $\sa{T}$, alors la règle
$\Gamma\vd\Exists{\und{y^1}}\,\Delta_1\vou \dots\vou \Exists{\und{y^m}}\,\Delta_m$ est \egmt valide.
\end{enumerate}
\end{enumerate}
 \end{rstra}
%----------- fin fact ----------------------------------------- 

%subsection{Collapsus}
%: Subsection{Collapsus}
\Subsection{Collapsus}\index{collapsus}\index{effondrement}\index{regle@règle!de collapsus}\label{NOTABot}

Une \rdy s'appelle une \textsl{règle de collapsus} ou \textsl{d'effondrement}\index{effondrement} lorsque le second membre est \gui{le $\Faux$}, que l'on note $\Bot$. 
On peut aussi voir $\Bot$ comme désignant la disjonction vide. Lorsque l'on a prouvé $\Bot$, l'univers du discours s'effondre, et toute formule atomique du langage est alors réputée \gui{vraie}, ou du moins \gui{valide}. C'est l'application de la règle \gui{ex falso quod libet}, qui est la signification intuitive pertinente du $\Faux$ en \coma. Ainsi dans les \tdys les règles 

\Regles{\lAb{Faux$_{P}$} $\,\,\Bot\vd P$ \quad (ex falso quod libet)} 

\noindent sont valides pour toutes les formules atomiques.

On donne aussi dans le langage la constante logique $\Top$ pour \gui{le $\Vrai$}, avec pour axiome la \ralg. \label{NOTATop} 

\Regles {\lAb{Vrai}$\,\,\vd \Top$}

On peut aussi voir $\Top$ comme désignant la conjonction vide\footnote{Quand il n'y a rien à démontrer, ne démontrons rien et tout sera OK. Par ailleurs, dans une \tdy comportant au moins une sorte $\iS$, $\Top$ est \eqv à $x=_\iS x$.}.
Les constantes \(\Bot\) et \(\Top\) sont (les seuls vrais) symboles logiques dans les \tdys.

Lorsqu'une \tdy ne comporte pas de règle de collapsus, elle admet toujours le modèle réduit à un point\footnote{S'il y a plusieurs sortes, chaque sorte est réduite à un point.} où toutes les formules atomiques sont évaluées vraies.
C'est l'objet final dans la catégorie des modèles de la théorie.

On peut dire qu'une \tdy sans règle de collapsus \textsl{s'effondre} si toutes les formules atomiques sont valides, à l'exception de $\Bot$.

On dit qu'une \tdy avec règle de collapsus \textsl{s'effondre} lorsque $\Bot$ est prouvable, et par suite toutes les \rdys aussi. Dans ce cas la théorie n'admet pas de modèle.
 
Considérer l'effondrement dans le sens du seul modèle réduit à un point, plutôt que dans le sens du pur néant, est seulement une affaire de goût qui ne change rien au fond des choses\footnote{En fait, un des auteurs doit avoir horreur du vide:-{\tiny)}: le silence de cet espace infini l'effraie :-{\tiny(}\,. Par ailleurs si la disparition totale dans le néant est la vraie signification du $\Faux$, il n'en reste pas moins que, avant même d'interdire l'existence des modèles, le $\Faux$ commence par les réduire à un seul point, qui satisfait tous les prédicats. Comme dit la chanson de Boris Vian: \gui{on est descendu chez Satan et en bas c'était épatant!}.}. 

Au lieu de dire qu'une \sad qui s'effondre n'a pas de modèle, on dit (sans négation) que tout modèle de cette \sad est trivial, réduit à un point, et que \gui{tout y est vrai}.

Pour concilier formellement ces deux points de vue, la meilleure solution semble être la suivante: chaque sorte $\iS$ introduite est accompagnée d'au moins deux constantes dans cette sorte, disons $0_\iS$ et $1_\iS$ pour fixer les idées, avec l'axiome \fbox{$0_\iS=_\iS 1_\iS\vd\Bot$}. Dans la suite c'est ce qui se produirait normalement pour la théorie des \trdis non triviaux et la théorie des anneaux commutatifs non nuls, ainsi que dans toutes leurs extensions. Mais nous préférons utiliser la convention suivante.

\Cadre{.9}{\noindent Dans tout ce mémoire, nous considérons, dans le cas d'un anneau ou d'un \trdi que le collapsus est toujours donné sous la forme $1=0$ ou une formule du même style, par exemple $0>0$ pour un corps ordonné. Un inconvénient de l'usage du $\Bot$ est qu'il fait sortir du cadre des \talgs quand on pourrait y demeurer.
\Llec qui le désire pourra ajouter un axiome du type $1=0\vd \Bot$.}

\Subsection{Classification des \tdys} 

\paragraph{Théories \agqs}~

\smallskip \noindent Une \rdy qui ne contient à droite du symbole $\vd\,$ ni $\vou $, ni $\,\Exists\,$, ni $\Bot$ est appelée une \textsl{règle algébrique}. Une \tdy est dite \textsl{\agq} lorsqu'elle ne comporte comme axiomes que des \ralgs. Dans la littérature anglaise, les \talgs ont pour nom \textsl{Horn theory}. Un cas particulier est fourni par les \tpes qui sont les \talgs avec une seule sorte et pour seul prédicat le prédicat 
d'\egt.%
\index{theorie@théorie!algébrique}\index{regle@règle!algébrique}%
\index{theory@Horn ---}

\smallskip
\Note Nous reprenons la terminologie suivante de \cite{CLR01}. Une \ralg est dite \textsl{directe} lorsque, à gauche du symbole $\vd\,$, on a uniquement des formules atomiques portant sur des variables différentes les unes des autres ou sur des constantes. Les autres \ralgs sont dites \textsl{règles de simplification}. Par exemple dans les deux \ralgs ci-dessous, la première est directe, la seconde non.\index{regle@règle!directe}\index{regle@règle!de simplification}

\DeuxRegles
{\labu $x=0\vet y=0\vd x+y=0$}
{\labu $x+y= 0\vet x=0\vd y=0$ }

\paragraph{Théories disjonctives}~

\smallskip \noindent Une \tdy est dite \textsl{\dij} si dans les axiomes, il n'y a pas de $\Exists$ à droite du~$\vd$\index{theorie@théorie!disjonctive}\index{regle@règle!disjonctive}. 

\paragraph{Théories existentielles}~

\smallskip \noindent Une \rdy est dite \textsl{existentielle simple} si le second membre (la conclusion) est de la forme $\Exists \uy\; \Delta$ où $\Delta$ est une liste finie de formules atomiques.\index{regle@règle!existentielle simple}

Une \tdy est dite \textsl{existentielle} si ses axiomes sont tous des \ralgs ou existentielles simples (une \ralg peut d'ailleurs être considérée comme un cas particulier de \rex). Une théorie existentielle typique est la théorie des \textsl{anneaux de Bézout} (tout \itf est principal). Dans la littérature anglaise concernant la logique catégorique
(étudiée dans le cadre des \clama), une théorie existentielle est
appelée une \textsl{regular theory}.%
\index{theorie@théorie!existentielle}

\paragraph{Théories existentiellement rigides, cartésiennes}~

\smallskip \noindent 
Les théories \textsl{existentiellement rigides} sont les \tdys dans lesquelles les axiomes existentiels sont simples et correspondent à des existences uniques. 
Cela généralise (de très peu) les théories \dijs.\index{theorie@théorie!existentiellement rigide} 

\smallskip 
Une théorie existentielle existentiellement rigide
est dite \textsl{cartésienne}.\index{theorie@théorie!cartésienne}
Cela généralise (de très peu) les \talgs.

\paragraph{Théories rigides}~

\smallskip \noindent Une \tdy est dite \textsl{rigide} si elle est existentiellement rigide
et si les disjoncts du second membre dans les axiomes disjonctifs sont deux à deux incompatibles. Par exemple la théorie des \cdis est rigide, mais pas celle des \alos. La théorie des \cdis réels clos peut être énoncée de manière rigide, pas celle des \cdacs. Voir \cite{Johnstone79}, qui utilise la terminologie \textsl{disjunctive theory}.

\paragraph{Théories propositionnelles}~

\smallskip \noindent
La logique (classique ou intuitionniste) des propositions a un caractère très abstrait, qui peut sembler inutile du point de vue dynamique, car elle est déjà présente sous la forme de certaines des règles structurelles admissibles \ref{rstr1}. 
Elle s'avère cependant utile pour la \dfn des \trdis associés à des \sads (section \ref{subsectrdisad}).
 
 La logique des propositions peut être présentée sous une forme minimale dynamique comme suit, sans aucune sorte, ce qui implique que les constantes doivent être interprétées comme de pures valeurs de vérité abstraites (en logique classique elles hésitent seulement entre $\Vrai$ et $\Faux$).

Les constantes sont donc \(\Bot\), \(\Top\), et des \textsl{constantes propositionnelles} ou \textsl{propositions}. 
 Pour définir une telle théorie \((\cL,\cA)\), on donne un ensemble \(G\) de propositions\footnote{Cela fixe le langage $\cL$ via la signature $\so{G,\Top,\Bot}$}
 et un ensemble \(\cA\) d'axiomes qui sont des \rdijs sur le langage \(\cL\)). 

On a tout d'abord les axiomes \(\Bot\vd p\) et \(p \vd \Top\), et les axiomes qui gèrent l'\egt dans \(G\): 
\(p\vd q\) chaque fois \(p=_G q\). Si \(G\) n'est pas un ensemble discret, ces axiomes reflètent la structure de l'ensemble \(G\) dans la catégorie informelle \sa{Set}.

Les axiomes \suls donnés dans \(\cA \) sont du type \(p_1\vet\dots\vet p_n\vd q_1 \vou \dots \vou q_m\) où les \(p_i\) et \(q_j\) sont des constantes dans \(G\) (avec possiblement \(m=0\) ou $n=0$).

Deux constantes \(p\) et \(q\) sont dites \textsl{opposées} ou \textsl{\cops} si elles satisfont les axiomes de la négation

\DeuxRegles
{\labu \(\vd p\vou q\)}
{\labu \( p\vet q\vd \Bot\)}

Si la logique classique des propositions peut être interprétée comme donnée par des \tdys sans sorte du type précédemment décrit, il n'en va pas de même pour la logique intuitionniste. En effet, le connecteur $\Rightarrow$ ne peut pas être décrit en se limitant aux \tdys. On peut évidemment introduire ce connecteur dans le langage, mais la règle structurelle externe permettant d'introduire le connecteur d'implication en déduction naturelle n'est pas formulable comme une \rdy. 

%%%%%%%%%%%%%%%%%%%%%%%%%%%%%%%%%%%%%%%%%
%: Subsection{Structures algébriques dynamiques}
%\vspace{-1em}
\section{Structures algébriques dynamiques}\label{subsubsecSAD}

Références: \cite{CLR01}, \cite{Lom06}, \cite{LM2022}. 

\smallskip Les \sads sont explicitement nommées dans~\cite{Lom06}. Dans \cite{CLR01} elles sont implicites, mais explicitées sous la forme de leurs présentations.
Elles sont \egmt implicites dans \cite{Lom02}, et, last but not least, dans \cite[{D5}, 1985]{D5}, qui a été une source d'inspiration essentielle: on peut calculer de manière sûre dans la clôture algébrique d'un corps discret, même quand il n'est pas possible de construire cette clôture algébrique. 
Il suffit donc de considérer la clôture \agq comme une \sad \gui{à la D5} plutôt que comme une structure \agq usuelle: \textsl{l'évaluation paresseuse à la D5 fournit une sémantique \cov pour la clôture \agq d'un \cdi}.

\Subsection{Définitions, exemples} 

%: Example{exaAc}
\begin{example} \label{exaAc} 
Notre premier exemple est la \tpe des \acos (avec une seule sorte, nommée $\Ac$) dans laquelle l'essentiel du calcul est confié à une machinerie extérieure à la théorie formelle. 
Cette possibilité repose sur le fait que les \elts de l'anneau $\ZZxn$ peuvent être réduits à une forme normale prédéfinie.
Cela implique que l'égalité de deux termes est \eqve à l'identité de leurs formes normales. Du coup le prédicat d'\egt binaire peut être remplacé par le prédicat d'\egt~à~$0$. 

\smallskip La \textsl{théorie \SA{Ac0} des anneaux commutatifs} est écrite sur la signature suivante. Il y a une seule sorte, nommée $\Ac$. 
\Sigt{\Ac}{\cdot=_\Ac0\mathrel{;}\cdot+\cdot,\cdot\times \cdot,-\,\cdot,0_\Ac,1_\Ac} \label{NOTASigAc}
\noindent Les seuls axiomes sont les suivants (ce sont des \reds)\footnote{Les noms des règles sont calligraphiés comme suit: pour les règles directes, tout en minuscule, pour les autres règles \agqs (les règles de simplification), la première lettre en majuscule, et enfin les autres \rdys, tout en majuscule.}: 

\DeuxRegles{
\Lab{ac0} $\vd 0_\Ac =_\Ac0$
\Lab{ac2} $\,\,x=_\Ac0\Vdi{x,y:\Ac} x\times y=_\Ac0 $
}
{
\Lab{ac1} $\,\, x=_\Ac 0\vet y=_\Ac 0\Vdi{x,y:\Ac} x+y=_\Ac0$
}

\smallskip 
Le terme \gui{$x-y$} est
une abréviation de \gui{$x+(-y)$} et le prédicat binaire \gui{$\cdot=\cdot$} est \textsl{défini}
par la convention: \gui{$x=y$} est une abréviation pour \gui{$x-y=0$}.

\smallskip \rdb 
On considère souvent la théorie $\SA{Ac}$ des \textsl{\acos \und{non triviaux}}, qui est obtenue à partir de $\sa{Ac0}$ en ajoutent l'axiome d'effondremment

\Regles{\lAb{CL$_{\Ac}$} $\,\,1=_\Ac0\vd \Bot$\label{AxCLnqAc}}

\smallskip \noindent \textbf{Explications.}\label{Ac-comments}

\noindent \textsl{1.} Les règles qui définissent la théorie \sa{Ac0} des anneaux commutatifs doivent être comprises \prmt 
comme ceci. Tout terme de la théorie peut-être vu comme un \pol à \coes entiers en les variables présentes. 
On utilise alors la machinerie calculatoire des \pols commutatifs à \coes entiers (\gui{extérieure} à la théorie), qui réécrit tout terme (formé sur les constantes et les variables) comme un \pol à \coes entiers sous une forme normale prédéfinie.

\noindent La règle de distributivité $x(y+z)=_\Ac xy+xz$, par exemple, est alors confiée à un calcul automatique qui réduit à $0$ le terme
 $x(y+z)-(xy+xz)$.

\noindent De même la transitivité de l'\egt binaire est gérée par la règle \Tsbf{ac1} et par le calcul automatique qui réduit à $(x-z)$ le terme $ (x-y)+(y-z) $.

\smallskip \noindent \textsl{2.} On reconnait dans les trois règles \Tsbf{ac0}, \Tsbf{ac1} et \Tsbf{ac2} les axiomes des \ids, qui permettent de créer une structure d'anneau quotient, et qui signifient la compatibilité
de l'\egt avec l'addition et la multiplication. Dans la théorie \sa{Ac0} toute formule atomique est de la forme \gui{$t(\xn)=_\Ac0$} où les $x_i$ sont des variables et $t$ un terme du langage. Toute formule atomique 
est donc \imdt \eqve à une formule atomique dans laquelle $t$ est un \elt de l'anneau $\ZZ[\xn]$, écrit sous la forme normale convenue. La théorie \sa{Ac0} est donc la \gui{théorie des \idas}, au sens ancien de l'expression.
\Prmt, on vérifie facilement que la validité d'une \ralg simple du style

\Regles{\labu $\,\,p_1=_\Ac0\vet\dots\vet p_m=_\Ac0\Vdi{\xr:\Ac} q=_\Ac0$}

\noindent signifie exactement que le \pol $q\in \ZZxr$ est dans l'\id engendré par les \pols $p_1,\dots,p_m$ de $\ZZxr$. Cette \prt est plutôt difficile à décider\footnote{Cela résulte par exemple du \tho VIII-1.5 dans \cite{MRR}.}. 

\noindent La \tdy d'une \tpe n'apporte aucun outil supplémentaire à une \tpe elle-même. Il en va de même pour les théories algébriques.
Il n'y a donc rien de vraiment \gui{dynamique} dans les \tdys algébriques. Les \tdys vraiment intéressantes sont obtenues en ajoutant des axiomes dynamiques à des \talgs.

\noindent Pour la théorie \sa{Ac0} la validité de règles plus compliquées que celles envisagées ci-dessus est gérée par la règle structurelle \ref{10rstr1} \paref{10rstr1}.

\smallskip \noindent \textsl{3.} La théorie \sa{Ac0} telle qu'elle est présentée ne semble pas \gui{purement équationnelle} au premier abord car les axiomes ne sont pas de simples \egts entre termes. Cela tient à notre parti pris de remplacer l'\egt par le prédicat unaire \gui{$\cdot=0$} accompagné de la machinerie calculatoire externe des \pols à \coes entiers. Ce parti pris a l'avantage, selon nous, de montrer la véritable
\gui{structure logique} de la théorie en la ramenant à trois axiomes très simples et en confiant à un calcul automatique ce qui peut lui être confié, qui n'a pas grand chose à voir avec la logique proprement dite.
La même remarque s'appliquera par la suite à de nombreuses théories que nous qualifierons de purement équationnelles. \eoe
\end{example}
%--------- fin example ---------------------------------------- 

%:     Definition{defiSAD}
\begin{definition} \label{defiSAD}
Si $\sa{T}=(\cL,\cA)$ est une \tdy, une \textsl{\sad de type~\sa{T}} est une extension simple de la théorie \sa T donnée par un 
ensemble $G$ de \textsl{générateurs} et un ensemble $R$ de \textsl{relations}. 
Une \gui{relation} est par \dfn une formule atomique $P(\und{t})$ construite sur le langage $\cL\cup G$ avec des termes $t_i$ clos dans ce langage. À une telle relation est associé l'axiome \gui{$\vdi P(\und{t})$} de la \sad. Cette \sad est donc la \tdy $(\cL\cup G,\cA\cup R)$, que l'on note aussi $\big((G,R),\sa{T}\big)$. 
\end{definition}

%: Example{exaSaCd}
\begin{example} \label{exaSaCd} 
 Par exemple on obtient une \sad de corps discret 

\snic{\gK=\big((G,R),\Sa{Cd}\big)}

\noindent en prenant~\hbox{$G=\so{a,b}$} et 
$R=\so{105=0,\,a^2+b^2-1=0}.$ Ce \cdi dynamique correspond à n'importe quel \cdi de caractéristique
$3$ ou $5$ ou $7$ engendré par deux \elts $\alpha$ et
$\beta$ vérifiant~\hbox{$\alpha^2+\beta^2=1$}. 

\noindent Outre les règles dynamiques valables dans tous les corps discrets, il y a maintenant celles que l'on obtient en élargissant le
langage avec les constantes prises dans $G$ et en ajoutant aux axiomes les 
relations prises dans $R$. Par exemple est valide la \rdij

\Regles{\labu $3=0 \vou 5=0  \vou 7=0$}

\noindent ainsi que la \rdy

\Regles{\labu $\Exists z \;15z=1 \vou \Exists z \; 21z=1  \vou \Exists z \;35z=1$ \eoe}

\end{example}
%--------- fin example ---------------------------------------- 

%d
%: Definition{defiFait}
\begin{definotas} \label{defiFait} \label{notasadreglevalide}
Soit $\gS=\big((G,R),\sa{T}\big)$ une \sad de type $\sa{T}=(\cL,\cA)$.
\begin{itemize}
\item Nous indiquerons que la règle \gui{$\,\Gamma\vd \dots$} est valide dans la \sad $\gS$ sous la forme abrégée suivante: \gui{$\, \Gamma\vdi_{\gS} \dots$}. On pourrait aussi utiliser la notation \gui{$\,R\vet \Gamma\vdi_{\sA{T}} \dots$}, qui signifie que la \demo peut utiliser une liste finie d'axiomes extraite de $R$.
\item L'ensemble des termes clos de $\gS$, \cad les termes construits sur $\cL\cup G$, se note $\Tcl(\gS)$. L'ensemble des formules atomiques closes se note~$\Atcl(\gS)$. 

\item Une \ralg~$\vd P$ pour~$P\in\Atcl(\gS)$ s'appelle un fait de $\gS$. L'ensemble des faits de $\gS$ valides dans $\gS$ se note $\Atclv(\gS)$.
Un fait concerne uniquement des objets définissables syntaxiquement dans la structure.
Il est clair $\gS$ prouve exactement les mêmes \rdys que la \sad $\wi\gS=\big((\Tcl(\gS),\Atclv(\gS)),\sa{T}\big)$.
\end{itemize}
%Soit $\gS=\big((G,R),\sa{T}\big)$ une \sad de type $\sa{T}$. Une \ralg sans hypothèse
%et sans variable s'appelle
%un \textsl{fait} (dans $\gS$). 
\end{definotas}
%----------- fin definition -------------------------------- 

 L'algèbre \gui{concrète} consiste très souvent à prouver des faits ou des 
\rdys dans des \sads particulières.
C'est un peu plus général que la théorie (inépuisable) des identités 
algébriques, \cad l'algèbre universelle,
à l'œuvre derrière une forte proportion des grands \thos d'algèbre abstraite.

\smallskip 
Dans le cas d'une \talg \sa{T}, une \sad de type \sa{T} donne une structure \agq usuelle, définie par \gtrs et relations, satisfaisant les \ralgs requises. 

\smallskip La méthode dynamique est souvent un moyen pratique de
construire des identités algébriques (des \gui{Positivstellens\"atze} par exemple), en suivant au plus près les pistes indiquées dans les preuves données en \clama. 
 
\smallskip Dans une \sad un fait $P(\und t)$ est \textsl{absolument vrai} s'il est prouvable 
(\cad si la règle \gui{$\vdi P(\und t)$} est valide). Il est \textsl{absolument faux}, ou plus 
justement \textsl{catastrophique} si \gui{$P(\und t)\vd \Bot$} est valide. 
Intermédiaires entre ces deux cas existent de nombreuses possibilités: une \sad n'a pas un modèle figé unique, mais représente à l'état potentiel toutes les réalisations éventuelles idéales de la structure
(cette notion reste volontairement floue).
Ajouter un fait catastrophique comme axiome revient à supprimer tous 
les modèles\footnote{Dans la variante où le collapsus réduit tout modèle à un singleton: \dots\ revient à n'autoriser que le modèle trivial.}.

%e
%: Example{exasdz}
\begin{example} \label{exasdz} 
Nous considérons une \pn $(G,R)$ dans le langage de la théorie des anneaux commutatifs \Sa{Ac}. Soit \sa T une \tdy qui étend la théorie \Sa{Ac} sans étendre le langage, 
par exemple la théorie \sa{Cd} des corps discrets. Tout terme clos de la \sad $\big((G,R),\sa T\big)$ se réécrit sous forme
d'un \pol $f(\ux)\in\ZG$ à \coes entiers en les \gui{constantes} $x_i\in G$ de la \sad.
Les \elts de $R$ sont des relations $f(\ux)=0$, de sorte que par un léger abus de langage, on peut considérer $R$ comme un ensemble d'\elts de $\ZG$. 

\noindent On est donc en train d'étudier l'anneau $\gA=\aqo\ZG R$, ou plus exactement ce que devient cet anneau lorsqu'on lui demande de satisfaire certains nouveaux axiomes. Nous noterons~$\sa T(\gA)$ la \sad
$\big((G,R),\sa T\big)$.

\noindent Dans beaucoup d'exemples, la théorie s'effondre \ssi $\gA$ est trivial. Par exemple la \sad $\Sa{Cd}(\gA) $ s'effondre \ssi $1=_\gA0$. En \clama on dit: en effet un anneau non trivial possède un \idep $\fp$, et le corps des fractions de l'anneau intègre~$\gA/\fp$ est un modèle non trivial de $\Sa{Cd}(\gA) $.
Plus simplement, sans utiliser la théorie des modèles ni l'axiome de l'\idep, on transforme une preuve de $1=0$ dans~$\Sa{Cd}(\gA) $
en une preuve \hbox{de $1=0$} dans~$\Sa{Ac}(\gA) $ (\demo analogue à celle de \cite[Theorem~2.4]{CLR01}).

\noindent Pour ce qui concerne les faits $\theta=0$ \footnote{Avec $\theta=t(\uxi)\in \gA$, où $t\in\ZZ[G]$ et les $\xi_k$ sont les $x_k$ vu dans le quotient $\gA$ de $\ZZ[G]$} valides dans la théorie $\sa T(\gA)$, la situation est un peu plus compliquée.\\
La théorie des anneaux locaux \Sa{Al} prouve $\theta=0$ exactement lorsque $\theta=_\gA0$, d'où la très grande importance des anneaux locaux en \alg commutative.
\\
La théorie \Sa{Cd} prouve $\theta=0$ exactement lorsque $\theta \in\!\sqrt[\gA]0$. cela correspond au quotient réduit de~$\gA$: si $\theta =0$
dans $\Sa{Cd}(\gA)$ alors $\theta$ est nilpotent dans $\gA$. 
C'est une forme abstraite (relativement faible) du \nst. La \demo est \elr. Naturellement, si deux théories prouvent les mêmes faits, elles peuvent se différencier au niveau des \rdys plus \gnles que les \ralgs. \eoe
\end{example}
%--------- fin example ---------------------------------------- 

\Subsection{Diagramme positif d'une structure algébrique.} 

%: Definition{defidiagramme}
\begin{definota} \label{defidiagramme}~
\begin{enumerate}
\item Soit $\sa T_1=(\cL_1,\cA_1)$ une \tdy et $\gA$ une structure \agq sur un langage $\cL\subseteq \cL_1$. On appelle \textsl{diagramme positif de $\gA$ pour le langage $\cL$}, une \pn $(G,R)$ de~$\gA$ comme structure \agq sur le langage $\cL$. Un tel diagramme est noté~$\Diag(\gA,\cL)$. \\
On note alors $\sa T_1(\gA)$ la \sad $\big(\Diag(\gA,\cL),\sab{T}_{\!1}\big)$.%
\index{diagramme positif!d'une structure algébrique sur un langage donné} 
\item Soient $\sa T=(\cL,\cA)$ une \tdy, $\gB$ un modèle de $\sa T$, et $\sab{T}_{\!1}=(\cL_1,\cA_1)$ une extension simple de $\sa{T}$. Considérons une \pn$(G,R)$ de $\gB$ comme \sad de type \sa{T}. Un tel diagramme est appelé \textsl{diagramme positif de $\gB$ pour la \tdy $\sa T$} et il est noté~$\Diag(\gB,\sa T)$. On note alors $\sab{T}_{\!1}(\gB)$ la \sad $\big(\Diag(\gB,\sa T),\sab{T}_{\!1}\big)$.%
\index{diagramme positif!d'un modèle d'une \tdy}
\end{enumerate}
\end{definota}
Le point 1 peut être vu comme un cas particulier du point 2, lorsque $\cA=\emptyset$.

%r
%: Remark{remdefidiagramme}
\begin{remarks} \label{remdefidiagramme}~

\noindent 1) Voici un exemple typique pour le point 1. La théorie $\sa T_1$ est une extension simple de la théorie~\sa{Ac} des anneaux commutatifs avec l'\egt pour seul prédicat. Notons $\cL$ le langage de~$\sa{Ac}$ et soit~$\gA$ un anneau commutatif. 
Pour \gtrs de $\Diag(\gA,\cL)$ on peut prendre les \elts de l'ensemble sous-jacent à $\gA$, et pour relations on 
peut se limiter aux \egts $0_\gA=0$, $1_\gA=1$, $-a=b$, $a+b=c$ et $ab=c$ lorsqu'elles sont satisfaites entre \elts de $\gA$. Ce diagramme positif ne contient aucune inégalité $a\neq b$ pour la simple raison qu'elles ne font pas partie du langage de \sa{Ac}. C'est pour cette raison que nous parlons de diagramme \gui{positif}.

\smallskip \noindent 2) 
Un \elt $a$ de $\gA$ n'a pas toujours de représentant canonique
dans un ensemble à la Bishop, même si l'ensemble est discret. Dans un tel cas, pour revenir à la \dfn de l'ensemble sous-jacent à $\gA$ selon Bishop, on peut prendre une constante $x_b$ différente pour chaque représentant $b$ 
de l'\elt~$a$. On trouve alors dans le diagramme positif de $\gA$ une relation $x_b=x_c$
chaque fois \hbox{que $b=_\gA c$}.
\eoe 
\end{remarks}
%----------- fin remark ---------------------------------- 

%: Subsubsection{Modèles d'une \sad}
\Subsection{Modèles constructifs versus modèles classiques}\label{subsecmodelescofsSAD}
%-----------
On considère une \sad $\gA=\big((G,R),\sa T\big)$ de type \sa{T} avec une ou plusieurs sortes. Pour simplifier les notations nous supposons une seule sorte. Un \textsl{modèle de~$\gA$} est une \salg usuelle (statique) $M$ décrite dans le langage 
associé à~$\gA$ et vérifiant les axiomes de $\gA$ (ceux de $\sa{T}$ et ceux donnés par la \pn de $\gA$).

Lorsque $\gA$ est défini par la \pn vide, on parle de \textsl{modèles de \sa T}.

\smallskip 
La notion de modèle est donc basée à priori sur une notion intuitive de \textsl{structure \agq} à la Bourbaki. 
Nous pouvons qualifier ces \salgs de \gui{statiques} par contraste avec les \sads \gnles. Notons qu'ici l'ensemble \gui{sous-jacent} à la structure est un ensemble \gui{naïf} (ou plusieurs ensembles naïfs s'il y a plusieurs sortes) structuré par la donnée de prédicats et de fonctions (au sens naïf)
soumis à certains axiomes.

D'un point de vue \cof, les modèles doivent satisfaire les axiomes en respectant le sens intuitif du \gui{ou} et du \gui{il existe}: pour prouver qu'une \salg particulière satisfait les axiomes, on autorise uniquement la logique intuitionniste. 
Notons aussi que la théorie des ensembles à laquelle nous nous référons est à priori celle, informelle, de Bishop.

%%%%%%%%%%%%%%%%%%%%%%%%%%%%%%%%%%%%%%%%%%%%%%%%%%%%%%%%%%%%%%%%%%%% 
%: Subsection{Morphismes entre \sads de même type}
\Subsection{Morphismes entre \sads de même type}
%-----------

Considérons une \tdij \sa{T}.
Dans ce cas, une notion naturelle possible de morphisme d'une 
\sad $\gA=\big((G, R), \sa{T}\big)$ vers une \sad \hbox{$\gA'=\big((G', R'), \sa{T}\big)$}
pour la \tdij \sa{T} est la suivante.

 Un \elt de \(\Sad_{\sA T}(\gA,\gA')\) est donné par une fonction \(\varphi\colon G\to \Tcl(\gA')\) qui interprète les \elts de \(G\) par des termes clos de \(\gA'\). Cette fonction s'étend de manière unique en une fonction \(\Tcl(\gA)\to \Tcl(\gA')\) en respectant la construction des termes au moyen des symboles de fonctions présents dans \(\cL\). En outre les \elts de $R$ doivent donner des faits valides dans $\gA'$ selon cette interprétation.
 
L'\egt entre deux \elts \(\varphi\) et \(\psi\) de l'ensemble \(M=\Sad_{\sA T}(\gA,\gA')\) est définie comme suit: on~a~\(\varphi=_M\psi \) \ssi pour tout \(x\in G\), l'\egt 
\(\varphi(x)=\psi(x) \) est valide dans \(\gA'\).
 
La composition des morphismes est définie de manière naturelle pour trois \sads \(\gA\), \(\gB\) et \(\gC\). 

\smallskip On obtient ainsi une catégorie (informelle) très intéressante.
Les objets sont les \sads de type \sa{T}. L'ensemble des flèches de \(\gA\) vers \(\gA'\) est \(\Sad_{\sA T}(\gA,\gA')\). 

Cette catégorie possède des limites et colimites arbitraires, construites très naïvement au niveau des \pns $(G,R)$, basées sur la théorie intuitive naïve des ensembles que l'on considère dans le monde \mathe ambiant. 

Par exemple le
produit dans la catégorie \(\Sad_{\sA T}\) de $\gA=\big((G, R), \sa{T}\big)$ et \hbox{$\gA'=\big((G', R'), \sa{T}\big)$} est la \sad de type \sa{T} dont la \pn est donnée par \((G\times G',R\times R')\). Lorsque~\sa{T} est la théorie des corps discrets, dans la version disjonctive où un \cdi est défini comme un anneau \zedr connexe (il faut introduire un symbole de fonction pour le quasi inverse). Le produit précédent est une nouvelle \sad de \cdi, et \(\gA\times \gA'\) est muni d'une structure \agq usuelle d'anneau \zedr.
Cette situation semble analogue à celle des faisceaux de \cdis (selon la sémantique de Kripke-Joyal%, voir la \dfn~\ref{contextBKJforSheaves}
), qui ne sont des \cdis que dans les fibres.

\smallskip \rem Si on a affaire à une \telri, on peut se ramener au cas d'une \tdij par skolémisation des \rexris. 
Il semble cependant que dans le cas d'une \tdy avec des axiomes existentiels non rigides, les choses ne soient pas très claires. \eoe

\smallskip 
Il arrive que l'on soit intéressé par une notion plus restrictive de morphisme
entre deux \sads $\gA$ et $\gA'$ de même type \sa{T}, par exemple la notion de morphisme local entre anneaux commutatifs, adaptée à un contexte précis.
 Dans un tel cas, on souhaite que les nouveaux morphismes de $\gA$ vers $\gA'$ puissent être traités comme donnés par des \sads pour une certaine \tdy définie à partir de \sa{T}, $\gA$ et $\gA'$ (comme cela peut être le cas pour les morphismes locaux).
 
\Subsection{Exemples}

\noindent 1) La \tdij \SA{Asdz0} (resp. \SA{Asdz}) des \textsl{anneaux sans diviseurs de zéro} est obtenue 
à partir de la théorie \Sa{Ac0} (resp. \Sa{Ac}) en ajoutant la \rdy%
\index{anneau!sans diviseur de zéro}%

\Regles {\laB{ASDZ} $\,\,xy=0\vd x=0 \vou y=0$}

\smallskip \noindent 2) Référence \cite[section VIII-3]{ACMC}. 
La \tel \SA{Alsdz0} (resp. \SA{Alsdz}) des \textsl{aneaux \lsdz} est obtenue en ajoutant à la théorie \Sa{Ac0} (resp. \Sa{Ac}) l'axiome \tsbf{LSDZ}:%
\index{anneau!localement sans diviseur de zéro}

\Regles {\Lab{LSDZ} $\,\,xy=0\vd \Exists u,v \;(ux=0\vet vy=0\vet u+v=1)$}

\smallskip \noindent 3) Théorie disjonctive \SA{Ai} des \textsl{anneaux intègres}. Avec la signature
\Sigt{\Ai}{\cdot=0,\cdot\neq0\mathrel{;}\cdot+\cdot,\cdot\times \cdot,-\,\cdot,0,1} \label{NOTASigAi}
\noindent la théorie des \textsl{anneaux intègres} est obtenue 
à partir de la théorie \Sa{Ac0} en ajoutant comme axiomes les \rdys suivantes%
\index{anneau!integre@intègre} 

\DeuxRegles 
{
\labu $\,\,x\neq 0\vet y= 0\vd x+y\neq0$
\labu $\,\,x\neq 0\vet y\neq 0 \vd xy\neq0$
\labu $\,\,x\neq 0\vet xy=0 \vd y=0$
\lAb{ED\inq} $\vd x=0\vou x\neq0$ \label{AxEdnq}
}
{
\labu $\vd 1\neq0$
\labu $\,\,xy\neq 0 \vd x\neq0$
\lab{} 
\lAb{col\inq} $\,\,0\neq 0 \vd 1=0$ \label{Axcolnq}
}

On notera la différence sensible entre les théories \sa{Asdz} et \sa{Ai}, ce qui correspond à une distinction importante en \coma mais invisible en \clama.

\smallskip \noindent 4) Théorie \SA{Al1} des anneaux locaux avec unités, basée sur la signature 
\Sigt{\Alu}{\,\cdot=0,\U(\cdot)\mathrel{;}\cdot+\cdot,\cdot\times \cdot,-\,\cdot,0,1\,}\label{NOTASigAlu}
\noindent Cette théorie est une extension de la théorie \Sa{Ac} des anneaux commutatifs.
Un prédicat~$\U(x)$ est introduit comme le prédicat d'inversibilité au moyen des deux axiomes convenables.

\DeuxRegles 
{\Lab{Uv} $\,\, xy=1\Vd \U(x)$
}
{\Lab{UV} $\,\, \U(x) \Vd \Exists y\;xy=1$
}

\noindent On~ajoute %l'axiome de collapsus \tsbf{cl$_{Al1}$} et
l'axiome \tsbf{AL1} des anneaux locaux proprement dit. 

\Regles 
{\Lab{AL1} $\,\, \U(x+y) \Vd \U(x) \vou \U(y)$}

\noindent La règle valide $\,\,\U(0)\Vd 1=0$ est vue comme une règle de collapsus.

\smallskip \noindent 5) Théorie \SA{Alrd} des \textsl{\alrds}. 

\noindent Un anneau local est dit \textsl{résiduellement discret} lorsque l'anneau résiduel est un corps discret. 
La notion d'anneau résiduel (le quotient de l'anneau par son radical de Jacobson) n'est pas évidente à introduire dans le cadre des théories dynamiques. Néanmoins pour les anneaux locaux résiduellement discrets, cela se passe convenablement.
La théorie dynamique des {anneaux locaux résiduellement discrets} est obtenue à partir de la théorie \Sa{Al1} en ajoutant un prédicat~$\Rn$ (pour les éléments résiduellement nuls) comme prédicat opposé au prédicat d'inversibilité au moyen des axiomes\index{anneau!local!résiduellement discret}

\DeuxRegles{
\Lab{Alrd} $\,\, \U(x)\vet \Rn(x)\vd 1=0$
}
{
\Lab{ALRD} $\vd \U(x) \vou \Rn(x)$
}

\noindent  
On a donc la signature suivante.
\Sigt{\Alrd}{\cdot=0,\U(\cdot),\Rn(\cdot)\mathrel{;}\cdot+\cdot,\cdot\times \cdot,-\,\cdot,0,1}
\label{NOTASigAlrd}

\noindent
En mathématiques classiques tout anneau local non trivial est résiduellement discret. Il est significatif que la différence entre les deux notions (présente en mathématiques constructives) soit naturelle au niveau des théories dynamiques.
 
\noindent 
La théorie \Sa{Alrd} peut aussi être décrite à partir de \Sa{Ac} en ajoutant  comme axiomes les règles directes  suivantes  

\vspace{-.2em}
\DeuxRegles{
\Lab{alrd0} $\,\,x=0\vd \Rn(x)$
\Lab{alrd2} $\,\,\Rn(x) \vd \Rn(xy)$
\Lab{alrd4} $\,\,\Rn(x)\vet \U(y) \vd \U(x+y)$
}
{
\Lab{alrd1} $\,\,\Rn(x)\vet \Rn(y) \vd \Rn(x+y)$
\Lab{alrd3} $\vd \U(1)$
\Lab{alrd5} $\,\,\U(x)\vet \U(y) \vd \U(xy)$
}

\noindent puis la \rsim $\,\,\U(xy) \vd \U(x)$ et les \rdys \Tsbf{IV}, \Tsbf{AL1},  \Tsbf{ALRD} (voir \cite{Lom1997}).
\eoe

%:Subsection{L'arithmétique primitive récursive}
\Subsection{L'arithmétique primitive récursive}
\label{PrimRec}
Cette section montre l'intérêt d'utiliser des sortes pour des fonctions
définissables \cot dans une structure \agq (par exemple ici le semi-anneau des entiers naturels) lorsque le langage où est défini la structure ne permet pas d'introduire des symboles fonctionnels correspondant à ces fonctions dans la \tgm correspondante.

Cet exemple nous motive pour introduire des sortes pour certaines fonctions \sagcs et leurs modules de continuité dans la théorie des \crc non discrets. 

\smallskip Bien avant la mise en place de la machinerie des \tdys, R.~L.~Goodstein a expliqué comment traiter une grande partie des \maths usuellement pratiquées au moyen de systèmes formels purement calculatoires, \textsl{sans logique}.

 Dans le livre \cite{Goo1957} l'auteur, suivant une suggestion de Skolem, montre comment un système de calcul \gui{sans logique, et sans quantificateurs} permet de développer une partie très importante de \gui{l'arithmétique}, comprise au sens d'une théorie formelle des entiers naturels.
 
Le seul \pb, mais c'est un \pb d'importance qui a de quoi rebuter plus
\stMF d'un\e mathématic\ien, c'est que les énoncés usuels des \maths doivent être codés sous forme de fonctions primitives récursives. Ceci peut sembler nous ramener au royaume de l'arithmétique du second ordre et des Reverse Mathematics. 

Dans un second ouvrage (\cite{Goo1961}) Goodstein étend son étude à l'analyse récursive. On recommande aussi vivement l'ouvrage \cite{Goo1979}.

%%%%%%%%%%%%%%%%%%%%%%%%%%%%%%%%%%%%%%%%%%%%%%%%%%%%%%%%%%%%%%%%%%%%
%:Le système formel de Goodstein
\Subsubsection {Le système formel de Goodstein}

Si on se limite à l'arithmétique primitive récursive\footnote{Le livre de Goodstein étudie des systèmes de calcul plus larges qui incluent les fonctions définies par récurrences multiples, comme la fonction d'Ackerman.}, on peut décrire le système formel proposé par Goodstein comme suit.

Comme dans la théorie formelle \sa{Peano}, les variables et les constantes représentent des entiers naturels. Il y a une seule constante, $0$, et un seul symbole de relation qui est l'égalité $x=y$. 

Pour toute fonction primitive récursive $f\colon \NN^{k+1}\to\NN$, définie par récurrence simple au moyen des équations
\[ 
\begin{array}{rcl} 
f(\xk,0) & = & g(\xk) \\[.3em] 
f(\xk,y+1) & = & h(\xk,y,f(\xk,y)) \\[.3em] 
 \end{array}
\] 
où $g$ et $h$ sont des fonctions primitives récursives préalablement définies,
on introduit un symbole fonctionnel correspondant à cette \dfn de $f$. 

De même, pour toute fonction primitive récursive définie par composition de fonctions primitives récursives préalablement définies, on introduit un symbole fonctionnel correspondant à cette \dfn.

Les symboles de fonctions qui \gui{initialisent} le système sont $\rS$ pour la fonction successeur, $0_1$ pour la fonction nulle à une variable \hbox{et $\pi_{n,k}$} pour la $k$-ème fonction coordonnée $\NN^n\to\NN$ ($k\in\lrbn$, $n\in\NN$).

On obtient de cette manière un symbole fonctionnel d'arité $r$ pour chaque \dfn de fonction primitive récursive $\NN^r\to\NN$.

Notons que l'on a un nom formel noté $\und n$, abréviation d'un $S(S(\dots(S(0))\dots))$,
pour chaque entier $n$.

Les calculs dans l'arithmétique primitive récursive à la Goodsein consistent à établir des \gui{\idts} entre deux fonctions ainsi définies correspondant à deux symboles fonctionnels\footnote{Par exemple en utilisant une numérotation convenable de tous les symboles fonctionnels créés.} $f_i$ et $f_j$ de même arité
$$ 
\forall \xk\;\; f_i(\xk)=f_j(\xk)
$$
ce que nous pouvons écrire sous forme d'une règle valide dans le système proposé: 

\Regles{ \lab{Eq$_{i,j}$} $\vd f_i(\xk)=f_j(\xk)$}

Il ne s'agit cependant pas d'une \tdy car les règles valides ne résultent pas d'un simple système axiomatique de \rdys. 

En effet, on doit bien évidemment prendre comme axiomes toutes les \egts évoquées préalablement qui servent à définir des fonctions primitives récursives arbitraires, mais cela n'est manifestement pas suffisant. 

En effet, si l'\egt $(\und m+\und n)+\und p=\und m+(\und n+\und p)$ peut être établie pour des entiers $m,n,p$ arbitraires par simple utilisation de la définition de $+$ par \recu, la règle correspondante

\Regles {\lab{~} $\vd (x+y)+z=x+(y+z)$}

\noindent ne résulte pas de manière purement finitaire des axiomes de définition de $+$.

La même remarque vaut par exemple pour une \egt $f_i(\und m,\und n)=f_j(\und m,\und n)$ qui pourrait être constatée pour tous entiers $m,n$ par simple application des axiomes définissant $f_i$ et $f_j$, alors que la règle correspondante \tsbf{Eq$_{i,j}$}
ne pourra en \gnl pas être établie de manière purement finitaire si les définitions de $f_i$ et $f_j$ utilisent le schéma de récurrence simple.

Le système de calcul permet de gérer la composition des fonctions\footnote{Par exemple, $g=f_1\circ (f_2\circ f_3)$ et $h=(f_1\circ f_2)\circ f_3$ qui correspondent à deux \dfns différentes, donnent bien lieu à l'identité $\vd g(x)=h(x)$ car d'après les \dfns $\vd g(x)=f_1(f_2(f_3(x)))$ et $\vd h(x)=f_1(f_2(f_3(x)))$.}, car il est soumis aux axiomes de l'\egt. 
Mais on a besoin pour compléter l'arithmétique primitive récursive des \textsl{règles externes} correspondant aux \dfns par récurrence. \Prmt, par exemple, à partir des deux règles valides suivantes 

\Regles {\lab{~} $\vd u(x,0)=v(x,0)$
\lab{~} $\,\,u(x,y)=v(x,y)\vd u(x,y+1)=v(x,y+1)$}

\noindent on déduit la validité de la règle

\Regles {\lab{~} $\vd u(x,y)=v(x,y)$} 
 
\noindent (ici $u$ et $v$ sont deux termes contenant $x$ et $y$ comme variables libres)

L'utilisation de ces règles externes permet d'éviter le recours à l'axiome correspondant qui peut être formulé au premier ordre, mais pas dans une \tdy finitaire:
$$
\big[\forall x\;u(x,0)=v(x,0)\,\vii\, \forall x,y\;\big((u(x,y)=v(x,y)\Rightarrow u(x,y+1)=u(x,y+1)\big)\big] \Rightarrow \forall x,y\;u(x,y)=v(x,y)
$$

%:subsubsection{PRA tgm finitaire} 
\Subsubsection{Une \tgm finitaire pour l'arithmétique primitive récursive} \label{secPRA}

Nous montrons dans cette sous-section comment traiter l''arithmétique primitive récursive à la Goodstein dans le cadre d'une \tdy.\index{arithmétique primitive récursive!à la Goodstein}

Nous expliquons maintenant comment l'utilisation de sortes pour les fonctions nous permet d'éviter le recours à ces règles externes et de définir un \sys formel simple (une \talg finitaire) pour l'arithmétique primitive récursive. 

\hum{Cela ressemble beaucoup à du $\lambda$-calcul. Si c'est bien le cas, il doit y avoir un truc tout à fait analogue pour traiter l'arithmétique primitive récursive en $\lambda$-calcul. Il faudrait donner des références.}

\begin{enumerate}
\item \textsl{Les sortes.} \\
Pour chaque entier $k$ on introduit la sorte $\iF_k$ des fonctions primitives récursives $f\colon \NN^k\to\NN$. La sorte $\iF_0$ est la sorte des entiers notée $\iN$.
\item \textsl{Les prédicats.} \\
Pour chaque sorte $\iF_k$, on a un symbole d'\egt correspondant $\cdot=_k\cdot$.
\item \textsl{Les constantes.}
\begin{enumerate}
\item 
Les constantes de base sont (des noms pour) 
\begin{itemize}
\item $0$ de sorte $\iN$,
\item $0_k$ de sorte $\iF_k$ (pour la fonction constante nulle, $k\geq 1$),
\item la fonction successeur $S$ de sorte $F_1$, 
\item pour $n\geq 1$, les $n$ fonctions coordonnées\footnote{Notons que $\pi_{1,1}$ est la constante qui désigne l'identité dans $\iF_1$} $\pi_{n,k}$ de sorte $\iF_n$ ($1\leq k\leq n$), 
\end{itemize}

\item Pour toute fonction primitive récursive $f\colon \NN^{k}\to\NN$ ($k\geq 1$), définie par récurrence simple au moyen des équations
\[ 
\begin{array}{rcll} 
f(\ux,0) & = & g(\ux) &g\colon \NN^{k-1}\to\NN\\[.3em] 
f(\ux,S(y)) & = & h(\ux,y,f(\ux,y))\qquad &h\colon \NN^{k+1}\to\NN %\\[.3em] 
 \end{array}
\] 
où $g$ et $h$ sont des fonctions primitives récursives préalablement définies,
on introduit un nom pour $f$ comme constante de sorte $\iF_{k}$\footnote{Pour $k=2$ par exemple cela pourrait être le nom $R2(G,H)$ si $G$ et $H$ sont des noms pour $g$ et $h$.}.

\item On introduit \egmt un nom pour toute fonction primitive récursive définie par composition de fonctions primitives récursives préalablement définies.
\end{enumerate}
\item \textsl{Les autres symboles de fonction.}\\
Les symboles de fonction donnés en 4b et 4c permettent d'éviter, si on le désire, de créer les constantes prévues en 3b et 3c.
\begin{enumerate}
\item \textsl{Les évaluations}. Pour chaque $\ell\geq 1$ on a symbole fonctionnel pour l'évaluation $\Ev_{\ell}$ de la constante $f\in F_\ell$ en un $k$-uplet d'entiers.
C'est un symbole du type $\iF_\ell\times \iN^\ell\to \iN$.
On abrège $\Ev_{\ell}(f,x_1,\dots,x_\ell)$ en $f(x_1,\dots,x_\ell)$. 
\item \textsl{La récurrence simple}. Pour $k\geq 1$ on a un symbole fonctionnel $\rR_{k}$
pour l'\elt $f\in \iF_k$ défini \gui{par \recu simple} à partir d'un \elt $g\in \iF_{k-1}$ et d'un \elt~\hbox{$h\in \iF_{k+1}$} (comme en 3b, mais ici, $f$, $g$, $h$ sont des variables). C'est un symbole du type~\hbox{$\iF_{k-1}\times \iF_{k+1}\to \iF_{k}$}.

\item \textsl{Les compositions}. Pour $k\geq 1$ et $\ell\geq 1$ on a un symbole fonctionnel $\rC_{\ell,k}$ pour la \gui{composition} de l'\elt $f\in \iF_\ell$ avec $\ell$ \elts $g_i\in \iF_k$. C'est un symbole du type $\iF_\ell\times \iF_k^{\,\ell}\to \iF_k$. On abrège $\rC_{\ell,k}(f,g_1,\dots,g_\ell)$ en $f\circ(g_1,\dots,g_\ell)$ ou $f(g_1,\dots,g_\ell)$.
On peut voir $\Ev_{\ell}$ comme le cas particulier $\rC_{\ell,0}$.
\end{enumerate}
\item \textsl{Les axiomes} sont les suivants.
\begin{enumerate}
\item Les axiomes d'\egt usuels pour les relations $=_k$.
\item Un axiome d'effondrement $\,\,S(0)=0\vd \Bot$. On note $1$ plutôt que $\und 1$ pour $S(0)$.
\item Les axiomes qui établissent les \egts liant les constantes et les autres symboles fonctionnels. Par exemple: 

\Regles{
\lab{~} $\Vdi{f:\iF_k} f\circ (\pi_{k,1},\dots,\pi_{k,k})=f$
\lab{~} $\Vdi{f_1,\dots,f_k:\iF_\ell} 0_k(f_1,\dots,f_k)=0_\ell$
\lab{~} $\Vdi{f_1,\dots,f_k:\iF_\ell} \pi_{k,i}(f_1,\dots,f_k)=f_i$
\lab{prod} $\Vdi{x,y:N} \Prod=\rR_2(0_1,\mathrm{\Sum}\circ (\pi_{3,3},\pi_{3,1}))$
}

La dernière règle donne la \dfn par \recu du produit (la constante $\Prod$ 
dans~$\iF_2$) à partir de l'addition (la constante $\mathrm{\Sum}$ dans $F_2$).

\hum{1. Il semble que la troisième règle résulte de la première et des axiomes pour la composition.\\
2. Donner aussi un exemple pour une \gui{fonction} définie par composition?}
\item Les axiomes pour l'associativité des compositions, y compris les cas des évaluations. \\
Par exemple:

\Regles{
\lab{asC$_{1,1,1}$} $\Vdi{f,g,h:\iF_1} f\circ (g\circ h)=(f\circ g)\circ h$ 
\lab{asC$_{1,1,0}$} $\Vdi{f,g:\iF_1;x:N} (f\circ g)(x)=f(g(x))$ 
\lab{asC$_{2,1,1}$} $\Vdi{f:\iF_2;g_1,g_2,h_1,h_2:\iF_1} f\circ (g_1\circ h_1,g_2\circ h_2)=(f\circ (g_1,g_2))\circ (h_1,h_2)$ 
}
\item Les axiomes pour les \dfns par \recu. Par exemple pour $\rR_{2}$:

\Regles{
\lab{Rec$_{2,\mathrm{ini}}$} $\,\,f=\rR_2(g,h)\Vdi{f:\iF_2;g:\iF_1;h:\iF_3} f\circ (\pi_{1,1},0_1)=g$
\lab{Rec$_{2,\mathrm{rec}}$}$\,\,f=\rR_2(g,h)\Vdi{f:F_2;g:F_1;h:F_3} f\circ (\pi_{2,1},S\circ \pi_{2,2})=h\circ (\pi_{2,1},\pi_{2,2},f)$
}
\item Les axiomes pour les \demos par \recu (un pour chaque arité). \\
Donnons par exemple $\tsbf{REC}_2$:

\Regles{
\lab{REC$_2$}
$ \left.
\begin{array}{rcll} 
f_1\circ (\pi_{1,1},0_1)&=&f_2\circ (\pi_{1,1},0_1)\vet \\[.3em] 
f_1\circ (\pi_{2,1},S\circ \pi_{2,2}) & = & h\circ (\pi_{2,1},\pi_{2,2},f_1)\vet \\[.3em] 
f_2\circ (\pi_{2,1},S\circ \pi_{2,2}) & = & h\circ (\pi_{2,1},\pi_{2,2},f_2) 
 \end{array} \right\}
\Vdi{f_1,f_2:F_2;h:F_3} f_1=f_2
$ 
}
\end{enumerate}
\end{enumerate}

Notez qu'en (e) les axiomes affirment que la fonction $f=\rR_2(g,h)$ vérifie les \prts attendues d'une \dfn par \recu tandis qu'en (f) l'axiome affirme l'unicité de la fonction vérifiant ces \prts. 

Voyons maintenant comment une \demo par \recu usuelle se traduit dans le \gui{langage des fonctions} que nous avons mis en place. Par exemple la distributivité de la multiplication sur l'addition. Dans la \demo usuelle on prouve l'\egt $x\times (y+z)=(x\times y)+(x\times z)$ par \recu sur $z$ comme suit:
\begin{itemize}
\item \textsl{Initialisation.} 
$$x\times (y+0)\,\eqdf{1}\, x\times y\,\eqdf{2}\, (x\times y)+0 \,\eqdf{3}\,(x\times y)+(x\times 0)
$$
avec: $1:$ initialisation de $a+\cdot$, $2:$ initialisation de $a+\cdot$, $3:$ initialisation de $a\times \cdot$
\item \textsl{Induction.} 
\[ 
\begin{array}{rccclll} 
x\times (y+S(z)) & \eqdf{4} & x\times S(y+z)& \eqdf{5} & (x\times (y+z))+x& \eqdf{6}\\[.3em] 
((x\times y)+(x\times z))+x & \eqdf{7} & (x\times y)+((x\times z)+x) 
 & \eqdf{8} & (x\times y)+(x\times S(z))
 \end{array}
\] 
avec $4:$ induction de $a+\cdot$, $5,8:$ induction de $a\times\cdot$, $6:$ \hdr,\\ $7:$ associativité de $+$ (démontrée auparavant),
%
%\item 
\end{itemize}

\smallskip Traduisons tout cela dans le langage des fonctions, pour les \elts de $F_3$ 
$$
f=\Prod(\pi_{3,1},\mathrm{\Sum}(\pi_{3,2},\pi_{3,3}))\hbox{ et } g=\mathrm{\Sum}(\Prod(\pi_{3,1},\pi_{3,2}),\Prod(\pi_{3,1},\pi_{3,3})).
$$

Pour valider $f=g$, on utilise le principe $\tsbf{REC}_3$. Pour cela on valide les trois hypothèses. 
Tout d'abord l'initialisation, qui est $f(\pi_{2,1},\pi_{2,2}, 0_1)=g(\pi_{2,1},\pi_{2,2}, 0_1)$. 

\begin{itemize}
\item On a $f(\pi_{2,1},\pi_{2,2}, 0_1)=\Prod(\pi_{2,1},\mathrm{\Sum}(\pi_{2,2},0_1))$.
Comme $\mathrm{\Sum}(\pi_{2,2},0_1)=\pi_{2,2}$ d'après l'initialisation dans la \dfn par \recu de $\mathrm{\Sum}$, on obtient 
$$
f(\pi_{2,1},\pi_{2,2}, 0_1)=\Prod(\pi_{2,1},\pi_{2,2}).
$$
\item On a $g(\pi_{2,1},\pi_{2,2}, 0_1)=\mathrm{\Sum}(\Prod(\pi_{2,1},\pi_{2,2}),\Prod(\pi_{2,1},0_1))$. Comme $\Prod(\pi_{2,1},0_1)=0_1$ d'après la \dfn par \recu de $\Prod$, on obtient 
$$
g(\pi_{2,1},\pi_{2,2}, 0_1)=\mathrm{\Sum}(\Prod(\pi_{2,1},\pi_{2,2}),0_1),
$$ puis $g(\pi_{2,1},\pi_{2,2}, 0_1)=\Prod(\pi_{2,1},\pi_{2,2})$ d'après la \dfn par \recu de $\Sum$.
\end{itemize}

\smallskip Ensuite nous devons valider \gui{le passage de $n$ à $S(n)$}, i.e. trouver un \elt $h\in F_4$ convenable, i.e. satisfaisant les \egts
\[ 
\begin{array}{rcl} 
f(\pi_{3,1},\pi_{3,2},S(\pi_{3,3})) & = & h(\pi_{3,1},\pi_{3,2},\pi_{3,3},f) \\[.3em] 
g(\pi_{3,1},\pi_{3,2},S(\pi_{3,3})) & = & h(\pi_{3,1},\pi_{3,2},\pi_{3,3},g) \end{array}
\]

\begin{itemize}
\item On a $f(\pi_{3,1},\pi_{3,2},S(\pi_{3,3})) = \Prod(\pi_{3,1},\Sum(\pi_{3,2},S(\pi_{3,3})))$
ce qui donne d'après l'induction dans la \dfn par \recu de $\Sum$
$$f(\pi_{3,1},\pi_{3,2},S (\pi_{3,3}))=\Prod(\pi_{3,1},S(\Sum(\pi_{3,2},\pi_{3,3}))),$$
puis d'après la \dfn par \recu de $\Prod$
$$f(\pi_{3,1},\pi_{3,2},S (\pi_{3,3}))=\Sum(\Prod(\pi_{3,1},\Sum(\pi_{3,2},\pi_{3,3})),\pi_{3,1})=\Sum(f,\pi_{3,1}).
$$ 
\item De même on aura (en utilisant notamment l'associativité de l'addition) 
$$
g(\pi_{3,1},\pi_{3,2},S (\pi_{3,3}))=\Sum(\Sum(\Prod(\pi_{3,1}, \pi_{3,2}),\Prod(\pi_{3,1}, \pi_{3,3})),\pi_{3,1})=\Sum(g,\pi_{3,1}).
$$
\item Nous avons donc validé les hypothèses avec l'\elt $h\in F_4$
défini par $$h=h(\pi_{4,1},\pi_{4,2},\pi_{4,3},\pi_{4,4}):=\Sum(\pi_{4,4},\pi_{4,1}).$$ 
\end{itemize}

Nous notons \SA{PRA} la \tdy de l'arithmétique primitive récursive que nous venons de définir. Cette \tdy démontre exactement les mêmes énoncés que le système développé par Goodstein.
 
%%%%%%%%%%%%%%%%%%%% SECTION %%%%%%%%%%%%%%%%%%%%%%%%%%%%
%%%%%%%%%%%%%%%%%%%%%%%%%%%%%%%%%%%%%%%%%%%%%%%%%%%%%%%%%%%%%%%%%%%%
\section{Extensions conservatives d'une \tdy} \label{secconservative}

\Subsection{Extensions \eseqs}

%: Definition{defithconserv}
\begin{definition} \label{defithconserv} On dit qu'une \tdy~$\sab{T}'$ est une \textsl{extension conservative de la théorie~\sa{T}}
si c'est un extension de \sab{T} et si les \rdys formulables dans \sab{T} et valides dans~$\sab{T}'$
sont valides dans \sa{T}\,\footnote{La réciproque est claire.}.\index{extension --- d'une théorie dynamique!conservative}
\end{definition}
%----------- fin definition -------------------------------- 

Deux \tdys sur le même langage sont dites \textsl{identiques} lorsqu'elles prouvent exactement les mêmes \rdys. Autrement dit les axiomes de l'une sont des \rdys valides dans l'autre. Chacune est évidemment une extension conservative de l'autre.

\smallskip Le cas le plus simple d'extension conservative lorsque le langage a augmenté est celui des extensions qui sont \eseqs.

%:     Definition{defi-exteseq}
\begin{definition} \label{defi-exteseq}
Une extension $\Tp$ de la \tdy \sa T est dite \textsl{\eseq}
si elle est obtenue, à renommages près, au moyen des procédés suivants, utilisés de manière itérative, chaque fois en donnant les axiomes convenables (voir les détails dans \cite[section 2.3]{LM2022}).

\noindent Les extensions \textsl{\esids} sont celles qui sont obtenues sans ajouter de nouvelles sortes.%
\index{extension!identique}%
\index{extension!essentiellement identique}%
\index{extension!essentiellement equi@essentiellement équivalente}
\begin{itemize}
\item Ajout d'abréviations.
\item Ajout de prédicats traduisant la conjonction ou la disjonction de prédicats déjà existants. Cela revient à accepter $\vii$ et $\vuu$ comme symboles logiques permettant de construire des formules composées, \cad à accepter un peu de logique intuitionniste dans le langage.\rdb\label{NOTAvii}\label{NOTAvuu} 
\item Ajout de prédicats traduisant une formule $\exists x P$ où $P$ est un prédicat déjà existant et $x$ une variable. Même commentaire que pour le point précédent.\rdb\label{NOTAexists} 
\item Ajout d'un symbole de fonction lorsqu'est valide une existence unique sous certaines hypothèses.\rdb\label{skolemunique} 
\item Ajout d'une sorte produit de plusieurs sortes.
\item Ajout d'une sorte réunion disjointe de plusieurs sortes.
\item Ajout d'une sous-sorte d'une sorte déjà existante, définie comme les \elts satisfaisant un prédicat existant. 
\item Ajout d'une sorte quotient d'une sorte déjà existante, définie par un prédicat binaire, qui est prouvablement une relation d'\eqvc, et qui définit la nouvelle égalité dans la sorte quotient.
\item Ajout d'une sorte dont les objets sont (certains) morphismes d'une sorte vers une autre qui partage une certaine structure \agq avec la premières. 
\end{itemize}
\end{definition}
%----------- fin definition -------------------------------- 

Les extensions \eseqs sont \inteqs  au sens de la \dfn \ref{defextintequiv}

%: definition{defextintequiv}
\begin{dfni} \label{defextintequiv}
On considère une \tdy \sa{T} et une extension $\Tp$ de \sa{T}. On dit que $\Tp$ est un extension  \textsl{intuitivement équivalente} de \sa T si sont vérifiées les trois \prts suivantes.\index{extension!intuitivement équivalente} 
\begin{enumerate}
\item La théorie $\Tp$ est une extension conservative de \sa{T}.
\item [\textsl{2}.] Toute \rdy formulée dans le langage de $\Tp$ est \eqve\footnote{L'\eqvc en question est une règle externe, à l'image des règles structurelles décrites précédemment. Elle peut dépendre de la logique utilisée dans le monde extérieur.} à une famille de \rdys formulées dans le langage de $\sa{T}$.
\item [\textsl{3}.] Pour toute \pn~$(G,R)$ dans le langage de \sa{T},
les \sads $\gA=\big((G,R),\sab{T}\big)$ \hbox{et $\Tp(\gA):=\big((G,R),\Tp\big)$} ont les mêmes modèles (en \coma comme en \clama).
\end{enumerate}
\end{dfni}
%--------- fin lemma ----------------------------------

\bigskip Nous donnons dans la suite de cette section deux cas très importants d'extensions conservatives
(\thos~\ref{thFond} et \ref{thFondExists}), relatifs à l'utilisation de la logique classique, qui ne relèvent pas du cas simple des théories \eseqs. 

Auparavant, nous donnons le \thref{thFond0} relatif à l'utilisation de la logique constructive (intuitionniste).

%: Subsection{Comparaison avec la logique intuitionniste}
\Subsection{Comparaison avec la logique intuitionniste}\label{subsubsecthFond0}

On peut considérer que les \tdys ne sont que des versions tronquées de la déduction naturelle intuitionniste, dans lesquelles on n'introduit ni le connecteur $\Rightarrow$ ni le quantificateur~$\forall$.

C'est ce qui fait précisément la force des \tdys: ne pas s'encombrer de formules \gui{compliquées} du style $(A\Rightarrow B)\Rightarrow C$, ou $\forall x\;\exists y\;\forall z \dots$, permet d'y voir plus clair et de simplifier un certain nombre de résultats non triviaux, lorsqu'ils peuvent se démontrer au niveau de base de la déduction naturelle, \cad avec le système \gui{sans logique} des preuves dynamiques.

%: --- Theorem{thFond0}------- 
\begin{theorem}[conservativité de la logique intuitionniste par rapport aux preuves dynamiques] 
\label{thFond0} 
On considère une \tdy finitaire \sa{T}. Si une formule cohérente sur le langage de~\sa{T} est démontrée en utilisant les axiomes de \sa{T} et la logique
intuitionniste, la \rdy correspondante peut être démontrée valide directement dans la \tdy \sa{T}.
\end{theorem}
%--- end-theorem-----------

%
\begin{proof}
Voir \cite[Coquand, 2005]{Coq2005}. 
\end{proof}

Thierry Coquand interprète intuitivement la démonstration comme la construction d'un certain type de modèle (un \textsl{modèle générique}) de la \sad considérée.
La démonstration dans \cite{Coq2005} est plus directe et plus intuitive que celle du résultat de conservativité un peu plus fort donné dans le \thref{thFond}.

%%%%%%%%%%%%%%%%%%%%%%%%%%%%%%%%%%%%%%%%%
%: Subsection{\Tho fondamental des \tdys}
\Subsection{\Tho fondamental des \tdys}\label{subsubsecthFond}

À une \tdy $\sa{T}$ correspond une théorie cohérente, ou \tgm du premier ordre, obtenue en remplaçant les \rdys par les formules correspondantes selon le schéma donné \paref{subsectdy} au début de la section \ref{subsectdy}. Cette théorie cohérente peut être traitée selon la logique classique ou selon la logique intuitionniste. Notons les respectivement $\sab{T}^{\rm c}$ et $\sab{T}^{\rm i}$. 
On a le \tho fondamental \ref{thFond} ci-après (cf. par exemple le Theorem 1 dans \cite{CLR01}).
Ce \tho est déjà donné pour les théories purement équationnelles
dans \cite[Prawitz, 1971]{pra1971}, et ce genre de résultat
est omniprésent dans la littérature contemporaine, sous des formes plus ou moins variées. 
Nous recommandons les progrès récents sur ce thème décrits dans les articles \cite{DN2015,Dyc2015} qui montrent que, convenablement traitées, les preuves classiques ne fournissent pas des preuves \covs sensiblement plus longues.
La \demo dans \cite{CLR01} est \cov et relativement intuitive, mais conduit à une explosion de la taille des \demos.

Elle repose sur le lemme suivant qui explique le caractère inoffensif, \und{dans} \und{certaines} \und{circonstances}, de la règle du tiers \hbox{exclu}.

%: Lemma{lemNegValide}
\begin{lemma}[élimination de la négation classique] \label{lemNegValide} ~ \\
Soit \sa{T} une \tdy finitaire, et $P(.,.)$
un prédicat faisant partie de la signature (on l'a pris ici d'arité 2 à titre d'exemple). Introduisons \gui{le prédicat opposé à $P$}, notons le $Q(.,.)$, avec les deux \rdys qui le définissent en \clama\footnote{La \dfn du prédicat opposé à un prédicat $P$ en \coma n'est pas la même, et elle ne se laisse pas traiter dans le cadre des \tdys, sauf dans le cas où le prédicat est décidable. La signification \cov de $\lnot P$ est $P\Rightarrow\Bot$ et l'implication \cov ne peut pas être traitée par la méthode dynamique toute seule.}:

\DeuxRegles{
\lab{In-non$_P$} $ \vd P(x,y) \;\vou\; Q(x,y)$
}{
\lab{El-non$_P$} $\,\, P(x,y) ,\, Q(x,y) \vd \Bot %1=0
$
}

\smallskip\noindent Alors, la nouvelle \tdy est une extension \cosv de \sa{T}. 
\end{lemma}
%--------- fin lemma -----------------------------------

On note que cette fois-ci certains modèles \cofs de la première théorie peuvent ne plus être des modèles \cofs de la seconde. Néanmoins, ce n'est pas trop grave, comme l'indique le lemme, et cela se généralise dans
le \tho fondamental suivant.

%: --- Theorem{thFond}------- 
\begin{theorem}[élimination des coupures] 
\label{thFond} ~ \\
Pour ce qui concerne les \tdys finitaires, la logique, y compris classique (et en particulier le principe du tiers exclu) ne sert qu'à raccourcir les preuves. 
Plus précisément:
une règle dynamique est valide dans une \tdy \sa{T} \ssi elle est valide dans la \tco classique correspondante (celle qui a la même signature et les mêmes axiomes que~\sa{T}): on utilise dans la \tco les connecteurs, les quantificateurs et la logique classique du premier ordre.
\end{theorem}
%--- end-theorem-----------

%r
%: Remark{remprogHilbert}
\begin{remark} \label{remprogHilbert} 
Le \tho précédent, et le suivant concernant la skolémisation, montrent que l'usage des \tdys permet de réaliser en partie le programme de Hilbert, en fournissant une sémantique \cov pour certains usages du princicpe du tiers exclu et de l'axiome du choix. 
\eoe
\end{remark}
%----------- fin remark ---------------------------------- 

%%%%%%%%%%%%%%%%%%%%%%%%%%%%%%%%%%%%%%%%%%%%%%%%%%%%%%%%%%%%%%%%%%%%
%: Subsection{Skolémisation}
\Subsection{Skolémisation}
%-----------
%%%%%%%%%%%%%%%%%%%%%%%%%%%%%%%%%%%%%%%%%%%%%%%%%%%%%%%%%%%%%%%%%%%%

Nous examinons maintenant le processus de skolémisation général
qui consiste à se débarrasser des $\,\Exists\,$ dans certaines règles valides 
d'une \tdy en remplaçant les \gui{existants} par des fonctions.

Nous avons déjà indiqué le cas où cette opération est inoffensive, 
selon la remarque informelle suivante: quand l'existant dans une règle valide est prouvablement unique, cela ne mange pas de pain de remplacer la variable muette qui désigne l'existant par un symbole de fonction. 

En revanche, remplacer la variable muette qui désigne l'existant par un symbole de fonction lorsque l'existant n'est pas prouvablement unique est plus problématique. C'est la skolémisation proprement dite. Certains modèles \cofs avant skolémisation peuvent ne plus convenir après skolémisation, et la nouvelle théorie peut ne plus avoir de modèle \cof connu. Et même en \clama,
si les modèles sont \gui{presque} les mêmes, c'est à condition de supposer l'axiome du choix.

%: --- Theorem{thFondExists}------- 
\begin{theorem}[skolémisation] 
\label{thFondExists}

On considère une \tdy \sa{T}. On note $\sab{T}'$ la théorie \gui{skolémisée}, où l'on a skolémisé tous les axiomes existentiels en remplaçant les $\Exists$ par l'introduction de symboles de fonctions. 
Alors~$\sab{T}'$ est une extension \cosv de \sa{T}.
\end{theorem}
%--- end-theorem-----------
 
%
\begin{proof}
Une preuve en \clama avec axiome du choix consiste à constater que les deux théories ont \gui{les mêmes modèles}. Une \demo syntaxique et \cov est donnée dans \cite[Bezem \& Coquand, 2019]{BC2019}.
\end{proof}

\section[Treillis distributifs associés à une \sad]{Treillis distributifs et espaces spectraux associés à une \sad}
\label{subsectrdisad}

Pour cette section, nous renvoyons à \cite[chapitres XI et XIII]{CACM} et \cite[sections 1 et~3]{LM2022}.

\Subsection{Treillis distributifs et relations implicatives}

Une règle particulièrement importante
pour les \trdis, appelée \textsl{coupure}, est la
suivante
%-----------------begin $$----------------
\begin{equation}\label{coupure1}
\hbox{si } \, x\vi a \leq b \, \hbox{ et } \, a \leq x\vu b 
\,\hbox{ alors } \, (a \leq b).
\end{equation}
%-----------------end $$------------------

Si $A\in\Pfe(\gT)$ (ensemble des parties finiment énumérées de~$\gT$) on notera
%-----------------begin $$----------------
\[\ndsp \Vu A:=\Vu _{x\in A}x\qquad {\rm et}\qquad \Vi A:=\Vi _{\!x\in A}x.
\]
%-----------------end $$------------------
On note $A \vda B$ ou $A \vdash_\gT B$ la relation définie comme suit sur l'ensemble $\Pfe(\gT)$:

\snic{A \vda B \; \; \equidef\; \; \Vi A\;\leq \;
\Vu B.}

Cette relation vérifie les axiomes suivants, dans lesquels on
écrit $x$ pour $\{x\}$ et $A, B$ \hbox{pour $ A\cup B$}.

\vspace{-.5em}
%--------------------begin array---------------
\[\arraycolsep3pt\begin{array}{rcrclll}
& & x &\vda& x &\; &(R) \\[1mm]
 \hbox{si } A \vda B & \hbox{alors} & A,A' &\vda& B,B' &\; &(M) \\[1mm]
\hbox{si } (A,x \vda B) \hbox{ et } (A \vda B,x) 
& \hbox{alors} & A &\vda& B &\;
&(T).
\end{array}
\]
%---------------------end array--------------
On dit que la relation est \textsl{réflexive}, \label{remotr} \textsl{monotone} et
\textsl{transitive}.
La troisième règle (transitivité) peut être vue comme une
réécriture de la règle (\ref{coupure1}) et s'appelle \egmt la
règle de \textsl{coupure}.

%: --- Definition{defEntrel}-------------
\begin{definition}
\label{defEntrel}
Pour un ensemble $S$ arbitraire, une relation binaire sur $\Pfe(S)$ qui est
réflexive, monotone et transitive est
appelée une {\sl \entrel} (en anglais, {\sl entailment relation}).
\end{definition}
%--- end-definition------------------------------------

Le \tho suivant est fondamental. Il dit que les
trois propriétés des \entrels sont exactement ce qu'il faut pour que
l'interprétation d'une relation implicative comme la trace de celle d'un \trdi soit adéquate.

%: Theorem{thEntRel1}----------------
\begin{theorem}[\tho fondamental des \entrels]
\label{thEntRel1} {\rm \cite[Satz 7]{Lor1951}, \cite{CC00}, \cite[\hbox{XI-5.3}]{ACMC}}.
Soit un ensemble~$S$ avec une \entrel
$\vdash_S$ sur $\Pfe(S)$. On considère le \trdi~$\gT$ défini par
\gtrs et relations comme suit: les \gtrs sont les
\elts de $S$ et les relations sont les

\snic {A\, \vdash_\gT \, B}

\noindent chaque fois que $A\, \vdash_S \, B$. Alors, pour tous $A$,
 $B$ dans $\Pfe(S)$, on a

\snic {\hbox{si }A\, \vdash_\gT \, B
\hbox{ alors } A\, \vdash_S \, B.}

\noindent 
En particulier, deux \elts $x$ et $y$ de $S$ définissent le même \elt de $\gT$ \ssi
\hbox{on a} $x \vdash_S y$ et $y \vdash_S x$.
\end{theorem}
%--- end-theorem------------------------------------

%:Subsubsection{Le spectre d'un \trdi}
\Subsection{Le spectre d'un \trdi} 

En \clama un {\sl \id premier} $\fp$ d'un \trdi $\gT\neq \Un$ est un \id dont
le complémentaire $\fv$ est un filtre (qui est alors un {\sl
filtre premier}). On~a alors $\gT/(\fp=0,\fv=1)\simeq\Deux$. Il
revient au même de se donner un \idep de $\gT$ ou un morphisme de
\trdis $\gT\rightarrow \Deux$.

%Nous noterons $\theta_\fp:\gT\to\Deux$ l'\homo
%associé à l'\idep $\fp$.

On vérifie facilement que si $S$ est une partie génératrice du
\trdi $\gT$, un \idep~$\fp$ de $\gT$ est complètement
caractérisé par sa trace sur $S$ (cf. \cite{CC00}).

%d
%: Definition{defiSpecTrdi}
\begin{definition} \label{defiSpecTrdi}
Le \textsl{spectre d'un \trdi $\gT$} est l'ensemble $\Spec\,\gT$ de ses \ideps, muni de la topologie suivante: une base
d'ouverts est donnée par les 
$$
\DT(a)\eqdefi\sotq{\fp\in\Spec
\,\gT}{a\notin\fp},\quad a\in \gT.
$$
 \end{definition}
%----------- fin definition -------------------------------- 
On vérifie en \clama que
%--- equation eqDa --------
\begin{equation} \label{eqDa}
\left.
\begin{array}{rclcrcl}
 \DT(a\vi b) & = & \DT(a)\cap \DT(b) ,&\quad & \DT(0) & = & 
\emptyset ,\\
 \DT(a\vu b) & = & \DT(a)\cup \DT(b) ,&& \DT(1) & = & 
\Spec\,\gT.
 \end{array}
\right\}
\end{equation}
%---------------------end equation--------------

Le complémentaire de $\DT(a)$ est un fermé que l'on note $\VT(a)$.

On étend la notation $\VT(a)$ comme suit: si $I\subseteq\gT$, on
pose $\VT(I)\eqdefi\bigcap_{x\in I}\VT(x)$. Si~$I$ engendre l'idéal $\fII$, on a $\VT(I)=\VT(\fII)$. On dit parfois que $\VT(I)$ est \textsl{la
variété associée~à~$I$}.

%d
%: Definition{defiEspaceSpectral}
\begin{definition} \label{defiEspaceSpectral}
 Un espace topologique homéomorphe à un espace $\Spec(\gT)$
est appelé un \textsl{espace spectral}. \end{definition}
%----------- fin definition -------------------------------- 

Les espaces spectraux proviennent de l'étude \cite[Stone, 1937]{Sto37}.

\cite{Joh1986} appelle ces espaces des \textsl{espaces cohérents}.
\cite{BW74} les appellent \textsl{Stone spaces}. Cette terminologie est obsolète, car depuis \cite{Joh1986} les espaces de Stone sont les espaces spectraux associés aux \trdis qui sont des \agBs.

Le nom \textsl{spectral space} est donné par \cite[Hochster, 1969]{Hoc1969}, qui les a popularisés dans la communauté \mathe après une hibernation prolongée depuis 1937.

Avec la logique classique et l'axiome du choix, l'espace $\Spec (\gT)$
a \gui{suffisamment de points}: on peut retrouver le treillis $\gT$
à partir de son spectre.

\smallskip 
On dit qu'un point $\fp$ d'un espace spectral $X$ est le
\textsl{point générique du fermé $F$} \hbox{si $F=\ov{\so{\fp}}$}. Ce point
(quand il existe) est nécessairement unique car les espaces spectraux
sont des espaces de Kolmogoroff. En fait, les fermés $\ov{\so{\fp}}$
sont exactement tous les fermés irréductibles de $X$. La relation
d'ordre $\fq\in\ov{\so{\fp}}$ sera notée $\fp\leq_X \fq$, et \hbox{l'on a} les équivalences
% equation label {eqOrdreSpec}
\begin{equation} \label {eqOrdreSpec}
\fp\leq_X \fq\;\Longleftrightarrow\; \ov{\so{\fq}}\subseteq \ov{\so{\fp}}\,.
\end{equation}
% end-equation
Les points fermés de $\Spec(\gT)$ sont les \idemas de $\gT$.
Lorsque $X=\Spec(\gT)$ la relation $\fp\leq_X \fq$ est simplement la
relation d'inclusion usuelle $\fp\subseteq \fq$ entre \ideps du \trdi $\gT$.

Dans la catégorie des espaces spectraux on définit les \textsl{morphismes (spectraux)} comme des applications telles que l'image réciproque de tout \oqc est un \oqc (en particulier elles sont continues).

\Subsection{L'antiéquivalence de Stone}

L'antiéquivalence de Stone affirme (en langage moderne) qu'en \clama la catégorie des \trdis est antiéquivalente à la catégorie des espaces spectraux.

Bien que les espaces spectraux aient envahi l'\alg abstraite contemporaine, c'est seulement en \coma qu'on accorde l'attention qu'elle mérite à cette antiéquivalence des \clama. 

L'objectif est de définir correctement les \trdis correspondant aux espaces spectraux de la littérature classique, et, si possible, de décrypter les discours classiques utilisant les espaces spectraux en discours \cofs concernant les \trdis correspondants.

\Subsection{Le treillis et le spectre de Zariski d'un anneau commutatif}

Le treillis de Zariski d'un anneau commutatif peut être obtenu à partir des règles valides dans différentes extensions de la théorie \sa{Ac} des anneaux commutatifs.

Nous choisissons la théorie des anneaux locaux en raison de leur rôle fondamental dans les schémas de Grothendieck.
On considère \prmt la \tdy \Sa{Al1} des  anneaux locaux avec unités.

Soit $\gA$ un anneau commutatif. Considérons la \entrel $\vdash_{\gA,\mathrm{Zar} }$ sur l'ensemble sous-jacent à $\gA$ définie par l'\eqvc suivante 
% equation label {eqZarclass}
\vspace{-.8em}
\begin{equation} \label {eqZarclass}
\begin{aligned} 
 a_1,\dots,a_n &\,\vdash_{\gA,\mathrm{Zar}} c_1,\dots,c_m 
 \qquad\quad \equidef \\[.3em] 
\U(a_1)\vet \dots\vet \U(a_n) & \Vdi{\sA{Al1}(\gA)} \U(c_1)\vou \dots\vou \U(c_m) 
 \end{aligned}
\end{equation}

On définit le \textsl{treillis de Zariski de $\gA$}, noté $\ZarA$ ou 
$\Zar(\gA)$,
 comme celui engendré par la \entrel $\vdash_{\gA,\mathrm{Zar}}$.
 
L'application correspondante $\DA:\gA\to\ZarA$ s'appelle le \textsl{support de Zariski de $\gA$}. Lorsque~$\gA$ est fixé par le contexte on note simplement $\rD$.

Le \textsl{spectre de Zariski} usuel est l'espace spectral dual de ce \trdi.

Notons que puisque $\rD(a_1)\vi\dots\vi\rD(a_n)=\rD(a_1\cdots a_n)$, les \elts de $\ZarA$ sont tous de la forme $\rD(c_1,\dots,c_m):= \rD(c_1)\vu\dots\vu\rD(c_m)$.

On démontre alors les \eqvcs suivantes. L'\eqvc des points (1) et (2) recopie pour l'essentiel la définition
de~$\ZarA$). L'\eqvc avec le point (3) est l'objet d'un \textsl{\nst formel}. Le \nst de Hilbert proprement dit est un résultat plus difficile. 
%t
%: Theorem{thNstFormel}
\begin{theorem}[\nst formel] \label{thNstFormel} ~\\
 Soit $\gA$ un anneau commutatif, et $\an,c_1,\dots c_m\in\gA$. \Propeq
\vspace{-.5em}
\[ 
\begin{aligned} 
(1)\qquad\quad\; \rD(a_1),\dots,\rD(a_n) &\,\vdash_{\ZarA} \rD(c_1),\dots,\rD(c_m) \\[.3em] 
(2)\qquad\quad \U(a_1)\vet \dots\vet \U(a_n) &\Vdi{\sA{Al1}(\gA)} \U(c_1)\vou \dots\vou \U(c_m) \\ 
(3)\qquad\; \exists k>0\;\;(a_1 \cdots a_n)^k&\,\in\gen{c_1,\dots,c_m} 
 \end{aligned} \label {eqNstFormel}
\] 
\end{theorem}
%----------- fin theorem ----------------------------- 

On peut donc identifier l'\elt $\rD(c_1,\dots,c_m)$ de $\ZarA$ 
à l'\id $\sqrt[\gA]{\gen{c_1,\dots,c_m}}$. Modulo cette identification, la relation d'ordre est la relation d'inclusion.

%c
%: Corollary{corZarA}
\begin{corollary} \label{corZarA}
Le treillis $\ZarA$ est engendré par la plus petite \entrel sur (l'ensemble sous-jacent~à)~$\gA$ satisfaisant les relations suivantes\footnote{En orthographe moderne recommandée, \gui{corolaire} prend un seul \gui{l}. Voir \url{http://www.renouvo.org/liste.php?t=3&lettre=c}.}.

\DeuxRegles
{
\labu $\;\;0\vdash 0_\gT$
\labu $\;\;ab\vdash a$
\labu $\;\;a+b\vdash a,b$
}
{
\labu $\;\; 1\vdash 1_\gT$
\labu $\;\;a,b\vdash ab$
} 

\noindent En d'autres termes, l'application $\rD:\gA\to\ZarA$ satisfait les relations 
\[
\rD(0)=0,\;\rD(1)=1,\;\rD(ab)=\rD(a)\vi\rD(b),\;\rD(a+b)\leq \rD(a)\vu \rD(b),
\]
et toute autre application $\rD':\gA\to \gT$ qui satisfait ces relations se factorise via $\ZarA$ avec un unique morphisme de \trdis $\ZarA\to \gT$.
\end{corollary}
%--------- fin corollary ------------------------------- 

\Subsection{Autres exemples}

Rappelons qu'une règle disjonctive est une \rdy où ne figure pas le symbole $\Exists$ et qu'une règle disjonctive simple est une \rdy de la forme suivante, avec $m,n\geq 0$.%
\index{regle@règle!disjonctive}\index{regle@règle!disjonctive!simple}
 
%---- equation {eqAxds} ----
\begin{equation} \label{eqAxds}
C_1\vet \ldots \vet C_n \vd D_1\vou \ldots \vou D_m\end{equation}
%------end equation----
où les $C_i$ et $D_j$ sont des formules atomiques.
Les \ralgs (simples) sont des cas particuliers de règles disjonctives (simples). 

\paragraph{Premier exemple.} Considérons une \sad $\gA=\big((G,R),\sa{T}\big)$ pour une \tdy $\sa{T}=(\cL,\cA)$.
 Si $P(x,y)$ est un prédicat binaire dans la signature, et si $Tcl=\Tcl(\gA)$ est l'ensemble des termes clos de~$\gA$, on obtient une \entrel $\vdash_{\gA,P}$ sur $Tcl \times Tcl$ en posant
 
\vspace{-.8em}
% equation label {eq1}
\begin{equation} \label {eq1}
\begin{aligned} 
 (a_1,b_1),\dots,(a_n,b_n) &\,\vdash_{\gA,P} (c_1,d_1),\dots,(c_m,d_m) 
 \qquad\quad \equidef \\[.3em] 
 P(a_1,b_1)\vet \dots\vet P(a_n, b_n) & \Vdi{\gA} P(c_1, d_1)\vou \dots\vou P(c_m, d_m) 
 \end{aligned}
\end{equation}
% end-equation

Intuitivement le \trdi engendré par cette \entrel est le treillis des \gui{valeurs de vérité} du prédicat~$P$ dans la \sad $\gA$.
 
\paragraph{Plus généralement.} Considérons une \sad $\gA=\big((G,R),\sa{T}\big)$ pour une \tdy $\sa{T}=(\cL,\cA)$.
Soit $S$ un ensemble de formules atomiques closes de $\gA$. On définit \textsl{la \entrel sur $S$ associée à $\gA$} comme suit: 
 
\vspace{-.5em}
% equation label {eq2}
\begin{equation} \label {eq2}
\begin{aligned}
 A_1,\dots,A_n &\,\vdash_{\gA,S} B_1,\dots,B_m 
 \qquad\quad \equidef \\[.2em] 
 A_1\vet \dots\vet A_n &\Vdi{\gA} B_1\vou\dots\vou B_m 
 \end{aligned}
\end{equation}
% end-equation
On pourra noter $\Zar(\gA,S)$ le \trdi engendré par cette \entrel.

En particulier, le treillis $\Zar(\gA,\Atcl(\gA))$ est appelé le \textsl{treillis de Zariski absolu de la \sad $\gA$}.\index{Treillis de Zariski!absolu d'une \sad}

\paragraph{Cas d'une extension \sa{T$_1$} qui reflète les \rdijs valides.} 
Soit \sa{T$_1$} une extension d'une \tdy \sa{T} qui prouve exactement les mêmes \rdijs (par exemple une extension conservative). Soient $\gA=\big((G,R),\sa{T}\big)$ et $\gA_1=\big((G,R),\sa{T$_1$}\big)$. Soit $S$ un ensemble de formules atomiques closes de $\gA$.
Alors les treillis de Zariski $\Zar(\gA,S)$ et $\Zar(\gA_1,S)$ sont isomorphes. 

En particulier, lorsque \sa{T$_1$} est une extension \eseq de \sa{T} les treillis de Zariski absolus de~$\gA$ et~$\gA_1$ sont isomorphes.

\paragraph{Les treillis de Zariski donnent cependant une image amoindrie} d'une \sad. 
D'une part, dans ces treillis de Zariski à priori rien n'est pris en compte qui corresponde aux \rdys valides lorsqu'elles ne sont pas disjonctives. 
D'autre part ajouter la logique classique et skolémiser une \tdy ne changent pas les treillis correspondant à $\Atcl(\gA)$, mais dans ce cas, le treillis de Zariski absolu de $\gA_1$ est l'\agB engendrée par 
$\Zar(\gA,\Atcl(\gA))$.
Pour retrouver la richesse des \tdys vues d'un point de vue \cof, il faut alors faire appel à la théorie des faisceaux ou des topos.

%%%%%%%%%%%%%%%%%%%%%%%%%%%%%%%%%%%%%%%%%%%%%%%%%%%%%%%%%%%%%%%%%%%%
%%%%%%%%%%%%%%%%%%%%%%%%%%%%%%%%%%%%%%%%%%%%%%%%%%%%%%%%%%%%%%%%%%%%
\section{Théorie des modèles}\label{subsecModeles}

%Dans cet article, les \thos ou lemmes de \clama qui n'ont pas de \demo \cov connue, et qui souvent ne peuvent pas en avoir, sont indiqués avec une étoile.

\Subsection{Théories qui partagent certaines règles}

On considère deux \tdys $\sa{T}_1$ et $\sa{T}_2$ dont les signatures étendent une même signature $\Sigma$. 

%d
%:     Definition{defiColSim}
\begin{definition} \label{defiColSim} On suppose chacune des deux théories a un axiome de collapus.
On dit que les \tdys $\sa{T}_1$ et $\sa{T}_2$ \textsl{s'effondrent simultanément} (pour $\Sigma$) lorsque, pour toute présentation $(G,R)$ sur $\Sigma$, les \sads $\gA_1=\big((G,R),\sa T_1\big)$ et $\gA_2=\big((G,R),\sa T_2\big)$ s'effondrent simultanément.  
\end{definition}
%----------- fin definition -------------------------------- 
  
%d
%:     Definition{defiMemesRalgs}
\begin{definition} \label{defiMemesRalgs}
On dit que les \tdys $\sa{T}_1$ et $\sa{T}_2$ \textsl{prouvent les mêmes \ralgs} (pour $\Sigma$) lorsque, pour tout présentation $(G,R)$ sur $\Sigma$, les \sads $\gA_1=\big((G,R),\sa T_1\big)$ et $\gA_2=\big((G,R),\sa T_2\big)$ prouvent les mêmes \ralgs. 
\end{definition}
%----------- fin definition -------------------------------- 

Il revient au même de dire que les théories prouvent les mêmes faits (\dfn \ref{defiFait}).

\Subsection{\Tho de complétude, effondrement simultané}
Voici tout d'abord le \tho de complétude sous sa forme minimale:
son interprétation intuitive en \clama est que la logique classique donne de manière exhaustive les règles de raisonnements conformes à la \gui{vérité absolue}, fondée sur un univers \mathe idéal dans lequel aucun doute n'est jamais permis, le principe du tiers exclu est absolument vrai et l'axiome du choix tout pareillement.
 
%: Theorem{thGodel1}
\begin{theoremc}[\tho de complétude de Gödel, première forme] \label{thGodel1} ~\\
Une \sad qui ne s'effondre pas admet un modèle non trivial. 
\end{theoremc}
%--------- fin theorem ----------------------------------- 

\comm Une forme équivalente du théorème de complétude est
le cas particulier suivant (lemme de Krull): \textsl{tout anneau commutatif non trivial possède un quotient intègre non trivial}.

\noindent
 La forme constructivement acceptable du lemme de Krull est le résultat facile suivant: 
 \textsl{lorsque l'on ajoute la règle dynamique \gui{$\,\,xy=0\vd x=0\vou y=0\,\,$} à la théorie des anneaux commutatifs, une structure algébrique dynamique s'effondre dans la première théorie \ssi elle s'effondre dans la seconde}.

\noindent
 Autrement dit, les théories \Sa{Ac} et \Sa{Asdz} s'effondrent simultanément. 
\eoe

%: Theorem{thGodel2}
\begin{theoremc}[\tho de complétude de Gödel, deuxième forme] \label{thGodel2} ~\\
On considère une \tdy \sa{T} et une \sad $\gA$ de type~\sa{T}. Un fait est valide dans $\gA$ \ssi il est satisfait dans tous les modèles de~$\gA$. 
\end{theoremc}
%--------- fin theorem ----------------------------------- 
%
Une \tdy qui en étend une autre (en ajoutant des sortes et/ou des prédicats et/ou des axiomes) prouve à priori plus de résultats. Un cas intéressant est lorsqu'elle prouve les mêmes résultats tout en offrant des facilités plus grandes pour les \demos. C'était l'essence des \thos fondamentaux 
\ref{thFond} et \ref{thFondExists}. Une variante en théorie des modèles, mais seulement en \clama,
est donnée par les \thos qui suivent.

%
%: Theorem{thcolsimcomp}
\begin{theoremc}[effondrement simultané et modèles non triviaux] \label{thcolsimcomp} ~\\
Soient \sa{T} une \tdy et $\sab{T}'$ une extension qui s'effondre simultanément avec \sa{T}. Si une \sad de type~\sa{T} admet un modèle non trivial, elle admet \egmt un modèle non trivial en tant que \sad de type~$\sab{T}'$. Plus \prmt si~$M$ est un modèle non trivial de \sa{T}, la \sad $(\Diag(M),\sab{T}')$ admet un modèle non trivial.
\end{theoremc}
%--------- fin theorem ----------------------------------- 

\comm
Une bonne version constructive du théorème de complétude est 
sans doute une pure tautologie: si une structure algébrique dynamique ne s'effondre pas, alors \dots\ elle ne s'effondre pas. Ou encore: si une théorie formelle du premier ordre ne s'effondre pas, alors \dots\ elle ne s'effondre pas. Et pour le \thref{thcolsimcomp}: si \sa{T} et $\Tp$ s'effondrent simultanément, alors \dots\ elles s'effondrent simultanément. Il en irait de même pour les théorèmes \ref{thGodel2} et~\ref{thcolsimralg}. 

\noindent
 En effet, ce que l'on appelle une version constructive d'un théorème classique \gui{douteux} est un énoncé, correct en mathématiques constructives, qui, en pratique, 
i.e. pour démontrer des résultats concrets, rend les mêmes services que le théorème classique. Or, en pratique, tous ces théorèmes \gui{abstraits} ne sont utilisés,
pour aboutir à des résultats concrets, que dans des raisonnements par l'absurde, lesquels utilisent des modèles fictifs pour conclure qu'ils ne peuvent pas exister. Le résultat concret, quant à lui, ressemble beaucoup plus à l'hypothèse du théorème abstrait qui a été invoqué. L'analyse détaillée de l'ensemble de la démonstration montre alors en \gnl que l'on a tautologisé en rond sans s'en rendre compte (voir par exemple \cite{Lom98} pour le 17\ieme\ problème de Hilbert). C'est une des raisons qui expliquent pourquoi les mathématiques classiques sont si souvent constructives, contrairement à l'apparence que donnent leurs démonstrations.
\eoe

%%%%%%%%%%%%%%%%%%%%%%%%%%%%%%%%%%%%%%%%%%%%%%%%%%%%%%%%%%%%%%%%%%%%
%: Subsection[Théorèmes de plongement] 
\Subsection{Théorème de plongement, théories qui prouvent les mêmes \ralgs}
%-----------

Le \tho suivant s'appelle aussi un \gui{\tho de représentation}.

%: Theorem{thcolsimralg}
\begin{theoremc}[\tho de plongement] \label{thcolsimralg} 
On considère une \tdy $\sab{T}'$ qui étend une \talg~\sa{T} et qui prouve les mêmes \ralgs. Toute structure \agq $\gA$ de type~\sa{T} est isomorphe\footnote{La terminologie anglaise est: any algebraic structure of type \sA{T} is a \textsl{subdirect product} of algebraic structures of type $\sAb{T}'$.} à une sous-\sa{T}-structure d'un produit de structures \agqs de type~$\sab{T}'$. 
\end{theoremc}
%--------- fin theorem ----------------------------------- 

Par exemple \sa{T} est la théorie des \grls et $\sab{T}'$
est la théorie des groupes abéliens totalement ordonnés. Le résultat important, \cof, est que ces deux théories prouvent les mêmes \ralgs. L'interprétation intuitive en \clama est que tout \grl est un sous-\grl d'un produit de groupes totalement ordonnés. Lorsque Paul Lorenzen démontre ce résultat il généralise le résultat analogue de Krull qui dit que la \cli d'un anneau intègre $\gA$ est l'intersection des \advs de son corps de fractions qui contiennent $\gA$. 

%%:Subsection{Théories qui prouvent les mêmes règles algébriques}
%\Subsection{Théories qui prouvent les mêmes règles algébriques}

Les \thos intuitifs suivants, qui nous seront utiles, concernent les extensions qui prouvent les mêmes \ralgs, ils sont démontrés dans \cite{Lom-tgac}. 

%: Theorem{thEseqMemesfaits}
\begin{theorem} \label{thEseqMemesfaits}
Soit $\sab{T}_{\!2}$ une \tdy extension simple d'une théorie $\sab{T}_{\!1}$ et qui prouve les mêmes \ralgs. On considère une extension \eseq $\sab{T}_{\!1}'$ de $\sab{T}_{\!1}$ obtenue sans ajout de prédicats existentiels. On suppose qu'il n'y a pas d'interférence syntaxique au niveau des extensions de langage entre $\sab{T}_{\!1}'$ et $\sab{T}_2$. On peut donc construire une extension \eseq $\sab{T}_{\!2}'$ de $\sab{T}_{\!2}$ en recopiant pour $\sab{T}_{\!2}$ ce qui a été fait pour $\sab{T}_{\!1}$. 
 Dans ces conditions,~$\sab{T}_{\!2}'$ prouve les mêmes \ralgs que 
 $\sab{T}_{\!1}'$.
\end{theorem}
%----------- fin theorem ----------------------------- 

%: Theorem{thEseqMemesfaitsbis}
\begin{theorem} \label{thEseqMemesfaitsbis}
Soit $\sab{T}_{\!2}$ une \tdy extension simple d'une théorie $\sab{T}_{\!1}$ et qui prouve les mêmes \ralgs. On considère une extension $\sab{T}_{\!1}'$ de $\sab{T}_{\!1}$ obtenue en ajoutant une famille d'axiomes qui sont tous des \ralgs. On considère l'extension $\sab{T}_{\!2}'$ de $\sab{T}_{\!2}$ obtenue en ajoutant les mêmes axiomes. 
 Dans ces conditions, $\sab{T}_{\!2}'$ prouve les mêmes \ralgs que 
 $\sab{T}_{\!1}'$.
\end{theorem}

\newpage \thispagestyle{empty}

%%%%%%%%%%%%%%%%% CHAPITRE %%%%%%%%%%%%%

\chapter{Théories \gmqs infinitaires}\label{subsecgeominfini}\label{chap-gmqinfini}
\index{theorie@théorie!géométrique infinitaire}\index{theorie@théorie!géométrique}
%-----------
\Today

\minitoc

%\section*{Introduction}
%\addtocontents{toc}{\vskip0.8em}
%\addcontentsline{toc}{section}{Introduction}
%\rdb

\section{Généralités}
Une notion \gnle très utile de \tgm est définie, qui ne s'exprime pas \ncrt de manière finitaire. On parle alors de \textsl{\tgm infinitaire}. Dans une \tgm infinitaire, on autorise des \rdys qui ont une disjonction infinie dans le second membre. Une restriction essentielle est à noter: les variables libres présentes dans une telle disjonction doivent être précisées d'avance et en nombre fini.

Intuitivement, on utilise
de telles règles dans le \sys de preuves des \tdys en \gui{ouvrant les branches de calcul correspondant à la disjonction infinie}. Qu'est-ce que cela signifie \prmt? Cela signifie qu'une conclusion sera déclarée valide si elle
est valide dans chacune des branches. 

Prenons un exemple simple, et indiquons ce qui se passe si on a dans les axiomes une règle infinitaire du type

\Regles{\lab {~} $\Vdi{x_1,\dots,x_k} \Vou_{i\in I}\; \Gamma_i$
}

\noindent avec un ensemble infini $I$ et les $\Gamma_i$ des listes de formules atomiques sans autres variables libres que celles mentionnées (i.e. $\xk$). Si pour chaque $i\in I$ on a une règle valide $\,\,\Gamma_i\vd B(\ux)$, alors on déclare valide la règle $\Vdi{x_1,\dots,x_k}B(\ux)$. 

Il intervient donc \ncrt une \demo
intuitive \und{externe} à la \tdy pour certifier que la conclusion souhaitée est valide dans chacune des branches. En effet le \sys de calcul \gui{sans logique} à l'œuvre dans la \tdy ne peut pas prendre en charge
la validité une telle infinité de règles de deduction. Un calcul purement mécanique ne saurait ouvrir une infinité de branches! Par exemple, avec $I=\NN$ la \demo intuitive externe pourra être une preuve par \recu.

Notons par contre que la \demo interne doit démontrer la validité de la conclusion souhaitée selon les règles de \demo \gui{sans logique} de la \tdy.

\smallskip La présentation qui précède n'est qu'une ébauche. Tout ceci mérite une \dfn plus formelle de ce qu'est le fonctionnement légal d'une \tgm infinitaire; même s'il y a un aspect informel inévitable dans le recours à des \demos \gui{externes} en mathématiques intuitives. 

On devra aussi toucher deux mots du fonctionnement des théories formelles intuitionniste et classique qui étendent la \tdy infinitaire (par l'ajout du connecteur $\Rightarrow$ et du quantificateur universel dans le cas intuitionniste, par ajout de l'axiome de tiers exclu dans le cas classique).

Comme dans le cas des \tgms finitaires nous réserverons l'appellation de \tdy aux démonstrations dont la partie interne est \gui{sans logique}.

%: Subsec{Exemple: les \elts nilpotents,
\Subsection{Exemple: les \elts nilpotents, la dimension de Krull} 

Un \elt $x$ d'un anneau est nilpotent s'il existe un $n\in\NN^+$
tel que $x^n=0$. Si l'on introduit un prédicat $\mathrm{Z}(x)$ pour \gui{$x$ est nilpotent}, il sera soumis aux axiomes naturels suivants:

\DeuxRegles{
\Lab {nil1} $\vd \mathrm{Z}(0)$
\Lab {nil2} $\,\, \mathrm{Z}(x),\,\mathrm{Z}(y)\vd \mathrm{Z}(x+y)$
\Lab {NIL1} $\,\,Z(x)\vd \Exists z\;z(1+x)=1 $
}{
\Lab {nil3} $\,\,\mathrm{Z}(x)\vd \mathrm{Z}(xy) $
\Lab {Nil} $\,\,\mathrm{Z}(x^2)\vd \mathrm{Z}(x) $
}

\smallskip Dans la \tdy correspondante, les seuls termes $t$ pour lesquels on pourra démontrer~$\mathrm{Z}(t)$ seront ceux pour lesquels on pourra démontrer $t^n=0$ pour un~$n>0$. Rien ne garantit cependant que dans un modèle de la théorie, le prédicat~$\mathrm{Z}(x)$ corresponde bien à \gui{$x$ est nilpotent}.

La seule manière de s'en assurer est d'introduire la \rdy infinitaire

\Regles{\Lab{NIL} $\,\,\mathrm{Z}(x)\vd\Vou_{n\in\NN^+} x^n=0$}

\smallskip Cette préoccupation est en relation directe avec la dimension de Krull des anneaux commutatifs.

La dimension de Krull d'un \trdi peut être formulée dans une \tdy comme suit.

\penalty-2500
%: Definition{defiDDKTRDI}
\begin{definition}\label{defiDDKTRDI}~
\begin{enumerate}
\item Deux suites $(\xzn)$ et $(\bzn)$ dans
un \trdi $\gT$ sont dites \textsl{complémentaires} si
%--- equation eqC2G --------
\begin{equation}\label{eqC2G}
\left.\arraycolsep3pt
\begin{array}{rcl}
 b_0\vi x_0& = & 0 \\
 b_1\vi x_1& \leq & b_0\vu x_0 \\
\vdots~~~~& \vdots &~~~~ \vdots \\
 b_n\vi x_n & \leq & b_{n -1}\vu x_{n -1} \\
 1& = & b_n\vu x_n
\end{array}
\right\}
\end{equation}
%---------------------end equation--------------
Une suite qui possède une suite \cop sera dite \textsl{singulière}.
\item Pour $n\geq0$ on dira que le \trdi $\gT$ est \textsl{de \ddk $\leq n$} si toute suite $(\xzn)$ dans $\gT$ est singulière.
Par ailleurs, on dira que le \trdi $\gT$ est de \ddk $-1$
s'il est trivial, \cad si $1_\gT=0_\gT$.
\end{enumerate}
\end{definition}

Par exemple pour $n=2$ les \egts et inégalités \pref{eqC2G} correspondent au dessin suivant dans~$\gT$.
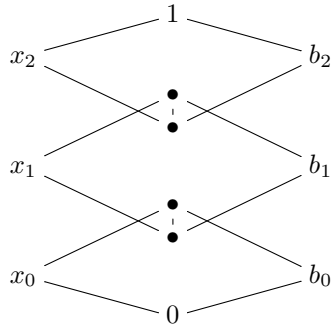
$$\SCO{x_0}{x_1}{x_2}{b_0}{b_1}{b_2}$$

Et la dimension $\leq 2$ correspond à l'axiome existentiel suivant.

\Regles{\Lab{KDIM2} $\,\,\vd \Exists b_0,b_1,b_2 \;(x_2\vu b_2=1,\,x_2\vi b_2\leq x_1\vu b_1,\,x_1\vi b_1\leq x_0\vu b_0, \,x_0\vi b_0=0 )$}

\smallskip Pour la dimension de Krull des anneaux, il faut faire intervenir le \trdi formé par les radicaux d'\itfs et l'on exprime par exemple
la dimension $\leq 2$ comme suit, en notant 
$\DA(x,y)=\sqrt{\gen{x,y}}=\sotq{z\in\gA}{\exists n\in\NN^+, z^n\in\gen{x,y}}$:

Pour tous $x_0$, $x_1$, $x_2 \in \gA$ il existe
$b_0$, $b_1$, $b_2\in \gA$ tels que
%--- equation eqCG --------
\begin{equation}\label{eqCG}
\left.\arraycolsep2pt
\begin{array}{rcl}
\DA(b_0x_0)& = &\DA(0) \\
\DA(b_1x_1)& \subseteq & \DA(b_0,x_0) \\
\DA(b_2 x_2 )& \subseteq & \DA(b_{1},x_{1}) \\
\DA(1)& = & \DA(b_2,x_2 )
\end{array}
\right\}
\end{equation}
%---------------------end equation--------------

Notons que $\DA(b_2 x_2 ) \subseteq \DA(b_{1},x_{1})$
signifie qu'il existe $a_1,y_1\in\gA$ et $n\in\NN^+$ tels que
$$
(b_2 x_2)^n=a_1b_1+y_1x_1.
$$ 

On pourra donc exprimer \gui{$\Kdim\gA\leq 2$} dans le cadre d'une \tgm infinitaire.
Et par exemple les \thos de Serre ou de Foster-Swan (\cite[Chapitre XIV]{CACM}) avec la dimension de Krull en hypothèse peuvent entièrement être traités dans le cadre de \tgms~(\cite{CL05}).

% Subsection{Le \tho de Barr}---- 
\section{Le \tho de Barr}

Le \tho fondamental des théories dynamiques (finitaires) \ref{thFond}
est une base solide pour le décryptage \cof des \demos classiques.
En \clama on démontre qu'une \tco prouve une \rdy en regardant ce qui se passe dans les modèles de la théorie, que l'on étudie avec des outils surpuissants mais douteux tels que le tiers exclu, l'axiome du choix et parfois même toute la puissance de \sa{ZFC}. Or le \thref{thFond}
nous assure que si la règle en question est prouvable dans la théorie formelle avec logique classique, elle est \egmt démontrable par les méthodes \elrs \gui{sans logique} que constituent les preuves dynamiques. 

L'essentiel du décryptage revient donc à vérifier que la \demo classique peut se formaliser en logique du premier ordre classique. Cela n'est pas toujours facile, car après tout la théorie~\sa{ZFC} peut être utilisée pour démontrer des résultats beaucoup plus \gui{étranges} que le \tho de complétude de Gödel, et pourquoi pas des résultats carrément faux si \sa{ZFC} est inconsistante. Mais en pratique, en \clama, même l'usage à outrance des ultrafiltres ou de l'hypothèse du continu semble toujours cacher des arguments plus simples.

\smallskip Un \textbf{\tho de Barr}, établi en \clama (et semble-t-il impossible à démontrer en \coma), dit que pour les \tgms, tout résultat démontré avec la logique classique peut \egmt être démontré avec la logique \cov. Il s'agit d'une \gnn du \thref{thFond} qui se trouve confirmée en pratique, même si l'on n'en a pas de certitude complète du point de vue \cof.
Une étude récente du problème est faite par Rathjen dans
l'article \cite{Rat2016} publié dans le livre~\cite{CPM2016}.

Le \tho de Barr nous donne une bonne raison de penser que le type de décryptage fourni par le \thref{thFond} s'applique aussi pour les \tgms infinitaires, avec les mêmes bémols que ceux que nous avons indiqués pour les théories finitaires. 

\Llec peut trouver des exemples de ce type dans \cite[sections XV-6 et XV-7]{CACM}.\label{inThGeomR}

\smallskip Nous illustrons à la fin de ce chapitre la manière dont il ne faut pas comprendre le \tho de Barr.

%%%%%%%%%%%%%%%%%%%%%%%%%%%%%%%%%%%%%%%%%%%%%%%%%%%%%%%%%%%%%%%%%%%%
\Subsection{Théorie \gmq infinitaire pour l'arithmétique primitive récursive}\label{PrimRecOmega}

On considère maintenant la \tgm infinitaire $\sa{PRA}\omega$ obtenue à partir de la \tgm finitaire \Sa{PRA} en ajoutant l'axiome suivant qui force la sorte $\iN$ à ne contenir que des entiers usuels.

\Regles{
\Lab{Nat} $\Vdi{x:\iN} \Vou_{n\in\NN}\; x=\und n$
}

Dans la théorie ainsi obtenue, on voit que pour démontrer $\vd f(x)=g(x)$ avec $f,g$ de sorte~$F_1$, il suffit de savoir démontrer que pour chaque $n\in\NN$ la règle
$\vd f(\und n)=g(\und n)$ est valide. En effet, on en déduit alors $\,\,x=\und n\vd f(x)=g(x)$ pour chaque entier concret $n$, et par utilisation de la règle~$\tsbf{Nat}$, on obtient $\vd f(x)=g(x)$.

Or pour des fonctions $f$ et $g$ définies dans \sa{PRA} (\cad deux fonctions primitives récursives arbitraires), $f(\und n)$ et $g(\und n)$ sont deux entiers usuels explicites. 
Ainsi deux fonctions primitives récursives sont \gui{prouvablement partout égales} dans la théorie \ssi elles prennent les mêmes valeurs en tout entier, \cad si elles sont concrètement égales, \cad si l'on dispose d'une \demo pour le fait qu'elles sont égales. 

Mais cette \demo n'est pas toujours formalisable à l'intérieur de la théorie (elle peut par exemple utiliser une \recu double). En tout cas, elle est censée être produite en \maths intuitives dans le monde \mathe intuitif des entiers naturels.

Un exemple est fourni par la fonction primitive récursive $C_{PRA}:\NN\to\NN$ qui s'annule partout \ssi la théorie \sa{PRA} et consistante (ce dont on est intimement convaincu). D'après le \tho d'incomplétude de Gödel, la théorie \sa{PRA} ne peut pas démontrer $\vd C_\PRA(x)=0$. Alors qu'elle démontre $\vd C_\PRA (\underline n)=0$ pour tout $n$.

De même il semble probable que la \tgm $\sa{PRA}\omega$, bien qu'elle soit capable de démontrer $\vd C_\PRA(x)=0$, ne puisse cependant pas démontrer $\vd C_\PRA=0_1$.

\smallskip L'axiome infinitaire \tsbf{Nat} que l'on ajoute autorise donc jusqu'à un certain point, mais jusqu'à un certain point seulement, l'utilisation de la $\omega$-règle pour l'\egt des fonctions primitives récursives.

\Subsection{Conclusion}

L'étude que l'on vient de faire jette une ombre sur le \tho de Barr, car il semblerait affirmer dans ce cas que toute fonction primitive récursive prouvée nulle en \clama serait prouvablement nulle en \coma. 

Il se peut cependant que la \demo intuitive \du mathematic\ien classique utilise des principes douteux, comme ceux formalisés dans la théorie \sa{ZF}, elle peut alors aboutir à des résultats tout à fait contestables. 
 
Par exemple considérons la fonction primitive récursive $\mathrm{consisZ}$ qui prend toujours la valeur~$0$ jusqu'au moment où éventuellement on trouve une preuve de $0=1$ dans \sa{Z} et la fonction prend alors la valeur~$1$.
 Ainsi $\mathrm{consisZ}$ est une constante de type $\iF_1$ bien définie de \sa{PRA}.

En \clama avec une théorie intuitive des ensembles suffisamment forte (par exemple \sa{ZF}), on peut démontrer que $\mathrm{consisZ}$ est identiquement nulle. Et cela démontre en \clama (en utilisant \sa{ZF} dans le monde \mathe intuitif extérieur) la règle
$\vd \mathrm{consisZ}(x)=0$ dans la \tgm $\sa{PRA}\omega$. 

Or ce résultat échappe manifestement à toute \prco. Et il ne peut y avoir de \prco de la règle $\vd \mathrm{consisZ}(x)=0$ dans la \tgm $\sa{PRA}\omega$. 

\smallskip En fait il faut préciser l'énoncé du \tho de Barr. Il ne dit pas que le cadre des \clama est conservatif pour les \prts \gmqs dans une \tgm. Il dit seulement que
lorsqu'on utilise les mêmes mathématiques à l'extérieur de la théorie formelle infinitaire, l'ajout des connecteurs, des quantificateurs et des règles de la logique classique à l'intérieur de la \tgm infinitaire, ne permet pas de démontrer de nouvelles \prts formulables comme règles \gmqs.

\smallskip \hum{Il faudra ici un énoncé précis et compréhensible du \tho de Barr.
Ce serait parfait si Thierry expliquait clairement ce que signifie vraiment le \tho de Barr et dans quelle mesure on peut lui accorder crédit.

Par ailleurs, La théorie formelle \sa{Peano} admet elle-même une version \tdy infinitaire. Il faudra en toucher deux mots, ainsi que pour des extensions où l'ordinal $\vep_0$ est remplacé par un ordinal dénombrbale plus grand.
}

%%%%%%%%%%%%%%%%%%%%%%%%%%%%%%%%%%%%%%%%%%%%%%%%%%%%%%%%%%%%%%%%%%%%
%%%%%%%%%%%%%%%%%%%%%%%%%%%%%%%%%%%%%%%%%%%%%%%%%%%%%%%%%%%%%%%%%%%%
%\section{Décryptage de \demos classiques}\label{secDecrypt}
%: Subsection{Décryptage de \demos classiques}

%\newpage \thispagestyle{empty}
%
%
%\chapter*{Conclusion}
%\addstarredchapter{Conclusion}

\part{Théories \gmqs finitaires pour l'algèbre réelle}

\chapter*{Introduction}
\addstarredchapter{Introduction}

Cette deuxième partie est consacrée à la mise au point d'une \tdy finitaire qui a pour ambition de décrire de manière exhaustive les \prts \agqs du corps des réels, et plus \gnlt d'un corps réel clos \textsl{non} discret\footnote{Dans ce texte, une négation est mise en italique lorsque l’affirmation correspondante, vraie en \clama, implique en \coma un principe non constructif bien répertorié, tel que \tsbf{LPO}
ou même \tsbf{MP}.}, du moins celles qui sont exprimables dans un langage restreint, proche du langage des anneaux ordonnés. Cela constitue un développement, avec quelques modifications terminologiques mineures, des idées données dans l'article \cite{LM2017}.
L'axiome d'archimédianité, introduit dans le dernier chapitre, nous fait quitter le domaine des \tgms finitaires.

\smallskip Le chapitre \ref{chapcoo} propose une première définition de la structure de corps ordonné en l'absence d'un test de signe. Il discute aussi la possibilité d'une axiomatique convenable pour les \crcs \textsl{non} discrets, comme le corps des nombres réels.

 La section \ref{seccodi} donne quelques rappels sur la \tdy des \codis et celle des \crcds.
 
 La section \ref{secPstFormels} discute les conséquences, décisives, des \Psts formels, dans notre cadre de travail.

 La section \ref{secConondisc} propose une axiomatique pour les corps ordonnés \textsl{non} discrets (\dfn \ref{defiConondiscret}). 

La section \ref{condna}
décrit un exemple de corps ordonné de Heyting \textsl{non} discret non archimédien.

 La section \ref{secCRCnondis} donne une première discussion sur les
\tdys acceptables pour $\RR$ en tant que corps réel clos \textsl{non} discret.

\smallskip Le chapitre \ref{chap-afr} traite des \tdys qui admettent des extensions \eseqs à la théorie \Sa{Co}. 

On commence (section \ref{sectrdisad}) par la théorie des \trdis (un corps ordonné \textsl{non} discret est un \trdi pour sa relation d'ordre). 

Dans la section \ref{secgrl} on traite les \grls (théorie \peq), valable pour l'addition sur les réels (ce sont les $\ell$-groups dans la littérature anglaise).

Ensuite (section \ref{secfrings}) nous passons aux \afrs ($f$-rings dans la littérature anglaise), théorie inspirée par les anneaux de fonctions réelles continues. 

La section \ref{secArftr} décrit des \tdys dans lesquelles on ajoute dans la signature la relation $\cdot>0$ (\asrs et variantes).

La section \ref{secCOG} propose un retour sur la théorie \Sa{Co} en la confrontant à des extensions convenables de la théorie des \asrs.

\smallskip Le chapitre \ref{chapreelclos} propose une \dfn la structure de corps ordonné réel clos en l'absence d'un test de signe.

La section \ref{subsecclotrlRR}
 explique comment introduire des racines carrées des \elts~\hbox{$\geq 0$} dans un corps ordonné \textsl{non} discret. Ceci à titre de mise en jambes pour introduire la notion plus \gnle de racines virtuelles.

 La section \ref{secCoVR} introduit les fonctions racines virtuelles et certaines \tdys correspondantes: notamment les \afrvrs et
les \covrs \textsl{non} discrets,

 La section \ref{secArc} traite les \arcs et la section \ref{secCrc2} propose une \dfn pour les \crcs \textsl{non} discrets comme \arcs locaux. La théorie des anneaux réels clos est ici présentée sous une forme \elr, \peq, dans le style de \cite{Tre2007}.

\smallskip Le chapitre \ref{secGeomReelsArchi} aborde une \tgm infinitaire où l'on ajoute l'axiome selon lequel le corps des nombres réels est \textsl{archimédien}.

\smallskip Ainsi, nous proposons pour la \tdy convoitée la structure d'anneau réel clos local (éventuellement archimédien si cela s'avérait utile).

\smallskip Dans un chapitre de conclusion, nous résumons la situation à laquelle nous avons abouti, en précisant les questions importantes, d'un point de vue \cof, auxquelles nous ne savons pas aujourd'hui répondre de manière satisfaisante. 
\newpage \thispagestyle{empty}

\chapter{Corps ordonnés} \label{chapcoo}
\Today
\minitoc
%%%%%%%%%%%%%%%%%%%%%%%%%%%%%%%%%%%%%%%%%%%%%%%%%%%%%%%%%%%%%%%%%%%%

\section*{Introduction}
\addtocontents{toc}{\vskip0.8em}
\addcontentsline{toc}{section}{Introduction}
\rdb

Ce chapitre propose une première approche \cov de la théorie classique des corps réels clos. En fait, la théorie classique ne s'applique qu'aux \crcs pour lesquels on dispose d'un test de signe sur n'importe quel \elt du corps, s'il est donné conformément à la \dfn. Autrement dit, la théorie classique usuelle est une théorie des \crcs \und{discrets}. 
Mais il est bien connu que l'analyse numérique matricielle, utilisée dans les applications de la théorie à des situations concrètes, n'utilise jamais un tel test de signe. Une approche \cov pour une théorie des \prts \agqs du corps des réels nécessite une \tdy des \crcs \und{\textsl{non} discrets}. 

La section \ref{seccodi} donne quelques rappels sur la \tdy des \codis et celle des \crcds.

La section \ref{secPstFormels} explique la grande utilité des \Psts formels.

La section \ref{secConondisc} propose une axiomatique pour les corps ordonnés \textsl{non} discrets (\dfn \ref{defiConondiscret}). Nous devons abandonner les axiomes d'ordre total dans leur formulation discrète usuelle et les remplacer par des \rdys pertinentes pour $\RR$. On constate alors que de nombreuses fonctions rationnelles continues bien définies sur $\QQ$, comme la fonction $\mathrm{sup}$, doivent être introduites dans le langage.

La section \ref{condna}
décrit un exemple de corps ordonné de Heyting \textsl{non} discret non archimédien.

La section \ref{secCRCnondis} donne une première discussion sur les
\tdys acceptables pour~$\RR$ en tant que corps réel clos \textsl{non} discret.
Nous sommes guidés par le principe de prolongement par continuité \ref{thRRcomplet} qui peut être vu comme une version \agq de la complétude de $\RR$. Dans ce cadre le \thref{thParamcontFsagc0} joue un rôle fondamental pour une \dfn pertinente des \fsagcs, en ramenant la \dfn au cas des \fsagcs dans le cadre discret de $\RRa$.

%\hum{développer un peu cette introduction}

\section{Rappels sur les corps ordonnés discrets}\label{seccodi}

Un \textsl{\codi} est un \cdi $\gR$ muni d'une relation d'ordre convenable, laquelle peut être décrite par la donnée de la partie $P$ formée par les \elts $\geq 0$. Cette partie $P$ doit satisfaire les \prts suivantes.
\begin{enumerate}
\item $P\cup -P=\gR$.
\item $P\cap -P =\so0$.
\item $P+P:=\sotq{x+y}{x,y\in P}\subseteq P$.
\item $P\cdot P:=\sotq{x\cdot y}{x,y\in P}\subseteq P$.
\item $\gK^2:=\sotq{x^2}{x\in \gR}\subseteq P$.
\end{enumerate}
La relation d'ordre est alors définie par  $x\geq y\equidef x-y\in P$, et la relation d'ordre strict est définie par $x>0\equidef x\geq0 \vii x\neq 0$. 

D'un point de vue \cof, le \gui{ou} caché dans le $\cup$ du point~1 doit être explicite. On a donc explicitement la trichotomie \fbox{$x=0\vuu x> 0\vuu x< 0$}.

En termes de \tdy, tout ceci est traduit par des prédicats et des lois données dans le langage de \Sa{Aso}. En vue de démontrer le \Pst \cot, on a choisi des axiomes convenablement mis en ordre, en commençant par les règles directes.

%:Subsection Une \tdy naturelle pour les \codis
\Subsection{Une \tdy naturelle pour les \codis}\label{subseccodi}

On rappelle ici la \tdy des \textsl{\codis} \SA{Cod} donnée dans \cite{CLR01}.\index{corps ordonné!discret}%

\vspace{-.2em}
\Sigt{\Aso}{\cdot=0,\cdot\geq 0,\cdot>0\mathrel{;}\cdot+\cdot, \cdot\times\cdot,-\cdot, 0,1} \label{NOTASigAso}

\vspace{-.2em}
Si l'on veut donner \textsl{un \codi dynamique}, i.e. une \sad de type \Sa{Cod}, on ajoute à la signature une présentation
par \gtrs et relations. 
%de la \sad considérée. 
Par exemple cela peut être la \pn vide, ou un ensemble dénombrable de \gtrs, sans aucune relation, ou encore cela peut être basé sur une structure \agq existante dans laquelle on demande de préserver certaines relations, par exemple toutes les relations d'égalité entre termes construits sur les \elts de la structure. Ainsi tout anneau définit un \codi dynamique.

\bni {\bf Abréviations}

\vspace{-1em} \DeuxCols{
\begin{itemize}
\itbu $x\neq0 $ signifie $x^2>0$
\itbu $x = y $ signifie $ x - y = 0$
\itbu $x > y $ signifie $ x - y > 0$
\itbu $x < y $ signifie $ y > x$
\end{itemize}
}
{\begin{itemize}
\itbu $x \geq y$ signifie $x-y\geq 0$ $\phantom {x^2>0}$
\itbu $x \neq y $ signifie $ x - y \neq 0$
\itbu $x \leq y $ signifie $ y \geq x$
\end{itemize}
} 

\bni {\bf Axiomes}

\medskip \noindent {\sl Règles directes}

\smallskip On a d'abord mis les axiomes des anneaux commutatifs, puis les règles qui concernent~\hbox{$\cdot=0$} \hbox{et $\cdot\geq 0$}, ensuite 
les règles qui font intervenir $\cdot>0$\footnote{Les règles \Tsbf{ga0} et \Tsbf{ga1} ont été introduites sous les noms \Tsbf{ac0} et \Tsbf{ac1} dans l'exemple \ref{exaAc}.}.

\DeuxRegles{
\laB{ga0} $\vd 0=0$
\laB{ac2} $\,\,x=0\vd xy=0 $
}
{
\laB{ga1} $\,\, x=0\vet y=0\vd x+y=0$
}

\DeuxRegles{
 \laB{gao0} $\vd 0 \geq 0$
 \Lab{ao1} $\vd x^2 \geq 0$
}
{
 \laB{gao1} $\,\, x \geq 0\vet y \geq 0 \vd x + y \geq 0$
 \Lab{ao2} $\,\, x \geq 0\vet y \geq 0 \vd x y\geq0$ 
}

\DeuxRegles{
 \Lab{aso1} $ \vd 1> 0$ \phantom{$x^2$}
 \Lab{aso2} $\,\, x> 0 \vd x \geq 0$
}
{
 \Lab{aso3} $\,\, x > 0\vet y \geq 0 \vd x + y > 0$
 \Lab{aso4} $\,\, x > 0\vet y > 0 \vd xy > 0$
}

\medskip \noindent {\sl Collapsus} 

\Regles{
\lAb{col\igt} $\,\,0> 0 \vd 1=0$ \label{Axcoligt}
}

%\qquad\qquad (i.e. $\,\,0\neq0\vd1=0$)

\medskip \noindent {\sl Règles de simplification}

\DeuxRegles{
\laB{Gao} $\,\, x\geq 0\vet x\leq 0 \vd x = 0$
}
{
\Lab{Iv} $\,\, xy = 1 \vd x\neq  0$
}

\medskip \noindent {\sl Règles dynamiques}

\DeuxRegles{
\Lab{IV} $\,\, x\neq 0 \vd \Exists y\; xy = 1$

\item [\Edinq]  $\vd x=0 \vou \,x\neq 0$
}
{
\Lab{OT} $ \vd x \geq 0 \vou x\leq 0$
}

\smallskip La \tdy \SA{Crcd} des \textsl{corps réels clos discrets} est obtenue à partir de la théorie \Sa{Cod} en ajoutant comme axiomes les \rdys \tsbf{RCF$_n$}\footnote{Un \tho essentiellement équivalent à ces règles est démontré par Bishop pour le corps $\RR$, mais en utilisant l'axiome du choix dépendant.}.\index{corps réel clos!discret}%

\Regles{
\lAb{RCF$_n$} {$ \,\, a< b\vet P(a)P(b)<0 \vd \Exists x\, \big(P(x)=0,\,a<x<b\big) $ \quad ($ P(x)=\sum_{k=0}^na_kx^k $) }\label{RCFn}}

\medskip Les règles \tsbf{gao0} et \tsbf{gao1} expriment, dans le contexte des groupes, la réflexivité et la transitivité de la relation d'ordre (compatible avec la loi de groupe). 
La règle \tsbf{Gao} correspond à l'antisymétrie pour la relation d'ordre. 

%%%%%%%%%%%
\smallskip Les règles \Edinq\ et \Tsbf{OT} expriment que l'\egt est discrète et l'ordre total. Elles ne sont pas satisfaites \cot pour~$\RR$. 
Pour les réels de Bishop, la règle \Ednq\ équivaut au principe d’omniscience \tsbf{LPO} et la règle \tsbf{OT} 
équivaut au principe~\tsbf{LLPO}. Notons aussi que le principe \gui{tout \elt régulier de $\RR$ est inversible} équivaut au principe de Markov\footnote{Équivalence suggérée par Fred Richman.} \tsbf{MP}.

\smallskip Vue la forme \gui{sans négation} adoptée ici pour le collapsus, l'anneau trivial est un corps ordonné discret, et l'axiome de collapsus \tsbf{col$_>$} est une conséquence de \tsbf{IV}.

\smallskip Au moyen des seules règles directes on voit que $\,1=0\vd (x=0\vet x\geq 0\vet x>0)$. Cela justifie de prendre $1=0$ comme un substitut de $\Bot$.

%%%%%%%%%%%%%%%%%%%%%%%%%%%%%%%%%%%%%%%%%%%%%%%%%%%%%%%%%%%%%%%%%%%%
\Subsection{Quelques règles dérivées dans \Sa{Cod}}
%: Subsection{Quelques règles dérivées}

 {\sl Quatre règles de simplification valides}
\label{regsimpcod}

\DeuxRegles{
\laB{Anz} $\,\, x^2= 0 \vd x = 0$
\Lab{Aonz} $\,\, c\geq 0\vet x(x^2+c)\geq 0 \vd x\geq 0$
}
{
\Lab{Aso1} $\,\, x> 0\vet xy\geq 0 \vd y\geq 0$
\Lab{Aso2} $\,\, x\geq 0\vet xy> 0 \vd y> 0$
}

\medskip \noindent {\sl Deux \rdys valides}

\DeuxRegles{
\Lab{OTF} $\,\, x+y> 0 \vd x > 0\vou y>0$
}
{
\lAb{OTF$\eti$} $\,\, xy< 0 \vd x <0 \vou y<0 $\label{AxOTFx}
}

\smallskip \noindent Notons que la règle \tsbf{Aso1} implique qu'un \elt $>0$ est régulier.

%l
%: Lemma{lemthRRco}
\begin{lemma} \label{lemthRRco} La règle suivante est prouvable
avec les axiomes directs.

\Regles{\Lab{Aonz2} $\,\, c\geq 0\vet x(x^2+c)\geq 0\vet x<0 \vd 0> 0$}
 
\end{lemma}
%----------- fin lemma ----------------------------------- 

%: Theorem{thRRco}
\begin{theorem} \label{thRRco}
Hormis les règles \Edinq\ et \Tsbf{OT}, toutes les règles énoncées précédemment sont valides \cot pour $\RR$, sans utilisation de l'axiome du choix dépendant.
\end{theorem}
%----------- fin theorem -----------------------------
%
\begin{proof}
Tout est clair sauf peut-être la règle \Tsbf{Aonz}.
Pour $x\in\RR$, on peut démontrer $x\geq 0$ en réduisant à l’absurde
$x<0$. Or, c'est ce que fait \Tsbf{Aonz2} (prendre $c=0$).
\end{proof}
%

%r
%: Remark{remOTFx}
\begin{remark} \label{remOTFx} 
 Les règles  \OTFx\ et \Edinq\ impliquent la règle
$\vd x=0 \,\vou\, x<0\,\vou\, x>0 $, et à fortiori~\tsbf{OT}. 
Voir aussi~\ref{lemAfrsdz}, \ref{lemAsrs} et \ref{lemArftr}. 
\eoe\end{remark}
%----------- fin remark ---------------------------------- 

Si $\gK$ est un \codi, on note $\Sa{Cod}(\gK)$ la \sad~de type \Sa{Cod} ayant pour \pn le \textsl{diagramme positif de $\gK$}. Un modèle non trivial de $\Sa{Cod}(\gK)$ est un \codi non trivial $\gL$ donné avec un morphisme $\gK\to\gL$.
De même, si $\gA$ est un anneau commutatif (ou un anneau ordonné), on note $\Sa{Cod}(\gA)$ la \sad~de type \Sa{Cod} ayant pour \pn le \textsl{diagramme positif de $\gA$}. Un modèle de $\Sa{Cod}(\gA)$ est un \codi $\gL$ donné avec un morphisme $\gA\to\gL$ (d'anneau commutatif, ou d'anneau ordonné).

\Subsection{Théories dynamiques plus faibles} 
%: Subsection{Théories dynamiques plus faibles}

La règle \Tsbf{Aonz} implique $x^3\geq 0\vd x\geq 0$, donc aussi,
sous \Tsbf{Gao},
$x^3= 0\vd x= 0$, et \textsl{à fortiori}~\Tsbf{Anz}.

%: Definition{defisaAor}
\begin{definition} \label{defisaAor}~
\\
Théories basées sur le langage des anneaux ordonnés\index{anneau!ordonné}
\sigt{\Ao}{\cdot=0,\cdot\geq 0\mathrel{;}\cdot+\cdot, \cdot\times\cdot,-\,\cdot,0,1}. \label{NOTASigAo}
\begin{enumerate}
\item La \talg \SA{Apo} des \textsl{anneaux préordonnés}.\index{anneau!préordonné}
Les axiomes sont ceux des anneaux commutatifs et les \reds \Tsbf{gao0}, \Tsbf{gao1}, \Tsbf{ao1}, \Tsbf{ao2}.%
\item La \talg \SA{Ao} des \textsl{anneaux ordonnés}.
Les axiomes sont ceux des anneaux préordonnés et la \rsim \Tsbf{Gao}.%
\item La\index{anneau!ordonné!strictement réduit} \talg \SA{Aonz} des \textsl{anneaux ordonnés strictement réduits\footnote{La règle \tsbf{Aonz} est plus forte que la règle \tsbf{Anz}, on a donc mis \gui{strictement réduit} plutôt que \gui{réduit}. Voir cependant le point \textsl{3} du lemme \ref{lemAtonz}. }} est obtenue en ajoutant la \rsim \Tsbf{Aonz} à la théorie \sa{Ao}.%
\item La \tdy \SA{Ato} des \textsl{anneaux totalement ordonnés} est obtenue en ajoutant la \rdy \Tsbf{OT} à la théorie \sa{Ao}.%
\item La \tdy \SA{Atonz} des \textsl{anneaux totalement ordonnés réduits} est obtenue en ajoutant la \rdy \Tsbf{Anz} à la théorie \sa{Ato}. Cette théorie
prouve les règles \tsbf{Aonz} et \tsbf{ASDZ}
(lemme~\ref{lemAtonz})%
\end{enumerate}
Théories basées sur le langage des anneaux \stm ordonnés:
on ajoute le prédicat $\cdot> 0$.
\begin{enumerate}\setcounter{enumi}{5}
\item La théorie directe \SA{Apro} des \textsl{anneaux proto-ordonnés} 
(cf. \cite{CLR01}).
Les axiomes sont ceux des anneaux commutatifs, toutes les \reds énoncées pour \Sa{Cod} (\Tsbf{gao0}, \Tsbf{gao1}, \Tsbf{ao1}, \Tsbf{ao2}, \Tsbf{aso1} à \Tsbf{aso4}) et le collapsus \coligt.%
\item La \talg \SA{Aso} des \textsl{anneaux \stm ordonnés} 
est la théorie \sa{Apro} à laquelle on ajoute les \rsims \Tsbf{Gao}, \Tsbf{Aso1} et \Tsbf{Aso2}. On peut la voir aussi comme construite à partir de \sa{Ao} en ajoutant le prédicat $\,\cdot> 0\,$ dans le langage, les \reds \Tsbf{aso1} à \Tsbf{aso4} et les \rsims \Tsbf{Aso1} et \Tsbf{Aso2}.\index{anneau!strictement ordonné}\index{anneau!ordonné!strictement ---}%
\item La \talg \SA{Asonz} des \textsl{anneaux \stm ordonnés réduits} (\gui{quasi-ordered rings} dans \cite{CLR01})
 est obtenue en ajoutant la \rsim \Tsbf{Aonz} à~\sa{Aso}. On peut aussi la voir comme la théorie \sa{Apro} à laquelle on ajoute les \rsims \Tsbf{Gao}, \Tsbf{Aonz}, \Tsbf{Aso1} et \Tsbf{Aso2}.%
\item La \tdy \SA{Asto} des \textsl{anneaux \stm totalement ordonnés} 
est la théorie \sa{Aso} à laquelle on ajoute la \rdy \Tsbf{OT}. On peut la voir aussi comme construite à partir de \sa{Ato} en ajoutant le prédicat $\cdot> 0$ dans le langage, les \reds \Tsbf{aso1} à \Tsbf{aso4} et les \rsims \Tsbf{Aso1} et \Tsbf{Aso2}.%
\item La \tdy \SA{Aito} des \textsl{anneaux intègres totalement ordonnés} 
 est obtenue en ajoutant la \rdy \ \tsbf{ED$_>$} \ à \sa{Asto}. Les règles 
 \tsbf{Aonz}, \tsbf{OTF} et \tsbf{OTF$\eti$} sont valides dans cette théorie.%
\end{enumerate} 
\end{definition}
%--------- fin definition --------------------------------

Dans les points 6, 7 et 8, la signification de $x>0$ n'est pas fixée à priori.
Cela peut aller de \gui{$x$ est régulier et $\geq 0$} jusqu'à \gui{$x$ est inversible et $\geq 0$}

La théorie directe \Sa{Apro} est celle dans laquelle le collapsus est le plus clair, directement donné par un certificat \agq, comme précisé dans le lemme qui suit. 

Rappelons que dans un anneau, un \textsl{cône}\index{cone@cône}
est une partie $C$ qui contient les carrés et qui est stable par addition et produit: $C+C\subseteq C$, $C\times C\subseteq C$.

%l
%: Lemma{lemColApo}
\begin{lemma}[certificat algébrique d'effondrement] \index{effondrement}\label{lemColApo}~\\
Soit $\gK$ une structure algébrique dynamique de type \Sa{Apro} donnée par une présentation $(G; R_{>0}, R_{\geq 0}, R_{=0})$ avec la signification suivante: 
$G$ est l'ensemble des générateurs de la structure, $R_{>0}$, $ R_{\geq 0}$ et $R_{=0}$ sont trois parties de $\ZZ[G]$, les éléments de $R_{>0}$ (resp. $R_{\geq 0}$, $R_{=0}$) sont supposés~$> 0$ (resp. $\geq 0$, $= 0$) dans $\gK$. \\
La structure algébrique dynamique $\gK$ s’effondre si, et seulement si, on a dans $\ZZ[G]$ une égalité
\[\fbox{$s+p+z=0$}\]
où 
 $s$ est dans le \mo multiplicatif engendré par $R_{>0}$,
 $p$ est dans le cône engendré par $R_{>0}\cup R_{\geq 0}$, et 
 $z$ dans l'idéal engendré par $R_{=0}$. 
\end{lemma}
%----------- fin lemma ----------------------------------- 

%l
%: Lemma{lemAtonz}
\begin{lemma}[voir aussi le lemme \ref{lemAfrsdz}] \label{lemAtonz} ~
\begin{enumerate}
\item La \tdy \sa{Aso} prouve la \rsim \Tsbf{Iv} et le collapsus \coligt.
\item \smallskip La \tdy \sa{Ato} prouve les \rsims suivantes. 

\DeuxRegles{
\Lab{Ato1} $\,\, x\geq 0 \vet xy=1\vd y\geq 0$
}
{
\Lab{Ato2} $\,\, c\geq 0\vet x(x^2+c)\geq 0 \vd x^3\geq 0$
}

\item La théorie \Sa{Atonz} prouve les règles 
 \Tsbf{Aonz} et \Tsbf{ASDZ}.
\end{enumerate}
 
\end{lemma}
%----------- fin lemma ----------------------------------- 
%
\begin{proof}
Les points \textsl{1} et \textsl{2} sont faciles (utiliser \Tsbf{Gao}). 
Pour le point \textsl{3}, \tsbf{Aonz} découle de~\tsbf{Ato2}. 
Voyons \tsbf{ASDZ}. Soient $a,b$ tels que $ab=0$; si $\abs a\leq \abs b$, 
on a $0\leq {\abs a}^2\leq \abs {ab}=0$ donc $a=0$, 
et dans le cas \hbox{où $\abs b\leq \abs a$} on obtient $b=0$.
\end{proof}

Dans la théorie \sa{Aso}, la règle \tsbf{Ato1} est une variante affaiblie de \Tsbf{Aso2}. 
Dans la théorie \sa{Ao}, la règle \tsbf{Ato2} est une variante affaiblie de \Tsbf{Aonz}.

\smallskip Dans un anneau totalement ordonné, si l'on définit \gui{$x>0$}
par \gui{$x$ est \ndz \hbox{et $\geq 0$}}, toutes les règles qui définissent \Sa{Aso} sont satisfaites (et $x\neq 0$ est à priori plus fort que la simple négation de $x=0$). Cela explique l'intérêt de la \talg \Sa{Aso}. On définira plus loin une variante réticulée, la théorie \Sa{Asr} des \asrs.

%%%%%%%%%%%%%%%%%%%%%%%%%%%%%%%%%%%%%%%%%
\Subsection{Un exemple avec des nilpotents}
%: Subsection{Un exemple avec des nilpotents}

%: Example{exatotordnonreduit}
\begin{example}[un anneau totalement ordonné non réduit] \label{exatotordnonreduit}
L'exemple que nous donnons maintenant est celui que l'on doit garder en tête pour bien comprendre la différence entre les anneaux totalement ordonnés et les anneaux totalement ordonnés intègres.

\noindent Il s'agit de l'anneau totalement ordonné $\QQ[\alpha]$ où $\alpha>0$ et $\alpha^6=0$ ($\alpha$ est un infinitésimal~\hbox{$>0$} nilpotent).
Soit $c$ un \elt tel que $c^2=0$ (par exemple $c=\alpha^5$). Le système de contraintes 

\snic{x^2=c^2,\;x\geq 0,}

\noindent qui pourrait être suggéré pour caractériser $\abs c$ sans utiliser le test de signe dans le cas d'un corps ordonné, admet maintenant une infinité de solutions: tous les $r\alpha^3+y\alpha^4$ où $r>0$ dans $\QQ$ et $y$ arbitraire dans $\QQ[\alpha]$. \eoe
\end{example}
%--------- fin example ---------------------------------------- 

%%%%%%%%%%%%%%%%%%%%%%%%%%%%%%%%%%%%%%%%%

%%%%%%%%%%%%%%%%%%%%%%%%%%%%%%%%%%%%%%%%%%%%%%%%%%%%%%%%%%%%%%%%%%%%
%:Subsection Ajout du symbole de fonction $\vu$
\Subsection{Ajout du symbole de fonction $\vu$ pour le sup}\label{secCodisup}

Dans un ensemble totalement ordonné et à fortiori dans un \codi toute paire d'\elts admet une borne supérieure: le plus grand des deux. On ne change donc rien d'essentiel à la théorie \sa{Cod} en ajoutant un symbole fonctionnel
$\cdot\vu\cdot$ soumis au trois axiomes qui définissent le sup de deux \elts, quand il existe pour une relation d'ordre donné.

%: Definition{deficodisup}
\begin{definition} \label{deficodisup}
 La \tdy des \textsl{\codis avec sup}, notée \SA{Codsup}, est la \tdy des \codis à laquelle on ajoute un symbole de fonction~$\cdot\vu\cdot$ et pour axiomes les \ralgs \tsbf{sup1}, \tsbf{sup2} et \tsbf{Sup} suivantes.\index{corps ordonné!discret avec sup}%

\DeuxRegles{
\laB{sup1} $ \vd x\vu y\,\geq x $
\laB{Sup} $ \,\, z\geq x\vet z\geq y\vd z\geq x\vu y$
}
{
\laB{sup2} $ \vd x\vu y\,\geq y $
}

\noindent On définit de la même manière à partir des théories \Sa{Ato}, \Sa{Asto}, \Sa{Aito} et \Sa{Crcd} les théories \SA{Atosup}, \SA{Astosup}, \SA{Aitosup} et \SA{Crcdsup}.
\end{definition}
%--------- fin definition --------------------------------

Dans le cas des \codis et des \crcds on aurait pu aussi remplacer la \rsim \Tsbf{Sup} par la \red \tsbf{sup} suivante, de manière à n’ajouter que des \reds à \Sa{Cod}.

\Regles{
\laB{sup} $ \vd \big((x\vu y)- x\big)\,\big((x\vu y)- y\big)=0$
}

 On peut aussi voir la théorie \Sa{Crcdsup} comme la théorie \Sa{Codsup} à laquelle on ajoute les axiomes de clôture réelle \RCFn.
 
%r
%: Remark{rem}Crdsup
\begin{remark} \label{remCrdsup} 
 Les théories \Sa{Cod} et \Sa{Codsup} sont \esids. Même chose pour les autres couples de théories dans la \dfn \ref{deficodisup}. 
\eoe\end{remark}
%----------- fin remark ---------------------------------- 

%%%%%%%%%%%%%%%%%%%%%%%%%%%%%%%%%%%%%%%%%%%%%%%%%%%%%%%%%%%%%%%%%%%%
\section{Positivstellensätze formels}\label{secPstFormels}
%: Subsection{Positivstellensätze formels}

Le Positivstellensatz formel des \clama (\cite[Theorem~4.4.2]{BCR}) admet la version \cov suivante (voir \cite{CLR01}).

%: pstf{Pst1}
\begin{pstf}[Positivstellensatz formel pour les corps ordonnés, 1] \label{Pst1}\index{Positivstellensatz!formel!pour les corps ordonnés, 1}
On considère les \sads sur le langage des anneaux strictement ordonnés. 
\begin{enumerate}
\item Les \tdys \Sa{Apro}, \Sa{Aso}, \Sa{Cod} et 
\Sa{Crcd} s'effondrent simultanément. 
\item Les \tdys \Sa{Asonz}, \Sa{Aito}, \Sa{Cod} et \Sa{Crcd} prouvent les mêmes \ralgs.
\end{enumerate}
\end{pstf}
%--------- fin pstf ----------------------------------- 

La force de ce \tho tient à ce que l'effondrement d'une \sad de type \Sa{Apro} est donné par un \textsl{certificat \agq d'effondrement} (lemme \ref{lemColApo}), que l'on appelle un \textsl{Positivstellensatz}.%
\index{Positivstellensatz}\index{certificat algébrique}

Comme cas particulier, pour un anneau commutatif $\gR$, la \sad $\Sa{Crcd}(\gR)$ s'effondre \ssi $-1$ est une somme de carrés dans $\gR$.

En \clama, grâce au \tho de complétude de Gödel \ref{thGodel1}, on déduit du point~\textsl{1} précédent le \textsl{Positivstellensatz formel abstrait} sous la forme suivante (voir~\cite{CLR01}).

\smallskip\noindent \textsl{Un \sys de conditions de signes imposées à des \elts d'un anneau $\gA$ admet un certificat \agq d'impossibilité \ssi le seul modèle de $\Sa{Cod}(\gA)$ est trivial, \ssi le seul modèle de $\Sa{Crcd}(\gA)$ est trivial.\index{Positivstellensatz!formel!abstrait}} 

\smallskip Le point \textsl{2} de \ref{Pst1} admet \egmt
des versions classiques \gui{abstraites} via la théorie des modèles, en application du \thref{thcolsimralg} (voir \cite{CLR01}).

\smallskip Nous examinons maintenant ce que deviennent les résultats précédents en l'absence du prédicat~\hbox{\gui{$\cdot>0$}} dans la \pn d'une \sad.

%: pstf{Pst1bis}
\begin{pstf}[Positivstellensatz formel, 1bis]\label{Pst1bis}%
\index{Positivstellensatz!formel!pour les corps ordonnés, 1bis}~
\\
On considère les \sads sur le langage des anneaux ordonnés.
\begin{enumerate}
\item Les \tdys \Sa{Ao}, \Sa{Apro}, \Sa{Ato}, \Sa{Cod} et 
\Sa{Crcd} s'effondrent simultanément. 
\item Les \tdys \Sa{Aonz}, \Sa{Atonz}, \Sa{Cod} et \Sa{Crcd} prouvent les mêmes \ralgs.
\end{enumerate}
NB. Le langage utilisé pour la \pn ne doit pas mentionner pas le prédicat $\cdot>0$. Il n'y a pas d'axiome de collapsus dans \Sa{Ao} et \Sa{Ato}, et l'on dit qu'une \sad de type \Sa{Ao} ou \Sa{Ato} s'effondre quand elle prouve $1=0$.

\smallskip \hum{On pourrait ajouter dans le point \emph{1} la théorie \Sa{Apo} en y définissant le collapsus par $-1\geq 0$.} 
\end{pstf}
%--------- fin pstf -----------------------------------
 
%
\begin{proof}

\noindent \textsl{1}.
On considère une \sad $\gA$ pour \Sa{Ao}. La même \pn donne une 
\sad~$\gA'$ pour \Sa{Apro}. Supposons que $\gA'$ prouve $1=0$.
Le collapsus pour~$\gA'$ (de type~\Sa{Apro}) a la forme d'un certificat \agq (un Positivstellensatz) bien précis. Ce certificat s'écrit $1+p=0$, où $p$ est $\geq 0$ en vertu de la \pn et des seuls axiomes de \Sa{Ao}. On a donc à la fois $1\geq 0$ et $1\leq 0$ dans $\gA$, donc $1=0$ en vertu de \Tsbf{Gao}.
Inversement, si $\gA$ prouve~\hbox{$1=0$}, alors à fortiori il en sera de même pour $\gA'$.
\\ Enfin, on applique le \pstref{Pst1} et l'on note que
 \Sa{Ato} est une théorie intermédiaire entre \Sa{Ao} et \Sa{Cod}.

\sni\textsl{2}. On considère une \sad $\gA$ pour \Sa{Aonz}. Il suffit de prouver le résultat pour un fait de la forme $x\geq 0$ (car $x=0$ équivaut à $x\geq 0$ et $x\leq 0$). Ce fait est valide dans~\Sa{Cod} \ssi le fait $x<0$ fait s'effondrer la \sad $\Sa{Cod}(\gA)$.
Cela correspond d'après le \pstref{Pst1} à un certificat \agq de la forme $x^{2n}+p=xq$, où~$p$ et~$q$ sont $\geq 0$
en vertu de la \pn et des seuls axiomes de \Sa{Ao}. Cela donne
$x(x^{2n}+p)\geq 0$, \hbox{puis $x^k(x^{2k}+p_1)\geq 0$} pour un $k$ impair convenable. La règle \Tsbf{Aonz} nous dit
que $x^k\geq 0$ dans $\gA$. Cette même règle montre \hbox{que $x^3\geq 0$} implique~\hbox{$x\geq 0$}, et par suite $x^k\geq 0$ implique $x\geq 0$ pour tout~$k$ impair.
\end{proof}

Une conséquence du point \textsl{1} en \clama (via le \thref{thcolsimcomp}) est qu'un corps~$\gK$ dans lequel $-1$ n'est pas une somme de carrés peut être ordonné. En revanche, la seule signification \und{calculatoire} connue de ce résultat des \clama est que la théorie~$\Sa{Cod}(\gK)$ s'effondre \ssi $-1$ est une somme de carrés dans $\gK$. 

D'un point de vue classique, comme la théorie $\Sa{Cod}(\RR)$ ne s'effondre pas, on peut munir $\RR$ d'une relation d'ordre qui prolonge la relation d'ordre usuel et qui est un ordre total. Mais la seule signification constructive de ce résultat des \clama est que $-1$ n'est pas une somme de carrés dans $\RR$.

Une autre conséquence du \pstref{Pst1bis} est le \corl suivant, pour lequel une \demo plus directe semble un défi. Une \demo directe sera donnée dans la théorie \Sa{Afrnz}.

%c
%:     Corollary{corPst1bis}
\begin{corollary} \label{corPst1bis}
Dans la théorie \sa{Aonz} la règle suivante est valide

\Regles{
\Lab{Aonz3} $\,\, a\geq 0\vet b\geq 0\vet a^2=b^2 \vd a=b$%
}
 
\end{corollary}
%--------- fin corollary ------------------------------- 

%%%%%%%%%%%%%%%%%%%%%%%%%%%%%%%%%%%%%%%%%%%%%%%%%%%%%%%%%%%%%%%%%%%%
\Subsection{La force démonstrative des Positivstellensätze formels}
%: Subsection{La force démonstrative des Positivstellensätze formels}

Les \tdys que nous explorons dans la suite pour décrire les \prts \agqs des nombres réels sont des extensions de \Sa{Asonz} (si le prédicat $\cdot > 0$ est présent) ou \Sa{Aonz} (dans le cas contraire). En outre les théories explorées sont toujours plus faibles que \Sa{Crcd}.
Et toute \ralg valide dans la \tdy \Sa{Crcd} est valide dans 
\sa{Asonz} (dans \sa{Aonz} si le prédicat $\cdot > 0$ est absent).

Or $\RR$ constitue un modèle \cof de la théorie \Sa{Asonz}
pour le langage basé sur la signature $(\cdot=0,\cdot>0,\cdot\geq0\cdot\mathrel{;}\cdot+\cdot, \cdot\times\cdot,-\,\cdot,0,1)$ (\thref{thRRco}).

Ainsi du point de vue des seules \ralgs, les Positivstellensätze formels nous disent que la théorie \Sa{Crcd} est entièrement satisfaisante, y compris pour $\RR$, qui ne satisfait pourtant
 ni~\coligt \, ni~\Tsbf{OT}. Cependant, pour tempérer cette déclaration optimiste, voici
le résultat précis. On notera aussi qu'il ne s'applique que pour les \ralgs,
pas pour les autres \rdys.

%t
%: Theorem{thRRtalg}
\begin{theorem} \label{thRRtalg}
Considérons une \ralg formulée dans la \sad $\gR=\Sa{Asonz}(\RR)$. Si les constantes qui interviennent dans la règle font partie d’un sous-corps discret $\gR_0$ de $\RR$, pour que la règle soit valide dans $\gR$, il suffit qu’elle soit valide dans $\Sa{Crcd}(\gR_0)$. 
\end{theorem}
%----------- fin theorem ----------------------------- 

Les règles existentielles satisfaites dans~$\RR$ et introduites dans les \tdys considérées seront autant que possible traitées dans le cadre d’existences prouvablement uniques et pourront donc être skolémisées sans dommage, fournissant des théories
sans axiomes existentiels \eseqs à celles qui auraient nécessité des axiomes existentiels.

%: Remark{remRRtdy0}
\begin{remark} \label{remRRtdy0} 
 On peut aussi appliquer le \thref{thFond} avec la \tdy $\Sa{Co--}(\RR)$
(\dfn \ref{defiCo0}). On introduira alors un prédicat $x\succeq y$ opposé à $x<y$. La nouvelle \tdy traitera $\RR$ comme un \codi et toute \rdy prouvée dans la nouvelle théorie mais n'utilisant pas $x\succeq y$
sera aussi valable dans $\RR$. L'inconvénient est naturellement que~$\RR$ n'est pas un modèle \cof de la nouvelle théorie. 
Un autre inconvénient est le statut mystérieux du nouveau prédicat $x\succeq y$, plus faible que $x\geq y$ dans la nouvelle \tdy.
En conclusion, l'avantage que semble procurer le \thref{thFond} (l'usage de la logique classique est inoffensif) ne semble pas aller au delà des considérations que l'on a développées sur le bon usage du Positivstellesatz formel.
\eoe\end{remark}
%----------- fin remark ---------------------------------- 

% subsubsection{Positivstellensatz concret}
\Subsection{Positivstellensatz concret}
%: Subsection{Positivstellensatz concret}

On rappelle tout d'abord le \tho fondamental de Tarski.
Pour une \demo simple, dite \gui{à la Cohen-Hormander}, voir \cite[Section 1.4]{BCR}, ou \cite[Lemme 3.12]{CLR01}. Quelques commentaires instructifs se trouvent dans \cite[\thos 10, 11, 12]{LR90}.

%t
%: Theorem{thTarski}
\begin{theorem}[élimination des quantificateurs] \label{thTarski}\label{NOTARa}
La théorie formelle intuitionniste du premier ordre associée à
la \tdy \Sa{Crcd} admet l'élimination des quantificateurs.
Elle est complète et décidable. En particulier elle décrit de manière exhaustive toutes les propriétés purement \agqs de $\RRa$ (celles formulées au premier ordre dans le langage des anneaux ordonnés). 
\end{theorem}
%----------- fin theorem ----------------------------- 

On donne dans ce paragraphe un \tho \eqv au Positivstellensatz de Krivine-Stengle, énoncé ici dans le langage des \sads.

Une \prco du \pstref{thPstStengle} se trouve dans \cite{CLR01} ou \cite{Lom91}. 
Elle est fondée sur Positivstellensatz formel d'une part et sur le \gui{lemme 3.12} de \cite{CLR01}, variante du \tho de Tarski.

Pour une approche plus conceptuelle et des bornes de complexité améliorées voir
\cite{LPR2015}. Pour la construction de la clôture réelle d’un \codi voir \cite{LR91,LR90}.

%: pst{thPstStengle}
\begin{pstc} \label{thPstStengle}%
\index{Positivstellensatz!concret} ~
\\
Soit $\gK$ un \codi et $\gR$ un corps réel clos discret contenant~$\gK$, (par exemple la clôture réelle de $\gK$).
Soit $\gA=\big((G,Rel),\sa{Cod}(\gK)\big)$ une \sad où $G=(\xn)$ et où $Rel$ est fini. 
\begin{enumerate}
\item La \sad $\gA$ s'effondre \ssi il est impossible de trouver un modèle de $\gA$ contenu dans $\gR$. 
\item L'effondrement s'il a lieu est donné par un certificat \agq conformément au point~1 du \thref{Pst1} et au lemme \ref{lemColApo}.
\item On a un \algo qui décide si $\gA$ s'effondre et qui en cas de réponse négative donne la description d'un \sys $(\xin)$ dans $\gR^n$ qui satisfait les contraintes données dans les relations $Rel$.
\end{enumerate} 
\end{pstc}
%--------- fin pst ----------------------------------- 

Cet énoncé n'est pas valable sous cette forme générale si l'on prend $\gK=\gR=\RR$ car il n'y a pas de test de signe dans $\RR$ et les \algos qui explicitent le \pstref{thPstStengle}\footnote{Ces \algos sont fournis par la \prco du \tho.} utilisent de manière cruciale ce test de signe. 

\smallskip Voici un petit exemple des \pbs auxquels on se heurte.
Sur $\RR$, comme sur un anneau local arbitraire\footnote{Pour le traitement constructif des anneaux locaux, du radical de Jacobson et des corps de Heyting voir par exemple \cite[section IX-1]{CACM}.} dans lequel $x\neq 0$ désigne le prédicat d'inversibilité, on a l'\eqvc
% equation label {eqTiersexclu1}
\begin{equation} \label {eqTiersexclu1}
\exists y\ \ x^2y=x \,\iff\;x=0 \;\; \mathrm{ ou }\;\; x\neq 0.
\end{equation}
% end-equation
En effet supposons $x(1-xy)=0$. Si $xy$ est inversible, alors $x$ est inversible, et si $1-xy$ est inversible, alors $x=0$. Cette \demo se traduit formellement dans la \tdy correspondante par l'établissement des trois règles valides suivantes: 

\Regles{
\itbu $ \,\, x^2y=x\vd x=0 \; \vou\; x\neq 0$, 
\itbu $\,\, x=0\vd\Exists y\ \ x^2y=x$,
\itbu $\,\, x\neq0\vd\Exists y\ \ x^2y=x$.
}

Ce cas simple d'élimination du quantificateur $\exists$ montre que l'on aboutit dans les calculs à des impasses du point de vue de la décidabilité,
puisque \gui{$x=0 \;\mathrm{ ou }\; x \hbox{ inversible}$} est indécidable dans~$\RR$.

\smallskip Néanmoins, dans la section finale de l'article \cite{GL93}, on trouve une forme \cov entièrement satisfaisante pour le 17\ieme\ \pb de Hilbert sur $\RR$. Et d'autres cas de Positivstellensätze \cot prouvables sur $\RR$ sont \egmt traités.

%

%%%%%%%%%%%%%%%%%%%%%%%%%%%%%%%%%%%%%%%%%%%%%%%%%%%%%%%%%%%%%%%%%%%%
\section{Corps ordonnés \textsl{non} discrets}\label{secConondisc}
%\Today

 En première approximation, et en suivant une suggestion de Heyting, on pourrait choisir comme théorie formelle du premier ordre pour les \prts \agqs de $\RR$ la théorie~\Sa{Asonz} (vue comme théorie formelle du premier ordre) à laquelle on ajoute les axiomes géométriques \Tsbf{IV} et \Tsbf{OTF} ainsi que l'axiome \tsbf{HOF} suivant, non géométrique donc indésirable (axiome de Heyting ordonné).\index{Heyting!corps de --- ordonné}\index{corps!de Heyting}

\Regles{
\Lab{HOF} $\quad (x> 0 \,\Rightarrow 1=0) \,\Rightarrow\, x\leq 0$\qquad (Heyting Ordered Field)
}

Cela revient à remplacer dans la théorie \Sa{Cod} les axiomes		
\coligt \ et \Tsbf{OT} par les axiomes~\Tsbf{OTF} et~\Tsbf{HOF}. 
 On a alors une structure d’anneau local, car les règles \Tsbf{Iv} et \Tsbf{IV} impliquent que \gui{$x\neq0 $} signifie \gui{$x$ est inversible}, donc \tsbf{OTF} implique que pour tout $x$, $x$ \hbox{ou $1-x$}
est inversible.

Dans ce cadre l'axiome~\tsbf{HOF} signifie que le radical de Jacobson est réduit à $0$.

%r
%: Remark{remHOH}
\begin{remark} \label{remHOF} 
Notons que l'axiome \tsbf{HOF}, formulable au premier ordre même s'il n'entre pas dans le cadre des \tdys, est satisfait de manière indirecte sous la forme suivante:
\textsl{dans une \sad de type \Sa{Asonz}, si un terme clos $t$ vérifie $\,\,t>0\vd\Bot$, alors il vérifie aussi $\vd t\leq 0$.} Cela résulte du Positivstellensatz formel. 
\\
À vrai dire on a même: \textsl{si un terme clos $t$ vérifie $\,\,t\geq 0\vd\Bot$, alors il vérifie aussi $\vd t< 0$}. Cela signifie que le principe de Markov, qui s'exprime sur $\RR$ par l'implication $\,\,\lnot (t\geq 0)\,\Rightarrow\,t<0 $ vaut comme une règle de déduction \textsl{externe}\footnote{Nous devons ajouter ici le mot \textsl{externe}, car il ne s'agit pas d'une règle valide
dans la \sad elle-même. Il s'agit du fait que l'on déduit la validité d'une règle de celle d'une autre règle.} admissible dans la \tdy $\sa{Asonz}(\RR)$.\\
Les mêmes remarques s'appliquent pour les \tdys qui étendent \sa{Asonz} tout en s'effondrant simultanément.
\eoe
\end{remark}
%----------- fin remark ---------------------------------- 

%%%%%%%%%%%%%%%%%%%%%%%%%%%%%%%%%%%%%%%%%%%%%%%%%%%%%%%%%%%%%%%%%%%%
%:Subsection
\Subsection{Une première \tdy}\label{subsecCo0}
 
 Outre le caractère indésirable de \Tsbf{HOF}, la théorie formelle envisagée au début de la section
présente un inconvénient majeur, qui est de ne pas pouvoir démontrer
l'existence de la borne supérieure de deux \elts: 
voir à ce sujet~\cite{Coq07}.

Il est donc légitime d'explorer les possibilités qu'offre l'ajout d'une loi pour cette borne supérieure, avec les règles adéquates.
 Nous proposons maintenant une \tdy minimaliste
 pour les corps ordonnés \textsl{non} discrets en introduisant le symbole de fonction $\vu$ dans le langage.
 
%d
%: Definition{defiCo0}
\begin{definition} \label{defiCo0}
Une première \tdy minimale pour les corps ordonnés \textsl{non} discrets, notée \SA{Co--}, est basée sur la signature suivante. Il y a une seule sorte, nommée $\CoO$.
\Sigt{\CoO}{\cdot=0,\cdot\geq 0,\cdot>0\mathrel{;}\cdot+\cdot, \cdot\times\cdot,\cdot\vu\cdot,-\,\cdot,0,1} \label{NOTASigCoO}
\noindent Les axiomes sont ceux de \Sa{Asonz}, les axiomes \Tsbf{IV} et \Tsbf{OTF}, et les axiomes naturels pour $\vu$: \Tsbf{sup1}, \Tsbf{sup2}, \Tsbf{Sup}, \Tsbf{grl} et~\Tsbf{afr}.
Ils sont tous énumérés ci-après ($x^+$ est une abréviation de $x\vu 0$).%

\DeuxRegles{
\laB{ga0} $\vd 0=0$
\laB{ac2} $\,\,x=0\vd xy=0 $
\laB{gao0} $\vd 0 \geq 0$
\laB{ao1} $\vd x^2 \geq 0$
}
{
\laB{ga1} $\,\, x=0\vet y=0\vd x+y=0$
\laB{}
\laB{gao1} $\,\, x \geq 0\vet y \geq 0 \vd x + y \geq 0$
\laB{ao2} $\,\, x \geq 0\vet y \geq 0 \vd x y\geq0$ 
}

\DeuxRegles{
\laB{aso1} $ \vd 1> 0$ \phantom{$x^2$}
\laB{aso2} $\,\, x> 0 \vd x \geq 0$
\laB{sup1} $ \vd x\vu y\,\geq x $
\laB{sup2} $ \vd x\vu y\,\geq y $
}
{
\laB{aso3} $\,\, x > 0\vet y \geq 0 \vd x + y > 0$
\laB{aso4} $\,\, x > 0\vet y > 0 \vd xy > 0$
\laB{grl} $\vd x+(y\vu z)=(x+y)\vu(x+z)$
\laB{afr} $\vd x^+\, (y\vu z)=(x^+\, y)\vu(x^+\, z)$
}

\DeuxRegles{
\laB{Gao} $\,\, x\geq 0\vet x\leq 0 \vd x = 0$
\laB{Anz} $\,\, x^2= 0 \vd x = 0$
\laB{Aonz} $\,\, c\geq 0\vet x(x^2+c)\geq 0 \vd x\geq 0$
\laB{Sup} $ \,\, z\geq x\vet z\geq y\vd z\geq x\vu y$
\laB{IV} $\,\, x> 0 \vd \Exists y\, xy = 1$
}
{
\laB{Iv} $\,\, xy = 1 \vd x^2> 0$
\laB{Aso1} $\,\, x> 0\vet xy\geq 0 \vd y\geq 0$
\laB{Aso2} $\,\, x\geq 0\vet xy> 0 \vd y> 0$
\lab{}
\laB{OTF} $\,\, x+y> 0 \vd x >0 \vou y>0 $
}

\end{definition}
%----------- fin definition -------------------------------- 

%r
%: Remark{remCo0}
\begin{remarks} \label{remCo0}~\\ 
 1) Notons que le collapsus \gui{$\,\,0>0\vd 1=0$\,\,} se déduit de \Tsbf{IV}. 

\smallskip \noindent 2)
Si l'on ajoute les axiomes \Tsbf{OT} et
\coligt \ à la théorie \Sa{Co--} on retrouve une théorie \esid à \Sa{Codsup} ou \Sa{Cod}. En fait, il suffit d'ajouter l'axiome \tsbf{ED$_>$}: voir le lemme \ref{lemAfrsdz}. 
\eoe\end{remarks}
%----------- fin remark ---------------------------------- 

%e
%: Example{exacorpsnondiscret}
\begin{examples} \label{exacorpsnondiscret} 
 De nombreux sous-corps \gui{naturels} de $\RR$ sont \textsl{non} discrets,
par exemple le corps énumérable $\RR_{\tsbf{PR}}$ des réels calculables en temps primitif récursif, ou le corps énumérable $\RR_{\tsbf{Ptime}}$ des réels calculables en temps polynomial, ou encore le corps \textsl{non} énumérable $\RR_{\tsbf{Rec}}$ des réels récursifs. 
Une \tdy satisfaisante pour les \prts \agqs des nombres réels devra accepter pour modèles ces sous-corps \gui{naturels} de $\RR$.\eoe
\end{examples}
%--------- fin example ---------------------------------------- 

Les sous-corps $\RR_{\tsbf{PR}}$ et $\RR_{\tsbf{Ptime}}$ peuvent être manipulés sur machine de manière plus sympathique que le corps $\RR_{\tsbf{Rec}}$. 
Ce sont en effet des corps énumérables (quoique \textsl{non} discrets),
dont les \elts n'ont pas besoin d'être accompagnés de \gui{certificats}
extérieurs à la \tdy considérée (une fonction récursive \gnle existe \cot seulement si elle est accompagnée d'un \gui{certificat}: une \prco du fait qu'elle est totale). 

Notons que le caractère \gui{complet} de $\RR$ semble relever plus de l’analyse que de l’algèbre. 
Notons aussi que le \gui{corps} des séries de Puiseux sur $\RR$ ne semble pas satisfaire \Tsbf{OTF} (pour n'importe quelle tentative de \dfn raisonnable pour la relation d'ordre).

%%%%%%%%%%%%%%%%%%%%%%%%%%%%%%%%%%%%%%%%%%%%%%%%%%%%%%%%%%%%%%%%%%%%
%:Subsection
\Subsection{L'axiome de convexité et la théorie \Sa{Co}
}\label{subsecCocvx}

Outre la fonction sup, d'autres fonctions \gui{rationnelles}
posent le même type de \pbs. 

\smallskip Dans la théorie des anneaux réels clos, en \clama, 
(voir les articles~\cite{Sch84,PN2002} et la section \ref{secArc}), l'axiome \gui{de convexité\footnote{Il y a deux utilisations fort distinctes du terme \gui{convexe} dans le texte présent. D'une part un anneau ordonné peut être déclaré convexe, comme ici.\index{convexe!anneau ordonné ---}\index{anneau!ordonné!convexe} D'autre part, un sous-groupe d'un groupe ordonné peut être déclaré convexe comme \paref{subsecgrlsolide}.}} suivant est satisfait

\Regles{\Lab{CVX} $\,\,0\leq a\leq b\vd \Exists z\; zb=a^2$ \quad (convexité)}

\noindent On note que si $zb=a^2$ alors $(z\vi a)b=zb\vi ab=a^2\vi ab=a(a \vi b)=a^2$. De même $(z\vu 0)b=a^2$. En conséquence, un axiome \eqv qui assure l’unicité de $z$ est donné par la \rdy suivante:

\Regles{
\Lab{FRAC} $\,\,0\leq a\leq b\vd \Exists z\; (zb=a^2\vet 0\leq z\leq a)$ 
}\label{RFRAC}

\noindent Cette règle est valide pour $\RR$ car l'\elt $z$ est une fonction continue de $(a,b)$ sur son domaine de définition.

%l
%: Lemma{lemUniqFRAC}
\begin{lemma} \label{lemUniqFRAC}
L’unicité de $z$ (lorsqu'il existe) dans la règle \Tsbf{FRAC} est assurée dans la \talg \Sa{Atonz} et dans~\Sa{Co--}. Le même calcul montre que l'unicité est assurée dans la théorie \agq \Sa{Afrnz} (lemme~\ref{lemAfrnzFRAC}).
\end{lemma}
%----------- fin lemma ----------------------------------- 
%
\begin{proof}
 En effet, supposons $yb=zb=a^2$, $0\leq z\leq a$ et $0\leq y\leq a$. On a ${(y-z)}\,b=0$, $\abs{y-z}\leq a\leq b$ et \hbox{donc $\abs{y-z}^2\leq\abs{y-z}\,b=0$}. 
\end{proof}
%

%: Lemma{lemCo0FRAC}
\begin{lemma} \label{lemCo0FRAC}
L'ajout de l'axiome \Tsbf{FRAC} à la théorie \Sa{Co--} peut être remplacé par l'introduction d'un symbole de fonction $\Fr$ avec les axiomes

\DeuxRegles{
\Lab{fr1} $ \vd \mathrm{Fr}(a,b)\, \abs b=(\abs a\!\vi\! \abs b)^2 $
}
{
\Lab{fr2} $ \vd 0\leq {\mathrm{Fr}(a,b)}\leq \abs a\!\vi\! \abs b \phantom{(\abs a\!\vi\! \abs b)^2}$
}

\noindent La même remarque s'applique aux théories \Sa{Atonz} et \Sa{Afrnz}.
\end{lemma}
%----------- fin lemma ----------------------------------- 

%
\begin{proof} La règle \tsbf{FRAC} est \eqve à la règle suivante

\Regles{\lab{~}$\vd \Exists z\; \big(z\abs b = (\abs a\vi\abs b)^2\vet\, 0\leq z\leq \abs a\vi\abs b\big)$}

\noindent On a donc simplement skolémisé une règle existentielle dans le cas d’existence unique (lemme~\ref{lemUniqFRAC}). Ceci nous donne une extension \esid (voir \paref{skolemunique}). Par ailleurs, une fois ajoutés le symbole de fonction~$\mathrm{Fr}$ et les axiomes \tsbf{fr1} et \tsbf{fr2}, la règle \tsbf{FRAC}
devient clairement inutile.
\end{proof}
%

%: Definition{defiConondiscret}
\begin{definition} \label{defiConondiscret} 
 La \tdy \SA{Co} des corps ordonnés \textsl{non} discrets\footnote{Dans \cite{LM2017}, on a donné une \dfn légèrement différente, un peu plus forte, voir la remarque \ref{remCo}.}
 est l'extension de la théorie \Sa{Co--} obtenue en ajoutant le symbole de fonction $\mathrm{Fr}$ et les axiomes \tsbf{fr1} et \tsbf{fr2}.\index{corps ordonné!(\textsl{non} discret)}
\end{definition}
%----------- fin definition -------------------------------- 

C'est une \tdy à une seule sorte avec la signature
\Sigt{\Co}{\cdot=0,\cdot\geq 0,\cdot>0\mathrel{;}\cdot+\cdot, \cdot\times\cdot,\cdot\vu\cdot,-\,\cdot,\Fr(\cdot,\cdot),0,1}
\label{NOTASigCo}
\noindent Les sous-corps \gui{naturels} de $\RR$ \textsl{non} discrets de l’exemple \ref{exacorpsnondiscret}
 sont des modèles de~\sa{Co}, et aussi de certaines extensions de \sa{Co} que nous définissons par la suite, comme~\Sa{Crc1} ou~\Sa{Corv}. 

%r
%: Remark{remPst1Co}
\begin{remark} \label{remPst1Co} 
Dans l'énoncé du \pstfref{Pst1}, on peut ajouter la théorie \sa{Codsup} qui est \esid à la théorie \sa{Cod}, et 
la théorie \sa{Co} qui est intermédiaire entre \sa{Asonz} et \sa{Codsup}. Pour plus de précision voir le \pstfref{Pst2}. \eoe
\end{remark}
%----------- fin remark ---------------------------------- 

\Subsection{Autres opérations \gui{rationnelles} continues}\label{subsecratcon}

 Voici un autre exemple paradigmatique avec une fonction continue partout définie
\vspace{-.5em}
% equation label {eqratcont}
\begin{equation} \label {eqratcont}
f(x,y)=\frac{(ax+by)xy}{x^2+y^2}
\end{equation}
Cette fraction rationnelle\footnote{La fonction est à priori définie pour $(x,y)\neq (0,0)$ et se prolonge par continuité en $f(0,0)=0$. On voit alors qu'elle est uniformément continue sur tout cube $[-a,+a]^4$.} est le prototype d'une famille, paramétrée par $a,b$, de fonctions réelles continues 
 $\RR^2\to \RR$ (ou d'une fonction réelle continue $\RR^4\to \RR$).

Une \rdy \gui{définit} cette fonction:
% equation label {eqratcont2}
\begin{equation} \label {eqratcont2}
\vd \Exists z\quad \big(z(x^2+y^2)=(ax+by)xy ,\;
\abs z\leq \abs{ax+by}\big)
\end{equation}
% end-equation
et elle ne semble pas valide dans la théorie de base \Sa{Co--}.

Dans cet exemple, si $a=b=1$, la fraction est du type $z=u/v$
avec $u^2\leq v^3$. Elle est caractérisée par les relations $zv=u$
et $\abs z ^2\leq \abs v$. Or les \rdys suivantes sont satisfaites pour $\RR$, et aussi bien pour les \crcds:

\Regles{\lAb{FRAC$_n$} $ \,\,\abs u^n\leq \abs v^{n+1}\vd \Exists z\; (zv=u\vet \abs z^{n}\leq \abs v) $ \quad ($ n\geq 1 $) }\label{AxFRACn}

\hum{On pourrait hésiter et préférer: $\,\,1\vi\abs u^n\leq \abs v^{n+1}\vd \Exists z\; (zv=u\vet 1\vi \abs z^{n}\leq 1\vi\abs v)$ car ce qui se passe lorsque $\abs v\geq 1$ n'a pas d'importance}

\smallskip Intuitivement cette règle signifie que la fraction $u/v$ est bien définie. Dans le cas d’un \codi on raisonne cas par cas: si $v\neq 0$ c'est clair, si $v=0$ la règle force $z=0$. 
Plus \gnlt on vérifie que l'existence, si elle est supposée, est prouvée unique dans la théorie \Sa{Afrnz} (\paref{theorieAfrnz}) comme suit.
\\
Si $zv=u=yv,\, \abs z^{n}\leq \abs v,\, \abs {y}^{n}\leq \abs v$, on pose $w=\abs{z-y}$ et l'on obtient 

\snic{w\abs v=0,\,\leq \abs z +\abs{y}\leq 2 {\abs v}^{\frac1 n},\,w^n\leq 2^n {\abs v},\,0\leq w^{n+1}\leq 2^n\abs v w=0
,}

\noindent donc $w^{n+1}=0$ et enfin $w=0$.

\noindent Un autre argument consiste à dire que la \rsim

\Regles {\lab{~} $\,\,zv=u\vet yv=u\vet \abs z^{n}\leq \abs v\vet \abs {y}^{n}\leq \abs v \vd z=y $}

\noindent est satisfaite dans \Sa{Cod}, et que \Sa{Afrnz} et \Sa{Cod} prouvent les mêmes \ralgs (\pstfref{Pst2}).

\smallskip Dans la \talg \Sa{Aonz} la règle \tsbf{FRAC} découle de \tsbf{FRAC$_1$} en posant~\hbox{$u=a^2$} \hbox{et $v=b$}. 

Inversement les règles $\tsbf{FRAC}_n$ se déduisent de la règle \tsbf{FRAC} dans un cadre assez \gnl (voir le lemme \ref{lemfracfrac2}). Nous utilisons donc par la suite uniquement la règle \tsbf{FRAC}.

%r
%: Remark{remCo}
\begin{remark} \label{remCo} 
Il aurait été plus logique de demander, dans la \dfn de la théorie \Sa{Co}, outre la validité de la règle \tsbf{FRAC}, celle des règles $\tsbf{FRAC}_n$
(c'est le choix qui était fait dans l'article \cite{LM2017}, \dfn 2.13). Plus \gnlt, pour toute fonction $f\colon \QQ^n\to \QQ$ qui prolonge par continuité\footnote{Plus \prmt: la fraction $h/p$ à \coes dans $\QQ$, définie pour $p(x)\neq 0$, doit se prolonger en une fonction continue $\QQ^n\to\QQ$.} une fraction $h/p$, où $h$ et $p$ sont des semipolynômes\footnote{Une sup-inf combinaison de \pols.}, les zéros de $p$ étant d'intérieur vide, on devrait demander que soit valide la règle
analogue à $\tsbf{FRAC}_n$ qui dit que \gui{$f$ existe}, et plus \prmt introduire un symbole de fonction correspondant avec les axiomes adéquats. 
Mais nous n'avons pas voulu trop compliquer la \dfn de la théorie \Sa{Co}
des corps ordonnés \textsl{non} discrets dans la mesure où nous avons essentiellement en vue la théorie des corps réels clos \textsl{non} discrets,
dans laquelle la règle \tsbf{FRAC} est suffisante.
\eoe\end{remark}
%----------- fin remark ---------------------------------- 

%: sous section supplémentaire
\section{Un corps ordonné \textsl{non} discret non archimédien}
\label{condna}

Nous décrivons dans cette section un exemple de corps ordonné de Heyting \textsl{non} discret non archimédien.

Soit $\vep$ une \idtr. On note $\gZ=\QQ[[\vep]]$ l'anneau des séries formelles à \coes rationnels et
$\gQ=\QQ((\vep)):=\QQ[[\vep]][1/\vep]$. En \clama $\gZ$ est un anneau local intègre hensélien et $\gQ$ est son corps de fractions. On munit $\gZ$
de la relation d'ordre pour laquelle $\vep$ est un infiniment petit $>0$
(i.e. $0<\vep$ et $\vep<r$ pour tout rationnel $r>0$). On notera encore $\gZ$ l'anneau totalement ordonné ainsi obtenu.
Le localisé $\gZ[1/\vep]$ muni de la relation d'ordre compatible avec celle de $\gZ$ sera encore noté~$\gQ$. C'est sans doute l'exemple le plus simple de corps ordonné non archimédien.

Ces affirmations, classiques, sont encore valables en \coma à condition d'utiliser des \dfns convenables. Par exemple nous allons voir que, modulo des \dfns convenables, $\gQ$ est un modèle de la théorie~\Sa{Co}.

D'un point de vue \cof, il n'y a pas de test de signe sur $\gZ$ ou sur $\gQ$.
Et il n'est pas immédiat de définir une structure d'ordre correspondant à l'intuition donnée par les \clama.

\smallskip Voici comment on peut traiter cet exemple \cot.
Un \elt de $\gZ$ est donné par une série formelle $\xi=\sum_{j=0}^{+\infty} x_j \vep^j$ avec les $x_j\in\QQ$.
Tout rationnel dans $\QQ$ peut être considéré comme un \elt de $\gZ$ selon le procédé usuel. 
Le \coe $x_j$ est noté $\rc_j(\xi)$. On pose conventionnellement $\rc_{j}(\xi)=0$ pour $j<0$ dans $\ZZ$ et pour tout $\xi\in\gZ$. 

Pour une série $\xi=\sum_{j=0}^\infty x_j \vep^j$ dans $\gZ$ on définit pour chaque $k\geq 0$ un \emph{signe potentiel en l'exposant $k$}, noté $\kappa_k(\xi)\in\so{-1,0,1}$ comme suit, par récurrence sur $k$: 
\begin{itemize}
\item $\kappa_{0}(\xi)$ est le signe de $x_{0}$;
\item si $\kappa_k(\xi)\neq 0$ alors $\kappa_{k+1}(\xi)=\kappa_k(\xi)$, sinon $\kappa_{k+1}(\xi)$ est le signe de $x_{k+1}$;
\item on pose conventionnellement $\kappa_{j}(\xi)=0$ pour $j<0$ dans $\ZZ$ et pour tout $\xi\in\gZ$.
\end{itemize}

\smallskip Au moins intuitivement on a le résultat suivant: si $\kappa_k(\xi)=1$, alors \hbox{$\xi>0$}; si $\kappa_k(\xi)=-1$, alors~\hbox{$\xi<0$}; si $\kappa_k(\xi)=0$, alors le signe de $\xi$ est à priori inconnu.

L'ensemble $\gZ$ est muni de la structure d'anneau usuelle (pour les séries formelles) et l'\egt $\xi=\zeta$ a lieu exactement lorsque les séries sont identiques. Cela équivaut à $\forall k\geq 0\; \kappa_k(\xi-\zeta)=0$.
Cet anneau est limite projective de la suite des morphismes surjectifs $\pi_{k+1,k}\colon\aqo{\QQ[\vep]}{\vep^{k+1}}\to \aqo{\QQ[\vep]}{\vep^k}$ ($k\in\NN$), via les morphismes naturels
$\pi_k\colon\gZ\to \aqo{\QQ[\vep]}{\vep^k}$ obtenus par la troncation des séries à l'ordre $k$. 

\smallskip
Les bases de la théorie \cov des anneaux locaux résiduellement discrets henséliens, y compris la construction du hensélisé d'un \alrd, sont traitées dans l'article \cite{ALP08}. 

\smallskip On introduit maintenant en détail tout ce qui est nécessaire
pour le traitement \cof de l'anneau ordonné~$\gZ$.

\begin{enumerate}
\item On définit $\xi>0$ par \fbox{$\exists k \;\kappa_k(\xi)=1$} et $\xi\geq 0$ par \fbox{$\forall k\;\kappa_k(\xi)\geq 0$}. %, $\xi>\zeta$ par $\xi-\zeta>0$, et $\xi\geq \zeta$ par $\xi-\zeta\geq 0$. 
\\On a alors:
\begin{itemize}
\item $\xi=0$ \ssi $\xi\geq 0$ et $\xi\leq 0$;
\item les règles \Tsbf{OTF} et \OTFx sont valides dans $\gZ$;
%
%\item 
%
\item l'axiome de Heyting \gui{$\lnot (x>0)\Rightarrow -x\geq 0$} est satisfait. 
%
%\item 
\end{itemize}

\item \textsl{Valeur absolue et fonction $\sup$}. On définit la fonction \gui{valeur absolue} $\xi\mapsto\abs\xi$ en posant \fbox{$\rc_k(\abs\xi)=\kappa_k (\xi) \rc_k(\xi)$} pour tout~$k$.
Enfin \fbox{$\xi\vu\zeta:=(\xi+\zeta+\abs{\xi-\zeta})/2$}.
\item On vérifie alors que $\gZ$ est un \asr (théorie \Sa{Asr}, chapitre \ref{chap-afr}), autrement dit, en ajoutant le fait que $\gZ$ est réduit, toutes les \ralgs valides dans la théorie \sa{Codsup} sont satifaites dans $\gZ$ (voir le \Pst formel \ref{Pst2}, point \ref{i4Pst2}). 
\item \textsl{Valuation}.\label{avalnd}\index{valuation} 
\begin{enumerate}
\item On définit $\rv(\xi)= k \equidef (\kappa_k(\xi)=\pm1,\; \kappa_{k-1}(\xi)=0)$.
\item On définit $\rv(\xi)> k \equidef \kappa_k(\xi)=0$.
\item On lit intuitivement $\rv(\xi)> k$ comme \gui{la valuation de $\xi$ est $>k$}. En fait, $\rv(\xi)$ n'est pas un \elt de $\NN$ mais %un \elt 
d'une compactification convenable de $\NN$ contenant $+\infty$\footnote{C'est l'espace métrique obtenu en prenant sur $\NN$ la métrique $d(k,\ell)=\sup(2^{-k},2^{-\ell})$ et en complétant (cela envoie $\infty$ sur $0$). }. 
\item On a \prmt la description suivante pour $\alpha<\beta$;
% equation label {eqcondna1}
\begin{equation} \label {eqcondna1}
\alpha<\beta\;\Longleftrightarrow\; \exists k\;
\formule{
(\kappa_k(\alpha)=-1,\; \kappa_k(\beta)\geq 0) & \vu\\
(\kappa_k(\alpha)\leq 0,\; \kappa_k(\beta)=+1) & \vu\\
(\rv(\alpha)=\rv(\beta)=k, \; \rc_k(\alpha)<\rc_k(\beta))}
\end{equation}
% end-equation

%
\item On a $\rv(\xi^2)=2\rv(\xi)$ (\egt définie par \gui{$\rv(\xi^2)>2k-1 \Leftrightarrow \rv(\xi^2)>2k \Leftrightarrow\rv(\xi)>k$}).
\item On en déduit que $\gZ$ est un anneau réduit (c'est donc un modèle \cof de la théorie \sa{Asrnz}).
\item Enfin, 
$\gZ$ est un \textsl{anneau de valuation} au sens suivant: c'est un \asr réduit dans lequel deux \elts $\alpha$ et $\beta$ strictement positifs sont toujours comparables pour la \dve. En d'autres termes, la règle suivante \tsbf{Val1} est valide\footnote{Nous donnons ici une \dfn pour le cas d'un \asr. On pourrait donner une \dfn plus \gnle pour un \alrd muni d'un prédicat $\cdot \neq 0$ convenable.}.\index{valuation!anneau de ---} 

\Regles
{\Lab{Val1}$ \,\,\alpha>0\vet \beta>0 \vd \Exists \xi\; \alpha\xi=\beta \vou \Exists \xi\; \beta\xi=\alpha$}

Une règle voisine \egmt valide dans $\gZ$ est la suivante.

\Regles
{\Lab{Val2} $ \,\, \beta\geq \alpha\geq 0\vet\beta>0 \vd \Exists \xi\; \alpha=\beta\xi$
}

Démontrons ces règles. De $\beta>0$ on déduit qu'il y a un $k$ tel que $\beta=\vep^k\gamma$ avec $\rv(\gamma)=0$. On va voir dans le point \ref{serconv} que $\gamma\in\gZ\eti$. Donc $\vep^k=\gamma^{-1}\beta$.
Pour \tsbf{Val1} on a aussi un $\ell$ tel que $\alpha=\vep^\ell\delta$ avec $\delta\in\gZ\eti$.
 D'où la disjonction selon que $\ell\geq k$ ou $k > \ell$. Pour \tsbf{Val2}, de $0\leq \alpha\leq \beta$ on déduit $\kappa_{k-1}(\alpha)=0$, donc $\alpha=\vep^k\delta=\beta\gamma^{-1}\delta$ pour un $\delta\in\gZ$.
 
 \smallskip \hum{On peut se demander si la règle \tsbf{Val2} ne serait pas valide dans la théorie \sa{Asrnz}. En effet, cette théorie prouve les mêmes \ralgs que \sa{Crcdsup} (\Pst formel \ref{Pst2} point \ref{i4Pst2}). Cela semble quand même improbable.}

\noindent Notons aussi que l'implication suivante est satisfaite.
 $$\rv(\alpha)=\rv(\beta)=k \Rightarrow \exists \xi\in\gZ\eti\;\; \alpha=\beta\xi$$
\item Le \textsl{groupe de valuation} est le groupe ordonné que l'on définit comme le symétrisé du monoïde de \dve $\gZ\eti/\gZ^+$ où $\gZ^+=\sotq{\alpha\in\gZ}{\exists k>0\;\rv(\alpha)=k}$. Ce groupe de valuation est isomorphe à $(\ZZ,+,\geq)$, il est engendré par (la classe de) $\vep$. Dans la terminologie usuelle, on dit que $\gZ$ est un \textsl{\adv discrète},
mais ici le mot discret n'a pas le sens habituel qu'on lui donne en \coma.%
\index{valuation!groupe de ---}\index{valuation!discrète} 
\end{enumerate}

\item \textsl{Séries convergentes}.\label{serconv} Si l'on a une suite infinie $(\xi_n)_{n\in\NN}$ dans $\gZ$ et si la suite des $\rv(\xi_n)$ tend vers $+\infty$, alors la somme infinie 
$\sum_{n\in\NN}\xi_n$ est bien définie. \\
En particulier, si $\rv(\xi)>0$
la somme $\zeta=1+\sum_{n\in\NN}\xi^n$ est bien définie et l'on a $\zeta(1-\xi)=1$. \\
De là on déduit les \prts suivantes.
\begin{enumerate}
\item\label{iaserconv} L'anneau $\gZ$ est un anneau local, dont le corps résiduel est discret, isomorphe à $\QQ$. 
\item\label{ibserconv} On a $\gZ\eti=\sotq{\xi\in\gZ}{\kappa_0(\xi)=\pm1}$ et $\Rad(\gZ)=\sotq{\xi\in\gZ}{\kappa_0(\xi)=0}$. 
\item\label{icserconv} Si $\rc(\beta)=k$, on écrit $\beta=\vep^k\gamma$ avec $\gamma=c_0(1- \alpha)$ et $\rv(\alpha)>0$, donc:\\ $\vep^k=\beta \gamma^{-1}=\beta c_0^{-1}(1+\sum_{n\in\NN}\alpha^n)$.
\item\label{idserconv} L'anneau $\gZ$ est hensélien. \Prmt, si $P\in\gZ[X]$ satisfait les conditions $\rv(P(0))>0$ et $\rv(P'(0))=0$, il existe un (unique) $\xi\in\gZ$ tel que $P(\xi)=0$ et $\rv(\xi)>0$. Pour construire la série $\xi$, on utilise la méthode de Newton.
\item\label{ieserconv} L'anneau $\gZ$ est le hensélisé de l'\alrd 
$(\QQ[\vep])_{1+\vep\QQ[\vep]}$. 
\end{enumerate}
\item Enfin, nous montrons que la règle \Tsbf{FRAC} est valide dans $\gZ$.
L'hypothèse est donnée par deux \elts $\xi,\zeta\in\gZ$ qui vérifient $0\leq \xi\leq \zeta$ et nous cherchons un $\rho$ tel que $0\leq \rho\leq \xi$ \hbox{et $\rho\zeta=\xi^2$}. 
D'après le lemme \ref{lemAfrnzFRAC}, l'unicité de~$\rho$
(si existence) est assurée dans les \afrs réduits, comme dans la théorie \sa{Co--} (lemme \ref{lemUniqFRAC}).
\\ 
Nous notons que $0\leq \xi\leq \zeta$ implique que $\rv(\xi)\geq \rv(\zeta)$. Nous définissons $c_k(\rho)$ comme suit. 
\begin{itemize}
\item Si $\kappa_k(\zeta)=0$, alors
$\kappa_k(\xi)=0$, ce qui force $\kappa_k(\rho)=0$, donc $\rc_k(\rho)=0$.
\item Si $\rv(\zeta)=k$, nous avons $\zeta= z_k \vep^k(1+\vep \alpha)$ avec $z_k>0$ et $\alpha\in\gZ$, et $\xi=\vep^k\beta$ avec $\beta\in\gZ$. 
Alors on doit avoir l'égalité \fbox{$\rho=z_k^{-1}\beta^2\vep^k(1+\vep \alpha)^{-1}$}. Comme cette \egt implique $\kappa_{k-1}(\rho)=0$, elle est compatible avec les \coes de $\rho$ calculés jusqu'à l'exposant $k-1$.
Cette \egt permet de définit $\rc_{k+m}(\rho)$ pour tous les $m>0$.
\item Enfin si $\kappa_{k-1}(\zeta)=\kappa_k(\zeta)=1$, on cherche l'exposant $\ell<k$ tel que $\kappa_{\ell-1}(\zeta)=0$ et $\kappa_\ell(\zeta)=1$, et l'on est ramené au cas précédent via le calcul de la série $\rho$. 
\end{itemize}
Donc $\rho$ est bien défini.
\end{enumerate}

Récapitulons les résultats obtenus.
%p
%: Proposition{propQ[[T]]}
\begin{proposition} \label{propQ[[T]]}
L'anneau $\gZ$ est un \alrd hensélien et un \asr réduit. En outre il satisfait les règles \tsbf{OTF}, \tsbf{FRAC}, \tsbf{Val1} et \tsbf{Val2}. 
\end{proposition}
%----------- fin proposition ----------------------------- 

Le test sur $\gZ$ pour $\forall \alpha (\alpha^2>0 \vuu \alpha=0)$ équivaut
à \tsbf{LPO}.

On ne peut pas non plus démontrer \cot que $\gZ$ est un anneau \sdz: l'hypothèse $\xi\zeta=0$ équivaut à $\abs\xi\vi\abs\zeta=0$, mais l'implication $\abs\xi\vi\abs\zeta=0\Rightarrow ((\abs\xi=0) \vuu(\abs\zeta=0))$ équivaut au principe \tsbf{LPPO}. 

Enfin, on ne peut pas démontrer que tout \elt régulier $\geq 0$ est strictement positif: c'est \eqv au principe de Markov \tsbf{MP}. 
L'anneau total des fractions de $\gZ$ est donc un objet un peu mystérieux, un anneau qui contient $\gQ$ et qui fort heureusement ne présente pas d'intérêt \mathe évident.

\smallskip \Note Les articles \cite{KL00,KLP03} proposent une théorie \cov des anneaux de valuation (sans relation d'ordre) uniquement dans le cas des anneaux intègres avec une relation de \dve décidable. Il serait utile de généraliser les résultats (obtenus \cot) à des \advs au sens donné pour $\gZ$, et à d'autres cas similaires (on a besoin d'une relation de séparation sur l'anneau)\footnote{On pourra s'appuyer sur l'article \cite{ALP08}.}. \eoe

\smallskip On déduit facilement de la proposition \ref{propQ[[T]]} le \tho suivant.
 
%p
%: theorem{propQ((T))}
\begin{theorem} \label{propQ((T))}
L'anneau $\gQ=\gZ[1/\vep]$ satisfait tous les axiomes de la théorie \Sa{Co} ainsi que l'axiome de Heyting ordonné. C'est un \alrd avec $\Rad(\gQ)=0$ (donc un corps de Heyting dans la terminologie de \cite{CACM} ou \cite{MRR}).
En bref, c'est un corps de Heyting non archimédien, \emph{non} discret et un corps ordonné au sens de la théorie \sa{Co}.
\end{theorem}
%----------- fin proposition ----------------------------- 

\noindent \textsl{Note.} Un \elt de $\gQ=\gZ[1/\vep]$ s'écrit $\vep^{j_0}\alpha$ avec $\alpha\in\gZ$ et $j_0\in\ZZ$, il peut être codé sous la forme $\gamma=(j_0,\alpha)$. 
L'\egt $(j_0,\alpha)=(j'_0,\alpha')$ est définie comme suit: 
\begin{itemize}
\item si $j_0\leq j'_0$, $\vep^{j'_0-j_0}\alpha=\alpha'$;
\item si $j'_0\leq j_0$, $\vep^{j_0-j'_0}\alpha'=\alpha$. 
\end{itemize}
On définit alors, pour $\gamma=(j_0,\alpha)$:
\begin{itemize}
\item $\kappa_k(\gamma)=\kappa_{k-j_0}(\alpha)$ (donc $\kappa_k(\gamma)=0$ pour $k<j_0$);
\item $\rc_k(\gamma)=\rc_{k-j_0}(\alpha)$ (donc $\rc_k(\gamma)=0$ pour $k<j_0$);
\item $\rv(\gamma)=\rv(\alpha)-j_0$ (donc $\rv(\gamma)\geq -j_0$). \eoe
\end{itemize}

%%%%%%%%%%%%%%%%%%%%%%%%%%%%%%%%%%%%%%%%%%%%%%%%%%%%%%%%%%%%%%%%%%%
%%%%%%%%%%%%%%%%%%%%%%%%%%%%%%%%%%%%%%%%%%%%%%%%%%%%%%%%%%%%%%%%%%%

% Subsection
\section{Corps réels clos \textsl{non} discrets: position du \pb}
\label{secCRCnondis}

Notre rêve est de répéter l'exploit que Artin, Schreier et Tarski ont réalisé pour la description des \prts \agqs de $\RR$ à travers la théorie des corps réels clos discrets, mais dans un cadre constructif, en logique intuitionniste sans principe du tiers exclu, en tenant compte du fait que $\RR$ \textsl{n'est pas} discret, et en évitant l'axiome du choix dépendant.

%r

%r
%: Remark{remRRtdy} 
\begin{remark} \label{remRRtdy} 
On peut considérer que notre quête est la suivante: fixer une signature~$\Sigma$ qui permette de décrire
 aussi précisément que possible la structure de corps réel clos \textsl{non} discret, décrire sur cette signature une \tdy
qui soit \eseq à une théorie plus faible que \Sa{Crcd}, tout en étant
la plus forte possible parmi les \tdys 
qui admettent $\RR$ comme modèle \cof, sans utiliser l’axiome du choix dépendant. Ce Graal semble hors de portée dans l'absolu, car on n'a pas de critère clair pour savoir si une \rdy est satisfaite \cot sur $\RR$, l'axiome du choix dépendant n'étant pas autorisé dans les \demos.
\eoe
\end{remark}
%----------- fin remark ---------------------------------- 

%%%%%%%%%%%%%%%%%%%%%%%%%%%%%%%%%%%%%%%%%%%%%%%%%%%%%%%%%%%%%%%%%%%%
%%%%%%%%%%%%%%%%%%%%%%%%%%%%%%%%%%%%%%%%%%%%%%%%%%%%%%%%%%%%%%%%%%%%
\Subsection{Le principe de prolongement par continuité}
%: Subsection{Le principe de prolongement par continuité}

La propriété de \gui{completion} de $\RR$ s’exprime naturellement sous la forme suivante, sans interférence avec l’axiome du choix dépendant.

%t
%: Theorem{thRRcomplet}
\begin{theorem} \label{thRRcomplet}
Si une fonction $f\colon \QQ^n\to\RR$ est uniformément continue sur tout borné
elle se prolonge de manière unique en une fonction $\wi f\colon \RR^n\to\RR$
uniformément continue sur tout borné. 
\end{theorem}
%----------- fin theorem ---------------- 

Ce \tho est un \tho d’analyse et ne peut pas s’exprimer directement
dans le cadre d’une \tdy qui vise les \prts \agqs de $\RR$, car la propriété \gui{être uniformément continue} n'est pas \gmq.

Néanmoins c’est essentiellement ce \tho qui nous guide dans notre quête exprimée dans la remarque \ref{remRRtdy}. Nous remplacerons pour cela la \prt \gui{être uniformément continue} par une formulation où la continuité uniforme est contrôlée à priori et ne cache plus de $\forall\exists$.

En outre, les seules fonctions que l'on puisse envisager d'introduire dans un cadre purement \agq sont les fonctions semi\agqs.

Nous devons donc nous fonder sur des propriétés pertinentes des \fsagcs, que nous
développons dans le paragraphe qui suit. 
 
%%%%%%%%%%%%%%%%%%%%%%%%%%%%%%%%%%%%%%%%%%%%%%%%%%%%%%%%%%%%%%%%%%%%
%:Subsection
\Subsection{Fonctions semialgébriques continues}

 Tout d’abord nous rappelons que la continuité uniforme sur tout borné d’une \fsagc $\gR^n\to\gR$, où $\gR$ est un \crcd, est contrôlée \gui{à la \L{ojasiewicz}}
\prmt comme suit.

%t
%: fact{factfsagcLoja}
\begin{lemma} \label{factfsagcLoja}
Soit $\gR$ un \crcd et $K\subseteq \gR^n$ un fermé semialgébrique borné. 
\begin{enumerate}
\item 
Soit $g\colon K\to\gR$ une \fsagc. Alors $g$ possède un module de continuité uniforme qui s'exprime à la \L ojasiewicz comme suit
(avec un $c\in\gR$ et $\ell$ entier~\hbox{$\geq 1$})
% equation label {eqfactfsagcLoja1}
\begin{equation} \label {eqfactfsagcLoja1}
\forall \uxi,\uxi'\in K\;\abs{g(\uxi)-g(\uxi')}^\ell
 \leq \abs c \, 
 \,\norm{\uxi-\uxi'}.
\end{equation}
\item 
Soit $f\colon \gR^n\to\gR$ une \fsagc. Alors $f$ possède un module de continuité uniforme sur tout borné qui s'exprime à la \L ojasiewicz comme suit
(avec \hbox{un $c\in\gR$} et des entiers $k,\ell$ $\geq 1$)
% equation label {eqfactfsagcLoja}
\begin{equation} \label {eqfactfsagcLoja}
\forall \uxi,\uxi'\in \gR^n\;\abs{f(\uxi)-f(\uxi')}^\ell
 \leq \abs c \, \big(1+\norm{\uxi}+\norm{\uxi'}\big)^k 
 \,\norm{\uxi-\uxi'}.
\end{equation}
% end-equation
%
\end{enumerate}
\end{lemma}
%----------- fin theorem ----------------------------- 
\begin{proof} C'est une conséquence du \tho 2.6.6 de \cite{BCR}
qui dit que sur un ensemble \sagq localement fermé, si l'on a deux 
\fsagcs~$F$ et $G$ telles que $G$ s'annule aux zéros de $F$, il existe un exposant $N$ et
une \fsagc $h$ telles que $G^N=hF$. Dans le cas compact on majore $h$ par une constante, dans le cas \gnl on majore $h$ par une fonction \pol.
On applique cela avec $F(\ux,\ux')=\norm{\ux-\ux'}$ \hbox{et $G(\ux,\ux')=\abs{f(\ux)-f(\ux')}$}. 
\end{proof}

%:paragraph{Paramétrisation continue des \fsagcs
\paragraph{Paramétrisation continue des \fsagcs}~

\noindent Nous présentons maintenant un résultat de paramétrisation disant que, {du point de vue des \fsagcs,
tout provient continument de ce qui se passe sur le sous-corps 
des réels algébriques~$\RRa$}. En d'autres termes, toute \fsagc~\hbox{$\gR^n\to\gR$} est un point à \coos dans $\gR$ d'une famille équicontinue définie sur $\RRa$.

L'idée est en fait une simple \gnn de la remarque suivante. La famille des \pols univariés $f(x)=ax^2+bx+c$ (famille paramétrée par $(a,b,c)\in\RR^3$) n'est jamais 
que le \pol en quatre variables $g(a,b,c,x)=ax^2+bx+c$ 
{défini sur $\QQ$}, où l'on prend $(a,b,c)$ pour paramètres et $x$ comme variable, tous dans $\RR$: plus trop de souci donc avec le caractère \textsl{non} discret de~$\RR$, puisque tout est défini sur $\QQ$.
Chaque fonction individuelle $f(x)$ (dépendant de paramètres pris dans $\RR^3$) est un point réel d'une
famille définie sur $\RRa$. Ce point réel provient du prolongement
par continuité à $\RR^4$ d'une fonction continue $\RRa^4\to\RRa$. 

%t
%: Theorem{thParamcontFsagc0}
\begin{theorem} \label{thParamcontFsagc0} 
Soit $\gR$ un \crcd et $f\colon \gR^n\to\gR$ une \fsagc. Il existe un entier $r\geq 0$, une \fsagc $g\colon \gR^{r+n}\to\gR$ définie sur~$\RRa$, et un \elt $\uy\in\gR^r$ tels que
$$
\forall \xn\in\gR\;\; f(\xn)=g(\yr,\xn).
$$ 
\end{theorem}
%----------- fin theorem ----------------------------- 

Ce résultat semble faire partie du folklore.
Nous donnons ici une démonstration inspirée par les conseils de Michel Coste
et Marcus Tressl. Elle n'est cependant pas entièrement \cov. Il faudrait par exemple pour cela faire une relecture \cov du chapitre 7 de \cite{BCR}. Voir la question \ref{questthParamcontFsagc0}.
\begin{proof}
La fonction $f$ a un graphe $F$ fermé \sagq réunion d'un nombre fini de fermés \sagqs de base $F_i=\sotq{(\ux,y)\in\gR^{n+1}}{p_{i}(\ux,y)\geq 0}$. Les \coes des $p_{i}$, qui sont dans $\gR$,
peuvent être vus comme des spécialisations de paramètres $c_{k}$ ($k\in\lrbm$) de sorte que l'on a des \pols $P_{i}(\uc,\ux,y)$ avec des paramètres $c_k$. Les inégalités $P_{i}(\uc,\ux,y)\geq 0$ définissent pour $i$ fixé un fermé \sagq $G_i\subseteq \gR^{m+n+1}$. La réunion des $G_i$, notée $G$, est un \sagq qui n'est pas suffisamment pertinent. On ajoute un paramètre $c_0$ et l'on va maintenant restreindre le domaine de variation des $c_k$ à un ensemble \sagq $S$ \gui{convenable}. On entend par \gui{convenable} le fait qu'est satisfaite la formule suivante
\[ 
\begin{array}{ccc} 
\hspace{-.5em} \forall \uxi,\uxi'\in \gR^n \,\forall \zeta,\zeta'\in \gR\;\;((\uc,\uxi,\zeta)\in G,\,(\uc,\uxi',\zeta')\in G) \;\Rightarrow\;\abs{\zeta-\zeta'}^\ell
 \leq \abs {c_0} \, \big(1+\norm{\uxi}+\norm{\uxi'}\big)^k 
 \,\norm{\uxi-\uxi'}.
% \\[.4em] 
% \forall \uxi\in \gR^n \,\forall \zeta \in \gR\;\;\;(\uc,\uxi,\zeta)\in G \;\Rightarrow\;\abs{\zeta}\leq \abs {c_0} \, \big(1+\norm{\uxi}\big)^k
 \end{array}
\]
où $k$ et $\ell$ sont des entiers pour lesquels la fonction $f$ satisfait ces inégalités (pour un certaine spécialisation de $c_0$). Notez que $S\subseteq \gR^{m+1}$. 
Il est clair que le \sagq $S$ est défini sur $\RRa$. 
Notons $H$ le sous-ensemble \sagq de $\gR^{m+n+2}$ formé par les points de $G$ dont les $m+1$ premières coordonnées (les paramètres) forment un point de $S$. L'ensemble \sagq $H$ est le graphe d'une fonction $h\colon S\times \gR^n\to \gR$, qui est vu comme une famille de fonctions $\gR^n\to\gR$ paramétrée par $S$.
Pour tout point $s\in S$
le graphe $H_s$ correspondant est celui d'une \fsagc $f_s:\gR^n\to\gR$
dont le module de continuité uniforme est contrôlé par $\abs{c_0}$,~$k$ et~$\ell$. La fonction $f$ de départ correspond à un point $s_0\in S$ à \coos dans $\gR$.
Au moyen d'une décomposition algébrique cylindrique de $\gR^{m+1}$ adaptée à $S$, on insère $s_0$ dans une cellule~$\Gamma$ définie sur $\RRa$ semialgébriquement homéomorphe à $\gR^q$ pour un $q\geq 0$ ($q=0$ implique que $s_0$ est à \coos dans $\RRa$). En outre l'homéomorphisme est défini sur $\RRa$.
On obtient alors une fonction semi\agq $\varphi\colon \gR^{q+n}\to\gR$ définie sur $\RRa$ qui a les \prts suivantes:
\begin{itemize}
\item Il y a un \elt $\gamma=(\gamma_1,\dots,\gamma_q)\in \gR^q$ tel que $\varphi(\gamma,\uxi)=f(\uxi)$ pour tout $\uxi\in\gR^n$.
\item Pour tout \elt $\alpha=(\alpha_1,\dots,\alpha_q)\in \gR^q$, la fonction
$\uxi\mapsto \varphi(\alpha,\uxi)$ est \sagc. 
%
%\item La fonction $\varphi$ est localement bornée. 
\end{itemize}
Sous ces hypothèses, la remarque 7.4.9 de \cite{BCR}, nous assure qu'il existe une partition semi\agq $A_1\cup\dots\cup A_k$ de l'espace des paramètres $\gR^q$, définie sur $\RRa$, telle que la fonction $\varphi$ restreinte à chacun des $A_i\times \gR^n$
est continue. Par exemple le point $\gamma=(\gamma_1,\dots,\gamma_q)$ appartient à $A_1$. Au moyen d'une décomposition algébrique cylindrique de $\gR^{q}$ adaptée à $A_1$, on insère
$\gamma$ dans une cellule $\Delta$ définie sur $\RRa$ semialgébriquement homéomorphe à $\gR^r$ pour un $r\geq 0$. En outre l'homéomorphisme est défini sur $\RRa$. Ceci fournit la \fsagc $g\colon \gR^{r+n}\to \gR$ définie sur $\RRa$ réclamée dans l'énoncé.
\end{proof}

\paragraph{Une définition raisonnable.}
Nous proposons donc la \dfn suivante en \coma, rendue légitime par le \thref{thParamcontFsagc0}.

%d
%: Definition{defiFSAGC2}
\begin{definota} \label{defiFSAGC2} 
Soit $\gR$ un 
sous-corps ordonné\footnote{\Prmt $\gR$ est un sous-objet de $\RR$ pour la structure de corps ordonné \textsl{non} discret définie par la théorie~\sA{Co}. En outre, pour la \rex \Tsbf{IV}, on demande qu'un \elt de $\gR$ inversible dans $\RR$ soit inversible dans~$\gR$.} de $\RR$ contenant le corps des réels \agqs~$\RRa$ et 
soit une fonction $f\colon \gR^n\to \gR$.%
\index{fonction semialgébrique continue!sur un 
sous-corps ordonné de $\RR$ contenant $\RRa$}
\begin{enumerate}
\item (cas \elr) La fonction $f$ est \sagc s'il existe une \fsagc $g\colon \RRa^{n}\to\RRa$ dont $f$ est le prolongement par continuité. \Prmt on doit avoir les deux \prts suivantes:
 $f$ prolonge $g$, et
 $f$ a le même module de continuité uniforme que $g$, donné dans le point \textsl{2} du lemme \ref{factfsagcLoja}.
\item (cas \gnl) La fonction $f$ est \sagc s'il existe un entier $r\geq 0$, des \elts $y_1,\dots,y_r\in \gR$ et une fonction $h\colon \gR^{r+n}\to \gR$ qui relève du cas \elr précédent tels que 
$$
\forall \xn\in\gR\;\; f(\xn)=h(\yr,\xn).
$$
\end{enumerate}
On note $\Sac_n(\gR)$ l'anneau de ces fonctions (c'est un \afr réduit pour la relation d'ordre naturelle). 
\end{definota}
%----------- fin definition -------------------------------- 

Quelques propriétés importantes de ces espaces de fonctions seront établies dans la section \ref{PropGenFsagcs}.

%%%%%%%%%%%%%%%%%%%%%%%%%%%%%%%%%%%%%%%%%%%%%%%%%%%%%%%%%%%%%%%%%%%%
%:Subsection
\Subsection{Théories dynamiques raisonnables pour l'\alg des nombres réels} \label{secCrcdesirable}

On rappelle que le corps $\RRa$ des réels algébriques 
 est un \crcd au sens \cof.
 
Dans le prolongement de la remarque \ref{remRRtdy} et de la \dfn \ref{defiFSAGC2},
voici maintenant les \prts que nous avons en vue pour une \tdy \sa{Crc} 
des corps réels clos (\textsl{non} discrets), décrites ici de manière plutôt informelle.%

%p
%: Proriétés{prptaCrcdesirable}
\begin{prpta} \label{prptaCrcdesirable}~ 
\begin{enumerate}
\item La théorie \sa{Crc} est une extension de \Sa{Co}.
\item Les corps $\RR$, $\RR_{\tsbf{PR}}$, $\RR_{\tsbf{Ptime}}$ et $\RR_{\tsbf{Rec}}$ (cf. exemple \ref{exacorpsnondiscret}) sont des modèles \cof de \sa{Crc} (sans utiliser l'axiome du choix dépendant).
\item La théorie \sa{Crc} devient \eseq à \Sa{Crcd}
lorsqu’on lui ajoute l'axiome~\tsbf{ED$_>$}.
\item \label{Ra} Les \fsagcs $\RRa^n\to\RRa$ sont définies de manière agréable dans le langage
de \sa{Crc} et les \ralgs qu’elles satisfont sont valides dans la théorie. 
\item \label{prolcont} Des principes de prolongement par continuité (les plus larges possibles) sont satisfaits sous une forme convenable dans la \tdy. 
\item \label{prirec} Des principes de recollement (les plus larges possibles)
pour des fonctions définies sur un recouvrement fini par des ouverts \sagqs, ou par des \fsas, sont
satisfaits sous une forme convenable dans la \tdy.
\end{enumerate}
Les points qui suivent sont sujets à discussion.
Le point \textsl{\ref{lsalg}} sera abandonné si l'on veut décrire plus de \gui{structure} sur $\RR$, par exemple pour une structure o-minimale. Le point \textsl{\ref{lprimrec}} sera abandonné par exemple si l'on désire introduire tous les réels comme constantes de la théorie: dans une \tdy générale $\sa T=(\cL,\cA)$, $\cL$ et $\cA$ sont seulement supposés être des ensembles naïfs (à la Bishop). 
\begin{enumerate}\setcounter{enumi}{6}
\item \label{lsalg} Tous les symboles de fonction de \sa{Crc} définissent sur $\RR$ des \fsagcs de leurs variables (\dfn \ref{defiFSAGC2}). 
\item \label{lprimrec} Le langage de \sa{Crc} est énuméré de manière naturelle et dans ce cadre les axiomes sont décidables
de manière primitive récursive. 
\end{enumerate}
 
\end{prpta}
%----------- fin proposition ----------------------------- 

Une manière un peu brutale d'obtenir une réponse relativement satisfaisante est de prendre au sérieux le point \textsl{\ref{Ra}} ci-dessus. 
%Le point \textsl{\ref{dur}} sera alors satisfait si l'on prend la \dfn \ref{defiFSAGC2} pour les \fsagcs. 
Voici ce que cela donne.

%d
%: Definition{defiCrc1}
\begin{definition} \label{defiCrc1} La \tdy \SA{Crc1} est obtenue
à partir de la \tdy \Sa{Co} en ajoutant un symbole de fonction et
des axiomes convenables pour chaque \fsagc $f\colon \RRa^n\to\RRa$.
\end{definition}
%----------- fin definition -------------------------------- 

\noindent \textsl{Explication.}
Plus \prmt on procède comme suit. On sait d'après le \tho de finitude
(\cite[theorem~2.7.1]{BCR}) que le graphe 
$G_f=\sotq{(\ux,y)}{\ux\in\gR^n,y=f(\ux)}$ de $f$ (supposée semialgébrique continue)
est un \fsa de~$\RRa^{n+1}$ qui peut être décrit comme l'ensemble des zéros d'une \textsl{fonction \spo} $F\colon\RRa^{n+1}\to \RRa$, \textsl{i.e.} une fonction qui s'écrit sous la forme
\[
\sup\nolimits_i\,(\inf\nolimits_{ij}\,p_{ij}) \quad \hbox{ où }p_{ij}\in\RRa[\xn,y]
\]
On sait décider si un fermé semialgébrique $G_f$ ainsi décrit est celui d'une \fsagc, et dans ce cas calculer un module de continuité uniforme à la {\L}ojasiewicz.
Chaque fois qu'une telle (description de) fonction \spo définit une \fsagc, nous introduisons un symbole de fonction $\fsac_F$ avec 
l'axiome correspondant\label{Notafsa}

\UneRegle{Df$_{F}$} {$\vd F(\ux,\fsac_F(\ux))=0$.}

Par ailleurs, pour un terme arbitraire $t(\ux)$ dans le langage ainsi défini, si ce terme définit une fonction partout nulle sur $\RRa^n$ ($n\geq 0$), on introduit l'axiome correspondant $\vd t(\ux)= 0$. 

Par exemple on aura pour la fonction $\fsac_F$ un axiome de continuité qui reprend celui qui est satisfait pour les réels \agqs.

\UneRegle{Cont$_F$}{$\vd\abs{\fsac_F(\ux)-\fsac_F(\ux')}^\ell
 \leq \abs c \, \big(1+\norm{\ux}+\norm{\ux'}\big)^k 
 \,\norm{\ux-\ux'}.$}

\noindent En effet, une inégalité entre deux termes, $t_1\leq t_2$, équivaut à l'annulation du terme $(t_2-t_1)^-$. 
\eoe 

\smallskip 
Naturellement, une telle \tdy est à priori frustrante, car elle n'est pas très naturelle et sans doute elle est difficile à pratiquer d'un point de vue concret.

Nous verrons cependant qu'une manière plus naturelle, avec l'ajout de peu de symboles de fonction,
que nous proposons par la suite, aboutit à la théorie \Sa{Corv}
\esid~à~\Sa{Crc1}.

Tout ceci est étroitement lié à la théorie des anneaux réels clos et à sa réécriture sous forme concrète dans \cite{Tre2007}.

%%%%%%%%%%%%%%%%%%%%%%%%%%%%%%%%%%%%%%%%%%%%%%%%%%%%%%%%%%%%%%%%%%%%
%section
\section{Propriétés générales des \fsagcs}\label{PropGenFsagcs}

Nous donnons dans cette section quelques \prts remarquables des anneaux
$\Sac_n(\gR)$ (\fsagcs $\gR^n\to\gR$ \textsl{selon la \dfn \ref{defiFSAGC2}}) sur le corps $\gR=\RR$.
Plus \gnlt on peut considérer un sous-corps ordonné $\gR$ de $\RR$ contenant $\RRa$ et dans lequel toute \fsagc définie sur $\RRa$ prend ses valeurs dans $\gR$ aux points dont les coordonnées sont dans $\gR$, par exemple $\RR_{\tsbf{PR}}$, $\RR_{\tsbf{Ptime}}$ ou $\RR_{\tsbf{Rec}}$.

Ces \prts connues pour les \crcs discrets se prolongent à $\RR$ car on prend la précaution de ne prendre que des \prts dont la formulation n'implique pas le caractère discret de l'ordre. 

\Subsection{Stabilité par composition} 

 Par exemple on compose $f,g\in\Sac_3(\gR)$ avec $h\in\Sac_2(\gR)$.
Supposons que 
\begin{itemize}
\item $f$ est donnée sous la forme $f(x,y,z)=\wi f(a,b,x,y,z)$, pour $a,b\in\gR$
et $\wi f\colon \gR^5\to\gR$ prolonge par continuité $\ov f\colon \RRa^5\to\RRa$,
\item $g$ est donnée sous la forme $g(x,y,z)=\wi g(c,x,y,z)$, pour $c\in\gR$
et $\wi g\colon \gR^4\to\gR$ prolonge par continuité $\ov g\colon \RRa^4\to\RRa$,
\item $h$ est donnée sous la forme $h(u,v)=\wi h(d,u,v)$, pour $d\in\gR$
et $\wi h\colon \gR^3\to\gR$ prolonge par continuité $\ov h\colon \RRa^3\to\RRa$,
\item alors $h\circ(f,g):\gR^3\to\gR$ est de la forme 
$k(x,y,z)=\wi k(a,b,c,d,x,y,z)$ pour $(a,b,c,d)\in\gR^4$
et $\wi k:\gR^7\to\gR$
 prolonge par continuité la fonction 
$\ov k:\RRa^7\to\RRa$ définie par $$\ov k(a,b,c,d,x,y,z)=h(d,f(a,b,x,y,z),g(c,x,y,z)).$$
\end{itemize}

\Subsection{Stabilité par borne supérieure}

 Par exemple on a $f\in \Sac_4(\gR)$ et l'on veut montrer qu'il y a
une $g\in \Sac_2(\gR)$ telle que $g(x,y)=\sup_{z,t\in[0,1]}f(x,y,z,t)$.
Si $f$ est donnée sous la forme $f(x,y,z,t)=\wi f(a,b,x,y,z,t)$, pour $a,b\in\gR$,
où $\wi f\colon \gR^6\to\gR$ prolonge par continuité $\ov f\colon \RRa^6\to\RRa$, on considère la \fsagc $\ov g\colon \RRa^4\to\RRa$ définie par\footnote{Ici, $a,b$ sont des variables.}
$\ov g(a,b,x,y)=\sup_{z,t\in[0,1]}\ov f(a,b,x,y,z,t)$. Elle se prolonge par continuité en une fonction $\wi g\colon \gR^4\to\gR$ et on définit $g(x,y)=\ov g(a,b,x,y)$. Le fait que $g$ est bien la borne supérieure voulue tient au fait que la borne supérieure sur un compact est une fonction continue des paramètres et que $\RRa$ est dense dans $\gR$. 

\Subsection{Propriétés de finitude}

 En \clama toute \fsagc $f\colon \RR\to\RR$ admet un tableau de signes et de variations fini. Mais ce tableau ne dépend pas continument des paramètres,
par exemple lorsque $f(x)=g(\an,x)$ pour des paramètres $\an\in\RR^n$, où $g$ est le prolongement par continuité d'une \fsagc $\ov g\colon \RRa^{n+1}\to\RRa$.

Néanmoins lorsque $f$ est un \pol \unt, on a une approche \cov de la question en utilisant les fonctions \ravs. Par exemple on peut voir les points \textsl{\ref{ivrBudan}}, \textsl{\ref{ivrTVI}}, \textsl{\ref{ivrExtrema}} et \textsl{\ref{ivrMinAbs}} du \thref{thVirtualRoots} ainsi que la proposition \ref{factTableaucomplet}. 

On doit établir des résultats analogues en \coma pour les \fsagcs arbitraires, au moins pour la proposition \ref{factTableaucomplet}, mais en se limitant aux fonctions sur l'intervalle $[-1,1]$ (sur $\RR$ ce ne serait pas possible). L'utilisation d'un $\Vou$ infini sera peut-être \ncr.

%%%%%%%%%%%%%%%%%%%%%%%%%%%%%%%%%%%%%%%%%%%%%%%%%%%%%%%%%%%%%%%%%%%%
%%%%%%%%%%%%%%%%%%%%%%%%%%%%%%%%%%%%%%%%%%%%%%%%%%%%%%%%%%%%%%%%%%%%

%%%%%%%%%%%%%%%%%%%%%%%%%%%%%%%%%%%%%%%%%%%%%%%%%%%%%%%%%%%%%%%%%%%%
%%%%%%%%%%%%%%%%%%%%%%%%%%%%%%%%%%%%%%%%%%%%%%%%%%%%%%%%%%%%%%%%%%%%
%Subsection
\section{Quelques questions} 

%q
%: Question{questthParamcontFsagc0}
\begin{question} \label{questthParamcontFsagc0} 
Donner une \prco complète du \thref{thParamcontFsagc0}.
\end{question}
%----------- fin question ----------------------------- 

%q
%: Question{questeqratcont}
\begin{question} \label{questeqratcont} Montrer que la \ralg \Tsbf{FRAC}
n'est pas valide dans \Sa{Co--}.
Montrer de même que la \ralg \pref{eqratcont2} qui équivaut à l'existence de la fonction \pref{eqratcont} \paref{subsecratcon} n'est pas valide 
dans~\Sa{Co--}. 
\end{question}
%----------- fin question ----------------------------- 

%q
%: Question{quest17thRR}
\begin{question} \label{quest17thRR} 
Déterminer quelles \prts \agqs de $\RR$ permettent de démon\-trer
les formes \cot satisfaisantes de Positivstellensätze prouvées pour~$\RR$ dans~\cite{GL93}. Voir en particulier la question \ref{quest17H}. 
\end{question}
%----------- fin question ----------------------------- 

\Subsection{Variations continues}

On aurait pu définir comme suit les \fsagcs
$\RR^n\to\RR$. C'était la \dfn adoptée dans \cite[Définition 3.3]{LM2017}.

%d
%: Definition{defiAlgAfrnz}
\begin{definition} \label{defiAlgAfrnz}
Soit $\gR$ un anneau commutatif. Une fonction $f\colon \gR^n\to\gR$ est dite \textsl{\agq sur $\Rxn=\Rux$} si l'on a un \pol \hbox{$g(\ux,y)=\sum_{k=0}^mg_k(\ux)y^k\in\gR[\ux,y]$},
avec au moins un des \coes d’un $g_k(\ux)$ \iv,
tel que $g(\uxi,f(\uxi))=0$ pour tout $(\uxi)\in\gR^n$.
\end{definition}
%----------- fin definition -------------------------------- 

%d
%: Definition{defiFSAGC}
\begin{definition}[\dfn alternative à \ref{defiFSAGC2}] \label{defiFSAGC} 
Soit $\gR$ un 
%\gui{corps ordonné}, en tout cas un modèle de la théorie~\Sa{Co--}.
sous-corps de $\RR$.\\
Une fonction $f\colon \gR^n\to \gR$ est dite \textsl{\sagc} si, d'une part, elle est \agq sur $\Rux$ et si, d'autre part, elle possède un module de continuité uniforme sur tout borné à la {\L}ojasiewicz, donné par une inégalité~\pref{eqfactfsagcLoja} comme dans le lemme \ref{factfsagcLoja}.\index{fonction semialgébrique continue!sur un sous-corps ordonné de $\RR$, \dfn alternative}
\end{definition}
%----------- fin definition -------------------------------- 

%\vspace{-.3em}
Cette \dfn est légitime pour les sous-corps de $\RR$ car 
%i
\begin{itemize}
\item elle est valable en \clama, 
\item elle a une signification \cov claire,
\item les fonctions 
\sagcs au sens de la \dfn \ref{defiFSAGC2} le sont aussi au sens de la \dfn \ref{defiFSAGC}.
\end{itemize}

%q
%: Question{questRRfsagc1}
\begin{question} \label{questRRfsagc1} 
 
Si une fonction $f\colon \RR^n\to \RR$ est \agq sur $\RR[\ux]$ et si elle est uniformément continue sur tout borné, est-ce qu’elle possède un module de continuité uniforme à la {\L}ojasiewicz, comme dans le lemme \ref{factfsagcLoja}?
\\
NB: la réponse est positive en \clama, mais elle semble 
nettement plus délicate en \coma.
\end{question}
%----------- fin question -----------------------------

%q
%: Question{questRRfsagc2}
\begin{question} \label{questRRfsagc2} Une fonction 
\sagc au sens de la \dfn \ref{defiFSAGC} l'est-elle aussi au sens de la \dfn \ref{defiFSAGC2}? 
\\
Oui en \clama, mais le \pb se pose en \coma, et semble fort difficile. 
Il se peut que nous soyons là, en préférant la \dfn \ref{defiFSAGC2} à la \dfn~\ref{defiFSAGC}, dans une situation analogue à celle qui a conduit Bishop à définir la continuité d'une fonction $\RR\to\RR$ comme signifiant la continuité uniforme sur tout intervalle borné.
\end{question}
%----------- fin question -----------------------------

\hum{La règle \Tsbf{FRAC} traite un cas de prolongement par continuité pour des fractions $g/p$ 
où les valeurs mal définies (aux zéros de $p$) sont en fait nulles. Pour cette règle, le dénominateur est n'importe quelle fonction définissable, et la fonction obtenue par prolongement par continuité est nulle aux points litigieux. C'est assez différent du cas du prolongement par continuité d'une fraction dont le dénominateur est un polynôme, comme dans la remarque \ref{remCo}. Le prolongement fonctionne alors, sous condition de continuité uniforme. avec des valeurs arbitraires aux points litigieux, car les zéros de $p$ forment une partie d'intérieur vide. 
}

\newpage \thispagestyle{empty}

%%%%%%%%%%%%%%%%% CHAPITRE %%%%%%%%%%%%%

\chapter{Anneaux \ftm \rtls}\label{chap-afr}
\Today
\minitoc

\section*{Introduction}
\addtocontents{toc}{\vskip0.8em}
\addcontentsline{toc}{section}{Introduction}
\rdb

Ce chapitre reprend la problématique des corps ordonnés \textsl{non} discrets
en repartant de zéro. 

Toutes les théories introduites admettent des extensions \eseqs à la théorie \Sa{Co} des corps ordonnés \textsl{non} discrets (\dfn \ref{defiConondiscret}). 

On commence (section \ref{sectrdisad}) par la théorie des \trdis (un corps ordonné \textsl{non} discret est un \trdi pour sa relation d'ordre). 

Dans la section \ref{secgrl} on traite les \grls (théorie \peq), valable pour l'addition sur les réels (ce sont les $\ell$-groups dans la littérature anglaise).

Ensuite (section \ref{secfrings}) nous passons aux \afrs ($f$-rings dans la littérature anglaise), théorie inspirée par les anneaux de fonctions réelles continues. 

La section \ref{secArftr} décrit des \tdys dans lesquelles on ajoute la relation $\cdot>0$ (\asrs et variantes).

La section \ref{secCOG} propose un retour sur la théorie \Sa{Co} en la confrontant à des extensions convenables de la théorie des \asrs.

\smallskip Dans ce chapitre on dit \gui{groupe} pour \gui{groupe abélien}. Et les anneaux sont commutatifs unitaires comme dans tout le mémoire.

\section{Treillis distributifs}\label{sectrdisad}

Références: \cite{CC00,CL05,CLQ2006,Lor1951}

\Subsection{Théorie des treillis distributifs}\label{subsecTrdi}

La théorie \SA{Tr0} des treillis est une \tpe basée sur la signature suivante\footnote{Pour plus de précision, on peut préférer $\cdot=_\TR\cdot$, $\cdot\vi_\TR\cdot$, $\cdot\vu_\TR\cdot$, $0_\TR$ et $1_\TR$.}, avec la seule sorte $\TR$.
\Sigt{\TR}{\cdot=\cdot\mathrel{;}\cdot \vi \cdot, \cdot \vu \cdot, 1, 0}
\label{NOTASigTr}

\vspace{-1em}

\noindent Les symboles $\vu$ et $\vi$ utilisés pour la \bsp et la \bif binaires
 ne doivent pas être confondu avec les symboles~$\vuu$ et $\vii$ de la disjonction et de la conjonction logiques.

 Outre les axiomes de l'\egt on a les axiomes suivants

\DeuxRegles{
 \labu $\vd 0\vi x=0$
 \labu $\vd x\vi x=x$
 \labu $\vd x\vi y=y\vi x$ 
 \labu $\vd (x\vi y)\vi z=x\vi (y\vi z)$ 
 \labu $\vd (x\vi y)\vu x=x$ 
}
{
 \labu $\vd 1\vu x=1$
 \labu $\vd x\vu x=x$
 \labu $\vd x\vu y=y\vu x$ 
 \labu $\vd (x\vu y)\vu z=x\vu (y\vu z)$ 
 \labu $\vd (x\vu y)\vi x=x$ 
} 

\noindent On définit $x\geq y$ comme une abréviation de $x=x\vu y$. On obtient ainsi une extension de la théorie des ensembles ordonnés.

\smallskip On démontre facilement les règles suivantes

\DeuxRegles{
 \labu $\vd 0\leq  x\leq 1$
}
{
 \labu $\,\,1=0\vd  x=0$
 }

\smallskip La théorie \SA{Tr} des \textsl{treillis non triviaux} s'obtient en ajoutant l'axiome de collapsus

\Regles {\lAb{CL$_{\Tr}$} $\,\,1=_{\Tr} 0\vd \Bot$}

\smallskip La théorie \SA{Trdi} des \textsl{\trdis} est obtenue en ajoutant l'axiome de distributivité suivant
(l'axiome dual, \cad obtenu en renversant la relation d'ordre, s'en déduit)

\Regles { \labu $\vd (x\vu y)\vi z=(x\vi z)\vu (y\vi z)$ }

\smallskip 
La théorie \SA{Etob} des ensembles totalement ordonnés bornés s'obtient à partir de \sa{Trdi} en ajoutant l'axiome d'ordre total 

\Regles {\labu $\vd x=x\vi y\vou y=x\vi y$} 

\smallskip  
La théorie \SA{AgB} des \textsl{algèbres de Boole} s'obtient à partir de \sa{Trdi} en ajoutant l'axiome 

\Regles {\labu $\vd \Exists y\,(x\vi y=0\vet x\vu y=1)$}

\smallskip  
La théorie \SA{ABo} des \textsl{anneaux de Boole} s'obtient à partir de la théorie des anneaux commutatifs non triviaux \sa{Ac} en ajoutant l'axiome 
 
\Regles {\labu $\vd x^2=x$}

Il est bien connu que les théories \sa{AgB} et \sa{ABo} sont \esids.

\smallskip 
La théorie \SA{AgBo} de l'algèbre de Boole \hbox{$\BB=\so{\zB,\uB}$} s'obtient à partir de \sa{Trdi} en ajoutant l'axiome 

\Regles {\labu $\vd x=0\vou x=1$}

%:subsection{Idéaux et filtres dans un \trdi}---- 
\Subsection{Idéaux et filtres dans un \trdi}\label{subsecTrdiIdeFi}
%-----------
 Un \textsl{idéal} $\fb $ d'un \trdi $(\gT,\vi,\vu,0,1)$ est une
partie qui satisfait les contraintes:
%--- equation eqIdeal --------
\begin{equation}\label{eqIdeal}
\left.
\begin{array}{rcl}
 & & 0 \in \fb \\
x,y\in \fb & \Longrightarrow & x\vu y \in \fb \\
x\in \fb ,\; z\in \gT& \Longrightarrow & x\vi z \in \fb \\
\end{array}
\right\}
\end{equation}
%---------------------end equation--------------
On note $\gT/(\fb=0)$ le treillis quotient obtenu en forçant les \elts de $\fb$ à être nuls. On peut aussi définir les \ids comme les noyaux de morphismes etre \trdis.

Un \textsl{\idp} est un \id engendré par un seul \elt $a$, il est noté $\dar a$.
On a $\dar a=\sotq{x\in\gT}{x\leq a}$. 
L'\id $\dar a$, muni des lois $\vi$ et $\vu$ de $\gT$ est un \trdi
dans lequel l'\elt maximum est $a$. L'injection canonique $\dar
a\rightarrow \gT$ n'est pas un morphisme de \trdis parce que
l'image de $a$ n'est pas égale à $1_\gT$. Par contre
l'application $\gT\rightarrow \dar a,\;x\mapsto x\vi a$
est un morphisme surjectif, qui définit donc $\dar a$ comme une
structure quotient $\gT/(a=1)$.

La notion de \textsl{filtre} est la notion opposée (obtenue en renversant la relation d'ordre) de celle d'idéal.

Soient $\fa$ un \id et $\ff$ un filtre de $\gT$.
On dit que $(\fa,\ff)$ est un \textsl{couple saturé} dans~$\gT$ si
$$
(g\in \ff,\; x\vi g \in \fa) \Longrightarrow x\in \fa,
\hbox{ et }(a\in \fa,\; x\vu a \in \ff) \Longrightarrow x\in \ff.
$$
Un \textsl{couple saturé} est un couple $(\varphi^{-1}(0),\varphi^{-1}(1))$ pour un morphisme $\varphi\colon \gT\to\gT'$ de \trdis.
Lorsque $(\fa,\ff)$ est un couple saturé, on a
les \eqvcs
\[
1\in \fa\; \;\Longleftrightarrow\;\; 0\in \ff
\;\; \Longleftrightarrow\;\; (\fa,\ff)=(\gT ,\gT )
\]
Si $A$ et $B$ sont deux parties de $\gT$ on note
%--- equation eqvuvi --------
\begin{equation}\label{eqvuvi}
A\vu B=\sotq{ a\vu b} {a\in A,\,b\in B} \; \hbox{ et } \; A\vi
B=\sotq{ a\vi b} {a\in A,\,b\in B}.
\end{equation}
%---------------------end equation--------------
Alors l'\id engendré par deux \ids $\fa$ et $\fb$ est égal à
 $\fa\vu\fb$.
L'ensemble des \ids de $\gT$ forme lui-même un \trdi\footnote{En fait
il faut introduire une restriction pour obtenir vraiment un ensemble, de fa\c{c}on à ce que l'on ait un procédé bien défini de construction des \ids concernés.
Par exemple on peut considérer l'ensemble des \ids obtenus à partir des \idps par certaines opérations prédéfinies, comme les réunions et intersections dénombrables.} pour
l'inclusion, avec pour borne inférieure de $\fa$ et $\fb$ l'\id 
 $\fc=\fa\vi\fb$.
Ainsi les opérations $\vu$ et $\vi$ définies en (\ref{eqvuvi})
correspondent au sup et au inf dans le treillis des \ids.

Quand on considère le treillis des filtres, il faut
faire attention à ce que produit le renversement de la relation
d'ordre: $\ff\cap\ffg=\ff\vu\ffg$ est le inf des filtres~$\ff$ et $\ffg$,
tandis que leur sup est le treillis engendré par $\ff\cup \ffg$, égal à $\ff\vi\ffg$.

%:--- SUBsection{Quotients}-----------
\Subsection{Quotients}

Un \textsl{\trdi quotient $\gT'$ de $\gT$} est donné par une relation 
binaire
$\preceq$ sur $\gT$ vérifiant les propriétés suivantes:
%--- equation eqPreceq --------
\begin{equation} \label{eqPreceq}
\left.
%--------------------begin array---------------
\begin{array}{rcl}
a\leq b& \Longrightarrow & a\preceq b \\
a\preceq b,\,b\preceq c& \Longrightarrow & a\preceq c \\
a\preceq b,\,a\preceq c& \Longrightarrow & a\preceq b\vi c \\
b\preceq a,\,c\preceq a& \Longrightarrow & b\vu c\preceq a
\end{array}
%---------------------end array--------------
\right\}
\end{equation}
%---------------------end equation--------------

%--- Proposition{propIdealFiltre}----
\begin{proposition}
\label{propIdealFiltre} Soit $\gT$ un \trdi et
$(J,U)$ un couple de parties de $\gT$.
On considère le quotient $\gT'$ de $\gT$ défini par les
relations $x=0$ pour les $x\in J$ et $y=1$ pour les $y\in U$. Alors
 on a $a\leq_{\gT'}b$ \ssi
il existe une partie finie $J_0$ de $J$ et une partie finie $U_0$ de
$U$ telles que:
%--- equation eqpropIdealFiltre --------
\begin{equation} \label{eqpropIdealFiltre}
a \vi \Vi U_0 \; \leq_\gT\; b \vu \Vu J_0
\end{equation}
%---------------------end equation--------------
Nous noterons $\gT/(J=0,U=1)$ ce treillis quotient $\gT'.$
\end{proposition}
%--- end-proposition----------------------------------------

\Subsection{Théorème de représentation}

Le théorème constructif suivant donne en \clama le \tho de représentation qui dit que tout \trdi est un sous-\trdi de l'\agB des parties d'un ensemble. Ou encore celui qui dit que tout \trdi est un sous-objet d'un produit d'ensembles totalement ordonnés bornés.

%:     Theorem{thTrdiAgb}
\begin{theorem} \label{thTrdiAgb}
La théorie \Sa{Trdi} des treillis distributifs, la théorie \Sa{Etob} des ensembles totalement ordonnés bornés, la théorie \Sa{AgB} des algèbres de Boole et la théorie \Sa{AgBo} 
prouvent les mêmes règles algébriques. 
\end{theorem}
%----------- fin theorem ----------------------------- 
\begin{proof} On doit démontrer qu'un fait dans une \sad de type \sa{Trdi}, \cad dans un \trdi $\gT_0$, est prouvable dans une des 4 théories
\ssi il est prouvable dans \sa{Trdi}.  
Il suffit de le voir pour \sa{Trdi} et \sa{AgBo} car les autres sont des théories intermédiaires. Pour une structure algébrique dynamique de type $\sa{AgBo}$ une démonstration utilise un arbre de calcul gouverné par l'axiome \sul $\vd x=0\vou x=1$. A un nœud $\nu$ d'un tel arbre, les hypothèses accumulées donnent un treillis distributif $\gT_\nu$ quotient de $\gT_0$ obtenu en forçant certains \elts à être égaux à $0$ ou à $1$. 
Si $\nu$ n'est pas une feuille, en dessous du nœud $\nu$ on a deux branches: un certain \elt $y\in\gT_0$ définit un élément $\ov y\in\gT_r$, dans la première branche on force $\ov y=0$ et dans l'autre on force $\ov y=1$.
Si un fait $a=b$ de $\gT_\nu$ est prouvé dans les deux nouvelles branches, il est prouvable aussi dans $\gT_\nu$ car le morphisme canonique $\gT_\nu\to\gT_\nu/(\ov y=0)\times \gT_\nu/(\ov y=1)$ est injectif.
\end{proof}

Un corolaire est la proposition suivante.

%p
%:     Proposition{propTriTrdi}
\begin{proposition} \label{propTriTrdi} Soient $\xn$ dans un \trdi.\\
On définit $\Tri(\ux)=[\Tri_1(\ux), \Tri_2(\ux),\ldots,\Tri_n(\ux)]$, où 
\[\ndsp
\Tri_k(\xn) =\Vi_{I\in \cP_{k,n}}\big(\Vu_{i\in I}x_i\big) \quad(k\in\lrbn)\footnote{On note $\cP_{k,n}$ l'ensemble des parties à $k$ \elts de $\cP_{n}:=\so{1,\dots,n}$.}.
\]
\noindent 
On a les résultats suivants.
\begin{enumerate}
\item $\Tri_k(\xn) =\Vu_{J\in \cP_{n-k+1,n}}\big(\Vi_{j\in J}x_j\big)$, $(k\in\lrbn)$.
\item $\Tri_1(\ux)\leq \Tri_2(\ux)\leq \cdots \leq \Tri_n(\ux).$
\item Si les $x_i$ sont deux à deux comparables, la liste $\Tri(\ux)$
est la liste des~$x_i$ rangée en ordre croissant (il n'est pas \ncr que
le treillis soit discret).
\end{enumerate} 
\end{proposition}
%----------- fin proposition ----------------------------- 
%
\begin{proof}
C'est clair dans le cas d'un ensemble totalement ordonné borné.
\end{proof}

\Subsection{Recouvrements et recollements de \trdis}

En \alg commutative, si $\fa$ et $\fb$ sont deux \ids d'un anneau 
$\gA$
on a une \gui{suite exacte} de \Amos (avec $j$ et $p$ des \homos 
d'anneaux)
$$0\to\gA/(\fa\cap\fb)\vers{j}(\gA/\fa) \times 
(\gA/\fb)\vers{p}\gA/(\fa+\fb)\to 0$$
qu'on peut lire en langage courant: le système de congruences 
$x\equiv
a\;\mod\;\fa$, $x\equiv b\;\mod\;\fb$ admet une solution \ssi $a\equiv
b\;\mod\;\fa+\fb$ et dans ce cas la solution est unique modulo 
$\fa\cap\fb$.
Il est remarquable que ce \gui{\tho des restes chinois} se 
généralise à un
système \textsl{quelconque} de congruences \ssi l'anneau est
\textsl{arithmétique} (\cite[Theorem XII-1.6]{CACM}), \cad si le treillis des \ids est distributif.
Le \tho des restes chinois \gui{contemporain} concerne le cas 
particulier d'une
famille d'\ids deux à deux comaximaux, et il fonctionne sans hypothèse sur
l'anneau de base.

D'autres épimorphismes de la catégorie des anneaux commutatifs
sont les localisations. Et il~y~a un principe de recouvrement analogue au \tho 
des restes chinois pour les localisations, extrêmement fécond (le principe 
local-global).

\smallskip De la même manière on peut récupérer un \trdi à partir
d'un nombre fini de ses quotients,
si l'information qu'ils contiennent est \gui{suffisante}. 
On peut voir ceci au choix comme une procédure de 
recouvrement (de passage du local au global), ou comme une version du \tho des restes chinois pour les \trdis. Voyons les choses plus précisément.

%--- Definition{defRecouvTD}---------
\begin{definition}
\label{defRecouvTD}
Soit $\gT$ un \trdi, $(\fa_i)_{i=1,\ldots n}$ (resp. 
$(\ff_i)_{i=1,\ldots n}$)
une famille finie d'\ids (resp. de filtres) de $\gT$. On dit que 
les \ids
$\fa_i$ \textsl{recouvrent $\gT$} si $\Vi_i\fa_i=\so{0}$. De 
même on dit
que les filtres $\ff_i$ \textsl{recouvrent $\gT$} si 
$\Vu_i\ff_i=\so{1}$.
\end{definition}
%--- end-definition----------------------------

%--- theorem{propRecouvTD}---------------
\begin{theorem}[recouvrements d'un \trdi]
\label{propRecouvTD}~
\begin{enumerate}
\item Soit $\gT$ un \trdi, $(\fa_i)_{i=1,\ldots n}$ une famille finie 
d'\ids principaux de~$\gT$ et $\fa=\Vi_i\fa_i$.
%
%-----------------begin enum------------------
\begin{enumerate}
\item Si $(x_i)$ est une famille d'\elts de $\gT$ telle que pour 
chaque $i,j$ on
a $x_i\equiv x_j\;\mod\;(\fa_i\vu\fa_j=0)$, alors il existe un unique $x$ 
modulo
$\fa=0$ vérifiant: $x\equiv x_i\;\mod\;(\fa_i=0)$ pour chaque $i$.
\item Notons $\gT_i=\gT/(\fa_i=0)$, 
$\gT_{ij}=\gT_{ji}=\gT/(\fa_i\vu\fa_j=0)$,
$\pi_i:\gT\to\gT_i$ et $\pi_{ij}:\gT_i\to\gT_{ij}$ les projections 
canoniques.
Si les $\fa_i$ recouvrent $\gT$, alors $(\gT,(\pi_i)_{i=1,\ldots n})$ est 
la limite
projective du diagramme $((\gT_i)_{1\leq i\leq n},(\gT_{ij})_{1\leq 
i<j\leq
n};(\pi_{ij})_{1\leq i\neq j\leq n})$ (voir la figure ci-après).
\end{enumerate}
\item Soit maintenant $(\ff_i)_{i=1,\ldots n}$ une famille finie de 
filtres principaux %($\ff_i=\uar s_i$), 
et $\ff=\Vu_i\ff_i$ 
\begin{enumerate}
\item Si $(x_i)$ est une famille d'\elts de $\gT$ telle que pour 
chaque $i,j$ on
a $x_i\equiv x_j\;\mod\;(\ff_i\vi\ff_j=1)$, alors il existe un unique $x$ 
modulo $\ff=1$ vérifiant: $x\equiv x_i\;\mod\;(\ff_i=1)$ pour chaque $i$.
\item Notons $\gT_i=\gT/(\ff_i=1)$, 
$\gT_{ij}=\gT_{ji}=\gT/(\ff_i\cup\ff_j=1)$,
$\pi_i:\gT\to\gT_i$ et $\pi_{ij}:\gT_i\to\gT_{ij}$ les projections 
canoniques. Si les $\ff_i$ recouvrent $\gT$, alors $(\gT,(\pi_i)_{i=1,\ldots n})$ est 
la limite
projective du diagramme $((\gT_i)_{1\leq i\leq n},(\gT_{ij})_{1\leq 
i<j\leq
n};(\pi_{ij})_{1\leq i\neq j\leq n})$ (voir la figure ci-après).
\end{enumerate}

\end{enumerate}
%-----------------end enum------------------
\end{theorem}
%--- end-proposition-----------------------------------------
 
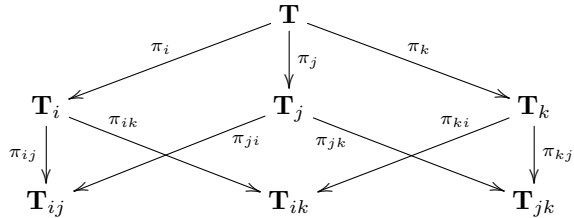
{\hspace*{10em}{
\xymatrix @R=2em @C=7em{
 & \gT \ar[rd]^{\pi _{k}}\ar[d]^{\pi _{j}}\ar[ld]_{\pi _{i}}\\
 \gT _i\ar[d]_{\pi _{ij}}\ar@/-0.75cm/[dr]^{\pi _{ik}} &
 \gT _j\ar@/-1cm/[dl]^{\pi _{ji}}\ar@/-1cm/[dr]_{\pi _{jk}} &
 \gT _k\ar@/-0.75cm/[dl]_{\pi _{ki}}\ar[d]^{\pi _{kj}} &
\\
 \gT _{ij} & 
 \gT _{ik} & 
 \gT _{jk} 
}
}}

\medskip Il y a aussi une procédure de recollement proprement dit (voir \cite[proposition 1.2.7]{CLQ2006}).

%--- theorem{thRecolTD}-------
\begin{theorem}[recollement de \trdis]
\label{thRecolTD}
Supposons donnés un ensemble fini totalement ordonné~$I$ et dans la catégorie des \trdis  un diagramme

\snic{\big((\gT_i)_{i\in I},(\gT_{ij})_{i<j\in I},(\gT_{ijk})_{i<j<k\in I};
(\pi_{ij})_{i\neq j},(\pi_{ijk})_{i< j, j\neq k\neq i}\big)}

\noindent 
comme dans la figure ci-après, 
ainsi qu'une famille d'\elts 

%$$
\snic
{(s_{ij})_{i\neq j\in I}\in \prod\nolimits_{i\neq j\in I}\gT_{i}}
%$$

\noindent satisfaisant les conditions suivantes:
\begin{itemize}
\item le diagramme est commutatif ($\pi_{ijk}\circ \pi_{ij}=\pi_{ikj}\circ \pi_{ik}$ pour tous $i$, $j$, $k$ distincts), 
\item pour $i\neq j$, $\pi_{ij}$ est un morphisme de passage au quotient par l'\id $\dar s_{ij}$,
\item pour $i$, $j$, $k$ distincts, $\pi_{ij}(s_{ik})=\pi_{ji}(s_{jk})$ et  $\pi_{ijk}$ est un morphisme de passage au quotient par \hbox{l'\id $\dar\pi_{ij}(s_{ik})$}.   
\end{itemize}

\smallskip 
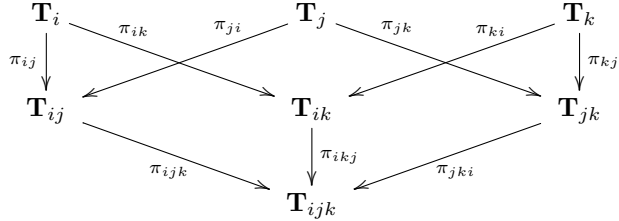
{\hspace*{10em}
\xymatrix @R=2em @C=7em{
 \gT_i\ar[d]_{\pi _{ij}}\ar@/-0.75cm/[dr]^{\pi _{ik}} &
     \gT_j\ar@/-.8cm/[dl]_{\pi _{ji}}\ar@/-.8cm/[dr]^{\pi _{jk}} &
        \gT_k\ar@/-0.75cm/[dl]_{\pi _{ki}}\ar[d]^{\pi _{kj}} &
\\
 ~\gT_{ij}~ \ar[rd]_{\pi _{ijk}} & 
    ~\gT_{ik}~  \ar[d]^{\pi _{ikj}} & 
      ~\gT_{jk}~  \ar[ld]^{\pi _{jki}} 
\\
   &  ~\gT_{ijk}~ 
\\
}
}

\smallskip \noindent Alors si $\big(\gT\,;\,(\pi_i)_{i\in I}\big)$ est la limite projective du diagramme, les~\hbox{$\pi_i:\gT\to \gT_i$} forment un recouvrement de $\gT$ par des quotients par des idéaux principaux, et le diagramme est isomorphe à celui obtenu
dans le \thref{propRecouvTD}.
Plus précisément, il existe des $s_i\in\gT$ tels que chaque~$\pi_i$ est un morphisme de passage au quotient par l'\id $\dar s_i$ et $\pi_i(s_j)=s_{ij}$ pour tous $i\neq j$.

\noindent Le résultat analogue est valable pour les quotients par des filtres principaux.
\end{theorem}
%--- end-proposition----------------------------------------

Les propositions précédentes de recouvrement et recollement ont des versions analogues pour la catégorie des modules sur un anneau commutatif et pour celle des \grls de \ddk $\leq 1$. Par contre dans la catégorie des anneaux commutatifs, seule la procédure de recouvrement (par la localisation en des \eco) est valable, et il a fallu que Grothendieck invente la catégorie des schémas pour avoir une procédure de recollement: un schéma quasi-compact quasi-séparé n'est autre que le recollement  d'un nombre fini de schémas affines le long d'ouverts quasi-compacts, ce qui correspond à un recollement \gui{abstrait} d'anneaux commutatifs dans le style du \thref{thRecolTD} pour les épimorphismes de localisation (voir \cite{CLS2009}).

\section{Groupes réticulés ($\ell$-groups)}\label{secgrl}

\Subsection{Groupes abéliens ordonnés et totalement ordonnés}

La théorie \SA{Gao} des \textsl{groupes (abéliens) ordonnés} est définie comme suit. Il y a une seule sorte, nommée $\Gao$.\index{groupe (abélien)!ordonné}

\vspace{-.5em}
\Sigt{\Gao}{\cdot=0,\cdot\geq 0\mathrel{;}\cdot+\cdot,-\,\cdot,0}
\label{NOTASigGao}

\noindent {\bf Abréviations}

\smallskip \noindent \textsl{Prédicats}

\vspace{-1em} \DeuxCols{
\begin{itemize}
\itbu $x = y $ signifie $ x - y = 0$
\itbu $x \leq y $ signifie $ y\geq x$
\end{itemize}
}
{\begin{itemize}
\itbu $x \geq  y $ signifie $ x - y \geq  0$
\end{itemize}
} 

\medskip \noindent {\bf Axiomes}

\smallskip \noindent {\it Règles directes des groupes abéliens}

\DeuxRegles{
\Lab{ga0} $\vd 0=0$
\Lab{ga2} $\,\,x=0\vd -x=0 $
}
{
\Lab{ga1} $\,\, x=0\vet y=0\vd x+y=0$
}

\smallskip \noindent NB. Les règles \tsbf{ga0}, \tsbf{ga1} et \tsbf{ga2} définissent la théorie purement équationnelle \SA{Ga} des groupes abéliens. On doit alors remplacer, dans l'explication donnée \paref{Ac-comments} pour les anneaux commutatifs (exemple \ref{exaAc}), la machinerie calculatoire des anneaux commutatifs $\ZZxn$ (librement engendrés par les $x_i$)
par celle des groupes abéliens $\ZZ^{\{\xn\}}$ (librement engendrés par les $x_i$).

\smallskip \noindent {\it Règles pour la relation d'ordre}

\DeuxRegles{
\Lab{gao0} $\vd 0 \geq 0$
\Lab{Gao} $\,\, x\geq 0\vet x\leq 0\vd x=0$
}
{
\Lab{gao1} $\,\, x\geq 0\vet y\geq 0\vd x+y\geq 0$
}

La réflexivité et transitivité de $\cdot\geq \cdot$, et la règle \tsbf{gao2} $\,\,x\geq y\vd x+z\geq y+z\,$ résultent directement des \dfns. La règle \tsbf{Gao} donne l'antisymétrie de $\cdot\geq \cdot$.

La théorie \SA{Gto} des \textsl{groupes (abéliens) totalement ordonnés}
est obtenue à partir de la théorie \Sa{Gao} en ajoutant comme axiome la \rdy \Tsbf{OT}.

\Subsection{\Dfn de la théorie purement équationnelle des \grls}

La théorie \SA{Grl} est définie comme suit. Il y a une seule sorte, nommée $\GRL$.\index{groupe!réticulé ($\ell$-group)}

\vspace{-.5em}
\Sigt{\GRL}{\cdot=0\mathrel{;}\cdot+\cdot,-\,\cdot,\cdot\vu\cdot,0}
\label{NOTASigGrl}
 
\medskip \noindent {\bf Abréviations}

\smallskip \noindent \textsl{Symboles fonctionnels}

\vspace{-1em} \DeuxCols{
\begin{itemize}
\itbu $x\vi y$ signifie $ - (-x\vu -y)$
\itbu $\abs{x}$ signifie $ x \vu -x$
\end{itemize}
}
{
\begin{itemize}
\itbu ${x}^+$ signifie $ x \vu 0$
\itbu ${x}^-$ signifie $ -x \vu 0$
\end{itemize}
}

\medskip \noindent \textsl{Prédicats}

\vspace{-1em} \DeuxCols{
\begin{itemize}
\itbu $x = y $ signifie $ x - y = 0$
\itbu $x \perp y $ signifie $ \abs x \vi \abs y =0$
\end{itemize}
}
{\begin{itemize}
\itbu $x \geq y $ signifie $ x \vu y = x$
\itbu $x \leq y $ signifie $ y\geq x$
\end{itemize}
}

\medskip \noindent {\bf Axiomes}

\smallskip Après les axiomes des groupes abéliens, on ajoute les axiomes suivants.

\smallskip\noindent {\it Règle pour la compatibilité de $\vu$ avec l'\egt}

\Regles{
\lAb{sup$_=$} $\,\, x=0\vet y=0\vd (u+x)\vu (v+y)= u\vu v$ \label{Axsup=}
}

\smallskip \noindent {\it Axiomes équationnels}

\smallskip 
Les identités suivantes expriment 
le fait que~$\vu$ définit un sup-demi treillis non borné ainsi que la compatibilité de $\vu$ avec $+$ (le fait que les translations sont des morphismes pour la loi $\vu$).

\DeuxRegles{
\Lab{sdt1} $\vd x\vu x=x $
\Lab{sdt2} $\vd x\vu y=y\vu x$
}
{
\Lab{sdt3} $\vd (x\vu y)\vu z=x\vu (y \vu z)$
\Lab{grl} $\vd x+(y\vu z)=(x+y)\;\vu\;(x+z)$
}

\smallskip On obtient ainsi un \grl (abélien) avec toutes les règles \gmqs afférentes (voir \cite{BKW}, 
\cite[Chapter 2]{Zaa97}, [Bourbaki, Algèbre, Chapitre 6], et \cite[Section XI-2, résultats 2.2-2.5, 2.11 et 2.12]{CACM}). En voici quelques unes.

\smallskip Une théorie \esid à \sa{Grl} est obtenue à partir de la théorie \Sa{Gao} en ajoutant les axiomes qui disent que deux \elts $x$ et $y$ du groupe ordonné ont toujours une borne supérieure, que l'on note $x\vu y$\footnote{L'existence de la borne supérieure est en effet une existence unique en raison de l'antisymétrie.}.\label{GrlGao} 

\DeuxRegles{
\Lab{sup1} $ \vd x\vu y\,\geq x $
\Lab{Sup} $ \,\, z\geq x\vet z\geq y\vd z\geq x\vu y$
}
{
\Lab{sup2} $ \vd x\vu y\,\geq y $
}

%%%%%%%%%%%%%%%%%%%%%%%%%%%%%%%%%%%%%%%%%%%%%%%%%%%%%%%%%%%%%%%%%%%%
\Subsection{Quelques règles dérivées dans \sa{Grl}}
%: Subsection{Quelques règles dérivées dans \sa{Grl}}

\DeuxRegles{
\Lab{grl1} $\cramped{\vd x\vu(y_1\vi y_2)= (x\vu y_1)\vi (x\vu y_2)}$
\Lab{grl2} $\cramped{\vd x\vi(y_1\vu y_2)= (x\vi y_1)\vu (x\vi y_2)}$
\Lab{grl3} $\vd (x\vi y)\vu x=x$
\Lab{grl4} $\vd (x\vu y)\vi x=x$
\Lab{grl5} $\vd (x\vi y)+(x\vu y) = x+y$
\Lab{grl6} $\vd x=x^+-x^-$
\Lab{grl7} $\vd \abs{x}=x^++x^-=x^+\vu x^-$
}
{
\laB{Sup} $ \,\, z\geq x\vet z\geq y\vd z\geq x\vu y$
\laB{Gao} $\,\, x\geq 0\vet x\leq 0 \vd x=0$

\Lab{Grl1} $\cramped{\,\, y\geq 0\vet z\geq 0\vet y\perp z \vd (y-z)^+=y}$
\Lab{Grl2} $\,\, x\leq z \vd (x\vi y)\vu z=x\vi(y\vu z) $
\lAb{Grl3$_n$} \label{AxGrl3} $\,\, nx\geq0\vd x\geq0\quad (n\in\NN, \,n>1) $
\lAb{Grl4$_n$} $\,\, nx\geq \Vi_{k=1}^n(ky+(n-k)x) \vd x\geq y $ \label{AxGrl4}
}

\vspace{-.7em}
\Regles{\Lab{Gauss} $\,\, x\geq 0\vet y\geq 0\vet z\geq 0\vet x\perp y\vet x\leq y+z \vd x\leq z $
}

%%%%%%%%%%%%%%%%%%%%%%%%%%%%%%%%%%%%%%%%%
% subsubsection{Structures quotients}
\Subsection{Structures quotients}\label{subsecgrlsolide}
%: Subsection{Structures quotients}

 Les noyaux des morphismes de groupes (abéliens) ordonnés sont les \textsl{sous-groupes convexes}: un sous-groupe $H$ est convexe \ssi il vérifie la \prt 
$${ (x\in H,\,y\in G,\, 0\leq y\leq x)\Rightarrow y\in H.}$$
Si un sous-groupe est convexe, la relation d'ordre \gui{passe au quotient}
dans $G/H$\index{convexe!sous-groupe --- (dans un groupe ordonné).}.

Les noyaux des morphismes de \grls sont les
\textsl{sous-groupes solides}. Un sous-groupe est solide \ssi il est un sous-\grl convexe, ou encore convexe et stable par $x\mapsto \abs x$ (\cite[\tho 2.2.1]{BKW})\index{solide!sous-groupe --- (dans un \grl)}.

\smallskip Le sous-groupe solide engendré par un \elt $a$ est $\cC(a):=\sotq{x\,}{\,\exists n\in\NN,\,\abs x\leq n\abs a}$.

Les sous-groupes solides \tf sont tous principaux: $\cC(\abs a + \abs b)=\cC(\abs a \vu \abs b)$ est le sous-groupe solide engendré par $a$ et $b$. Les sous-groupes solides principaux forment un \trdi (avec $\cC(a)\cap\cC(b)=\cC(\abs a \vi\abs b)$), à ceci près qu'il manque un \elt maximum, que l'on peut ajouter formellement.

La dimension de Krull de ce \trdi 
%\hum{il faudrait lui donner un nom,}
est appelée la \textsl{dimension, ou hauteur, du \grl}. 
Il s'agit d'une \dfn \cov \eqve à la \dfn classique en \clama, mais qui ne nécessite pas l'existence de sous-groupes convexes \gui{premiers} (voir \cite[section XIII-6]{CACM} % ou la \dfn \ref{defiDDKTRDI}
 pour la dimension de Krull des \trdis).
Dans le cas des groupes totalement ordonnés, cela correspond au rang du groupe.

%%%%%%%%%%%%%%%%%%%%%%%%%%%%%%%%%%%%%%%%%
\Subsection{Théorème de représentation} 
%: Subsection{Théorème de représentation}
En \clama tout \grl est un sous-objet d'un produit de groupes totalement ordonnés dans la catégorie des \grls.

La méthode de démonstration expliquée dans \cite[Principe XI-2.10]{CACM} donne un \gui{équivalent \cof} de cette \prt:
pour prouver un \gui{fait concret} dans un \grl, on peut toujours faire comme si l'on était en présence d'un produit de groupes totalement ordonnés.

\smallskip En fait, nous avons une {meilleure} formulation (plus formelle) dans le langage des \tdys: \textsl{les deux \tdys (avec et sans l'axiome de l'ordre total) prouvent les mêmes \ralgs.} Voyons ceci plus précisément.

%: Definition{defiGto}
\begin{definition} \label{defiGtosup}
La \tdy \SA{Gtosup} des \textsl{groupes totalement ordonnés avec sup} est la \tdy des \grls à laquelle on ajoute la \rdy \tsbf{OT} disant que l'ordre est total.

\Regles{
\laB{OT} {$\vd x\geq 0\;\vou\;x\leq 0 $}} 
\end{definition}
%--------- fin definition -----------------

\vspace{-.5em}
Notez que par rapport à la théorie usuelle des groupes totalement ordonnés \Sa{Gto} 
nous avons introduit dans la signature la loi $\cdot\vu\cdot$ qui est bien définie. La théorie \Sa{Gtosup} est \esid à la théorie \Sa{Gto}.

%: theorem thfairecommesi0
\begin{pstf}[pour les \grls] \label{thfairecommesi0} \index{Positivstellensatz!formel!pour les \grls} ~\\
Les \tdys \Sa{Grl} et \Sa{Gtosup} prouvent les mêmes \ralgs. 
\end{pstf}
\begin{proof}
\Llec peut se reporter à la \demo du \Pst formel~\ref{thfairecommesi}, et changer le tout petit peu qu'il y a à changer.
\end{proof}

\Llec pourra par exemple démontrer facilement les règles \Tsbf{Grl2} et \Grlqn\ en utilisant le \pstref{thfairecommesi0}, ce qui serait
nettement moins simple sinon.

Comme \corl du \pstref{thfairecommesi0} on obtient en \clama le \tho de plongement suivant (comme cas particulier du \thref{thcolsimralg}).

%: Corollary{corthfairecommesi0}
\begin{corollaryc}[\tho de plongement]
 \label{corthfairecommesi0} 
\emph{Voir \cite[Lorenzen, 1939]{Lor1939}, et les développements \cite{Lor1950,Lor1953} commentés dans \cite{CLN2019b}.}
Tout \grl $G$ est un sous-objet d'un produit de groupes totalement ordonnés\footnote{Tout \grl $G$ est une sous-{\usefont{T1}{pzc}{m}{it}{Grl}}-structure d'un produit de groupes totalement ordonnés quotients de~$G$. 
En d'autres termes, il y a un sous-\grl d'un produit de groupes totalement ordonnés qui, en tant que \grl, est isomorphe au \grl de départ. La terminologie anglaise est: any lattice group is a \textsl{subdirect product} of linearly ordered groups.} quotients de $G$.\end{corollaryc}
%--------- fin corollary ------------------------------- 

%r
%: Remark{remLLM01}
\begin{remark} \label{remLLM01} 
 La théorie de la complexité \algq dans l'espace des fonctions réelles continues sur
l'intervalle $[0,1]$ utilise de manière naturelle la structure de \grl divisible (la 2-divisibilité suffit). Cet espace de fonctions est vu essentiellement comme un espace de Riesz, et la multiplication des fonctions est reléguée au second plan. Voir par exemple \cite[définition 3.2.1]{LLM2001}. À noter \egmt que
dans cette théorie les formules du langage sont remplacées par des circuits (un circuit court peut coder une formule très longue). Nous sommes dans ce cas en analyse plutôt qu'en \alg abstraite. 
\eoe\end{remark}
%----------- fin remark ---------------------------------- 

Un exemple d'applications du \Pst formel pour les \grls est donné dans \cite[Fact XI-2.12]{CACM} que nous reproduisons ci-dessous.

%%: Fact{factGpRtcl}
\begin{fact}[autres identités dans les \grls]\label{factGpRtcl} ~\\
Soient $x$, $y$, $x'$, $y'$, $z$, $t\in G$, $n\in\NN$, $x_1$, \dots, $x_n\in G$.

\vspace{-.10em} 
\begin{enumerate}%\itemsep=1pt
\item \label{i1factGpRtcl} $x+y =\abs{x-y} +2(x\vi y)$
\item $(x\vi y)^+=x^+\vi y^+$, $(x\vi y)^-=x^-\vu y^-$, \\ $(x\vu y)^+=x^+\vu y^+$, $(x\vu y)^-=x^-\vi y^-$.
\item $2(x \vi y)^+ \leq (x+y)^+ \leq x^++y^+$.
\item $\abs{x+y} \leq \abs{x}+\abs{y}\;:\;$
$\abs{x}+\abs{y}=\abs{x+y}+2(x^+\vi y^-) +2( x^-\vi y^+)$.
\item $\abs{x-y} \leq \abs{x}+\abs{y}\;:\;$
$\abs{x}+\abs{y}=\abs{x-y}+2(x^+\vi y^+) +2( x^-\vi y^-)$.
\item $\abs{x+y}\vu\abs{x-y}=\abs{x}+\abs{y}$.
\item $\abs{x+y}\vi\abs{x-y}=\abS{\abs{x}-\abs{y}}$.

\item $\abs{x-y}=(x\vu y)-(x\vi y)$.
\item $\abs{(x\vu z)-(y\vu z)}+\abs{(x\vi z)-(y\vi z)}= \abs{x-y}.$
\item $\abs{x^+ - y^+} + \abs{x^- - y^-} = \abs{x-y}$.
\item \label{i11factGpRtcl} $x\leq z \;\Longrightarrow\; (x\vi y)\vu z= x\vi (y\vu z)$.
\item $x+y=z+t \;\Longrightarrow\; x+y=(x\vu z)+(y\vi t)$.
\item \label{i13factGpRtcl} $n\, x\geq \Vi_{k=1}^n (k y+(n-k)x) \,\Longrightarrow\,x\geq y$.
\item $\Vu_{i=1}^nx_i = \sum_{k=1}^n(-1)^{k-1}
 \big(\sum_{I\in \cP_{k,n}}\Vi_{i\in I}x_i\big)$\,\, ($n=2$, $x\vi y + x\vu y =x+y$).
\item $x\perp y\,\Longleftrightarrow\, \abs{x+y}=\abs{x-y}
\,\Longleftrightarrow\, \abs{x+y}=\abs{x}\vu \abs{y}$.
\item $x\perp y\,\Longrightarrow\, \abs{x+y}=\abs{x}
+\abs{y}=\abs{x}\vu \abs{y}$.
\item \label{i15bisfactGpRtcl} $(x\perp y,\,x'\perp y,\,x\perp y',\,x'\perp y',\,x+y=x'+y') \,\Longrightarrow\, (x=x', \, y=y')$.

\end{enumerate}
Supposons $u$, $v$, $w\in G^+$.
\begin{enumerate}\setcounter{enumi}{17}\itemsep=1pt
\item \label{i18factGpRtcl} $u\perp v\,\Longleftrightarrow\, u+v=\abs{u-v}$.
\item $(u+v)\vi w \leq (u\vi w)+(v\vi w)$.
\item $(x+y)\vu w \leq (x\vu w)+(y\vu w)$.
\item $v\perp w \,\Longrightarrow\,(u+v)\vi w = u\vi w$.
\item $u\perp v \,\Longrightarrow\,(u+v)\vi w = (u\vi w)+(v\vi w)$.
\end{enumerate}
\end{fact}
\begin{proof}
Tout ceci est à peu près \imd dans un groupe totalement ordonné,
en raisonnant cas par cas. On conclut avec le \Pst formel.
\end{proof}

\Subsection{Groupe réticulé engendré par un groupe ordonné} 
 
On a un morphisme naturel entre \talgs, de $\Sa{Gao}$ vers $\Sa{Grl}$, car on a défini $x\geq 0$ comme une simple abréviation dans la théorie $\sa{Grl}$.

Si $(G,\geq)$ est un groupe ordonné, comme $\sa{Grl}$ est une \talg, la \sad $\sa{Grl}(G)$ définit une structure algébrique usuelle $H$ de groupe réticulé. C'est le \grl \textsl{engendré par} $G$ au sens usuel.

Une description précise de la structure algébrique usuelle $\sa{Grl}(G)$ a été donnée par Lorenzen. Voir \cite{Lor1951,CLN2019b}. Rappelons ces résultats.

\begin{theorem}\label{thmgogrlfree}\emph{(\cite[Theorem 4.15]{CLN2019b})}
Si~\(G\) est un groupe ordonné, on peut construire un \grl~\(H\) avec un morphisme~\(\varphi\colon G\to H\) tel que 
\[
0\leq_H\varphi(a) \iff \exists n\in \NN\etl, \;0\leq_G\varphi(na)
\]
Plus \prmt, \(H\) est le \grl librement engendré par~\(G\) (au sens du foncteur adjoint à gauche au foncteur d'oubli).\\ 
Voici la construction. On considère le \trdi $\gT$  engendré par la \entrel  sur~$G$ définie comme suit. 
\begin{enumerate}
\item $\varphi(a)\vdash \varphi(b) $ pour les $(a,b)$ tels que $a\leq b$,
\item $\varphi(a_1),\dots,\varphi(a_k)\vdash \varphi(b_1),\dots,\varphi(b_\ell)$
\ssi il existe des entiers~\(n_1,\dots,n_k,m_1,\dots,\alb m_\ell\geq 0\) tels que
\begin{itemize}
\item \(n_1+\cdots+n_k=m_1+\dots+m_\ell\geq1\) et
\item \(n_1a_1+\dots+n_ka_k\leq_Gm_1b_1+\dots+m_\ell b_\ell\). 
\end{itemize}
\end{enumerate}
Sur ce \trdi $\gT$, la loi de groupe de $G$ se prolonge de manière unique, ce qui définit le \grl $H$.\\
Ainsi $\varphi(a_1),\dots,\varphi(a_k) \vdash  \varphi(b_1),\dots,\varphi(b_\ell)
$ signifie \fbox{$\varphi(a_1)\vi\dots\vi\varphi(a_k)\leq_H\varphi(b_1)\vu\dots\vu\varphi(b_\ell)$}. 
\end{theorem}

Ainsi le morphisme $\varphi$ est injectif \ssi pour tout $a\in G$ et tout entier $n\geq 1$, l'inégalité $na\geq 0$ implique $a\geq 0$.
Ce \tho \cof implique en \clama le \tho de Lorenzen-Clifford-Dieudonné selon lequel un groupe abélien peut être muni d'une structure de groupe ordonné \ssi l'\egt $na= 0$ implique $a= 0$.

En outre Lorenzen a mis en évidence une \rsim cruciale valide dans la théorie \sa{Grl} qui est la \textsl{régularité} (voir \cite{Lor1953} et \cite[Section 2]{CLN2019b})   

\Regles{\Lab{Reg}$\,\,x+a\geq 0\vet y+b\geq 0\vd x+b\geq 0\vou y+a\geq 0 $}

Le résultat de Lorenzen est qu'un groupe ordonné $G$ vérifie la \prt de régularité \ssi il s'identifie à un sous-groupe ordonné d'un \grl (ou, ce qui revient au même, du \grl qu'il engendre).

\medskip \noindent {\bf Une autre manière de définir le \grl engendré par un groupe ordonné}

\smallskip La structure de \grl sur un groupe abélien $H$ peut être définie à partir de la loi unaire $x\mt\abs{x}$. En effet 
$2x^+=x+\abs x$, $a\vu b=b+(a-b)^+$, et $x\geq 0$ équivaut à $x=\abs x$.

Nous indiquons maintenant quels axiomes doit satisfaire cette loi unaire pour qu'elle définisse bien une loi $a\vu b$ de \grl.

Tout d'abord le groupe doit être \textsl{$2$-divisible}, \cad vérifier les axiomes suivants.\label{grl-abs}

\DeuxRegles
{
\Lab{2div1} $\,\,2y=0\vd y=0$
}    
{
\Lab{2div2} $\vd \Exists x\;2y=x $
}

Auquel cas on peut définir sans ambigüité $\frac 1 2\,x$.
Ensuite la loi $x\mt \abs x$ doit satisfaire les axiomes suivants. 

\DeuxRegles
{
\lAb{abs$_=$}$\,\, x=0\vd \abs {y+x}=\abs{y}$
\Lab{abs1}$\vd \abs x=\abs{-x}$
\Lab{abs3}$\vd\abs x+ x= \abS{\abs{x}+ x}$
}    
{
\Lab{abs0}$\vd \abs{0}=0$
\Lab{abs2}$\vd\abs x + \abs y= \abS { \abs x + \abs y}$
\Lab{Abs1}$\,\,z\geq x\vet z\geq -x\vd z\geq\abs x$
}

Lorsque ces conditions sont remplies on utilise les abréviations suivantes.
\begin{itemize}
\item $x\geq 0$ signifie $\abs x=x$.
\item $a\geq b$ signifie $a-b\geq 0$.
\item $x^+:=\frac 1 2\,(x+\abs x) $. 
\item $a\vu b:=b+(a-b)^+$. 
\end{itemize}

\smallskip On doit vérifier la validité des règles suivantes

\DeuxRegles{
\laB{gao0} $\vd 0 \geq 0$
\laB{Gao} $\,\, x\geq 0\vet x\leq 0\vd x=0$
\laB{sup1} $ \vd x\vu y\,\geq x $
\laB{Sup} $ \,\, z\geq x\vet z\geq y\vd z\geq x\vu y$
}
{
\laB{gao1} $\,\, x\geq 0\vet y\geq 0\vd x+y\geq 0$
\item[ ]
\laB{sup2} $\vd x\vu y\,\geq y $
\laB{grl} $\vd x+(y\vu z)=(x+y)\;\vu\;(x+z)$
}

\begin{enumerate}
\item Un petit calcul donne  \fbox{$a\vu b = \frac{a+b}2+\frac{\abs{a-b}}2 $}. Avec \tsbf{abs1} cela montre que $a\vu b =b\vu a$.
\item On a $2\abs u= \abs{2u}$ d'après \tsbf{abs2} en prenant $u=x=y$. On a $2\abs u=2u\vd \abs u=u$ d'après \tsbf{2div1}. Donc  \fbox{$\abs {2u}=2u\vd \abs u=u$}. Par suite $\abs {z}=z\vd \abs {\frac z 2}=\frac z 2$ et donc 
\fbox{$\abs {\frac z 2}=\frac{\abs z} 2 $}. On a donc aussi $2u\geq 0\vd u\geq 0$, et \fbox{$2a\geq 2b\vd a\geq b$}.
\item Notons que la compatibilité avec l'\egt des  prédicats $x\geq 0$, $a\geq b$ et des lois $x\mt x^+$ et $(x,y)\mt x \vu y$ est une conséquence automatique du fait que la loi $\abs x$ respecte l'\egt.
\item \tsbf{gao0} signifie $\abs0=0$ qui est \tsbf{abs0}.
\item \tsbf{gao1} signifie $\abs x= x\vet\abs y =y\vd \abs{x+y}=x+y$. Cela résulte de \tsbf{abs2} car les hypothèses impliquent qu'on peut y remplacer $\abs x$ et $\abs y$ par  $x$ et~$y$.  
\item \tsbf{Gao} signifie $\abs x= x\vet\abs x =-x\vd x=0$. En utilisant  \tsbf{abs1} et la symétrie et la transitivité de l'\egt,
on déduit $x=-x$. On conclut avec~\tsbf{2div1}. 
\item \tsbf{sup1} signifie $y+(x-y)+-x\geq 0$, \ie $(x-y)^+-(x-y)\geq 0$. Il suffit donc de démontrer $a^+-a\geq 0$, \ie $\frac 1 2\,(a+\abs a)-a\geq 0$, \ie $\frac 1 2\,(\abs a -a)\geq 0$. On conclut avec \tsbf{abs3} et $\abs {\frac z 2}=\frac{\abs z} 2 $. 
\item \tsbf{sup2} s'en déduit parce que $a\vu b=b\vu a$.  
\item  \tsbf{Sup} signifie $\abs{z-x}=z-x,\vet \abs{z-y}=z-y\vd z\geq \abs{x+y}$. \\
En utilisant $\abs u = \abs {-u}$ les hypothèses donnent $2z-(x+y)=\abs{z-x}+\abs{y-z}$. Comme $\abs a \geq a$, on obtient avec \tsbf{gao1} $2z-(x+y)\geq y-x$. Symétriquement $2z-(x+y)\geq x-y$. Donc  \tsbf{Abs1} donne $2z-(x+y)\geq \abs{y-x}$, \ie $2z\geq 2 (y\vu x)$ et on conclut avec $2a\geq 2b\vd a\geq b$.
\item \tsbf{grl} signifie $2x+(y+z)+\abs{y-z}=(x+y)+(x+z)+\abs{(x+y)-(x+z)}$, ce qui est clair.  
\end{enumerate}

%%%%%%%%%%%%%%%%%%%%%%%%%%%%%%%%%%%%%%%%%%%%%%%%%%%%%%%%%%%%%%%%%%%%
\section{Anneaux fortement réticulés ($f$-rings)}\label{secfrings}

Références: \cite{BKW}, \cite[Section V-4]{Joh1986}, \cite{BP56,DM1995,Mad10}. 

\smallskip La terminologie de \cite{BKW} est \gui{f-anneau} ou \gui{anneau de fonctions}. Cet ouvrage étudie le cas des anneaux non commutatifs et non unitaires, pour lesquels les résultats sont plus subtils et délicats que ceux que nous donnons ici dans le cas commutatif unitaire.
La terminologie \gui{\afr} pour les \gui{$f$-rings} de la littérature anglaise se trouve dans les exercices de Bourbaki (Algèbre, chapitre VI, Exercices, \S2, exercice 5). 

\smallskip Les \afrs (commutatifs unitaires) sont définis par une théorie purement équationnelle.\index{anneau!fortement reticule@fortement réticulé}\index{f-ring@$f$-ring!anneau fortement réticulé} 
Les axiomes sont ceux des anneaux commutatifs, ceux des \grls pour l'addition, et enfin la règle équationnelle \Tsbf{afr} qui exprime une forme de compatibilité de $\vu$ avec la multiplication\footnote{Par rapport à la théorie {\usefont{T1}{pzc}{m}{it}{Grl}}, 
on a ajouté la loi \gui{$\cdot\times \cdot$} et les règles \tsbf{ac2} et \tsbf{afr}. En outre, la machinerie calculatoire qui réduit tout terme sur les variables $\Xn$ à son écriture canonique dans le groupe abélien libre $\ZZ^{\{\Xn\}}$ a été remplacée par la machinerie calculatoire qui réduit tout \elt de $\ZZXn$ à une forme normale.}.
Voici tout en détail.

%%%%%%%%%%%%%%%%%%%%%%%%%%%%%%%%%%%%%%%%%%%%%%%%%%%%%%%%%%%%%%%%%%%%
\Subsectio{\Dfn de la \tpe \sa{Afr}}{Théorie purement équationnelle des \afrs}

La théorie \SA{Afr} est définie comme suit.

\vspace{-.5em}
\Sigt{\AfR}{\cdot=0\mathrel{;}\cdot+\cdot, \cdot\times\cdot,\cdot\vu\cdot,-\,\cdot,0,1}
\label{NOTASigAfr}

\noindent {\bf Abréviations} (comme pour les \grls)

\smallskip \noindent \textsl{Symboles fonctionnels}

\vspace{-1em} \DeuxCols{
\begin{itemize}
\itbu $x\vi y$ signifie $ - (-x\vu -y)$
\itbu $\abs{x}$ signifie $ x \vu -x$
\end{itemize}
}
{
\begin{itemize}
\itbu ${x}^+$ signifie $ x \vu 0$
\itbu ${x}^-$ signifie $ -x \vu 0$
\end{itemize}
}

\medskip \noindent \textsl{Prédicats}

\vspace{-1em} \DeuxCols{
\begin{itemize}
\itbu $x = y $ signifie $ x - y = 0$
\itbu $x \perp y $ signifie $ \abs x \vi \abs y =0$
\end{itemize}
}
{\begin{itemize}
\itbu $x \geq y $ signifie $ x \vu y = x$
\itbu $x \leq y $ signifie $ y\geq x$
\end{itemize}
}

\medskip \noindent {\bf Axiomes}

\smallskip\noindent \textsl{Règles des anneaux commutatifs}

\DeuxRegles{
\laB{ga0} $\vd 0=0$
\laB{ga1} $\,\, x=0\vet y=0\vd x+y=0$
}
{
\laB{ac2} $\,\, x=0\vd xy=0$
}

\smallskip \noindent \textsl{Règles de compatibilité de $\vu$ avec l'\egt}

\Regles{
\lab{\supieq} $\,\, x=0\vet y=0\vd (u+x)\vu (v+y)= u\vu v$ 
}

\smallskip \noindent \textsl{Règles équationnelles}

\DeuxRegles{
\laB{sdt1} $\vd x\vu x=x $
\laB{sdt2} $\vd x\vu y=y\vu x$
\laB{sdt3} $\vd (x\vu y)\vu z=x\vu (y \vu z)$
}
{
\laB{grl} $\vd x+(y\vu z)=(x+y)\vu(x+z)$
\Lab{afr} $\vd x^+\, (y\vu z)=(x^+\, y)\vu(x^+\, z)$
}

%:     note \label{noteAfrp}
\begin{notE} \label{noteAfrp}
 Si l'on utilise la signature

\Sigt{{\AfR'}}{\cdot=0,\cdot\geq 0\mathrel{;}\cdot+\cdot, \cdot\times\cdot,\cdot\vu\cdot,-\,\cdot,0,1}
\label{NOTASigAfr'}

\vspace{-.8em}
\noindent on donne les trois règles \Tsbf{sup1}, \Tsbf{sup2} et \Tsbf{Sup} (\paref{Axsup1})
 pour relier $\cdot\geq 0$ \hbox{et $\cdot\vu\cdot$}. On définit ainsi la théorie \SA{Afr'} \esid à \sa{Afr}. \eoe
\end{notE}

%%%%%%%%%%%%%%%%%%%%%%%%%%%%%%%%%%%%%%%%%%%%%%%%%%%%%%%%%%%%%%%%%%%%
\Subsection{Note sur les anneaux réticulés ($\ell$-rings)}
%: Subsection{Note sur les anneaux réticulés ($\ell$-rings)}

La théorie \SA{Arl} des \textsl{anneaux réticulés} ($\ell$-rings dans la littérature anglaise) est définie en remplaçant la règle \Tsbf{afr}
par les règles \Tsbf{ao1} et \Tsbf{ao2} des anneaux ordonnés, valides dans~\Sa{Afr}.\index{anneau!reticule@réticulé}\label{ao1}

\DeuxRegles{
\Lab{ao1} $\vd \,x^2\geq 0$
}
{
\Lab{ao2} $\,\,x\geq 0\vet y\geq 0\vd xy\geq 0$ \phantom{$,a^2\geq 0$}
}

%: Lemma{lemaoafr}
\begin{lemma} \label{lemaoafr}
Dans la théorie des anneaux \rtls les règles suivantes sont toutes équivalentes.

\DeuxRegles{
\laB{afr} $\vd a^+\, (b\vu c)=(a^+\, b)\vu(a^+\, c)$
\Lab{afr'} $\vd a^+\, (b\vi c)=(a^+\, b)\vi(a^+\, c)$
\Lab{afr0} $\vd b^-\vi a^+b^+= 0 $
\Lab{afr1} $\vd a^+\ a^-= 0 $
\Lab{afr2} $\vd \abs a\,\abs b= \abs{ab} $
\Lab{afr3a} $\vd (ab)^+=a^+b^++a^-b^-$
\Lab{afr3b} $\vd (ab)^-=a^+b^-+a^-b^+$
\Lab{afr4} $\vd c^+\abs a= \abs{c^+a}$
\Lab{afr5} $\vd (a\vi b)(a\vu b)=ab $
}
{
\Lab{Afr} $\,\,a \geq 0\vd a(b\vu c)=ab\vu ac\phantom{a^+}$
\Lab{Afr'} $\,\,a \geq 0\vd a(b\vi c)=ab\vi ac\phantom{a^+}$
\Lab{Afr0} $\,\,b\vi c = 0\vet a \geq 0 \vd b\vi ac = 0\phantom{a^+}$
\Lab{Afr1} $\,\, a \vi b =0\vd ab=0 $ 
\Lab{afr6a} $\vd a^2=(a^+)^2+(a^-)^2$
\Lab{afr6b} $\vd a^2= {\abs a}^2$
\Lab{sup} $ \vd \big((x\vu y)- x\big)\,\big((x\vu y)- y\big)=0$
\Lab{Afr2} $\,\, b\perp c\vd ab\perp ac$
%\Lab{Afr3} $\,\,b\vi c=0\vet x\geq 0\vd b\vi xc=0$
}
\end{lemma}
%--------- fin lemma ----------------------------------- 

\vspace{-.5em}
Autrement dit chacune de ces règles peut servir à définir les \afrs
en l'ajoutant à la théorie \Sa{Arl}. Le fait que \tsbf{afr} implique
\tsbf{ao1}, \tsbf{ao2} et les règles signalées en \ref{lemaoafr} résulte du \Pst formel \ref{thfairecommesi}. 

Dans le cas d'un anneau réticulé non unitaire
la règle \tsbf{afr0} est plus forte que les autres (voir \cite{BKW}, proposition 9.1.10\footnote{Le livre traite plus \gnlt les anneaux ordonnés non nécessairement commutatifs, ni unitaires. La condition \tsbf{afr0} doit alors être dédoublée pour tenir compte de la non commutativité.}).

\smallskip \noindent\hum{Donner une démonstration de \ref{lemaoafr} ne semble pas indispensable. On est surtout intéressé par le fait que les \afrs vérifient toutes ces \prts, et cela résulte du \ref{thfairecommesi}. Il faudrait cependant une référence pour ce lemme \ref{lemaoafr}. Je ne sais plus où je l'ai trouvé. Peut-être dans un brouillon à moi? L'exercice de Bourbaki ne le donne qu'en partie. Même chose pour \cite{Joh1986}.}

%%%%%%%%%%%%%%%%%%%%%%%%%%%%%%%%%%%%%%%%%%%%%%%%%%%%%%%%%%%%%%%%%%%%
\Subsectio{Quelques règles dérivées dans la théorie \sa{Afr}}{Quelques règles dérivées}
%: Subsection{Quelques règles dérivées dans \sa{Afr}}

 Outre les règles dérivées pour les \grls et celles signalées dans le lemme \ref{lemaoafr}, voici des règles classiques fort utiles dans lesquelles intervient la multiplication.

\DeuxRegles{
\laB{Ato1} $\,\, b\geq 0 \vet ab=1\vd a\geq 0 \vphantom{\abs{a}^2} $}
{
\laB{Ato2} $\,\, c\geq 0\vet a(a^2+c)\geq 0\vd a^3\geq 0 \vphantom{\big)}$
}

\vspace{-1em}
\Regles{
\Lab{afr7} {$\vd ab^+ = (ab \vi (a^2 +1)b) \vu (-(a^2 +1)b\vi 0)$} 
}

\smallskip \rem La règle \tsbf{afr7}
sert à démontrer la possibilité d'écrire sous une forme simplifiée les termes
dans un \afr libre: voir le lemme \ref{lemAfrReecriture}. 
\eoe

%%%%%%%%%%%%%%%%%%%%%%%%%%%%%%%%%%%%%%%%%
\Subsection{Structures quotients}
%: Subsection{Structures quotients}

\paragraph{Idéaux solides}~

\smallskip Par \dfn, les noyaux des morphismes d'\afrs sont appelés les
\textsl{idéaux solides}\footnote{L'ouvrage \cite{BKW} dit \textsl{un $\ell$-idéal}, conformément à la terminologie anglaise, et il ajoute une précision \ncr dans la cas non commutatif.} (un $\ell$-idéal dans la littérature anglaise)\index{idéal!solide (dans un \afr)}\index{solide!idéal --- (dans un \afr), ou encore un $\ell$-idéal.}. 

Un idéal est solide \ssi il est solide en tant que sous-groupe. 

L'idéal solide engendré par un \elt $a$ est 
$$\cI(a):=\sotq{x\,}{\,\exists y,\,\abs x\leq \abs {ya}}.$$ 
On a $\cI(a)=\cI(\abs a)$ et $\cI(a)\cap\cI(b)=\cI(\abs a\vi \abs b)$.
Enfin, l'idéal solide engendré par $\an$ est 
$$\cI(\an)=\cI(\abs {a_1}+\dots+\abs {a_n})= \cI(\abs {a_1}\vu \cdots \vu \abs {a_n}).$$ 

\paragraph{Idéaux solides irréductibles}~

\smallskip On dit qu'\textsl{un idéal solide $\fa $ d'un \afr $\gA$ est \textsl{irréductible}} si l'\afr quotient est totalement ordonné. Autrement dit, pour tout $x\in\gA$, $x^+\in \fa$ ou $x^-\in \fa$.\index{idéal!solide irréductible (d'un \afr)}\index{irréductible!idéal solide --- dans un \afr)} 

D'après le lemme \ref{lemAfrsdz}, tout idéal solide premier est irréductible.

Par ailleurs un idéal premier convexe (en tant que sous-groupe additif) $\fp$ est solide: on doit voir qu'il est stable par $\vu$. Si $a,b\in \fp$
on a $(a-b)^+$ ou $(a-b)^-\in\fp$. Et les identités $b+(a-b)^+=a\vu b=a+(a-b)^-$ sont valides dans les \grls (et à fortiori dans les \afrs) parce qu'elles le sont dans les groupes totalement ordonnés (\Pst formel~\ref{thfairecommesi0}).

%\paragraph{Corps résiduel d'un \afr local}~

%%%%%%%%%%%%%%%%%%%%%%%%%%%%%%%%%%%%%%%%%
\Subsectio{\Pst formel et théorème de plongement pour les \afrs}{\Pst formel et théorème de plongement} 
%: Subsection{\Pst formel et théorème de plongement pour les \afrs}

%\noindent 
Rappelons que la \tdy \Sa{Atosup} des anneaux totalement ordonnés avec sup est la \tdy des anneaux totalement ordonnés à laquelle on ajoute un symbole de fonction~$\cdot\vu\cdot$ qui doit satisfaire les \ralgs suivantes.%

\DeuxRegles{
\laB{sup1} $ \vd x\vu y\,\geq x $
\laB{Sup} $ \,\, z\geq x\vet z\geq y \vd z\geq x\vu y$
}
{
\laB{sup2} $ \vd x\vu y\,\geq y $
}

\smallskip\noindent On peut aussi voir \sa{Atosup} comme la théorie des \afrs à laquelle on ajoute comme axiome la \rdy \tsbf{OT} (disant que l'ordre est total).

\Regles{
\laB{OT} $\vd x\geq 0\;\vou\;x\leq 0 $
}

Vu l'existence unique du sup dans un anneau totalement ordonné, les théories \Sa{Ato} et \Sa{Atosup} sont \esids.
En particulier, elles prouvent les mêmes \rdys
(lorsqu'elles sont formulées sans utiliser $\vu$). 
%\end{lemma}
%--------- fin lemma ----------------------------------- 

Le \tho pour les \afrs analogue au \Pst formel~\ref{thfairecommesi0} pour les \grls est le suivant.
Il s'agit d'un résultat du même type que point \textsl{2} du \Pst formel \ref{Pst1bis}.

%: positivstellensatz thfairecommesi
\begin{pstf}[pour les \afrs] 
\index{Positivstellensatz!formel!pour les \afrs} \label{thfairecommesi}~\\
We consider \sads on the signature\\ 
\centerline{\sIgt{\AfR'}{\cdot=0,\cdot\geq 0\mathrel{;}\cdot+\cdot, \cdot\times\cdot,\cdot\vu\cdot,-\,\cdot,0,1} }\\
Les théories \Sa{Afr}, \Sa{Ato} et \Sa{Atosup} prouvent les mêmes \ralgs.
\end{pstf}
\begin{proof} Tout d'abord \Sa{Ato} et \Sa{Atosup} sont \esids.
Considérons maintenant une \ralg prouvée dans la \tdy $\sa{Atosup}\!$.
On peut supposer \spdg que la conclusion de la règle est une \egt $t=0$ pour un terme $t$ convenable.
 Dans le calcul correspondant, en présence d'un terme $u$, on est autorisé par \Tsbf{OT} à ouvrir deux branches. L'une où $u\geq 0$, l'autre où $u\leq 0$. À chaque nœud de la preuve dynamique, on travaille en fait dans un \afr défini par \gtrs et relations: les \gtrs sont donnés dans la présentation et dans les hypothèses de la règle \agq à démontrer; pour les relations de même, avec en outre celles que l'on a ajouté, dans la branche où l'on se trouve, aux embranchements qui précèdent le nœud.
Supposons qu'à un moment donné, pour deux termes~$a$ et $b$, on ait ouvert une branche où $a\geq b$ et une autre où $a\leq b$. 
Posons $c=b-a$. Dans la première branche on a ajouté l'hypothèse $c^-=0$, dans la seconde l'hypothèse $c^+=0$. 
Si dans chacune des branches on peut prouver $t=0$,
cela veut dire que l'on a, dans l'\afr correspondant au nœud
en question, d'une part $t\in \cI(c^-)$, et d'autre part $t\in \cI(c^+)$. Or dans un \afr $\cI(c^+)\cap \cI(c^-)= \cI(c^+\vi c^-)=\so 0$.
\end{proof}

Démontrons par la méthode du \Pst~\ref{thfairecommesi} la validité des règles \agqs \Tsbf{Ato1} et \Tsbf{Ato2} dans \sa{Afr}.

\DeuxRegles{
\lab{Ato1} $\,\, y\geq 0 \vet xy=1\vd x\geq 0\phantom{x^3}$
}
{
\lab{Ato2} $\,\, c\geq 0\vet x(x^2+c)\geq 0\vd x^3\geq 0$
}

\smallskip Dans les deux cas, on ouvre deux branches, l'une où $x\geq 0$, et le résultat est clair, l'autre \hbox{où $x\leq 0$}. Pour \Tsbf{Ato1}
on en déduit que $1\leq 0$, puis $1=0$, puis $x=0$.
Pour \Tsbf{Ato2}
on en déduit \hbox{que $x^3\geq -xc\geq 0$}.

\smallskip De même, on démontre \Tsbf{afr7} en examinant séparément les cas \gui{$b\geq 0$}, \gui{$b\leq 0,\,a\geq 0$}
\hbox{et \gui{$b\leq 0,\,a\leq 0$}}.

Comme conséquence du \Pst formel \ref{thfairecommesi} on obtient en \clama le \tho de plongement suivant (comme cas particulier du \thref{thcolsimralg}).
%: Corollary{corthfairecommesi}
\begin{corollaryc}[\tho de plongement] \label{corthfairecommesi}
Tout \afr $\gA$ est un sous-objet d'un produit d'anneaux totalement ordonnés quotients de $\gA$. 
\end{corollaryc}
%--------- fin corollary ------------------------------- 

\noindent Le \tho suivant est du même type que point \textsl{1} du \Pst formel~\ref{Pst1bis}. %On dit qu'\textsl{une \sad de type \sa{Afr} s'effondre lorsque la règle $\vd 1=0$ est valide}. 
Ce résultat peut être vu comme une deuxième forme du \Pst formel \ref{thfairecommesi} pour les \afrs. 
%: Theorem{thColsimafr}
\begin{theorem}[simultaneous collapse for \afrs]~\label{thColsimafr}\\
We consider \sads on the signature\\ 
\centerline{\sIgt{\AfR'}{\cdot=0,\cdot\geq 0\mathrel{;}\cdot+\cdot, \cdot\times\cdot,\cdot\vu\cdot,-\,\cdot,0,1} }\\  
The theories \Sa{Afr}, \Sa{Crcdsup} and all intermediate theories collapse simultaneously. 
\end{theorem}
\begin{proof}~\\
Les théories \Sa{Afr} et \Sa{Atosup} 
s'effondrent simultanément d'après le \Pst~\ref{thfairecommesi}.
\\
Les théories \Sa{Ato} et \Sa{Crcd} s'effondrent simultanément d'après
le point \textsl{1} du \Pst \ref{Pst1bis}.
\\
Enfin les théories \Sa{Ato} et \Sa{Crcd} sont \esids respectivement aux théories \Sa{Atosup} et \Sa{Crcdsup}.
\end{proof}
%

%%%%%%%%%%%%%%%%%%%%%%%%%%%%%%%%%%%%%%%%%%%%%%%%%%%%%%%%%%%%%%%%%%%%
\Subsection{Localisations des \afrs}\label{locafr}
%: Subsection{Localisation d'un \afr}

\paragraph{Généralités}~

\smallskip On considère un \mo $S$ dans un \afr $\gA$ et l'on construit la solution du \pb universel (dans la catégorie des \afrs) consistant à inverser les \elts de~$S$. 

Pour cela, il suffit de considérer le localisé usuel $S^{-1}\gA$ et de définir correctement la loi $\vu$. Comme inverser $s$ ou inverser $s^2$ revient au même, on peut ne considérer que des fractions à dénominateur $\geq 0$. On définit alors 

\snic{ \dsp\frac{\,a\,}{s}\vu\frac{\,b\,}{t}:=\frac{\,at\vu bs\,}{st} \qquad (s,t\geq 0).}
 
\Note On n'a pas le choix, car puisque $s,t\geq 0$, on doit avoir $st\,(\frac{a}{s}\vu\frac{b}{t})=st\,\frac{a}{s}\vu st\, \frac{b}{t} =at\vu bs$ dans~$S^{-1}\gA$. 
Il reste à voir que la loi est bien définie et qu'elle continue de satisfaire les axiomes requis. Vérifions par exemple qu'elle est bien définie. Supposons que $\frac{a_1}{s_1}=\frac{a_2}{s_2}$,
i.e. que $a_1s_2s_3=a_2s_1s_3$ pour un~$s_3\geq 0$ dans~$S$.
Alors on vérifie facilement que les deux \elts $\frac{a_i}{s_i}\vu\frac{b}{t}$
donnés selon la \dfn ci-dessus ont bien égaux dans $S^{-1}\gA$. 
Il s'agit du même calcul que celui qui a été fait pour justifier l'addition dans $S^{-1}\gA$ quand on était petit\footnote{Quand on est tombé d'admiration devant Claude Chevalley qui osait inverser des diviseurs de zéro, et rien d'affreux n'en résultait, bien au contraire.}. On remplace seulement $+$
par $\vu$, avec la précaution d'avoir des dénominateurs $\geq 0$.\eoe

%l
%:     Lemma{lemAfrqlg}
\begin{lemma} \label{lemAfrqlg}
Un \afr $\gA$ peut toujours être considéré comme plongé dans une \QQlg \frle. 
\end{lemma}
%----------- fin lemma ----------------------------------- 
%
\begin{proof}
En effet, d'après \Tsbf{Grl3}$_n$, les \gui{entiers} $n.1_\gA$ sont réguliers et donc $\gA$ s'injecte dans la \QQlg $\QQ\otimes_\ZZ \gA$ qui est \frle en tant que localisée de $\gA$\footnote{C'est vrai même si $\gA$ est trivial: le seul cas où la \QQlg en question ne contient pas $\QQ$ comme sous-anneau.}.
\end{proof}
%

%%%%%%%%%%%%%%%%%%%%%%%%%%%%%%%%%%%%%%%%%%%%%%%%%%%%%%%%%%%%%%%%%%%%
\paragraph{Recollement de structures d'\afrs}
%: paragraph{Recollement d'\afrs}

\smallskip On suppose donnée sur un anneau commutatif $\gA$, de manière locale, une structure d'\afr. Si ces structures définies localement coïncident deux par deux sur les ouverts où elles sont définies en commun, alors on peut les recoller en une structure d'\afr définie sur tout l'anneau $\gA$. Cela est à rapprocher du principe de recouvrement pour les \trdis (principe \ref{propRecouvTD}) et s'énonce d'une manière précise similaire comme suit.

%: Plgc {plcc.frings}
\begin{plcc}[recollement concret de structures d'\afrs sur un anneau]\label{plcc.frings}  
Soient $S_1$, $\dots$, $S_n$ des \moco d'un anneau $\gA$. On note $\gA_i$ pour $\gA_{S_i}$, $\gA_{ij}$ pour~$\gA_{S_iS_j}$, et l'on suppose donné sur chaque $\gA_i$ une structure d'\afr
avec une loi $\vu\!_i$. On suppose en outre que les images dans $\gA_{ij}$ des lois $\vu\!_i$ et $\vu\!_j$ coïncident. Alors il existe une unique structure d'\afr sur~$\gA$ qui induit par \lon en chaque $S_i$ la structure définie sur~$\gA_i$. Cet \afr s'identifie à la limite projective du diagramme
\[
\big(( \gA_i)_{i\in\lrbn},( \gA_{ij})_{i<j\in\lrbn};(\alpha_{ij})_{i\neq j\in\lrbn}\big),
\]
où les $\alpha_{ij}$ sont les morphismes de \lon, dans la catégorie des \afrs.
\end{plcc}
%--- end-plgc-----------------------------------------
 
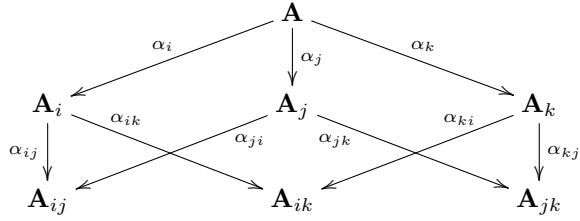
{\hspace*{10em}{
\xymatrix @R=2em @C=7em{
 & \gA \ar[rd]^{\alpha _{k}}\ar[d]^{\alpha _{j}}\ar[ld]_{\alpha _{i}}\\
 \gA _i\ar[d]_{\alpha _{ij}}\ar@/-0.75cm/[dr]^{\alpha _{ik}} &
 \gA _j\ar@/-1cm/[dl]^{\alpha _{ji}}\ar@/-1cm/[dr]_{\alpha _{jk}} &
 \gA _k\ar@/-0.75cm/[dl]_{\alpha _{ki}}\ar[d]^{\alpha _{kj}} &
\\
 \gA _{ij} & 
 \gA _{ik} & 
 \gA _{jk} 
}
}}

\begin{proof}
L'anneau $\gA$ est la limite du système projectif formé par les $\gA_i$ et $\gA_{ij}$ dans la catégorie des anneaux commutatifs, donc aussi dans la catégorie des ensembles. 
%Par suite l'ensemble~\hbox{$\gA\times \gA$} est la limite du système projectif formé par les $\gA_i\times \gA_i$ et $\gA_{ij}\times \gA_{ij}$ dans la catégorie des ensembles. 
Il s'ensuit qu'il existe une unique loi $\vu$ sur $\gA$ qui donne les $\vu_i$ sur les $\gA_i$ par les applications canoniques $\gA\to\gA_i$. Il reste à vérifier qu'elle satisfait les axiomes de la loi $\vu$ pour un \afr. Cela résulte du fait que ces axiomes sont donnés par des \egts entre termes, et du fait que l'application naturelle  $\varphi\colon \gA\to \prod_i\gA_i$, d'une part préserve les lois de la structure d'\afr, et d'autre part est injective.
\end{proof}
%

%\paragraph{Localisation en un filtre premier}~ 

\paragraph{Schémas réels}~ 

%: Remark{remplcc.frings}
\begin{remark} \label{remplcc.frings} 
Un \corl du principe de recollement \ref{plcc.frings} est que la notion de \gui{schéma (de Grothendieck) fortement réticulé} est bien définie. Un \gui{schéma fortement réticulé} semble d'ailleurs être la définition la plus naturelle pour la notion de schéma réel. En effet, cela autorise les nilpotents et donc une bonne théorie des multiplicités dans les schémas réels.
\\ 
Plus \prmt on pourrait définir un \gui{schéma réel} comme un schéma obtenu en recollant un nombre fini de schémas affines définis par des \afrs qui sont des \RRlgs \pf, ou plus \gnlt des \Rlgs \pf pour un \covr $\gR$ donné.
\\
Mais cela reste apparemment inexploré.
\eoe\end{remark}
%----------- fin remark ---------------------------------- 
\hum{ajouter deux mots sur les fibres d'un tel schéma}

%%%%%%%%%%%%%%%%%%%%%%%%%%%%%%%%%%%%%%%%%%%%%%%%%%%%%%%%%%%%%%%%%%%%
\Subsection{Réécriture de termes dans les \afrs}
%: Subsection{Réécriture de termes dans les \afrs}

Référence: \cite{Del86}. 

Contrairement à la théorie des anneaux commutatifs dans laquelle les termes se réécrivent sous une forme normale unique, on n'a pas pour les \afrs un résultat aussi satisfaisant. On a néanmoins une réécriture sous une forme simplifiée, du même style que la forme normale conjonctive dans les \trdis.

%: Lemma{lemAfrReecriture}
\begin{lemma} \label{lemAfrReecriture}
Soit $\gA$ un \afr et $t$ un terme écrit sur des \idtrs $\xn$
et des constantes dans $\gA$. Ce terme se réécrit sous forme 
$$\sup\nolimits_{i\in I}\!\big(\inf\nolimits_{j\in J_i}(f_{i,j}(\ux))\big)
$$ 
pour une famille finie convenable de \pols $f_{i,j}\in\AXn$.
%indexée par une somme disjointe $\bigcup_{i\in I}J_i$. 
\end{lemma}
%--------- fin lemma -----------------------------------
%
\begin{proof} Vu les réécritures usuelles dans les \trdis et vu que $x\mapsto -x$
échange~$\vu$ et $\vi$,
il suffit de savoir réécrire $a+(b\vu c)$ et $a(b\vu c)$ sous la forme voulue. Cela résulte des règles équationnelles \Tsbf{grl}, \Tsbf{grl6} et \Tsbf{afr7}.
\end{proof}
%

%n
%: definotation{notaAFR}
\begin{definota} \label{notaAFR}~
\begin{enumerate}
\item  Soit $\gB=\big((G,R),\Sa{Afr}\big)$ un \afr dynamique. Comme la théorie \Sa{Afr} est \agq, $\gB$ admet un \hyperref[Modelegnq]{\textsl{modèle générique}}, noté $\AFR(\gB)$. C'est l'\afr défini par les \gtrs $G$ et les relations $R$. 
\item Lorsque $G=\so{1,\dots,n}$ et $R$ est vide, on obtient l'\afr $\AFR(\ZZ[\Xn])$ librement engendré par $n$ \elts.
\end{enumerate}
\end{definota}
%--------- fin notation ------------------------------ 

Comme la théorie \Sa{Afr} est \peq, le lemme \ref{lemAfrReecriture} équivaut 
à son énoncé restreint aux cas particuliers où $\gA$ est un \afr libre sur un ensemble fini.
 
%l
%: Lemma{lemafrgenerique}
\begin{lemma} \label{lemafrgenerique}~
\begin{enumerate}
\item Les \elts de l'anneau $\AFR(\gB)$ peuvent tous s'écrire sous la forme donnée dans le lemme \ref{lemAfrReecriture} avec les $f_{ij}\in\ZG$. 

\item Si $\gC$ est un anneau commutatif, prenons pour $(G,R)$
le diagramme positif de $\gC$. Alors $\AFR(\gC)$ est l'\afr librement engendré par l'anneau commutatif $\gC$, et les \elts
 de $\AFR(\gC)$ s'écrivent sous forme $\sup\nolimits_{i\in I}\!\big(\inf\nolimits_{j\in J_i}a_{ij})\big)$ avec des \elts~$a_{ij}$ de~$\gC$. 

\end{enumerate}
 \end{lemma}
%----------- fin lemma ----------------------------------- 

%%%%%%%%%%%%%%%%%%%%%%%%%%%%%%%%%%%%%%%%%%%%%%%%%%%%%%%%%%%%%%%%%%%%
\Subsection{Anneaux de fonctions \frls, semi\pols}
%: Subsection{Anneaux de fonctions \frls, semi\pols}

Pour tout ensemble $E$ et tout \afr $\gA$ l'anneau des fonctions \hbox{$f\colon E\to\gA$} est muni d'une structure naturelle d'\afr (c'est la structure produit).

%: Definota{defiSIPD}
\begin{definota} \label{defiSIPD}
Soit $\varphi\colon \gA\to\gB$ un morphisme d'\afrs. L'anneau des \textsl{$\gA$-\spos en~$n$ variables\footnote{Dans la littérature anglaise, les \spos sont appelés les \gui{SIPD} ou \gui{sup-inf-polynomially-defined functions}.} sur $\gB$} est le sous-\afr de fonctions $f\colon \gB^n\to\gB$ engendré par les constantes dans $\varphi(\gA)$ et les fonctions \coos. 
On le notera $\SIPD_n(\gA,\gB)$. On abrège $\SIPD_n(\gA,\gA)$ en $\SIPD_n(\gA)$.
\\
La \dfn s'étend au cas où $\gA$ et/ou $\gB$ sont des anneaux totalement ordonnés, que l'on considère comme des \afrs.
\end{definota}
%--------- fin definition -------------------------------- 

Notons que ce n'est pas vraiment restrictif de supposer que $\varphi$ est injectif, ce qui permet de regarder~$\gA$ comme un sous-\afr de $\gB$.

%: Lemma{lemSIPD}
\begin{lemma} \label{lemSIPD}
On suppose que $\gA\subseteq \gB$. Tout \elt de $\SIPD_n(\gA,\gB)$ se réécrit sous forme $\sup_{i\in I}\!\big(\inf_{j\in J_i}(f_{i,j})\big)$ pour une famille finie convenable de \pols $f_{i,j}\in\Axn$.
%indexée par une somme disjointe $\bigcup_{i\in I}J_i$. 
\end{lemma}
%--------- fin lemma -----------------------------------
%
\begin{proof} C'est à très peu près le lemme \ref{lemAfrReecriture}.
\end{proof}
%

%e
%: Example{exaSIPD}
\begin{examples} \label{exaSIPD}~

\noindent 1. 
Les deux \elts $x\vu (1-x)$ et $1\vu x\vu (1-x)$ définissent la même fonction dans $\SIPD_1(\ZZ)$, mais pas dans $\SIPD_1(\QQ)$.

\smallskip \noindent 2. Soit $\gK=\QQ(\epsilon)$ avec $\epsilon$ infinitésimal positif et $\gR$ la clôture réelle de $\gK$. 
\\
Le \spo $f=x^+\vi-(x^2-\epsilon)(x^3-\epsilon)$ définit la fonction nulle sur $\gK$ mais ne définit pas une fonction nulle sur $\gR$: l'intervalle $[\epsilon^{1/2},\epsilon^{1/3}]$ est invisible sur $\gK$. On peut \gui{simplifier} cet exemple
en prenant $\gK=\QQ[\epsilon]$ avec $\epsilon>0$ nilpotent convenable.
\eoe
\end{examples}
%--------- fin example ---------------------------------------------- 

%:\newpage
%\newpage
\section{Au delà des \tpes} \label{secArftr}
%%%%%%%%%%%%%%%%%%%%%%%%%%%%%%%%%%%%%%%%%%%%%%%%%%%%%%%%%%%%%%%%%%%%

%%%%%%%%%%%%%%%%%%%%%%%%%%%%%%%%%%%%%%%%%%%%%%%%%%%%%%%%%%%%%%%%%%%%
\Subsection{Anneaux \frls \sdz} 
%:Subsection{Anneaux \frls \sdz}

%: Lemma{lemAfrsdz}
\begin{lemma} \label{lemAfrsdz} 
Un \afr \sdz est totalement ordonné. En d'autres termes, si l'on ajoute l'axiome \Tsbf{ASDZ} à la théorie \sa{Afr},
 la règle \Tsbf{OT} est valide. Autrement dit encore, la théorie \SA{Afrsdz} que l'on obtient 
est \esid à la théorie \Sa{Atonz} des anneaux totalement ordonnés \sdz (voir le point 3 du lemme \ref{lemAtonz}).

\DeuxRegles{\lab{ASDZ} {$\,\,xy=0\vd x=0 \;\vou\;y=0$}}
{\lab{OT} {$\vd x\geq 0 \;\vou\;-x\geq 0$}}
 
\end{lemma}
%----------- fin lemma ----------------------------------- 
%
\begin{proof}
En effet comme $x^+x^-=0$ (\Tsbf{afr1}), on obtient la règle valide

\UneRegle{~} {$\vd x^+=0 \;\vou\;x^-=0$}

\vspace{-1em}
\end{proof}

%%%%%%%%%%%%%%%%%%%%%%%%%%%%%%%%%%%%%%%%%
\Subsection{Anneaux fortement réticulés locaux}
%: Subsection{Anneau fortement réticulé local}

%l
%: Lemma{lemAfrLoc0}
\begin{lemma} \label{lemAfrLoc0}
Soit $\gA$ un \afr local et $x\in\Ati$, alors $x$ est $\geq 0$ ou $\leq 0$.
\end{lemma}
%----------- fin lemma ----------------------------------- 
%
\begin{proof}
Soit $x\in\Ati$, on écrit $x=x^+-x^-$, donc $x^+\in\Ati$ ou $x^-\in\Ati$.
Or $x^+x^-=0$. Dans le premier cas on obtient $x^-=0$, dans le second cas $x^+=0$. 
\end{proof}
%

%%%%%%%%%%%%%%%%%%%%%%%%%%%%%%%%%%%%%%%%%%%%%%%%%%%%%%%%%%%%%%%%%%%%
%:Subsection{Anneaux strictement réticulés}
\Subsection{Anneaux strictement réticulés}
\rdb

La théorie suivante \gui{fusionne} les théories \sa{Afr} et \sa{Aso}.
Cette théorie est \esid à celle définie dans l'article \cite{LM2017}.
%: Definition{defiAsr}
\begin{definition} \label{defiAsr}
Le langage de la \talg \SA{Asr} des \textsl{\asrs} est donné par la signature suivante.
\Sigt{\AsR}{\cdot=0,\cdot\geq 0,\cdot>0\mathrel{;}\cdot+\cdot, \cdot\times\cdot,\cdot\vu\cdot,-\,\cdot,0,1} \label{NOTASigAsr}
\noindent Les axiomes sont les suivants.\index{anneau!strictement réticulé}
%i
\begin{itemize}
\item les règles de la théorie purement équationnelle \Sa{Afr},
\item les \reds de \Tsbf{aso1} à \Tsbf{aso4}, (\paref{Axaso1})
\item les \ralgs
 \coligt, \Tsbf{Iv}, \Tsbf{Aso1} et \Tsbf{Aso2}, (\paref{AxAso1})
\item enfin, on a les trois règles \Tsbf{sup1}, \Tsbf{sup2} et \Tsbf{Sup} (\paref{Axsup1})
 pour relier $\cdot\geq 0$ \hbox{et $\cdot\vu\cdot$}.
\end{itemize}

\end{definition}
%--------- fin definition --------------------------------

Nous avons mis le prédicat~\hbox{\gui{$\cdot\geq 0$}} directement dans le langage plutôt que le définir à partir de la \hbox{loi $\cdot\vu\cdot$}.
 
La signification de $x>0$ n'est pas fixée à priori par les axiomes.
Cela peut aller de \gui{$x$ est régulier et $\geq 0$} jusqu'à \gui{$x$ est inversible et $\geq 0$}

% En passant de \Sa{Afr} à \Sa{Asr} on sort du cadre des théories purement équationnelles. % déjà annoncé dans le titre de la secion

%: Lemma{lemAsrs}
\begin{lemma} \label{lemAsrs}~
%
%\begin{enumerate}
%%
%\item 
On considère la \talg des \asrs à laquelle on ajoute l'axiome \Tsbf{OTF}.
Alors est \egmt valide la règle \OTFx (on rappelle ces règles ci-dessous). 

\DeuxRegles{
\lab{OTF} $\,\, x+y> 0 \vd x >0 \;\vou\; y>0 $
}
{
\lab{OTF$\eti$} $\,\, xy< 0 \vd x <0 \;\vou\; y<0 $
}

%
%\end{enumerate}
\end{lemma}
%--------- fin lemma -----------------------------------
%
\begin{proof} 
On suppose $xy<0$, d'où $x^2y^2>0$, d'où, par \Tsbf{Aso2}, $x^2>0$.
\\
On note que $x^2=(x^+)^2+(x^-)^2$. Donc par \Tsbf{OTF}, il suffit de traiter séparément les \hbox{cas $(x^+)^2>0$} et $(x^-)^2>0$. \\
Si $(x^+)^2>0$, on a $xx^+=(x^+)^2>0$, donc par \Tsbf{Aso2}, $x>0$. Et de nouveau par \Tsbf{Aso2}, on obtient~\hbox{$y<0$}.\\
Si $(x^-)^2>0$, on a $-xx^-=(x^-)^2>0$, donc par \Tsbf{Aso2}, $x<0$.
\end{proof}
%

%%%%%%%%%%%%%%%%%%%%%%%%%%%%%%%%%%%%%%%%%%%%%%%%%%%%%%%%%%%%%%%%%%%%
%:Subsection Théorie des anneaux fortement réticulés réduits
\Subsection{Anneaux fortement réticulés réduits}
\rdb

Nous examinons ici la \talg \SA{Afrnz} des \textsl{\afrs réduits}. On ajoute donc à \Sa{Afr}
l'axiome \Tsbf{Anz} des anneaux réduits, qui est une \rsim.

\Regles{
\laB{Anz} $\,\, a^2=0 \vd a =0$}%

De même la \talg \SA{Asrnz} des \textsl{\asrs réduits} est la théorie obtenue à partir de la théorie \Sa{Asr} en ajoutant la \ralg \tsbf{Anz}.

%%%%%%%%%%%%%%%%%%%%%%%%%%%%%%%%%%%%%%%%%%%%%%%%%%%%%%%%%%%%%%%%%%%%
\paragraph{Quelques règles dérivées}~

\smallskip Démontrons dans la théorie \sa{Afrnz} les quatre règles règles \Tsbf{Afrnz1}, \Tsbf{Afrnz2}, \Tsbf{Aonz}, et \Tsbf{Aonz3}
(les deux dernières ont été introduites pages \pageref{AxAonz} et \pageref{AxAonz3}).

\Regles{
\Lab{Afrnz1} $\,\, x^{3}\geq 0 \vd x\geq 0$%
}

\noindent On écrit $x=x^+-x^-$. Comme $x^+\,x^-=0$ on a $x^{3}=(x^+)^{3}-(x^-)^{3}\geq 0$. On multiplie par $x^-$, il vient $(x^-)^{4}\leq 0$, donc $(x^-)^{4}= 0$.
Or l'anneau est réduit: $x^-=0$ et $x\geq 0$.\eop

\smallskip Notez que de \tsbf{Afrnz1} on déduit la même règle pour un exposant impair arbitraire qui remplace l'exposant 3.

\smallskip On a aussi la réciproque suivante de la règle \Tsbf{Afr1}.

\Regles{
\Lab{Afrnz2} $\,\, ab=0 \vd \abs a \vi \abs b =0$%
}

\noindent En effet si $ab= 0$, alors $(\abs a\vi\abs b)^2\leq \abs a\,\abs b=0$, donc $\abs a\vi\abs b=0$. \eop

\smallskip Ainsi, pour $a,b\geq 0$, $ab=0$ équivaut à $a\vi b=0$.

\Regles{
\laB{Aonz} $\,\, c\geq 0\vet x(x^2+c)\geq 0 \vd x\geq 0\,$ 
}

\noindent En effet, par \Tsbf{Ato2} nous avons $x^3\geq 0$, d'où $x\geq 0$ par \Tsbf{Afrnz1}. 
\smallskip  Et enfin: 

\Regles{
\laB{Aonz3} $\,\, a\geq 0\vet b\geq 0\vet a^2=b^2 \vd a=b$%
}

\noindent En effet $ \abS{a-b}^2\leq \abS{a-b}\, \abS{a+b} =\abs{a^2-b^2}$.
\eop

\medskip On a maintenant facilement le résultat suivant. 

%: Lemma{lemAfr-Afrnz}
\begin{lemma} \label{lemAfr-Afrnz}
Dans la théorie \Sa{Afr} les règles \Tsbf{Afrnz1}, \Tsbf{Afrnz2}, \Tsbf{Aonz3},
\Tsbf{Aonz} et~\Tsbf{Anz} sont \eqves. 
\end{lemma}
%--------- fin lemma ----------------------------------- 

Voici un exemple simple d'application du \pstfref{Pst2}.

%l
%: Lemma{lemsupdansAonz}
\begin{lemma} \label{lemsupdansAonz}
Dans un \afr réduit, l'\elt $c=a\vu b$ est caractérisé par les \egts et inégalités suivantes
$$
c\geq a,\;c\geq b, \;(c-a)(c-b)=0.
$$
Plus \prmt, la théorie \sa{Afrnz} prouve la \ralg suivante

\Regles{\lab{~}$\,\,x\geq a\vet x\geq b\vet (x-a)(x-b)=0\vd x=a\vu b$.} 
\end{lemma}
%----------- fin lemma ----------------------------------- 
%
\begin{proof}
En effet, comme la règle est valide pour la théorie \Sa{Codsup}, cela résulte du point~\textsl{3} du \pstfref{Pst2}.
\end{proof}
%

%%%%%%%%%%%%%%%%%%%%%%%%%%%%%%%%%%%%%%%%%%%%%%%%%%%%%%%%%%%%%%%%%%%%
\paragraph{La règle \tsbf{FRAC} dans \sa{Afrnz}}~

\smallskip Rappelons les règles \Tsbf{FRAC} et \FRACn.

\Regles{\lab{FRAC} $\,\,0\leq a\leq b\vd \Exists z\; (zb=a^2\vet\,0\leq z\leq a)$ 
\lab{FRAC$_n$} $\,\,\abs u^n\leq \abs v^{n+1}\vd \Exists z\; (zv=u\vet\, \abs z^{n}\leq \abs v)$ \quad ($n\geq 1$)}

Notons que la règle \tsbf{FRAC}, appliquée avec $a=1$ implique l'inversibilité de tout \elt $\geq 1$.

\smallskip On rappelle maintenant pour la théorie \Sa{Afrnz} l'analogue du lemme \ref{lemCo0FRAC} pour la théorie \Sa{Co--}. 
%: Lemma{lemAfrnzFRAC}
\begin{lemma} \label{lemAfrnzFRAC}
L'ajout de l'axiome \tsbf{FRAC} à la théorie \sa{Afrnz} peut être remplacé par l'introduction d'un symbole de fonction $\Fr$ avec les axiomes
\tsbf{fr1} et \tsbf{fr2} que nous rappelons ci-dessous

\DeuxRegles{
\lab{fr1} $ \vd \Fr(a,b)\, \abs b=(\abs a\!\vi\! \abs b)^2 $
}
{
\lab{fr2} $ \vd 0\leq {\Fr(a,b)}\leq \abs a\!\vi\! \abs b \phantom{(\abs a\!\vi\! \abs b)^2}$
}
\end{lemma}
%----------- fin lemma ----------------------------------- 
%
\begin{proof} Comme dans le lemme \ref{lemCo0FRAC}, on constate qu'il s'agit de skolémiser une règle existentielle. Cela donne une théorie \esid lorsqu'il s'agit d'une règle existentielle à existence unique.
Pour le vérifier, la \demo est celle du lemme \ref{lemUniqFRAC}. On suppose $yb=zb=a^2$, $0\leq y\leq a$ et $0\leq z\leq a$. On a ${(y-z)}\,b=0$, $\abs{y-z}\leq a\leq b$\footnote{On est dans un \grl pour l'addition, on peut raisonner cas par cas, séparément avec $0\leq y\leq z\leq a$ et $0\leq z\leq y\leq a$.
On obtient dans les deux cas $0\leq \abs{z- y}\leq a$.} et \hbox{donc $\abs{y-z}^2\leq\abs{y-z}\,b=0$}. 
\end{proof}
%

%l
%: Lemma{lemfracfrac2}
\begin{lemma} \label{lemfracfrac2}
Dans la théorie \sa{Afrnz} la règle $\tsbf{FRAC}_2$ se déduit de la règle \tsbf{FRAC} et de la règle affirmant l'existence de la racine sixième $\geq 0$ d'un \elt $\geq 0$.
\end{lemma}
%----------- fin lemma ----------------------------------- 
%
\begin{proof}
On suppose $u^2\leq v^3$ et on veut trouver un $z$ tel que $zv=u$ et $z^2\leq v$. 

\smallskip \noindent 
Supposons d'abord $u\geq 0$ et montrons que l'on a un $z$ tel que $zv=u$. La règle \tsbf{FRAC} implique que la fraction $t=\frac{u^4}{v^3}$ est bien définie avec $0\leq t\leq u^2$. De nouveau \tsbf{FRAC} donne le fait que la fraction $w=\frac{t^2}{u^2}$ avec $0\leq w\leq t\leq u^2$ est bien définie. On obtient alors $u^2wv^6=t^2v^6=u^8$. Donc 
$u^2(wv^6-u^6)=0$. Or $w\leq u^2$, donc $\abs{wv^6-u^6}\leq u^2(v^6+u^4)$,
d'où $\abs{wv^6-u^6}^2\leq \abs{wv^6-u^6} u^2(v^6+u^4)=0$. Vu la règle \tsbf{Anz}, on obtient $wv^6=u^6$. On prend $z=w^{\frac1 6}$ et l'on a $zv=u$. En outre $z^6\leq u^2\leq v^3$ implique $z^2\leq v$.

\smallskip \noindent Pour un $u$ arbitraire on écrit $u=u^+-u^-$; on a $u^2=(u^+)^2+(u^-)^2\leq v^3$. On obtient un $z_1$ tel que $z_1v=u^+$ et un 
$z_2$ tel que $z_2v=u^-$, on pose $z=z_1-z_2$ et l'on a $zv=u$.
On obtient en outre $z^2\leq 2v$. Comme $z^2v^2=u^2\leq v^3$, on a $v^2(z^2-v)\leq 0$. Les inégalités $0\leq z^2\leq 2v$ impliquent $\abs{z^2-v}\leq v$, d'où $(z^2-v)^2\leq v^2$ et $(z^2-v)^3\leq v^2(z^2-v)\leq 0$, ce qui implique $z^2-v\leq 0$.
\end{proof}
%

%

%%%%%%%%%%%%%%%%%%%%%%%%%%%%%%%%%%%%%%%%%%%%%%%%%%%%%%%%%%%%%%%%%%%%
%:Subsection Retour sur les corps ordonnés
\section{Retour sur les corps ordonnés}
\label{secCOG}

%%%%%%%%%%%%%%%%%%%%%%%%%%%%%%%%%%%%%%%%%%%%%%%%%%%%%%%%%%%%%%%%%%%%
%:Subsection Anneaux réticulés \ftm réels 
\Subsection{Anneaux (réticulés) \ftm réels}

Le corps des réels vérifie toutes les \ralgs de la théorie des corps réels clos discrets, mais pas toutes les \rdys.
Rappelons les \rdys suivantes satisfaites par les \codis.

\DeuxRegles{
\laB{IV} $\,\, x> 0 \vd \Exists y\, xy = 1$
\laB{OTF} $\,\, x+y> 0 \vd x > 0\;\vou\;y>0$
\laB{FRAC} $\,\,0\leq a\leq b\vd \Exists z\; (zb=a^2\vet 0\leq z\leq a)$
}
{
\lab{\Edinq} $ \vd x^2> 0 \;\vou\; x= 0$
\laB{OT} $ \vd x \geq 0 \;\vou\; x\leq 0$
}

Le corps des réels vérifie \tsbf{IV}, \tsbf{OTF} et \tsbf{FRAC} mais ni \Ednq, ni \tsbf{OT}.

\smallskip Le lemme suivant prépare la \dfn de la \tdy \Sa{Aftr}.
%: Lemma{lemArftr0}
\begin{lemma} \label{lemArftr0}~ 
\label{i1lemArftr0} Dans la théorie \Sa{Afrnz} à laquelle on ajoute la règle \Tsbf{FRAC}, si l'on définit \gui{$x>0$} comme abréviation de \gui{$x\geq 0 \vii \exists z\; zx=1$}, le prédicat $\cdot>0$ satisfait tous les axiomes de \Sa{Asrnz} où il est présent ainsi que la règle \Tsbf{IV}.
\end{lemma}
%----------- fin lemma ----------------------------------- 

\begin{proof} Les règles \Tsbf{IV}, \Tsbf{aso1}, \Tsbf{aso2}, \Tsbf{aso4}, \CLinqAc\ sont trivialement valides. Voyons la règle \Tsbf{aso3}:
un \elt supérieur ou égal à un \elt \iv positif est \iv. Cela résulte de la règle \tsbf{FRAC} car si $b\geq a>0$, on a un~$z$ tel \hbox{que $zb=a^2$}, donc $b$ est \iv. Pour \Tsbf{Aso1}, \Tsbf{Aso2} et \Tsbf{Iv}, on commence par valider la règle
$\,\,x>0\vet xu=1\vd u\geq 0$: en effet $u=xu^2\geq 0$. Le reste suit.\end{proof}

%: Definition{defiAftr}
\begin{definition} \label{defiAftr}~ 
Nous définissons maintenant la \tdy \SA{Aftr} des \textsl{\arftrs} (ou pour abréger, des \textsl{\aftrs}) sur la signature suivante.
\Sigt{\AftR}{\cdot=0,\cdot\geq 0,\cdot>0\mathrel{;}\cdot+\cdot, \cdot\times\cdot,\cdot\vu\cdot,-\,\cdot,\Fr(\cdot),0,1} \label{NOTASigAftr} 
\noindent Les axiomes de la théorie \sa{Aftr} sont les suivants:
\begin{itemize}
\item les axiomes de \Sa{Afrnz};
\item les axiomes qui définissent $x>0$ comme une abréviation de \gui{$x\geq 0\vii \exists z\,xz=1$}; 
\item les axiomes \Tsbf{fr1} et \Tsbf{fr2}. 
\end{itemize}
\index{anneau!fortement reel@fortement réel}%
%
%\end{enumerate}
\end{definition}
%--------- fin definition --------------------------------

Cette \dfn est justifiée par le fait que la théorie des \aftrs \und{locaux} est \esid à la théorie \sa{Co} des corps ordonnés \textsl{non} discrets: point~\emph{3} du lemme~\ref{lemArftr}.
%l
%: Lemma{lemArftr01}
\begin{lemma} \label{lemArftr01}\label{i2lemArftr0} La \tdy \sa{Aftr} est \esid à la théorie \sa{Asrnz} à laquelle on ajoute les axiomes \Tsbf{IV} et \Tsbf{FRAC}\footnote{Cela implique que la théorie {\ssa{Aftr}} définie ici est légèrement plus forte que celle définie dans \cite{LM2017}.}.
\end{lemma}
%----------- fin lemma ----------------------------------- 
%
\begin{proof} 
La théorie \sa{Asrnz} contient dans sa signature le prédicat~\gui{$\cdot>0$}
qui n'est pas présent dans \sa{Afrnz}. Quand on ajoute l'axiome \tsbf{IV}, on a $x>0$ \ssi $x$ est $\geq 0$ et inversible.
Le lemme \ref{lemArftr0} implique donc que la théorie \sa{Afrnz} à laquelle on ajoute l'axiome \tsbf{FRAC} est \esid à théorie \sa{Asrnz} à laquelle on ajoute les axiomes \tsbf{IV} et \tsbf{FRAC}. 
Enfin le lemme \ref{lemAfrnzFRAC} montre qu'ajouter l'axiome \tsbf{FRAC} à \sa{Afrnz} revient à ajouter la fonction $\Fr$ avec les axiomes \tsbf{fr1} et~\tsbf{fr2}.
\end{proof}

Un \aftr est donc une \QQlg \ftm réticulée réduite dans laquelle tout \elt supérieur à un \elt positif \iv est lui même \iv.
En outre, la validité de la règle~\tsbf{FRAC} ajoute un petit quelque chose.

%: Lemma{lemArftr}
\begin{lemma} \label{lemArftr}~
\begin{enumerate}
\item 
 La \talg \Sa{Asrnz} à laquelle on ajoute comme axiomes les \rdys \Tsbf{IV} et \Tsbf{OTF} est \esid à la théorie \Sa{Co--}.
\item La \tdy \Sa{Aftr} à laquelle on ajoute l'axiome \Edinq\  
est \esid à la théorie \Sa{Codsup} des \codis avec sup (ou encore à \Sa{Cod}).
\item La \tdy \Sa{Aftr} à laquelle on ajoute l'axiome \Tsbf{OTF}
est \esid à la théorie \Sa{Co} des corps ordonnés \emph{non} discrets:
 un \aftr local est un corps ordonné \emph{non} discret\footnote{On aurait pu éviter d'introduire le prédicat $x>0$ dans \ssa{Aftr} car il est défini comme une abréviation. Il faudrait remplacer la règle \Tsbf{OTF} par la règle \Tsbf{AFRL} dans le point 2. La théorie serait alors une \talg. Cela ne saurait vraiment nous étonner puisque la théorie des \arcs est \peq et qu'un \crc \textsl{non} discret est un \arc local.}.
% dans lequel on n'utilise pas le prédicat $\cdot>0$.
%\item 
%
\end{enumerate}
\end{lemma}
%----------- fin lemma ----------------------------------- 
%
\begin{proof} \textsl{1.} En comparant les théories \sa{Co--} (\dfn \ref{defiCo0}) et \sa{Asrnz}, on trouve dans la première les
axiomes supplémentaires \tsbf{Aonz}, \tsbf{IV} et~\tsbf{OTF}
et il manque l'axiome de collapus. 
Mais le collapus résulte de \tsbf{IV} et \tsbf{Aonz} se déduit de \tsbf{Anz} (lemme \ref{lemAfr-Afrnz}). \\
NB: L'axiome~\tsbf{IV} implique que $x>0$
équivaut à \gui{$x$ est $\geq 0$ et inversible}. L'axiome \tsbf{OTF} ajoute le fait que l'anneau est local. 

\smallskip \noindent \textsl{2}.
La nouvelle théorie est une extension de \sa{Co--} d'après le point \textsl{1}. 
Pour passer à \sa{Codsup} il suffit d'ajouter \Edinq\  et \Tsbf{OT}. Comme tout \elt est nul ou inversible, on a affaire à un corps discret, et la règle \tsbf{ASDZ} est valide.
Le lemme~\ref{lemAfrsdz} dit alors que la règle \tsbf{OT} est \egmt valide.

\smallskip \noindent \textsl{3}. Résulte du point \textsl{1} et des lemmes \ref{lemAfrnzFRAC},
\ref{lemArftr0} et \ref{lemAftrloc}.
\end{proof}
%

%l
%: Lemma{lemAftrloc}
\begin{lemma} \label{lemAftrloc}
Dans un \afr réduit où la règle \tsbf{FRAC} est valide, la règle suivante \tsbf{AFRL} remplace la règle \Tsbf{OTF} lorsqu'on définit le prédicat $a >0$ comme abréviation de $a\geq 0\vii \exists z\;az=1$.

\Regles{\Lab{AFRL} $\;\;z(x+y)=1\vet x+y\geq 0\vd \Exists u\;(ux=1\vet x\geq 0 )\;\vou\; \Exists v\;(vy=1\vet y\geq 0 )$ }

\end{lemma}
%----------- fin lemma -----------------------------------
%
\begin{proof} \textsl{Implication directe}. On suppose $\gA$ local et on démontre la règle \tsbf{AFRL}. Puisque $x+y$ est \iv, $x$ ou $y$ est \iv. Par exemple~\hbox{$x\in\Ati$}. D'après le lemme \ref{lemAfrLoc0} on a $x\geq 0$ ou $x\leq 0$.
Si $x\leq 0$, alors $y\geq x+y\geq 0$, donc $y\in \Ati$ (lemme \ref{lemArftr0}).\\
\textsl{Réciproque}. Si $x+y$ est \iv, alors $x+y\geq 0$ ou $x+y\leq 0$ (lemme \ref{lemAfrLoc0}).
Dans le premier cas, \tsbf{AFRL} montre que $x$ ou $y$ est \iv. Dans le second cas,
$-x$ ou $-y$ est \iv.
\end{proof}

\Subsection{Posivitstellensätze formels avec sup}
% section{Posivitstellensatz formel avec sup}
\rdb

Pour ce qui concerne l'effondrement simultané, on a déjà donné le \thref{thColsimafr}.

%: pstf {Pst2}
\begin{pstf}[Positivstellensatz formel, 2] \label{Pst2} ~%
\index{Positivstellensatz!formel!pour les corps ordonnés, 2}\\ 
Les \tdys suivantes prouvent les mêmes \ralgs.
\begin{enumerate}
\item \label{i1Pst2} Sur la signature \sigt{\Ao}{\cdot=0,\cdot\geq 0\mathrel{;}\cdot+\cdot, \cdot\times\cdot,-\,\cdot,0,1}: la théorie \Sa{Aonz} des anneaux ordonnés \stm réduits (\dfn \ref{defisaAor}), la théorie \Sa{Crcdsup} des \crcds avec sup et les théories intermédiaires 
(par exemple \Sa{Afrnz}, \Sa{Atonz}, \Sa{Cod} ou \Sa{Co--}).

\item \label{i2Pst2} Sur la signature \sigt{\Aso}{\cdot=0,\cdot\geq 0,\cdot>0\mathrel{;}\cdot+\cdot, \cdot\times\cdot,-\cdot, 0,1}: la théorie \Sa{Asonz} des anneaux \stm ordonnés réduits, la théorie \Sa{Crcdsup} et les théories intermédiaires 
(par exemple \Sa{Asrnz}, \Sa{Cod} ou~\Sa{Co--}).
\item \label{i3Pst2} Sur la signature \sigt{\AfR'}{\cdot=0,\cdot\geq 0\mathrel{;}\cdot+\cdot, \cdot\times\cdot,\cdot\vu\cdot,-\,\cdot,0,1}:
la théorie \Sa{Afrnz} des \afrs réduits, la théorie \Sa{Crcdsup} et les théories intermédiaires (par exemple \Sa{Asrnz} ou  \Sa{Co--}). 
\item \label{i4Pst2} Sur la signature \sigt{\AsR}{\cdot=0,\cdot\geq 0,\cdot>0\mathrel{;}\cdot+\cdot, \cdot\times\cdot,\cdot\vu\cdot,-\,\cdot,0,1}: la théorie \Sa{Asrnz} des \asrs réduits, la théorie \Sa{Crcdsup} et les théories intermédiaires 
(par exemple~\Sa{Co} ou \Sa{Codsup}).
\item \label{i5Pst2} Sur la signature \sigt{\Co}{\cdot=0,\cdot\geq 0,\cdot>0\mathrel{;}\cdot+\cdot, \cdot\times\cdot,\cdot\vu\cdot,-\,\cdot,\Fr(\cdot,\cdot),0,1}: la \tdy \Sa{Aftr} des \aftrs, la théorie \Sa{Crcdsup} et les théories intermédiaires (par exemple \Sa{Co} ou \Sa{Codsup}).
\end{enumerate}
\end{pstf}
%--------- fin pstf ---------------------------------------------- 

%
\begin{proof} \textsl{1} et \textsl{2}. La théorie \sa{Crcdsup} est \esid à \sa{Crcd}. Donc
les Positivstellensätze \ref{Pst1bis} et \ref{Pst1} donnent les points \textsl{1} et \textsl{2}. 

\smallskip \noindent 
\textsl{3}. Les théories \sa{Atonz} et \sa{Crcd} prouvent les mêmes \ralgs (\ref{Pst1bis}). Les théories \sa{Ato} et \sa{Afr} prouvent les mêmes \ralgs (\pstref{thfairecommesi}), donc aussi les théories \sa{Atonz} et \sa{Afrnz} (\thref{thEseqMemesfaitsbis}).
Enfin  \sa{Crcd} et  \sa{Crcdsup}  sont \esids.

\smallskip \noindent 
\textsl{4}. Même raisonnement qu'au point précédent.

\smallskip \noindent 
\textsl{5}. On note que \sa{Codsup} valide les règles \tsbf{IV} et \tsbf{FRAC}, laquelle peut être remplacée par l'introduction de $\Fr$
avec ses deux axiomes (lemme \ref{lemAfrnzFRAC}). 
Par ailleurs la \tdy \sa{Aftr} est \esid à la théorie \sa{Asrnz} à laquelle on ajoute les axiomes \tsbf{IV} et \tsbf{FRAC} (point \textsl{2} du lemme~\ref{lemArftr0}). Le point \textsl{5} résulte donc du point \textsl{4}.
\end{proof}
%

%r
%: Remark{remPst2}
\begin{remark} \label{remPst2} 
Dans les théories \Sa{Aftr} et \Sa{Co}, on a ajouté le symbole de fonction
$\Fr$ avec les axiomes \Tsbf{fr1} et \Tsbf{fr2}, ce qui augmente les \ralgs formulables dans ces théories.
Néanmoins l'usage du symbole $\Fr$ peut être éliminé dans les \ralgs
au profit l'axiome \Tsbf{FRAC} (voir l'ajout d'un symbole de fonction \paref{skolemunique} et le lemme \ref{lemCo0FRAC}). Comme cet axiome est satisfait dans la théorie plus forte \Sa{Cod}, le 
\pstfref{Pst2} n'est pas affecté par la présence de $\Fr$. 
On pourrait d'ailleurs accepter la présence du symbole de fonction $\Fr$ avec les axiomes \Tsbf{fr1} et \Tsbf{fr2} dans la théorie \Sa{Codsup}.
\eoe\end{remark}
%----------- fin remark ---------------------------------- 

%c
%: Corollary{corPst2}
\begin{corollary} \label{corPst2} On considère la \tdy $\Sa{Asrnz}= \Sa{Asrnz}(\QQ)$ des \asrs réduits.
\begin{enumerate}
\item Soit $\gK$ un corps ordonné discret et $\gR$ sa clôture \agq. Soit $\gA=\big((G,Rel),\Sa{Asrnz}(\gK)\big)$ une \sad avec \hbox{$G=(\xn)$} et $Rel$ fini. 
On a un \algo qui décide si $\gA$ s'effondre et qui en cas de réponse négative donne la description d'un \sys $(\xin)$ dans $\gR^n$ qui satisfait les contraintes données dans les relations $Rel$.
\item On a un \algo qui décide si une \ralg de $\sa{Asrnz}$ est valide. En cas de réponse négative l'\algo donne la description d'un \sys $(\xin)$ dans $\gRa^n$ qui contredit la \ralg.
\item Les mêmes résultats sont valables avec $\Sa{Aftr}= \Sa{Aftr}(\QQ\cap[0,1])$
à la place de $\sa{Asrnz}$. Dans ce cas on ajoute à \Sa{Crcdsup} le symbole de fonction $\Fr$ avec les deux axiomes qui l'accompagnent.
\end{enumerate}
\end{corollary}
%--------- fin corollary -------------------------------
%
\begin{proof} 
On vient de voir (\pstfref{Pst2}) que la \talg \Sa{Asrnz} (resp. \Sa{Aftr}) prouve les mêmes \ralgs que \Sa{Crcdsup} (resp. en ajoutant $\Fr$). Par ailleurs on constate qu'une \ralg de \Sa{Crcdsup} (éventuellement en ajoutant $\Fr$) est \eqve à une famille de \ralgs de \Sa{Crcd}. On peut donc conclure par le Positivstellensatz concret~\ref{thPstStengle}.
\end{proof}

L'anneau $\SIPD_n(\gA)$ des \spos sur $\gA$ en $n$ variables est défini en \ref{defiSIPD}, l’\afr $\AFR(\gA)$ engendré par $\gA$ est défini en \ref{notaAFR}.

%: Theorem{thSIPDreduit}
\begin{theorem} \label{thSIPDreduit} On fixe $n$ et on note $\Kux=\Kxn$.
\begin{enumerate}
\item Soit $\gK$ un \codi et $\gR$ sa clôture réelle. L'anneau $\SIPD_n(\gK,\gR)$ s'identifie à l'\afr engendré par $\Kux$.
Plus \prmt:
la structure de~$\gK$ confère~à~$\Kux$ une structure d'\afr dynamique et l'unique~\hbox{$\gK$-morphisme} d'\afrs de $\AFR(\Kux)$ \hbox{vers $\SIPD_n(\gK,\gR)$} est un \iso. 
\item Soit $\gK$ un \codi et $\gR$ sa clôture réelle. Si tout ouvert \sagq de~$\gR^n$ contient des points de $\gK^n$, l'anneau $\SIPD_n(\gK)$ s'identifie à $\AFR(\Kux)$.
\item (\demo incomplète) Si $\gK$ est une \QQlg contenue dans $\RR$, l'anneau $\SIPD_n(\gK)$ s'identifie à $\AFR(\Kux)$.
\end{enumerate}
\end{theorem}
%--------- fin theorem -----------------------------------
%
\begin{proof} Il faut montrer que si une expression de la forme donnée dans le lemme~\ref{lemSIPD} définit la fonction identiquement nulle, cela peut être démontré en utilisant uniquement
les \ralgs des \afrs réduits.

\sni \textsl{1.}
Par le Positivstellensatz, le fait qu'un \spo est nul en tout point de $\gR^n$ possède un certificat \agq sur $\gK$. Or \Sa{Afrnz} et \Sa{Crcdsup} prouvent les mêmes \ralgs. (Pour plus de détails sur ce genre de sujet on peut voir \cite{GL93}). 

\sni \textsl{2.} Résulte du point précédent car sous l'hypothèse envisagée, un $\gK$-\spo non partout nul sur $\gR^n$ est non partout nul sur $\gK^n$. 

\sni \textsl{3.} Si $\gK$ est discret, cela résulte du point \textsl{1},
car un $\gK$-\spo nul sur $\gK$ est nul sur~$\QQ$ donc sur~$\RR$ et à fortiori sur la clôture réelle du corps de fractions de $\gK$. Apparemment, il faut se fatiguer un peu pour obtenir \cot le résultat en toute généralité, alors qu'il est clair en \clama. C'est le même genre de gymnastique 
que pour la \demo \cov complète de la solution 17\ieme\ \pb de Hilbert
sur $\RR$, donnée dans \cite{GL93}. Le bonus est que la solution est alors complètement explicite, \cad qu'elle n'utilise pas de test de signe (ni d'axiome du choix dépendant) sur $\RR$.
\end{proof}

\hum{En l'état actuel la \demo du point \textsl{3} est incomplète.
Elle est sans doute délicate. À mettre en annexe?}

%section{Treillis reel} 
\section{Le treillis et le spectre réels d'un anneau commutatif} 
\label{subsecTRR}
%-----------------------------------------

Un cône premier de l'anneau commutatif $\gA$, \cad un \elt du spectre réel 
$\Sper\,\gA$, peut être donné sous forme d'un anneau quotient intègre non 
trivial $\gA/\fP$ muni d'une structure d'anneau totalement ordonné, autrement dit comme un modèle minimal 
de la théorie $\Sa{Aito}(\gA)$ des anneaux totalement ordonnés intègres non triviaux basée sur l'anneau $\gA$ 
(voir la \dfn \ref{defidiagramme}).

Comme la théorie des anneaux totalement ordonnés intègres non triviaux est une théorie dynamique sans axiomes existentiels, le spectre réel s'identifie au spectre du \trdi $\Reel(\gA)$ obtenu en \gui{recopiant\footnote{Comme on le fait ci-dessous, dans la \dfn \ref{defZarR}.}} les axiomes 
de la théorie dynamique $\sa{Aito}(\gA)$.

Pour retrouver sur l'ensemble $\Spec(\Reel\gA)$ la topologie usuelle\footnote{Selon la tradition établie lors de l'invention du spectre réel.} du spectre réel $\Sper\gA$, nous devons considérer le treillis basé sur le seul 
prédicat $x>0$. Ceci donne la définition \ref{defZarR}, qui correspond aux \rdys valides suivantes dans \sa{Cod}

\DeuxRegles
{
\lab{\coligt} $\,\,0> 0 \vd \Bot$
\laB{aso3} $\,\, x > 0\vet y > 0 \vd x +y > 0$
\laB{aso4} $\,\, x > 0\vet y > 0 \vd x y > 0$
}
{
\laB{aso1} $ \vd 1 > 0$ 
\laB{OTF} $\,\, x+y> 0 \vd x > 0\vou y>0$
%\lAb{\OTFx} $\,\, xy< 0 \vd x <0 \vou y<0 $
\lab{OTFx} $\,\, xy< 0 \vd x <0 \vou y<0 $
}

%r
%: Remark{remdefZarR}
\begin{remark} \label{remdefZarR} 
Si l'on se base sur le seul prédicat $x>0$ et si l'on introduit le prédicat $x\geq 0$ comme opposé au prédicat $-x>0$, et le prédicat $x=0$ comme la conjonction $x\geq 0 \vii -x\geq 0$, on obtient sur la base des seuls axiomes précédents une extension conservative qui satisfait tous les axiomes de \sa{Aito}. Les modèles minimaux de la \tdy décrite par les 6 axiomes précédents sont donc bien (en \clama) les quotients intègres de $\gA$ munis d'une relation d'ordre total. Ceci justifie la \dfn qui suit: le spectre du treillis $\Reel\gA$
s'identifie bien au spectre réel de $\gA$ (en \clama).
\eoe
\end{remark}
%----------- fin remark ---------------------------------- 

%--Definition{defZarR}---- 
\begin{definition} 
\label{defZarR} 
Le \textsl{treillis réel}, noté $\Reel\gA$, d'un anneau commutatif $\gA$ 
est le \trdi engendré par $(\gA,\vdash)$ o\`u $\vdash$
est la plus petite \entrel vérifiant
\begin{equation} %\label{} 
\left.
\begin{array}{rclcrcl} 
 0 & \vdash & &\qquad \qquad &
 & \vdash & 1 \\
 x,\,y & \vdash & x+y &\quad & x+y & \vdash & x ,\, y \\
x,\,y & \vdash & xy & & -x y & \vdash & x ,\, y 
 \end{array}
\quad \right\}
\end{equation}
\end{definition}
%--- end-definition------------------------------------

Cette manière simple et constructive de définir le treillis réel remonte à \cite{CC00}, qui s"inspire de \cite[Section V-4.11]{Joh1986}. 

\smallskip \hum{Cette section devrait être développée, notamment en introduisant des faisceaux convenables sur le spectre réel d'un anneau commutatif.}

\newpage \thispagestyle{empty}

%%%%%%%%%%%%%%%%%%%%%%%%%%%%%%%%%%%%%%%%%%%%%%%%%%%%%%%%%%%%%%%%%%%%
%%%%%%%%%%%%%%%%%%%%%%%%%%%%%%%%%%%%%%%%%%%%%%%%%%%%%%%%%%%%%%%%%%%%
\chapter{Corps réels clos \textsl{non} discrets} \label{chapreelclos}
\Today
\minitoc

\section*{Introduction}
\addtocontents{toc}{\vskip0.8em}
\addcontentsline{toc}{section}{Introduction}
\rdb

La section \ref{subsecclotrlRR}
 explique comment introduire des racines carrées des \elts~\hbox{$\geq 0$} dans un corps ordonné \textsl{non} discret. Ceci à titre de mise en jambes pour introduire la notion plus \gnle de racines virtuelles.
Le cas des \codis a été traité dans \cite[section 3.2]{LR91}. La morale de l'affaire est que l'on n'a pas besoin de savoir si une racine carrée est déjà présente dans le corps ordonné pour l'introduire formellement sans risque de contradiction. On voit ici la supériorité du point de vue \cof sur le point de vue classique (qui usuellement utilise le principe du tiers exclu pour décider si la racine carrée convoitée est déjà présente ou non).

La section \ref{secCoVR} explique comment ajouter les fonctions racines virtuelles dans la théorie des \afrs. 
Les racines virtuelles ont été introduites dans \cite{GLM98}. Le but était de disposer, pour un \pol \unt réel, de fonctions continues des \coes qui recouvrent les racines réelles. Cela permettait notamment d'avoir une version \cov du \tho des valeurs intermédiaires dans laquelle aucun test de signe n'était utilisé. En fait, le travail analogue peut être réalisé sur n'importe quel \afr. 
D'où notamment les théories des \afrvrs et des \covrs \textsl{non} discrets,

 La section \ref{secArc} traite les \arcs et la section \ref{secCrc2} propose une \dfn pour les \crcs \textsl{non} discrets comme \arcs locaux. La théorie des anneaux réels clos est ici présentée sous une forme \elr, \peq, dans le style de \cite{Tre2007}. 

%%%%%%%%%%%%%%%%%%%%%%%%%%%%%%%%%%%%%%%%%%%%%%%%%%%%%%%%%%%%%%%%%%%%
%%%%%%%%%%%%%%%%%%%%%%%%%%%%%%%%%%%%%%%%%%%%%%%%%%%%%%%%%%%%%%%%%%%%
%:Subsection  Corps ordonné 2-clos (ou euclidien)
\section{Corps ordonné 2-clos (ou euclidien)} 
\label{subsecclotrlRR}

À titre de mise en jambes, voyons la question de l'introduction des racines carrées des \elts~\hbox{$\geq 0$}.
Nous sommes intéressés par la \rdy suivante qui dit que dans un \afr les \elts $\geq 0$ sont des carrés d'\elts $\geq 0$.

\Regles{
\Lab{sqr} $\vd \Exists z\;\;(x^+=z^2\vet z=z^+)$
}

Une autre manière d'affirmer la même chose est de postuler

\Regles{
\Lab{sqa} $\vd \Exists z\;\;(\abs x=z^2\vet z=\abs z)$
}

Bien que dans un \grl on ait $x^+\neq \abs x$ en \gnl, les conditions $x=x^+$ et $x=\abs x$ sont \eqves (elles signifient $x\geq 0$).

\vspace{-1em}
\Subsection{Le cas d'un \codi} 

Référence: \cite[section 3.2]{LR91}.

%d
%:     Definition{deficdireel}
\begin{definition} \label{deficdireel}
Un corps discret est dit \textsl{réel} si $\sum_{i=1}^nx_i^2=0$ implique $x_i=0$ ($1\leq i\leq n$). 
\end{definition}
%----------- fin definition -------------------------------- 

Un \codi est réel. En \clama, tout \cdi réel peut être ordonné. En \coma cela signifie simplement que si on ajoute formellement une structure d'ordre (avec les axiomes des \codis) à un \cdi, et si la théorie formelle prouve $1=0$, alors on avait déjà $1=0$ dans le \cdi lui-même. 

%:     Definition{deficdidc}
\begin{definition} \label{deficdidc}
Un \codi est dit \textsl{$d$-clos}  si tout polynôme $P$ de degré $\leq d$ qui change de signe entre $a$ et~$b$ possède une racine sur l'intervalle d'extrémités $a$ et $b$. Un \codi $2$-clos est aussi appelé un \textsl{\cdi euclidien}. 
\end{definition}
%----------- fin definition -------------------------------- 

La méthode usuelle de résolution des équations du second degré donne le lemme suivant. Mais cela ne fonctionne que pour un corps discret.

%:     lemma{lemmacodi2clos}
\begin{lemma} \label{lemmacodi2clos}  Soit $\gK$ un \codi.
\Propeq
\begin{enumerate}
\item Tout \elt $\geq 0$ est un carré.
\item $\gK$ est un \codi $2$-clos.
\end{enumerate}
 
\end{lemma}
%----------- fin lemma ----------------------------- 

%p
%:     Proposition{propcodi2clos}
\begin{proposition} \label{propcodi2clos} Soit $\gK$ un \cdi réel\footnote{On utilise l'axiome de réalité uniquement pour le cas particulier $\,\,x^2+y^2+z^2=0\vd x=y=z=0$.}. \Propeq
\begin{enumerate}
\item $\gK$ peut être muni d'un ordre (au sens des \codis) pour lequel tout \elt~$\geq 0$ est un carré. Auquel cas la relation d'ordre possible est unique.
\item Tout carré est une puissance $4$.
\end{enumerate}
 \end{proposition}
%----------- fin proposition ----------------------------- 
%
\begin{proof}
\textsl{1}  $\Rightarrow$ \textsl{2}. Soit $x$ égal à un carré $y^2$. Si $\gK$ est un \codi, on a $y\geq 0$ ou $y\leq 0$. Donc il existe $z$ tel que $y=z^2$ ou $-y=z^2$. Dans les deux cas $y^2=z^4$.

\smallskip \noindent \textsl{2} $\Rightarrow$ \textsl{1}. De $y^2=z^4$ on déduit que $y=z^2$ ou $-y=z^2$. Si $P$ est l'ensemble des carrés, on a donc $P\cup -P=\gK$.\\
Montrons ensuite que pour tous $x,y$, il existe $z$ tel que $x^2+y^2=z^2$.
En effet $x^2+y^2=\pm z^2$ et si $x^2+y^2+ z^2=0$ alors $x=y=z=0$ et  $x^2+y^2=z^2$. On a donc $P+P\subseteq P$. On a aussi $P\cap -P=\so 0$. En effet si $y^2=-y^2$ alors $y^2+y^2=0$ donc $y=0$. Enfin $P\cdot P\subseteq P$ car $x^2y^2=(xy)^2$.
Ainsi le corps peut être ordonné, de manière unique.
Et tout positif est un carré.
\end{proof}

Dans \cite{LR91} si $\gK$ est un \codi et $a\geq 0$, les auteurs introduisent formellement une racine carrée $\alpha\geq 0$ de $a$ et  démontrent que $\gK[\alpha]$ peut être muni d'une structure de \codi sans qu'on ait besoin de savoir si $\alpha\in\gK$. Autrement dit, on connait la structure de $\gK[\alpha]$ comme \codi, mais à priori on ne sait pas si $\gK[\alpha]$ est de dimension~$1$ ou~$2$ comme \Kev.

Cette construction \elr est la pierre de base pour construire la $2$-clôture d'un \codi. 

Les choses s'avèrent nettement plus difficiles pour un corps ordonné \textsl{non} discret.

\Subsection{Le cas d'un \afr} 

%: Remark{rem2closreduit}
\begin{remark} \label{rem2closreduit} 
Dans un \afr, en présence de nilpotents, deux \elts $z$ et $y \geq 0$ qui ont le même carré ne sont pas \ncrt égaux, mais si l'anneau est réduit, ils sont égaux, en vertu de la \rsim \Tsbf{Aonz3}. Dans la théorie \Sa{Afrnz} la règle \Tsbf{sqr} est donc une \rex à existence unique. Si on la skolémise  on obtient donc une théorie \esid.
\eoe
\end{remark}
%----------- fin remark ---------------------------------- 

On présente maintenant une version dans laquelle une racine carrée positive d'un \elt positif est donnée comme une loi unaire dans la \tdy qui étend \sa{Afr}

\smallskip 
\centerline{$
\Sqr\,:\gR\to\gR,\;x\mapsto \sqrt{ x^+}.$} 

\smallskip \noindent Cette loi unaire doit obéir à certaines règles naturelles.

\DeuxRegles{
\lAb{sqr$_=$} $\,\, y= 0 \vd \Sqr(x+y)=\Sqr(x)$\label{Axsqreg}
\Lab{sqr1} $ \vd \Sqr(x)= {\Sqr(x)}^+$
\Lab{sqr3} $ \vd \Sqr(x)^2=x^+$
}
{
\Lab{sqr0} $ \vd \Sqr(0)=0$
\Lab{sqr2} $ \vd \Sqr(x)=\Sqr(x^+)$
\Lab{sqr4} $ \vd \Sqr(x^+y^+)=\Sqr(x)\Sqr(y)$
}

Notons que $\Sqr(x)=0$ lorsque $x\leq 0$ et $\Sqr(x)=\sqrt x$ lorsque $x\geq 0$.

%d
%: Definition{defiAfr2c}
\begin{definition} \label{defiAfr2c}\label{defiAsr2c}~
\begin{itemize}
\item Nous notons \SA{Afr2c} la \tpe des \textsl{\afrs $2$-clos}: elle est obtenue à partir de \Sa{Afr} en ajoutant le symbole de fonction
unaire $\Sqr$ avec les six axiomes précédents:\index{anneau!fortement reticule@fortement réticulé!2-clos}
\sIgt{\AfRdc}{\cdot=0,\mathrel{;}\cdot+\cdot, \cdot\times\cdot,\cdot\vu\cdot,-\,\cdot,\Sqr(\cdot),0,1}. \label{NOTASigAfr2c}

\item Nous notons \SA{Asr2c} la \tdy des \textsl{\asrs 2-clos} obtenue à partir de \Sa{Asr} de la même manière que \Sa{Afr2c} a été obtenue à partir de \Sa{Afr}:\index{anneau!strictement réticulé!2-clos}
\sIgt{\AsRdc}{\cdot=0,\cdot\geq 0,\cdot>0\mathrel{;}\cdot+\cdot, \cdot\times\cdot,\cdot\vu\cdot,-\,\cdot,\Sqr(\cdot),0,1}.\label{NOTASigAsr2c} 

\item Nous notons \SA{Aftr2c} la \tdy des \textsl{\aftrs 2-clos} obtenue à partir de \Sa{Aftr} de la même manière que \Sa{Afr2c} a été obtenue à partir de \Sa{Afr}:\\
\sIgt{\AftRdc}{\cdot=0,\cdot\geq 0,\cdot>0\mathrel{;}\cdot+\cdot, \cdot\times\cdot,\cdot\vu\cdot,-\,\cdot,\Fr(\cdot),\Sqr(\cdot),0,1}.\label{NOTASigAftr2c}\index{anneau!fortement reel@fortement réel!2-clos}
\item Nous notons \SA{Co2c} la \tdy des \textsl{corps ordonnés 2-clos} (\textsl{non} discrets) obtenue à partir de \Sa{Co} de la même manière: même signature que \Sa{Aftr2c}.\label{NOTASigCo2c}\index{corps ordonné!2-clos} 
\end{itemize}
\end{definition}
%----------- fin definition -------------------------------- 

%:     Remark{remSqa}
\begin{remark} \label{remSqa} 
Nous avons choisi la fonction $\Sqr(x)$ qui vérifie  \fbox{$\Sqr(x)\geq 0$ et $\Sqr(x)^2=x^+$} parce qu'elle correspond à la deuxième racine virtuelle du \pol $Y^2-x$ dans le cas d'un \codi. Mais on pourrait aussi utiliser \fbox{$\Sqa(x):=\Sqr(\abs x)$} qui satisfait l'\egt
\fbox{$\Sqa(x^+)=\Sqr(x)$} et qui peut être caractérisée par \fbox{$\Sqa(x)\geq 0$ et $\Sqa(x)^2=\abs x$}.
On vérifie facilement que la fonction $\Sqa$ peut être introduite avec les axiomes suivants.

\DeuxRegles{
\lAb{sqa$_=$} $\,\, y= 0 \vd \Sqa(x+y)=\Sqa(x)$\label{Axsqaeg}
\Lab{sqa1} $ \vd \Sqa(x)= \abS{\Sqa(x)}$
\Lab{sqa3} $ \vd \Sqa(x)^2=\abs x$
}
{
\Lab{sqa0} $ \vd \Sqa(0)=0$
\Lab{sqa2} $ \vd \Sqa(x)=\Sqa(\abs x)$
\Lab{sqa4} $ \vd \Sqa(x y)=\Sqa(x)\Sqa(y)$
}

\end{remark}
%----------- fin remark ---------------------------------- 

%l
%: Lemma{lemAr2creduit}
\begin{lemma} \label{lemAr2creduit}
Un \afr $2$-clos est réduit. 
\end{lemma}
%----------- fin lemma ----------------------------------- 
%
\begin{proof}
Si $x^2= 0$, alors $0=\Sqr(x^2)=\Sqr(x)^2=x^+=x$. Comme $(-x)^2=0$, on a aussi $x^-=0$. Et $x=x^+-x^-$. 
\end{proof}

\Subsubsection{Quelques règles dérivées dans \sa{Afr2c}}

\DeuxRegles{
\Lab{Aonz0} $\,\, x^2+y^2+z^2=0 \vd  x=0$
\Lab{Afr21} $\,\, (x^2+y^2)^2= z^4 \vd  x^2+y^2 = z^2$
}
{
\Lab{AFR2} $ \vd \Exists y\; x^2=y^4$
}

Notons que la règle \tsbf{Aonz0} est valide dans \Sa{Aonz} (anneaux ordonnés strictement réduits).

\Subsubsection{Théories \esids à $\sa{Afr2c}$}

%l
%: Lemma{lemSqr}
\begin{lemma} \label{lemSqr}
Les cinq extensions suivantes de la théorie \sa{Afr} sont \esids.
\begin{enumerate}
\item La théorie \peq \sa{Afr2c}.
\item On ajoute le symbole de fonction $\Sqa$  et les 6 axiomes indiqués dans la remarque \ref{remSqa}.
\item On ajoute comme axiomes les \rdys \Tsbf{Anz} et \Tsbf{sqr}.
\item On ajoute comme axiomes les \rdys \Tsbf{Anz} et \Tsbf{sqa}.
\item On ajoute comme axiomes les \rdys \Tsbf{Anz} et \Tsbf{AFR2}.
\end{enumerate}
\end{lemma}
%----------- fin lemma -----------------------------------

Dans les points \textsl{3}, \textsl{4}, \textsl{5}, on ne change pas le langage. 
\begin{proof} On démontre que la théorie du point \textsl{3} est \esid à \sa{Afr2c}. Tout d'abord, dans la théorie \sa{Afrnz} la \rex \tsbf{sqr} est à existence unique en vertu de la remarque \ref{rem2closreduit}. Ensuite, on vérifie que la fonction $\Sqr$ obtenue en skolémisant l'axiome existentiel~\tsbf{sqr} satisfait les 6 axiomes voulus.
\\
Montrons aussi que la théorie du point \textsl{5} est \esid à celle du point \textsl{4}. Dans cette dernière on a $\Sqa(x)\geq 0$ et $\Sqa(x)^4= {\abs x} ^2=x^2$, donc l'axiome \tsbf{AFR2} est vérifié. Dans la théorie du point \textsl{5}, on a $x^2={\abs x} ^2=y^4$, donc par \Tsbf{Aonz3} $\abs x=y^2={\abs y}^2$. Ainsi $\abs y$ vaut pour le $z$ dans \Tsbf{sqa}.   
\\
Le reste est laissé \alec.
\end{proof}

Le lemme suivant peut être vu comme un \gnn au cas \textsl{non} discret du fait que sur un corps réel $2$-clos discret, la structure d'anneau commutatif détermine complètement la structure d'ordre.

%: Lemma{lemAr2cunique}
\begin{lemma} \label{lemAr2cunique}
Sur un anneau commutatif, s'il existe une structure d'\afr 2-clos, celle-ci est unique. Plus \gnlt, un morphisme d'anneaux entre deux \afrs 2-clos est un morphisme d'\afrs 2-clos. 
\end{lemma}
%----------- fin lemma ----------------------------------- 
\begin{proof} Soit $\varphi\colon \gA\to\gB$
un morphisme d'anneaux où $\gA$ et $\gB$ sont des \afrs 2-clos. Comme les $x\geq 0$
sont les carrés, la relation d'ordre est respectée. Maintenant dans un \afr 2-clos (ou plus \gnlt dans un \afr réduit) l'\elt $c=a\vu b$ est caractérisé par les \egts et inégalités 
$c\geq a$, \hbox{$c\geq b$} et \hbox{$(c-a)(c-b)=0$}
(lemme \ref{lemsupdansAonz}). On en déduit que le morphisme d'anneau est aussi un morphisme pour les lois $\vu$. Enfin, comme dans un \afr réduit,
deux \elts $\geq 0$ qui ont même carré sont égaux (remarque \ref{rem2closreduit}), la loi $\Sqr$ est elle aussi respectée par le morphisme d'anneaux.
\end{proof}

 Comme la théorie \Sa{Afr2c} est \peq, tout \afr $\gA$ engendre
librement un \afr 2-clos: sa 2-clôture $\AFRdC(\gA)$. 
Se pose alors la question: à quoi ressemble la 2-clôture d'un \afr?
Voici un premier \elt de réponse.\index{2-clôture!d'un \afr}\label{notaAFR2C} 

%l
%: Lemma{lemAfr2ccloture}
\begin{lemma} \label{lemAfr2ccloture}
Tout \afr réduit s'injecte dans 
\hbox{sa 2-clôture}. 
\end{lemma}
%----------- fin lemma ----------------------------------- 
%
\begin{proof}
La théorie \Sa{Afr2c} prouve les mêmes \ralgs que \Sa{Afrnz}: cela résulte du point \textsl{3} du \pstfref{Pst2}, car la fonction~$\Sqr$ ajoutée à la théorie \Sa{Crcd} donne une théorie \esid. On n'obtient donc pas de nouvelle \egt entre \elts de l'\afr de départ après avoir ajouté formellement les racines carrées des \elts $\geq 0$. 
\end{proof}
Cela généralise le fait qu'un \codi s'injecte dans \hbox{sa 2-clôture} (\cite{LR91,LR90}), laquelle est un \codi.
Pour le cas \textsl{non} discret se pose la question naturelle \ref{quest2cloture}.

\Subsubsection{Axiomes pour qu'un anneau commutatif soit un \afrdc}

Ce paragraphe éclaire le lemme \ref{lemAr2cunique}. Il généralise au cas d'un anneau commutatif ce qui était fait pour un \codi en vue de le rendre $2$-clos. Il suffisait d'imposer les axiomes \tsbf{Aonz0} (réalité) et \tsbf{AFR2} (tout carré est une puissance 4, voir la proposition \ref{propcodi2clos}).
Nous proposons le système d'axiomes suivants, qui ajouté à la théorie des anneaux commutatifs, donne une théorie \esid à \Sa{Afr2c}. Il faut prendre $x\geq 0$ comme une abréviation de $\Exists z \; x=z^2$, et $x\geq y$ comme une abréviation de $x-y\geq 0$.

\DeuxRegles{
\laB{Aonz0} $\,\, x^2+y^2+z^2=0 \vd  x=0$
\laB{Afr21} $\,\, (x^2+y^2)^2= z^4 \vd  x^2+y^2 = z^2$
\Lab{Afr23} $\,\,  z\geq x\vet z\geq -x\vet x^2=y^4\vd z\geq y^2 $
}
{
\laB{AFR2} $ \vd \Exists y\; x^2=y^4$
\Lab{Afr22} $\,\, x^2= y^4\vet (y^2+x)^2=z^4 \vd  y^2+x = z^2$
}

Voici quelques explications. L'axiome \tsbf{Afr21} admet comme cas particulier $u^4=v^4\vd u^2=v^2$. Cela permet de voir que dans \tsbf{AFR2} l'\elt $y^2$ est uniquement déterminé à partir de $x$ et peut donc être skolémisé sous le nom de $\abs x$. On voit que la définition pour $x\geq 0$ est \eqve à $x=\abs x$. Ensuite nous devons voir que cette fonction $\abs {\,\cdot\,}$ satisfait les axiomes que nous avons proposés pour la \dfn d'une structure de \grl sur un groupe abélien donné \paref{grl-abs}.
Nous devons introduire une constante $\frac 1 2$ pour avoir les axiomes \Tsbf{2div1} et \Tsbf{2div2}. Nous obtenons ainsi une structure de \grl sur le groupe additif de l'anneau\footnote{En fait $\gA$ est remplacé par $\gA/\sqrt{\gen{0}}$ et si $nx=0$ pour un entier $n\geq 1$, il faut aussi annuler $x$. Une situation plus confortable serait de supposer qu'on part d'une \QQlg réduite. Dans ce cas $\gA$ s'injecte dans l'\afrdc que l'on construit.}. Il suffit ensuite de vérifier la validité des axiomes \Tsbf{ao1}, \Tsbf{ao2} et \Tsbf{afr6b} (\paref{ao1}), ce qui est immédiat.

\smallskip \rem On aimerait bien pouvoir démontrer que les axiomes \tsbf{Afr22} et \tsbf{Afr23}, qui correspondent à \Tsbf{abs3} et \Tsbf{Abs1}, sont conséquences des autres axiomes.\eoe  

\Subsection{Le cas d'un \ndsof}

%l
%: Lemma{lemdefiAsr2c}
\begin{lemma} \label{lemdefiAsr2c} La théorie \Sa{Co2c} est \esid aux deux théories suivantes.
\begin{enumerate} 
\item Sur la signature
$(\cdot=0\mathrel{;}\cdot+\cdot, \cdot\times\cdot,\cdot\vu\cdot,-\,\cdot,\,\Fr(\cdot),\,\Sqr(\cdot),0,1)$ la théorie obtenue à partir de \Sa{Afr2c} en ajoutant le symbole de fonction $\Fr$ et les axiomes \Tsbf{fr1}, \Tsbf{fr2} et \Tsbf{AFRL}.
\item Sur la signature 
$\sIgt{\AsR}{\cdot=0,\cdot\geq 0,\cdot>0\mathrel{;}\cdot+\cdot, \cdot\times\cdot,\cdot\vu\cdot,-\,\cdot,0,1}
$ 
la théorie obtenue en ajoutant à \Sa{Asr} les règles \Tsbf{IV}, \Tsbf{OTF}, \Tsbf{FRAC}, \Tsbf{Anz} et \Tsbf{sqr}.
\end{enumerate}
\end{lemma}
%----------- fin lemma ----------------------------------- 

%%%%%%%%%%%%%%%%%%%%%%%%%%%%%%%%%%%%%%%%%%%%%%%%%%%%%%%%%%%%%%%%%%%%

\section{Racines virtuelles} \label{secCoVR}

%:Subsection Définition et premières propriétés
\Subsection{Définition et premières propriétés}\label{subsecrappelsvr}

Références: \cite{GLM98,CLLR06,AG2012,BG2011,Gal2013}.

L'idée qui a guidé l'introduction des racines virtuelles était de disposer, pour un \pol \unt réel, de fonctions continues des \coes qui recouvrent les racines réelles.
Quand une racine réelle s'évanouit dans le plan complexe, on peut la relayer par la racine de la dérivée qui coïncide avec la racine réelle double au moment où elle disparait.

Par exemple les racines carrées virtuelles d'un réel arbitraire $a$ (à savoir $-\Sqr(a)$ et $+\Sqr(a)$)
sont égales à $- \sqrt a$ et $ \sqrt a$ lorsque $a\geq 0$, et dans le cas contraire, elles sont nulles: c'est la valeur qu'elles avaient au moment de disparaitre (imaginer le \pol $X^2-a$ variant continument avec $a\in\gR$).

\smallskip Rappelons d'abord le \tho \agq des accroissements finis.

%: Lemma{lemaccrfinis}
\begin{lemma}[\tho \agq des accroissements finis] \emph{\cite{LR90,LR91}} \label{lemaccrfinis} \\
Il existe deux familles $(\lambda_{i,j})_{1\leq i\leq j\leq n}$ et 
$(r_{i,j})_{1\leq i\leq j\leq n}$ dans $\QQ\;\cap\; ]0,1[\,$ 
avec $\sum_{i=1}^nr_{i,n}=1$ pour tout $n\geq 1$ et telles 
que, pour tout \pol $f \in \QQ [X]$ de degré $ \leq n$, on ait 
dans~$\QQ[a,b]$:
 $$
 f(b)-f(a)=(b-a)\times\som_{i=1}^nr_{i,n}\cdot f'(a+\lambda_{i,n}(b- a)).
 $$ 
Le résultat s'applique à toute \QQlg $\gA$ (en particulier aux corps ordonnés \emph{non} discrets). 
Si $\gA$ est une \QQlg strictement ordonnée, cela montre qu'un \pol dont la dérivée \hbox{est $>0$} sur un intervalle $\,]a,b[\,$ est une fonction strictement croissante sur $[a,b]$. 
\end{lemma}
%--------- fin lemma ----------------------------------- 

%e
%: Example{exaACF}
\begin{example} \label{exaACF}
Par exemple pour les \pols de degré $\leq 4$ on a 
\[\ndsp 
\frac {f(1)-f(-1)}2= \frac1 3\,f'(-\frac2 3 )+
\frac1 6\,f'(-\frac1 3 )+
\frac1 6\,f'(\frac1 3 )+ 
\frac1 3\,f'(\frac2 3 ), 
\] 
et plus \gnlt avec $\Delta=b-a$
$$\ndsp
f(b)-f(a)=\Delta\cdot\big(
\frac1 3\,f'(a+\frac1 6 \Delta)+
\frac1 6\,f'(a+\frac1 3 \Delta)+
\frac1 6\,f'(a+\frac2 3 \Delta)+ 
\frac1 3\,f'(a+\frac5 6 \Delta)\big).
$$
\eoe
\end{example}
%--------- fin example ---------------------------------------------- 

%: lemBasicVirtualRoots
\begin{lemma}[légère variation sur {\cite[Proposition 1.2]{GLM98}}] \label{lemBasicVirtualRoots}~
\begin{enumerate}
\item Soit $\sigma=\pm1$ et $f\colon [a,b]\to\RR$ ($a\leq b \in \RR$)
une fonction continument dérivable telle que $\sigma\,f'(x)>0$ sur $\;]a,b[\,$. Alors $f$
atteint son minimum en valeur absolue en un unique $x\in[a,b]$. 
%:2025 $\rR (a,b,f,\sigma)$ remplace $\rR (a,b,f)$
Nous notons $\rR (a,b,f,\sigma)$ ce réel $x$. \\
On a $(x-a)(x-b)f(x)=0$ et $x$ est l'unique réel vérifiant le système d'inégalités suivant. 
%-----------------begin item------------------

\DeuxRegles{
\labu $a \leq x \leq b$
\labu $ \sigma\, (x - a) f(a) \leq 0$
\labu $ \sigma\, (b-x) f(b) \geq 0$
}
{
\labu $ \sigma\, (x-a) f(x)\leq 0 $
\labu $ \sigma\, (b-x) f(x) \geq 0$
}

\item %[2.] 
Soit $\sigma=\pm1$ et $f\colon [a,+\infty[\ \to\RR$ une fonction continument dérivable telle que \hbox{$\sigma\,f'(x)>0$} sur $\;]a,+\infty[\,$. On suppose qu'il y a un $b>a$ tel que~\hbox{$\sigma\,f(b)>0$}. Alors $f$ atteint son minimum en valeur absolue en un unique $x\in[a,+\infty[$. Nous notons $\rR (a,+\infty,f,\sigma)$ ce réel.\\ 
On a $(x-a)f(x)=0$ et $x$ est l'unique réel vérifiant le système d'inégalités suivant: 
%-----------------begin item------------------

\DeuxRegles{
\labu $a \leq x $
\labu $\sigma\,(x - a) f(a) \leq 0$
}
{
\labu $ \sigma\,(x - a) f(x) \leq 0$
\labu $ \sigma\, f(x) \geq 0$
}

\item %[3.] 
Un énoncé analogue au précédent, laissé \alec, pour l'intervalle $\;]-\infty,a]$. 

\item %[4.]%
Ce lemme est également valable pour un corps réel clos discret
$\gR$
si $f$ est une \fsagc continument dérivable sur un intervalle $[a,b]$.
\end{enumerate}

%-----------------end item------------------
\end{lemma}

%r
%: Remark{remlemBasicVirtualRoots}
\begin{remark} \label{remlemBasicVirtualRoots} 
 1) Dans l'article \cite{GLM98}, lorsque $f$ est un \pol unitaire de degré $d$, l'hypothèse est formulée sous la forme 
%:  $\sigma f'(x)\geq 0$
 $\sigma f'(x)\geq 0$ sur $[a,b]$, ce qui fait que l'ensemble des paramètres ($a$,~$b$ et les \coes de $f$) satisfaisant l'hypothèse est un \fsa  $F\subset \RR^{d+2}$. On démontre alors que la fonction $\rR (a,b,f,\sigma)$
est \sagc sur ce fermé $F$. 
%:2025 ajout ci dessous
Ce sera le cas ici pour \ref{prdfVirtualRoots} et \ref{thVirtualRoots}.
\\
2) On note que les points \textsl{2} et \textsl{3} sont décalqués du point \textsl{1}.\\
3) 
%:2025  on doit pouvoir
On doit pouvoir expliciter un module de continuité uniforme pour $\rR$ si l'on donne certaines précisions sur la \fsagc $f'$ (précisions disponibles lorsque que $f$ est un \pol unitaire). 
\eoe
 \end{remark}
%----------- fin remark ---------------------------------- 

\smallskip À partir de ce lemme on obtient la construction des \textsl{racines virtuelles} pour un \pol \unt de degré $d$: d'une part elles \gui{couvrent} toutes les racines réelles, d'autre part elles varient continument en fonction des \coes du \pol.

Pour un \pol $f$ \unt de degré $d$, \textsl{nous notons $f^{[k]}$ la dérivée $k$-ème de~$f$ divisée par son \coe dominant} $(0\leq k< d)$: c'est un \pol \unt de degré $d-k$.

%: PropDef{prdfVirtualRoots}
\begin{propdef} \label{prdfVirtualRoots} Soit $\gR$ un \crcd ou le corps~$\RR$.
Pour tout \pol \unt 

%\vspace{-.2em}
\snic{f(X) = X^{d} - ( a_{d-1} X^{d-1} + \cdots +a_1X+ a_0) \quad (d\geq 1)}

\noindent  
%:2025 ajout  il est légitime de 
il est légitime de définir les fonctions \emph{racines virtuelles de $f$} 

\snic{\rho_{d,j}(f) = \rho_{d,j}(a_{d-1}, \ldots, a_0)}

\noindent pour $1 \leq j \leq d$ par \recu sur $d$ comme suit: (on abrège ci-dessous $\rho_{\delta,j}(f^{[d-\delta]})$ en $\rho_{\delta,j}$) 
%-----------------begin item------------------

\vspace{.1em}
\begin{enumerate}
\item $\rho_{1,1}(X - a) = \rho_{1,1}(a) := a$;
\item $\rho_{d,j} := \rR (\rho_{d-1,j-1}, \rho_{d-1,j}, f,(-1)^{d-j})$ 
pour $1 \leq j \leq d$ (lorsque $d\geq 2$).
\end{enumerate}

\vspace{.2em}
\noindent NB. On a posé par convention $\rho_{\delta,0}=(-1)^d\infty$ et $\rho_{\delta,\delta+1}=+\infty$ pour tout $\delta\geq 1$, et la fonction $\rR$ est celle définie dans le lemme \ref{lemBasicVirtualRoots}.

%:2025 ajout d'une phrase
\smallskip \noindent En d'autres termes, la récurrence  est correcte parce que lorsque $\rho_{d-1,j-1}< \rho_{d-1,j}$ on a $(-1)^{d-j}\,f'(x)>0$ sur l'intervalle ouvert correspondant. 
\end{propdef}
%--------- fin propdef ------------------------------

Cette proposition se démontre simultanément avec les points
\textsl{\ref{ivrvariation}} et \textsl{\ref{ivrsigne}} du \tho qui suit, en utilisant les
lemmes \ref{lemaccrfinis} et \ref{lemBasicVirtualRoots}.

%: Theorem{thVirtualRoots}
\begin{theorem}[quelques \prts des racines virtuelles] \label{thVirtualRoots} \emph{\cite{GLM98,CLLR06}}\\
Soit $\gR$ un \crcd ou le corps $\RR$. On note $\xi$ un \elt arbitraire du corps. On considère un \pol unitaire $f$ donné de degré $d$.

\begin{enumerate}
%:2025  \delta  plutôt que  k dans tout le théorème et la suite
\item Les $\frac{d(d+1)}2$ \elts $\rho_{\delta,j}(f^{[d-\delta]})$, pour $1\leq j\leq \delta\leq d$ introduits par \recu dans la proposition \ref{prdfVirtualRoots} sont caractérisés par un \sys d'inégalités larges.
\item Chaque fonction $\rho_{d,j}:\gR^d\to\gR$ est uniformément continue sur toute boule\footnote{$\mathrm{B}_{d,M}:=\sotQ{(a_{d-1},\dots,a_0)\,}{\,\sum_ia_i^2\leq M}$, ($M>0$). La continuité peut être donnée sous forme complètement explicite à la {\L}ojasiewicz.} $\mathrm{B}_{d,M}$.
\item 
On note \fbox{$\wi f=\prod_{j=1}^d(X-\rho_{d,j})$} et \fbox{$f\sta=\prod_{j=0}^{d-1} f^{[j]}$}, 
 \fbox{$\rho_{\delta,j}=\rho_{\delta,j}(f^{[d-\delta]})$} , et l'on utilise les conventions \fbox{$\rho_{\delta,0}=(-1)^\delta\infty$} et \fbox{$\rho_{\delta,\delta+1}=+\infty$}.
\begin{enumerate}
\item On a $\rho_{d,1}\leq \rho_{d-1,1}\leq \dots\leq \rho_{d-1,j}\leq \rho_{d,j+1}\leq \rho_{d-1,j+1}\leq \dots \leq \rho_{d-1,d-1}\leq \rho_{d,d}$.
\item Si $d\geq 2$ et $f=X^d-a$, alors $\rho_{d,d}=\sqrt[d]{a^+}$, 
$\rho_{d,j}=0$ pour $1<j<d$, $\rho_{d,1}+\rho_{d,d}=0$ si $d$ pair et $\rho_{d,1}+\rho_{d,d}=\sqrt[d]{a}$ si $d$ impair.
%:2025  Vi et Vu qaund nécessaire
\item \label{ivrsupinf} Si $f=\prod_{i=1}^d(X-\xi_i)$ pour des $\xi_i\in\gR$, alors $\wi f=f$. 
En conséquence $\rho_{d,1}=\Vi_{i}\,\xi_i$, 
$\rho_{d,d}=\Vu_{\!i}\,\xi_i$ et 
$\rho_{d,k}=\Vi_{J\subseteq \lrb{1..d}, \#J=k}(\Vu_{\!i\in J}\,\xi_i)$.
\item \label{ivrvariation} 
\begin{itemize}
\item Si $\rho_{d-1,j}<\rho_{d-1,j+1}$, $(0\leq j\leq  d-1)$, alors $f$ est \stm monotone sur l'intervalle correspondant, croissante si $d-j$ impair, décroissante sinon
%:2025 footnote améliorée
\footnote{Ceci implique que dans le système d'inégalités larges qui définit les $\rho_{\delta,j}$ pour  $1\leq j\leq \delta\leq d$, le signe $\sigma$ dans le lemme \ref{lemBasicVirtualRoots}  peut être donné directement, comme dans l'exemple qui suit le \tho, ce qui simplifie un peu les choses: les $\sigma$ \gui{disparaissent}.}.
\item Pour $0\leq j\leq d-1$, on a $(-1)^{d-j}\big(f(\rho_{d-1,j+1})-f(\rho_{d-1,j})\big)\leq 0$.
\end{itemize}
\item \label{ivrsigne} 
 Si $\rho_{d,j}<\xi<\rho_{d,j+1}$, $(0\leq j\leq d)$, alors $(-1)^{d-j}f(\xi)>0$.
\item Les zéros de $f$ sont des zéros de $\wi f$, avec une multiplicité supérieure ou égale dans $\wi f$. Plus \prmt: 
%i
\begin{itemize}\itemsep.2em
\item Si $f(\xi)=0$, alors $\wi f(\xi)=0$;
\item Si $\abS{\wi f(\xi)} > 0$, alors $\abs{f(\xi)} > 0$;
\item Si $f^{[j]}(\xi)=0$ pour $j\in\lrbk$, alors 
${\wi f}^{[j]}(\xi)=0$ pour $j\in\lrbk$;
\item Si $f^{[j]}(\xi)=0$ pour $j\in\lrbk$ et $\abS{{\wi f}^{[k+1]}(\xi)}> 0$, alors $\abs{f^{[k+1]}(\xi)}> 0$.
\end{itemize}
En outre si les multiplicités sont connues, la différence des multiplicités pour $f$ et $\wi f$ est paire (par exemple, un $\rho_{d,j}$ non zéro de $f$ est de multiplicité paire dans $\wi f$). 
\item Les zéros réels de $f\sta$ sont exactement les $\rho_{\delta,j}$ ($1\leq j\leq \delta\leq d$). Plus \prmt:
\begin{itemize}
\item chaque $\rho_{\delta,j}$ est un zéro de $f\sta$,
\item si tous les $|\xi-\rho_{\delta,j}|$ sont $>0$, alors $|f\sta(\xi)|>0$, 
\item le \pol $\wi f$ divise $(f\sta)^d$. 
\end{itemize}
\item \label{ivrBudan} \emph{(Compte de Budan Fourier)} 
Soit $a\in \gR$ tel que les $\abs{f^{[\delta]}(a)} > 0$ pour $0\leq \delta\leq d$, et soit $r$
le nombre de changements de signes dans la suite des $f^{[\delta]}(a)$ $(\delta=d,\dots,0)$, $(r\in\lrb{0..d})$. \\
 Alors $\rho_{d,d-r}<a<\rho_{d,d-r+1}$. 
\item \label{ivrTVI} \emph{(\Tho de la valeur intermédiaire)}\\ 
Si $a<b$ et $f(a)f(b)<0$, on a
\smash{$\prod_{j=1}^{d}f(\mu_j)=0$, où \fbox{$\mu_j=a\vu(b\vi\rho_{d,j})$}.} \\
Cas particuliers.
\vspace{.1em}
\begin{itemize}
%
%\item 
%
\item Si $d$ est impair, alors $\prod_{j=1}^{d}f(\rho_{d,j})=0$.
\item Si $0\leq k<\ell\leq d$ et $f(\rho_{d-1,k})f(\rho_{d-1,\ell})<0$, alors
$\prod_{j=k}^{\ell-1}f(\rho_{d,j})=0$.
\vspace{.2em}
\item Si, selon le point \ref{ivrBudan} %précédent, 
on a $\rho_{d,k}<a<\rho_{d,k+1}<b<\rho_{d,k+2}$, alors $f(\rho_{d,k+1})=0$. 
\end{itemize}
\item \label{ivrExtrema} \emph{(\Tho des valeurs extrema)} Le \pol \unt $f$ \gui{atteint sa borne supérieure et sa borne inférieure sur tout intervalle fermé borné} au sens précis suivant: si $a<b$, on~a
\vspace{-.9em} 
\[ 
\begin{array}{rcl} 
\sup_{\xi\in[a,b]}f(\xi) & = & f(a)\vu f(b) \vu \sup_{j=1}^{d-1}f(\nu_j) \quad \hbox{où}\;\; \fbox{$\nu_j=a\vu(b\vi\rho_{d-1,j})$}\,, \\[.3em] 
\inf_{\xi\in[a,b]}f(\xi) & = & f(a)\vi f(b) 
\vi \inf_{j=1}^{d-1}f(\nu_j) \,. 
 \end{array}
\]

\vspace{-.6em} 
Si $f$ est de signe strict constant $\varepsilon=\pm1$ sur $[a,b]$, on a $\inf_{\xi\in[a,b]}\big(\varepsilon\,f(\xi)\big)>0$.
%$\varepsilon f(\xi)\geq \alpha>0$ sur~$[a,b]$ pour $\alpha=\inf_{\xi\in[a,b]}\varepsilon{f(\xi)}
%$. 
%
\item \label{ivrMinAbs} \emph{(\Tho du minimum en valeur absolue)}\\ Si $a<b$, on~a

\vspace{-.2em}
\snic{\inf_{\xi\in[a,b]}\abs{f(\xi)}=
\abs{f(a)} \vi \abs{f(b)} \vi \inf_{j=1}^{d}\abs{f(\mu_j)}
.\qquad\qquad \phantom{a}}

\noindent En outre, si le second membre est $>0$, alors $f$ est de signe constant sur $[a,b]$.

\item \emph{(Une borne)} Si $f(x)=x^d+\sum_{k=0}^{d-1}a_kx^k$
on a $\abs{\rho_{d,j}}\leq \sup_{k=0}^{d}(1+\abs{a_k})$ ($1\leq j\leq d$). 
\item \emph{(Changement de variable)} Soit $f(x)=x^d+\sum_{k=0}^{d-1}a_kx^k$ et
$g(x)=x^d+\sum_{k=0}^{d-1}c^{d-k}a_kx^k$ (formellement $g(x)=c^df(x/c)$).\begin{itemize}
\item Si $c\geq 0$, on a $\rho_{d,j}(g)=c\rho_{d,j}(f)$ ($1\leq j\leq d$).
\item Si $c\leq 0$, on a $\rho_{d,j}(g)=c\rho_{d,d+1-j}(f)$ ($1\leq j\leq d$). 
\item Dans tous les cas, $\prod_{1\leq j\leq d}(x-c\rho_{d,j}(f))=\prod_{1\leq j\leq d}(x-\rho_{d,j}(g))$.
\end{itemize}

\end{enumerate}
\end{enumerate} 
\end{theorem}
%--------- fin theorem ----------------------------------- 

\hum{Si \ncr on peut aussi donner le \cdv par translation.}

\hum{Il faut expliciter le module de continuité uniforme dans le point 2}
%: Example{exavr}
\begin{example} \label{exavr} Nous explicitons ici les inégalités
évoquées dans le point~\textsl{1} du \thref{thVirtualRoots} conduisant à $\rho_{4,3}(f)$ 
pour un \pol $f(X)=X^4-(a_3X^3+a_2X^2+a_1X+a_0)$, écrites ici sous forme de \reds sans hypothèses.
On reprend les conventions du point \textsl{3} du \thref{thVirtualRoots}. 
Ainsi, on pose 
$\rho_{1,1}=\frac {a_3} 4$,
$\rho_{2,j}=\rho_{2,j}(\frac {a_3}2, \frac {a_2}6)$,
$\rho_{3,j}=\rho_{3,j}(\frac {3a_3}4, \frac {a_2}2, \frac {a_1}4)$,
$\rho_{4,j}=\rho_{4,j}(a_3,{a_2},{a_1},{a_0})$. 
 On donne les inégalités qui caractérisent $\rho_{1,1}$, $\rho_{2,2}$, $\rho_{3,2}$ et $\rho_{4,3}$. On notera que dans la \dfn des \ravs,
le signe $\sigma$ devant $x-a$ ou $b-x$ dans le 
%:2025  point \textsl{1} du    supprimé 
lemme \ref{lemBasicVirtualRoots} est connu
%:2025 \ref{ivrvariation}  plutôt que  3d
en raison du point~\textsl{\ref{ivrvariation}} du \thref{thVirtualRoots}, ceci explique pourquoi ce signe n'apparaît pas dans les inégalités ci-dessous.
 
\Regles{\lab{vr$_{1,1}$} $\vd \rho_{1,1} = \frac {a_3} 4$}

\vspace{-.8em}
\DeuxRegles{
\lab{vr$_{2,1,0} $} $ \vd \rho_{2,1}\leq \rho_{1,1}
\phantom{(x^2-a^2)} $ 
\lab{vr$_{2,1,1}$} $ \vd (\rho_{2,1} - \rho_{1,1}) \, f^{[2]}(\rho_{1,1}) \leq 0 $ 
}
{
\lab{vr$_{2,1,2}$} $\vd (\rho_{2,1} - \rho_{1,1})\, f^{[2]}(\rho_{2,1}) \geq 0$
\lab{vr$_{2,1,3}$} $\vd f^{[2]}(\rho_{2,1}) \geq 0$
}

\vspace{-.8em}
\DeuxRegles{
\lab{vr$_{2,2,0} $} $ \vd \rho_{1,1}\leq \rho_{2,2}
\phantom{(x^2-a^2)} $ 
\lab{vr$_{2,2,1}$} $ \vd (\rho_{2,2} - \rho_{1,1}) \, f^{[2]}(\rho_{1,1}) \geq 0 $ 
}
{
\lab{vr$_{2,2,2}$} $\vd (\rho_{2,2} - \rho_{1,1})\, f^{[2]}(\rho_{2,2}) \leq 0$
\lab{vr$_{2,2,3}$} $\vd f^{[2]}(\rho_{2,2}) \geq 0$
}

\vspace{-.8em}
\DeuxRegles{ 
\lab{vr$_{3,3,0}$} $\vd \rho_{2,2}\leq \rho_{3,3}
\phantom{f^{[1]}(\rho_{3,3}}$
\lab{vr$_{3,3,1}$} $\vd (\rho_{3,3} - \rho_{2,2})\, f^{[1]}(\rho_{1,1})\geq 0$
}
{
\lab{vr$_{3,3,2}$} $\vd (\rho_{3,3} - \rho_{2,2})\, f^{[1]}(\rho_{3,3}) \leq 0$
\lab{vr$_{3,3,3}$} $\vd f^{[1]}(\rho_{3,3}) \geq 0$
}

\vspace{-.8em}
\DeuxRegles{ 
\lab{vr$_{3,2,0}$} $\vd \rho_{2,1}\leq \rho_{3,2}\leq \rho_{2,2} 
\phantom{f^{[1]}(\rho_{3,2})}$
\lab{vr$_{3,2,1}$} $\vd (\rho_{3,2} - \rho_{2,1})\, f^{[1]}(\rho_{2,1})\, \geq 0$
\lab{vr$_{3,2,2}$} $\vd (\rho_{3,2} - \rho_{2,2})\, f^{[1]}(\rho_{2,2})\, \geq 0$
}
{
\lab{vr$_{3,2,3}$} $\vd (\rho_{3,2} - \rho_{2,1})\, f^{[1]}(\rho_{3,2})\, \geq 0$
\lab{vr$_{3,2,4}$} $\vd (\rho_{3,2} - \rho_{2,2})\, f^{[1]}(\rho_{3,2})\, \geq 0$
}

\vspace{-.8em}
\DeuxRegles{
\lab{vr$_{4,3,0}$} $\vd \rho_{3,2}\leq \rho_{4,3}\leq \rho_{3,3}
\phantom{f(\rho_{4,3})}$
\lab{vr$_{4,3,1}$} $\vd (\rho_{4,3} - \rho_{3,2})\, f(\rho_{3,2})\, \geq 0$
\lab{vr$_{4,3,2}$} $\vd (\rho_{4,3} - \rho_{3,3})\, f(\rho_{3,3})\, \geq 0$
}
{
\lab{vr$_{4,3,3}$} $\vd (\rho_{4,3} - \rho_{3,2})\, 
f(\rho_{4,3})\geq 0$
\lab{vr$_{4,3,4}$} $\vd (\rho_{4,3} - \rho_{3,3})\, f(\rho_{4,3})\,\geq 0$
}
\eoe
\end{example}
%--------- fin example ---------------------------------------------

\vspace{-1em}

% subsubsection{Un résultat à la Pierce-Birkhoff}
\Subsection{Un résultat à la Pierce-Birkhoff}
%: Subsection{Un résultat à la Pierce-Birkhoff}

On appelle \textsl{fonction polyracine} une fonction $\gR^m\to \gR$
qui peut s'écrire sous la forme $\rho_{d,j}(f_1, \dots, f_d)$ pour des entiers $1\leq j\leq d$ et des \pols $f_j\in\gR[\xm]$.%
\index{fonction polyracine}

\smallskip Le \tho à la Pierce-Birhoff suivant mérite d'être signalé.
Cela ressemble à un Nusllstellensatz: il exprime qu'il y a une raison purement \agq au fait qu'une fonction
soit \sagc et entière sur l'anneau des \pols. 

%: Theorem{thPBpolyroots}
\begin{theorem} \label{thPBpolyroots} \emph{(\cite[Theorem 6.4]{GLM98})}
Soit $\gR$ un \crcd et soit $g\colon \gR^m\to\gR$ une \fsagc
 entière sur l'anneau $\Rxm$ (vu comme un anneau de fonctions). Alors $g$ est une combinaison par $\vu$, $\vi$ et $\,+\,$ de fonctions polyracines $\gR^m\to\gR$. Plus \prmt,
si $g(\ux)$ annule le \pol $Y$-unitaire $P(Y,\ux)$ de degré $d$, elle s'exprime comme sup-inf combinaison de fonctions de la forme

\vspace{-.8em}
\begin{equation} \label {eqthPBpolyroots}
\rho_{d,j}(P)+\sqrt[r]{ R_\ell^+ \cdot\bigg(1+{\som_{i=1}^nx_i^2}
\bigg)^s}
%= \rho_{d,j}(P)+\rho_{2r,2r}{\big(R_\ell^2\cdot(1+{\norm\ux}^2)^{2s}\big)}
\end{equation}

\vspace{-.4em}
\noindent pour des $R_\ell\in\Rxm$ (le deuxième terme dans la somme \pref{eqthPBpolyroots} est aussi une fonction polyracine, voir le point {3b} du \thref{thVirtualRoots}). 
\end{theorem}
%--------- fin theorem -----------------------------------

\rem
Lorsque la fonction $g$ est \polle par morceaux, elle annule un \pol \unt~$P(Y)=\prod_{i=1}^d(Y-f_i)$ pour des $f_i\in\Rxm$. Dans l'expression obtenue en~\pref{eqthPBpolyroots} pour~$g$, c'est l'inégalité de {\L}ojasiewicz qui est responsable de l'extraction de racine $r$-ème dans la formule. Quant aux $\rho_{d,j}(P)$ ce sont des sup-inf combinaisons des $f_i$ (point \textsl{\ref{ivrsupinf}} de \ref{thVirtualRoots}).
\eoe

%%%%%%%%%%%%%%%%%%%%%%%%%%%%%%%%%%%%%%%%%%%%%%%%%%%%%%%%%%%%%%%%%%%% 
\Subsection{\Afrs avec \ravs}\label{subsecafrvr}
%: Subsection{\Afrs avec \ravs}

%: Example{exavr2}
\begin{example} \label{exavr2} 
Nous reprenons l'exemple \ref{exatotordnonreduit} de la \QQlg totalement ordonnée $\QQ[\alpha]$, \hbox{avec $\alpha>0$} et $\alpha^6=0$.
Nous allons voir que les contraintes imposées pour~$\rho_{2,2}(f)$, lorsque~\hbox{$f=X^2-a^2$} et~\hbox{$a\geq 0$}, n'impliquent pas \ncrt \hbox{que $\rho_{2,2}=a$}. Les contraintes sont les suivantes \hbox{pour $x=\rho_{2,2}$}
 (notez \hbox{que $\rho_{1,1} = 0$}):

\DeuxRegles{
\lab{vr$_{2,2,0}$} $\vd 0\leq x
\phantom{(x^2-a^2)}$
\lab{vr$_{2,2,1}$} $\vd - x \, a^2 \leq 0$
}
{
\lab{vr$_{2,2,2}$} $\vd x \, (x^2-a^2) \leq 0$
\lab{vr$_{2,2,3}$} $\vd (x^2-a^2) \geq 0$
}

\smallskip\noindent 
Si nous prenons $a=\alpha$, tous les $x\geq 0$ tels que $x^2=a^2$ conviennent, et donc tous les $\alpha+y\alpha^5$ pour $y\in\QQ[\alpha]$ sont solutions. Si nous prenons $a^2=0$ les contraintes équivalent à \gui{$x\geq 0$ et $x^3\leq 0$} et tout \elt de l'intervalle $[0,\alpha^2]$ est solution, y compris $\zeta=\alpha^2$ alors que $\zeta^2> 0$.
\eoe
\end{example}
%--------- fin example ---------------------------------------------- 

%: Lemma{lem2afrvrreduit}
\begin{lemma} \label{lem2afrvrreduit}
Sur un \afr réduit, le \sys d'inégalités que satisfont les \ravs $\rho_{k,j}$ ($1\leq j\leq k\leq d$) pour un \pol unitaire de degré $d$ donné, si elles existent, définit ces \elts sans ambiguité. 
\end{lemma}
%--------- fin lemma ----------------------------------- 
%
\begin{proof}
L'unicité en question s'exprime au moyen de \ralgs. Or la théorie \Sa{Afrnz}
prouve les mêmes \ralgs que la théorie des corps réels clos discrets avec sup (\pstfref{Pst2}). Dans cette dernière théorie, l'unicité est assurée (point \textsl{1} du \thref{thVirtualRoots}).
\end{proof}

L'exemple \ref{exavr2} et les lemmes \ref{lem2afrvrreduit} et \ref{lemafrvrreduit} justifient la \dfn qui suit. 
%: Definition{defiAfrvr}
\begin{definition} \label{defiAfrvr}\index{anneau!fortement reticule@fortement réticulé!avec racines virtuelles} ~
\begin{enumerate}
\item 
 La théorie purement équationnelle \SA{Afrrv} des \textsl{\afrs avec racines virtuelles}
est obtenue comme suit à partir de la \tpe \Sa{Afr}.

\vspace{.1em} \begin{itemize}
\item Pour $1\leq j\leq d$ dans $\NN$, on ajoute un symbole de fonction $\rho_{d,j}$ d'arité $d$;
\item on ajoute comme axiomes les inégalités décrites dans le point \textsl{1} du \thref{thVirtualRoots};
\item on ajoute la règle \tsbf{vrsup} suivante 

\regles{\Lab{vrsup} $\vd \rho_{2,2}(a+b,-ab)=a\vu b$.}
\end{itemize} 
La signature est donc la suivante: \sIgt{\AfRrv}{\cdot=0 \mathrel{;}\cdot+\cdot, \cdot\times\cdot,\cdot\vu\cdot,-\,\cdot,(\rho_{d,j})_{1\leq j\leq d},0,1}.\label{NOTASigAfrrv}
\item 
On définit de la même manière la \talg \SA{Asrrv} des \textsl{\asrs avec racines virtuelles} à partir de la \talg \Sa{Asr}:\\
\sIgt{\AsRrv}{\cdot=0,\cdot\geq 0,\cdot>0 \mathrel{;}\cdot+\cdot, \cdot\times\cdot,\cdot\vu\cdot,-\,\cdot,(\rho_{d,j})_{1\leq j\leq d},0,1}.\label{NOTASigAsrrv}
\end{enumerate}

\end{definition}
%--------- fin definition --------------------------------

%: Lemma{lemafrvrreduit}
\begin{lemma} \label{lemafrvrreduit}
Un \afrvr est réduit. 
\end{lemma}
%--------- fin lemma ----------------------------------- 
%
\begin{proof}
Vu la règle \Tsbf{vrsup}, si $a^2=0$, on a avec $b=-a$: 
\[0=\rho_{2,2}(0,0)=\rho_{2,2}(a+b,-ab)=\abs a.
\]
\end{proof}

\vspace{-1em}
%: paragraph{Variante intègre}
\paragraph{Variante intègre}

%: Definition{defiAfrvrbis}
\begin{definition} \label{defiAfrvrbis} \textsl{(Anneau fortement réticulé avec racines virtuelles, variante intègre)}
\\ La \talg \SA{Aitorv} des \textsl{anneaux intègres totalement ordonnés avec racines virtuelles}
est obtenue à partir de la \talg \Sa{Aito} des anneaux intègres totalement ordonnés en ajoutant les \ravs de la même manière que la théorie \sa{Afrrv} est obtenue à partir de la théorie~\sa{Afr} dans le \dfn \ref{defiAfrvr}. 
\end{definition}
%--------- fin definition -------------------------------- 

Notez qu'on n'a pas besoin de mettre la règle \tsbf{vrsup} dans les axiomes.

%l
%: Lemma{lemAitorv}
\begin{lemma} \label{lemAitorv}
Un anneau intègre totalement ordonné avec racines virtuelles est \icl et son corps de fractions est réel clos discret. Réciproquement un anneau intègre \icl dont le corps de fractions est réel clos discret est un anneau intègre totalement ordonné avec racines virtuelles.
\end{lemma}
%----------- fin lemma ----------------------------------- 
%
\begin{proof}
On note $\gA$ l'anneau intègre et $\gK$ son corps de fractions, qui est discret.\\
\textsl{Implication directe}. Un \pol unitaire $f\in\AX$ vérifie \RCFn\ en raison du point \textsl{3i} de \ref{thVirtualRoots} et du fait que $\gK$ est discret. Pour un \pol arbitraire de $\KX$ on utilise le \cdv dans le point \textsl{3m} pour se ramener à un \pol unitaire de $\AX$. Donc $\gK$ est réel clos discret. Enfin $\gA$ est \icl en raison du point \textsl{3f}.\\
\textsl{Implication réciproque}. L'ordre sur $\gK$ induit un ordre total sur $\gA$. Il faut montrer que pour un \pol unitaire $f\in\AX$ les $\rho_{d,j}(f)$ sont dans $\gA$. Or ce sont des zéros de $f\sta$, \pol unitaire de $\AX$, ils sont dans $\gK$, et $\gA$ est \icl, donc ils sont dans $\gA$. Ainsi les fonctions $\rho_{d,j}(a_0,\dots,a_n)$ définies de $\gK^n$ vers $\gK$ se restreignent en des fonctions $\gA^n\to\gA$. 
\end{proof}
%

%: paragraph{Anneaux de \fsagces}
\paragraph{Anneaux de \fsagces}

\smallskip Le \thref{thPBpolyroots} (pour les \crcs discrets) légitime la \dfn suivante.

%: Definota{defiSacem}
\begin{definota} \label{defiSacem}
Soit $\gR$ un \afrvr (cas particuliers: un \hyperref[defiCorv]{\covr} ou un \aRc). Les familles $\Sace_m(\gR)$ ($m\in\NN$) de \textsl{\fsagces} sont définies comme les familles de fonctions de $\gR^m$ dans $\gR$ stables par composition, contenant les fonctions \pols (à \coes dans~$\gR$) et les fonctions racines virtuelles.
En d'autres mots, un \elt de $\Sace_m(\gR)$ est une fonction de~$\gR^m$ dans $\gR$ définie par un terme du langage de $\Sa{Afrrv}(\gR)$ avec les $m$ variables $\xm$ (certaines peuvent être absentes). 
\end{definota}
%--------- fin definota -------------------------------- 

%%%%%%%%%%%%%%%%%%%%%%%%%%%%%%%%%%%%%%%%%%%%%%%%%%%%%%%%%%%%%%%%%%%%
\paragraph{Anneaux de Pierce-Birkhoff}
%: paragraph{Anneaux de Pierce-Birkhoff}

%: Definota{defiSacembis}
\begin{definotas} \label{defiSacembis}
Soit $\gA$ un anneau, ou plus \gnlt une \sad 
d'\afr.
\begin{enumerate}
\item L'anneau $\AFRNZ(\gA)$ est l'\textsl{\afr réduit engendré par $\gA$}. 
\item L'anneau $\AFRRV(\gA)$ est l'\textsl{\afr avec \ravs engendré par~$\gA$}.
\item L'anneau $\PPM(\gA)$ est défini comme le sous-\afr de $\AFRRV(\gA)$ formé par les \elts~$x$ qui annulent un \pol $\prod_{i=1}^k(X-a_i)$ pour des $a_i\in\gA$. 
\item Un anneau $\gA$ est appelé un \textsl{anneau de Pierce-Birkhoff}
lorsque le morphisme naturel $\AFRNZ(\gA)\to\PPM(\gA)$ est un \iso.
\end{enumerate} 
\end{definotas}
%--------- fin definition --------------------------------

Voir la question \ref{questpbring}.

%%%%%%%%%%%%%%%%%%%%%%%%%%%%%%%%%%%%%%%%%%%%%%%%%%%%%%%%%%%%%%%%%%%%
\section{Anneaux réels clos}\label{secArc}

%%%%%%%%%%%%%%%%%%%%%%%%%%%%%%%%%%%%%%%%%%%%%%%%%%%%%%%%%%%%%%%%%%%%
\Subsection{Définition constructive et variantes}\label{subsecArcconst}
%: Subsection{Définition constructive}

%: Definition{defiArc}
\begin{definition}[anneaux réels clos] \label{defiArc}\index{anneau!réel clos} 
 La \tpe \SA{Arc} des \textsl{anneaux réels clos},
est obtenue en ajoutant à la théorie \Sa{Afrrv} le symbole de fonction $\mathrm{Fr}$ et les axiomes \Tsbf{fr1} et \Tsbf{fr2}.
\\
La signature est donc la suivante:
\Sigt{\ArC}{\cdot=0\mathrel{;}\cdot+\cdot, \cdot\times\cdot,\cdot\vu\cdot,-\,\cdot,(\rho_{d,j})_{1\leq j\leq d},\mathrm{Fr}(\cdot,\cdot),0,1}\label{NOTASigArc}
\end{definition}
%--------- fin definition -------------------------------- 

\vspace{-1em}

%: Lemma{lemArcunique}
\begin{lemma} \label{lemArcunique}
Sur un anneau commutatif, s'il existe une structure d'\arc, celle-ci est unique. Plus \gnlt, un morphisme d'anneaux entre deux \arcs est un morphisme d'\arcs. 
\end{lemma}
%----------- fin lemma ----------------------------------- 
%
\begin{proof}
Résulte des lemmes \ref{lemAr2cunique} et \ref{lem2afrvrreduit}.
\end{proof}

 \CAdre{.9}{Dans la suite, lorsque nous ne précisons pas le contraire, un \gui{\arc} désigne toujours un anneau défini dans \ref{defiArc}.}

 %l
 %: Lemma{lemArcbis}
 \begin{lemma}[variantes pour \sa{Arc}] \label{lemArcbis} ~
\begin{enumerate}
\item La théorie \sa{Arc} peut aussi être obtenue à partir de la \talg \Sa{Aftr} en ajoutant les \ravs de la même manière que la théorie \sa{Afrrv} est obtenue à partir de la théorie~\sa{Afr} dans la \dfn \ref{defiAfrvr}.
En outre, vu le lemme \ref{lemafrvrreduit}, l'axiome \Tsbf{Anz} de la théorie \sa{Aftr} peut être omis. Un \arc peut donc être vu comme un \emph{\aftr avec racines virtuelles}.
\item La théorie \sa{Arc} est \esid à la théorie \Sa{Asrrv} des \asrs avec racines virtuelles à laquelle on ajoute le symbole de fonction $\mathrm{Fr}$ et les axiomes \Tsbf{fr1} et \Tsbf{fr2}.
 NB: le prédicat $x>0$ doit être ajouté à \sa{Arc} comme une abréviation de \gui{$x$ \hbox{est $\geq 0$} et \iv}.
\end{enumerate}
\end{lemma}
 %----------- fin lemma ----------------------------------- 
%
\begin{proof} Le point \textsl{1} est clair. On en déduit le point \textsl{2}
en rappelant le lemme \ref{lemArftr0}.
\end{proof}
%

%: paragraph{Anneaux de \fsagcs}
\paragraph{Fonctions \sagcs}

\smallskip 
On reprend maintenant la \dfn \ref{defiFSAGC2} (légitimée par le \thref{thParamcontFsagc0}) en l'étendant aux \arcs.
On note que tout \arc contient une copie conforme de $\RRa$.

On suppose aussi que l'on a démontré le \thref{thArc3} et ses \corls.
%d
%: Definota{defiFSAGC2+}
\begin{definota} \label{defiFSAGC2+} 
Soit $\gR$ un 
\arc et 
soit une fonction $f\colon \gR^n\to \gR$.\index{fonction semialgébrique continue!sur un anneau réel clos}
\begin{enumerate}
\item (Cas \elr) La fonction $f$ est dite \textsl{\sagc} (de façon \elr) s'il existe une \fsagc $g\colon \RRa^{n}\to\RRa$ exprimée par un terme $t(\xn)$ de \sa{Arc} et si $f$ coïncide avec la fonction définie par ce terme.
\item (Cas \gnl) La fonction $f$ est \textsl{\sagc} s'il existe un entier $r\geq 0$, des \elts $y_1,\dots,y_r\in \gR$ et une fonction $h\colon \gR^{r+n}\to \gR$ qui relève du cas \elr précédent tels que 
$$
\forall \xn\in\gR\;\; f(\xn)=h(\yr,\xn).
$$
\end{enumerate}
On note $\Sac_n(\gR)$ l'anneau de ces fonctions (c'est un \arc pour la relation d'ordre naturelle).\label{Sacn} Le \thref{thParamcontFsagc0} montre que pour un \crcd $\gR$ on retrouve la \dfn usuelle des \fsagcs.

\end{definota}
%----------- fin definota -------------------------------- 
%\hum{Cela semblerait plus naturel de se limiter aux \arcs?}

%Tout terme du langage $\sa{Arc}(\gR)$ dont les variables sont prises parmi $\xn$ définit un \elt de $\Sac_n(\gR)$. Pour une réciproque voir le \thref{thArc3}.

Pour la comparaison de $\Sac_n(\gR)$ avec $\Sace_n(\gR)$ voir la question \ref{Qu-polrootsafrvr}.

%%%%%%%%%%%%%%%%%%%%%%%%%%%%%%%%%%%%%%%%%%%%%%%%%%%%%%%%%%%%%%%%%%%%
%: paragraph{Un exemple}
\paragraph{Un exemple}
%p
%: Proposition{propfsagccovrsup}
\begin{proposition} \label{propfsagccovrsup}
Soit $\gR$ un \afrvr et soit \hbox{$f\colon \gR^n\times \gR^p\to \gR$} une \fsagc. La fonction 
$$
g\colon \gR^p\to\gR,\,\ux \mapsto \sup\nolimits_{\uz\in [0,1]^n}f(\uz,\ux)
$$ 
est bien définie et \sagc. 
\end{proposition}
%----------- fin proposition ----------------------------- 
%
\begin{proof}
Vu la \dfn ad hoc des anneaux $\Sac_{m}(\gR)$ on est immédiatement ramené au cas \hbox{où $\gR=\RRa$.}
\end{proof}

Voir aussi les questions \ref{questArc2}.
%%%%%%%%%%%%%%%%%%%%%%%%%%%%%%%%%%%%%%%%%%%%%%%%%%%%%%%%%%%%%%%%%%%% 
%:Subsection Corps ordonnés avec racines virtuelles
\Subsection{Corps ordonnés avec racines virtuelles}\label{subsecCo0rv}

%: Definition{defiCorv}
\begin{definition} \label{defiCorv} (Comparer avec la \dfn \ref{defiAfrvr}, voir le lemme \ref{lem2afrvrreduit}).
\begin{enumerate}
\item La \tdy \SA{Corv} des \textsl{corps ordonnés avec racines virtuelles} est obtenue comme suit à partir de la \tdy \Sa{Co} 
des \ndsofs.\index{corps ordonné!avec racines virtuelles}

\vspace{.1em} \begin{itemize}
\item Pour $1\leq j\leq d$ dans $\NN$, on ajoute un symbole de fonction $\rho_{d,j}$ d'arité $d$;
\item on ajoute comme axiomes les inégalités décrites dans le point \textsl{1} du \thref{thVirtualRoots}.

\end{itemize}
La signature est donc la suivante:
\Sigt{\CoRv}{\cdot=0,\cdot\geq 0,\cdot>0\mathrel{;}\cdot+\cdot, \cdot\times\cdot,\cdot\vu\cdot,-\,\cdot,(\rho_{d,j})_{1\leq j\leq d},\mathrm{Fr}(\cdot,\cdot),0,1}\label{NOTASigCorv}

\vspace{-1em}
\item La \tdy \SA{Co--rv} est obtenue de la même manière à partir
de la théorie \Sa{Co--}.
\item La \tdy \SA{Codrv} est obtenue de la même manière à partir
de la théorie \Sa{Cod}.
\end{enumerate}
\end{definition}
%----------- fin definition -------------------------------- 

%
\rem La théorie \Sa{Codrv} est \esid à la théorie 
obtenue en ajoutant à \Sa{Co0rv} l'axiome \gui{de tiers exclu} \Edinq\  (voir la remarque \ref{remCo0}).

%:Subsection Positivstellensatz formel
\Subsection{Positivstellensatz formel}\label{subsecPst3}

%: pstf {Pst3}
\begin{pstf}[Positivstellensatz formel, 3] \label{Pst3} ~ \index{Positivstellensatz!formel!pour les corps ordonnés, 3} 
\begin{enumerate}
\item Les théories \Sa{Codrv}, \Sa{Crcd} et \Sa{Crcdsup} sont \esids.
\item 
Les \tdys suivantes prouvent les mêmes \ralgs
(écrites dans le langage de \Sa{Afrrv}).
\begin{enumerate}
\item La théorie \Sa{Afrrv} des \afrs avec \ravs.
\item La théorie \Sa{Arc} des \arcs.
\item  La théorie \Sa{Codrv} des \codis avec \ravs.
\end{enumerate}
\item 
Les \tdys suivantes prouvent les mêmes \ralgs
(écrites dans le langage de \Sa{Asrrv}).
\begin{enumerate}
\item La théorie \Sa{Asrrv} des \asrs avec \ravs.
%
%\item La théorie \Sa{Co0rv}.
%
\item La théorie \Sa{Corv} des \covrs.
\item La théorie \Sa{Codrv} des \codis avec \ravs.
\end{enumerate}
\item \label{i4Pst3} Le \thref{thVirtualRoots}
est entièrement valable pour la théorie \sa{Asrrv} (donc aussi pour \sa{Corv}). 
Il en va de même pour la \tpe \sa{Afrrv} (donc aussi pour \sa{Arc}) si les points qui utilisent le prédicat $\cdot> 0$ sont supprimés ou convenablement reformulés \hbox{avec $\cdot\geq 0$}. 
\end{enumerate}
\end{pstf}
%%%%%%%%%%%%%%%%%%%%%%%%%%%%%%%%%%%%%%%%%%%%%%%%%%%%%%%%%%%%%%%%%%%% 
%
\begin{proof} Le premier point est clair. Les points \textsl{2} et \textsl{3} sont donc
des variantes du \pstfref{Pst2} en tenant compte du \thref{thEseqMemesfaits}. 
\\
Enfin, pour le point \textsl{4}, il résulte des points précédents car 
les affirmations des points \textsl{2} et \textsl{3} du \thref{thVirtualRoots} peuvent s'écrire sous forme de \ralgs.
\end{proof}

En \clama vu le \tho général de représentation \ref{thcolsimralg}, les Positivstllensätze formels énoncés jusqu'ici donnent les résultats suivants,
qui peuvent être vus en \clama comme caractérisant les \tdys considérées. 
%c
%: Corollary{corPst3}
\begin{corollaryc} \label{corPst3}
Sur leurs signatures respectives, les objets suivants sont tous isomorphes à des sous-\sa{T}-structures de produits de \crcds (considérés avec le prédicat $x>0$ et les fonctions $\sup$, $\mathrm{Fr}$ et~$\rho_{d,j}$).
\begin{itemize}
\item Un \afr réduit (théorie \Sa{Afrnz}).
\item Un \asr réduit (théorie \Sa{Asrnz}). 
\item Un \aftr (théorie \Sa{Aftr}). 
\item Un \aftr local (théorie \Sa{Co}). 
\item Un \afr avec \ravs (théorie \Sa{Afrrv}). 
\item Un \arc (théorie \Sa{Arc}). 
\item Un \asr avec \ravs (théorie \Sa{Asrrv}). 
\item Un \covr (théorie \Sa{Corv})
%%
%\item 
\end{itemize}

\end{corollaryc}
%--------- fin corollary ------------------------------- 

%

%

\Subsection{Quotient, localisation et recollement d'\arcs}
%: Subsection{Quotient, localisation et recollement d'\arcs}

%l
%: Lemma{lem1Afrc}
\begin{lemma}[structure quotient] \label{lem1Afrc} 
Soit $\gA$ un anneau réel clos et $I$ un idéal radical. Alors $\gA/I$
est un anneau réel clos. 
\end{lemma}
%----------- fin lemma ----------------------------------- 
%
\begin{proof} Montrons d'abord que l'\id radical $I$ est solide. On doit d'abord montrer si $x\in I$, alors $\abs x \in I$: en effet ${\abs x}^2=x^2\in I$. Ensuite si $0\leq \abs x\leq y$ \hbox{avec $y\in I$}, on doit montrer que $x\in I$. Or, par \tsbf{FRAC}, $y$ divise $\abs x^2=x^2$, donc $x^2\in I$, puis $x\in I$.
Le quotient~$\gA/I$ est donc un \afr réduit. 
Ensuite il faut voir que les fonctions \ravs $\rho_{d,j}$
et la fonction \gui{fraction} $\mathrm{Fr}$ \gui{passent au quotient}. Or ces fonctions, quand elles existent dans un \afr réduit, sont bien définies par les systèmes d'inégalités qu'elles satisfont (lemme \ref{lem2afrvrreduit}). Comme ces inégalités passent au quotient, tout est dans les clous. 
\end{proof}
%

%l
%: Lemma{lem2Afrc}
\begin{lemma}[localisation] \label{lem2Afrc}
Soit $\gA$ un anneau réel clos et $S$ un monoïde. Alors $S^{-1}\gA$
est un anneau réel clos. 
\end{lemma}
%----------- fin lemma ----------------------------------- 
%
%
\begin{proof} On sait déjà que $S^{-1}\gA$ est muni d'une loi $\vu$ qui en fait un \afr réduit. Voyons ce qu'il advient pour les \ravs. Prenons l'exemple \ref{exavr} avec un \pol $f(X)=X^4-(\frac{a_3}{s}X^3+\frac{a_2}{s}X^2+\frac{a_1}{s}X+\frac{a_0}{s})$ avec les $a_i\in\gA$ et $s\in S^+$. 
On a dans $S^{-1}\gA$ avec $Y=sX$
$$s^4f(X)=Y^4-(a_3Y^3+sa_2Y^2+s^2a_1Y+s^3a_0)=Y^4-(b_3Y^3+b_2Y^2+b_1Y+b_0)=g(Y)$$
et donc aussi $s^4g(\frac Y s)=f(X)$.
On considère les \ravs $\rho_{i,j}$ pour le \pol unitaire~$g$ de $\gA[Y]$
avec $(i,j)$ égal à $(1,1)$, $(2,1)$, $(2,2)$, $(3,3)$, $(3,2)$, $(4,3)$, et enfin les $\rho'_{i,j}=\frac{\rho_{i,j}}s$. On constate que les inégalités de 
l'exemple \ref{exavr}, comme elles sont satisfaites pour les $\rho_{i,j}$ relativement au \pol $g$ dans $\gA[Y]$, sont ipso facto satisfaites pour les $\rho'_{i,j}$ relativement au \pol $f$ dans $S^{-1}\gA[X]$. Or ces inégalités caractérisent complètement les racines virtuelles lorsqu'elles existent 
(lemme~\ref{lem2afrvrreduit}).
\\
Un raisonnement analogue fonctionne pour la fonction $\mathrm{Fr}(\cdot,\cdot)$.
\end{proof}
%

%: Plgc {plcc.arc}
\begin{plcc}[recollement concret d'\arcs]\label{plcc.arc} ~
\\
Soient $S_1$, $\dots$, $S_n$ des \moco d'un anneau $\gA$. On note $\gA_i$ pour $\gA_{S_i}$, $\gA_{ij}$ pour~$\gA_{S_iS_j}$, et l'on suppose donné sur chaque $\gA_i$ une structure de type \Sa{Arc}. On suppose en outre que les images dans $\gA_{ij}$ des lois de $\gA_i$ et $\gA_j$ coïncident. Alors il existe une unique structure d'\arc sur~$\gA$ qui induit par \lon en chaque $S_i$ la structure définie sur~$\gA_i$. Cet \arc s'identifie à la limite projective du diagramme
$$\big(( \gA_i)_{i\in\lrbn},( \gA_{ij})_{i<j\in\lrbn};(\alpha_{ij})_{i\neq j\in\lrbn}\big),
$$
où les $\alpha_{ij}$ sont les morphismes de \lon, dans la catégorie des \arcs.

\end{plcc}
%--- end-plgc-----------------------------------------
%
\begin{proof}
On recopie, mutatis mutandis, la \demo du \plgc~\ref{plcc.frings}
pour les \afrs.
\end{proof}

\rems 1) Cela implique que la notion de schéma réel clos est bien définie.

\noindent 2) Les \plgcs analogues, avec la même \demo, sont valables pour les \afrs réduits, pour les \aftrs et pour les \afrvrs.
\eoe

%%%%%%%%%%%%%%%%%%%%%%%%%%%%%%%%%%%%%%%%%%%%%%%%%%%%%%%%%%%%%%%%%%%%
%:Subsection
\Subsection{Comparaison avec la \dfn en \clama}

Références: \cite{Sch84,PN2002,Tre2007}.

\smallskip \hum{Il faudra dire deux mots de \cite[Joyal \& Reyes, 1986]{JR86} \gui{Separably real closed local rings}.}

\smallskip La structure d'\textsl{anneau réel clos} est définie par N. Schwartz de manière très abstraite dans sa thèse~\hbox{\cite[Schwartz, 1984]{Schthesis}}. Une axiomatisation en tant que théorie cohérente a été proposée dans \cite[Prestel \& Schwartz, 2002]{PN2002} (voir la \dfn \ref{defiPrScArc} et la proposition \ref{propPrScArc}). 

Le but de N. Schwartz était de donner une description abstraite des anneaux de \fsagcs sur les \fsas pour un corps réel clos fixé $\gR$, et de définir des \gui{espaces réels clos} abstraits.

%:paragraph{Une axiomatique de Niels Schwartz
\paragraph{Une axiomatique de Niels Schwartz}~\label{subsecAxScw}

\smallskip Voici la \dfn des anneaux réels clos en \clama donnée dans \cite[Schwartz, 1986]{Sch84}.

%d
%: Definition{defiArcclama}
\begin{definitionc} \label{defiArcclama} 
 Un anneau \textsl{réel clos} est un anneau réduit $\gA$ vérifiant les \prts suivantes.
\begin{enumerate}
\item L'ensemble des carrés de $\gA$ est l'ensemble des \elts $\geq 0$ d'un ordre partiel qui fait de~$\gA$ un \aFr.
\item Si $0\leq a\leq b$, il existe $z$ tel que $zb=a^2$ (\textsl{axiome de convexité}).
\item Pour tout \idep $\fp$, l'anneau résiduel $\gA/\fp$ est \icl
et son corps de fractions est un \crc.
\end{enumerate}
\end{definitionc}
%----------- fin definition -------------------------------- 

%On montre plus loin (proposition \ref{propAfrc}) qu'en \clama les \dfns
%\ref{defiArc} (point 1) et~\ref{defiArcclama}
%sont \eqves. 

%l

%: Proposition{propAfrc}
\begin{propositionc} \label{propAfrc} 
 En \clama les \dfns \ref{defiArc} et \ref{defiArcclama}
sont \eqves. 
\end{propositionc}
%----------- fin proposition -----------------------------
%
\begin{proof}
\textsl{Direct.} Pour un anneau réel clos $\gA$ de la \dfn \ref{defiArcclama}, les fonctions \ravs sont bien définies, car on sait que toutes les \fsagcs définies sur $\RRa$ sont définies sur $\gA$. Il en va de même pour la fonction \gui{fraction} $\mathrm{Fr}$.

\smallskip \noindent \textsl{Réciproque.} Pour un anneau réel clos $\gA$ de la \dfn \ref{defiArc}, il faut montrer que le point \textsl{3} de la 
\dfn~\ref{defiArcclama} est satisfait. Soit un \idep $\fp$, l'anneau résiduel $\gA/\fp$ est \sdz donc totalement ordonné (lemme \ref{lemAfrsdz}). C'est aussi un \arc d'après le lemme \ref{lem1Afrc}. Le lemme \ref{lemAitorv} nous dit que $\gA/\fp$ est \icl
et que son corps de fractions est un \crc.% , \ref{lem2Afrc} et .
\end{proof}
%

%%%%%%%%%%%%%%%%%%%%%%%%%%%%%%%%%%%%%%%%%%%%%%%%%%%%%%%%%%%%%%%%%%%%
%:paragraph{L'axiomatique de Prestel-Schwartz
\paragraph{L'axiomatique de Prestel-Schwartz}~\label{subsecAxPrSc}

\smallskip \noindent L'article \cite[Prestel \& Schwartz, 2002]{PN2002} montre en \clama que la structure d'anneau réel clos de la \dfn \ref{defiArcclama} est décrite par une théorie cohérente. Les axiomes existentiels proposés par les auteurs pour remplacer le point \textsl{3} de \ref{defiArcclama} sont très sophistiqués et la \demo est \egmt un tour de force étonnant.

%d
%: Definition{defiPrScArc}
\begin{definition} \label{defiPrScArc} \textsl{(Anneaux réels clos à la Prestel-Schwartz)}\\
Un anneau commutatif est dit réel clos s'il satisfait les axiomes suivants. 
 \begin{itemize}
\item [i-iv)] L'anneau commutatif $\gA$ est réduit, les \elts $\geq 0$ sont exactement les carrés et la relation d'ordre fait de $\gA$ un \afr convexe (axiome \Tsbf{CVX}) 
\item [v)] Pour chaque $d\geq 1$, soit $f(x)=x^{2d+1}+\sum_{k=0}^{2d}a_kx^k$, $\delta=\mathrm{discr}_x(P)$ son discriminant, et $g(x)=x^{2d+1}+\sum_{k=0}^{2d}\delta^{2(2d+1-k)}a_kx^k$, on pose l'axiome $$\vd \Exists z \;g(z)=0$$
\item [vi)] Pour chaque $d\geq 1$ on pose l'axiome
$$\,\, x^d+\som_{k=0}^{d-1}a_kx^ky^{d-k}=0\vd \Exists z_1,\dots,z_d \;y(x-z_1y)\cdots(x-z_dy)=0$$
\end{itemize}
\end{definition}
%p
%: Proposition{propPrScArc}
\begin{proposition} \label{propPrScArc}
Dans la théorie \sa{Arc} les axiomes de la \dfn \ref{defiPrScArc} sont des \rdys valides.
\end{proposition}
%----------- fin proposition ----------------------------- 
%
\begin{proof}
Rappelons que le \thref{thVirtualRoots} est entièrement valable dans la théorie \sa{Arc} (point~\textsl{\ref{i4Pst3}} de \ref{Pst3}).
 
\noindent Voyons l'axiome \textsl{vi)}. Posons $f=x^d+\som_{k=0}^{d-1}a_kx^k$ et $g=x^d+\som_{k=0}^{d-1}a_kx^ky^{d-k}$. 
\\
Nous notons $\wi f=\prod_{1\leq j\leq d}(x-\rho_{d,j}(f))$ et $\wi g=\prod_{1\leq j\leq d}(x-\rho_{d,j}(g))$.\\ Le point \textsl{3m} du \thref{thVirtualRoots} donne l'\egt 
$$\prod\nolimits_{1\leq j\leq d}(x-y\rho_{d,j}(f))=\prod\nolimits_{1\leq j\leq d}(x-\rho_{d,j}(g)).
$$
Par ailleurs le point \textsl{3f} du \thref{thVirtualRoots} pour le \pol $g$ donne
$$\,\, x^d+\som_{k=0}^{d-1}a_kx^ky^{d-k}=0\vd (x-\rho_{d,1}(g))\cdots(x-\rho_{d,d}(g))=0.$$
On obtient donc dans la théorie \sa{Arc}, en prenant $z_k=\rho_{d,k}(f)$, la règle valide 
$$\,\, x^d+\som_{k=0}^{d-1}a_kx^ky^{d-k}=0\vd \;(x-z_1y)\cdots(x-z_dy)=0.$$

\noindent Voyons l'axiome \textsl{v)}. 
Nous allons démontrer que l'\elt $z$
dont l'existence est affirmée peut être choisi comme une \fsagc des paramètres $a_k$. Comme cette fonction est annulée par le \pol \unt $Q$
on conclut alors par le \tho \gui{à la Pierce-Birkhoff}~\ref{thPBpolyroots}. 
Vu le Positivstellensatz formel \ref{Pst3} (point \textsl{2}) il nous suffit 
en effet de démontrer la validité de la règle dans la théorie \sa{Codrv}. Considérons donc un \crc discret et dans l'espace des paramètres une composante connexe de l'ouvert $\so{\delta\neq 0}$. Sur cette composante connexe, les zéros réels de $f$ sont simples (il en existe au moins un car le degré est impair) et varient continument en fonction des paramètres. 
Ceux de $g$ sont simplement multipliés par $\delta^2$. On a donc sur cette composante connexe l'\elt $z$ cherché comme une \fsagc des paramètres en choisissant le plus grand des zéros réels. Quand on s'approche d'un bord d'une composante connexe, ces zéros tendent vers $0$ (ce sont des zéros de $f$ multipliés par $\delta^2$). Donc ces \fsagcs se recollent en une \fsagc globale. 
\end{proof}

On démontre en \clama l'implication réciproque: les axiomes de Prestel-Schwartz impliquent l'existence des racines virtuelles (car ce sont des \fsagcs). Cela donne l'\eqvc en \clama de notre axiomatique et de celle de Prestel-Schwartz. 

%:paragraph{L'axiomatique de Tressl
\paragraph{L'axiomatique de Marcus Tressl}~\label{subsecAxTr}

\smallskip \noindent 
Une version plus \elr, semblable à celle que nous proposons, pour la théorie des \arcs se trouve dans \cite[Tressl, 2007]{Tre2007} (voir aussi \cite{Sch84,Sch84b,Sch97,Sch89}). Dans cet article, un \arc est un \afr $\gR$ sur lequel sont données toutes les \fsagcs définies sur $\RRa$, et pour lequel toutes les identités \agqs liant ces fonctions sur~$\RRa$ sont satisfaites dans $\gR$. 

Une bonne analyse des articles de \clama sur les \arcs devrait nous permettre de comprendre pourquoi il suffit d'ajouter les fractions autorisées par la règle \Tsbf{FRAC} à un \afrvr pour pouvoir capturer toutes les 
\fsagcs ${\RRa}^m\to\RRa$. C'est l'objet des résultats concrets suivants, valides en \clama, mais dont nous désirons une \prco. Voir notamment la question \ref{questArc3}.

Rappelons que d'après le \tho de finitude
(\cite[theorem~2.7.1]{BCR}) le graphe $G_f=\sotq{(\ux,y)}{\ux\in\gR^n,y=f(\ux)}$ d'une \fsagc $f\colon \RRa^n\to\RRa$
est un \fsa de~$\RRa^{n+1}$ qui peut être décrit comme l'ensemble des zéros d'une \textsl{fonction \spo} $F:\RRa^{n+1}\to \RRa$, \textsl{i.e.} une fonction qui s'écrit sous la forme
\[
\sup\nolimits_i\,(\inf\nolimits_{1\leq j\leq k_i}\,p_{ij}) \quad \hbox{ où }p_{ij}\in\RRa[\xn,y]
\]
On sait décider si un tel graphe est celui d'une \fsagc. Le \tho suivant revient à dire que dans un tel cas on sait prouver l'existence du $y$ dépendant des $x_i$ directement dans la théorie \peq \sa{Arc}. 

%: theorem{thArc3}
\begin{theorem} \label{thArc3} 
Toute \fsagc $\RRa^n\to\RRa$ peut être définie par un terme de la théorie~\sa{Arc}.
\end{theorem}
%----------- fin theorem ----------------------------- 

%c
%: Corollary{corthArc3}
\begin{corollary} \label{corthArc3}
L'axiomatisation proposée en \ref{defiArc} pour les \arcs est \eqve à celle proposée par Tressl \cite{Tre2007}. 
\end{corollary}
%--------- fin corollary ------------------------------- 

%: Corollary{cor2thArc3}
\begin{corollary} \label{cor2thArc3}
Soit $\gR$ un \arc. Toute \fsagc $\gR^n\to\gR$ (\dfn \ref{defiFSAGC2+}) est définie par un terme de $\sa{Arc}(\gR)$ avec $n$ variables libres. 
\end{corollary}
%--------- fin corollary ------------------------------- 

Lorsqu'on ajoute l'axiome \Tsbf{OTF}, le \thref{thArc3} donne le \corl qui suit.
%: corollary{corCrc1Corv}
\begin{corollary} \label{corCrc1Corv} 
Les théories \Sa{Crc1} et \Sa{Corv} sont \esids.
 \end{corollary}
%----------- fin corollary --------------

Le \corl suivant est plus problématique, peut-on se ramener au cas $\gR=\RRa$?

%: corollary{cor3thArc3}
\begin{corollary} \label{cor3thArc3} 
On considère un \arc $\gR$, une \fsagc \hbox{$g\colon \gR^n \to \gR$} (un \elt de $\SaC_n(\gR)$) et un \pol $p\in\gR[\xn]$
 avec au moins un \coe inversible. On suppose que, sur l'ensemble $\sotq{\uxi\in\gR^n}{\abS{p(\uxi)}> 0}$, la fraction $f=g/p$ satisfait un module de continuité uniforme sur tout borné à la \L ojasiewicz (comme dans le lemme \ref{factfsagcLoja}). Alors il existe une unique \fsagc \hbox{$h\colon \gR^n \to \gR$} tel que $hp=g$.
\end{corollary}
%----------- fin corollary ----------------------------- 
%
\begin{proof}
\fbox{??, à préciser.}
\end{proof}
%

%%%%%%%%%%%%%%%%%%%%%%%%%%%%%%%%%%%%%%%%%%%%%%%%%%%%%%%%%%%%%%%%%%%%

\section{Corps réels clos \textsl{non} discrets}\label{secCrc2}

%:Subsection Définition \label{subsecdefiCrc2}
\Subsection{Une définition raisonnable}\label{subsecdefiCrc2}

%: Lemma{lemArcloc}
\begin{lemma} \label{lemArcloc}
Un \arc est local \ssi il satisfait la règle \Tsbf{AFRL}.
\end{lemma}
%----------- fin lemma ----------------------------------- 
%
\begin{proof} Voir le lemme \ref{lemAftrloc}.
\end{proof}

On propose maintenant pour la théorie des \crc \textsl{non} discrets
une formulation \esid à \Sa{Corv}, mais presque \peq. La règle \tsbf{AFRL} est préférée à la règle \Tsbf{OTF} car on n'introduit pas le prédicat $\cdot>0$ qui nous ferait sortir des \tpes pour \Sa{Arc}.

%d
%: Definition{defiCrc2}
\begin{definition} \label{defiCrc2} 
 La \textsl{\tdy des corps réels clos (\emph{non} discrets)}, notée \SA{Crc2},
est l’extension de la \tpe
 \sa{Arc} obtenue en ajoutant
la règle \tsbf{AFRL}. Autrement dit un corps réels clos \textsl{non} discret n'est rien d'autre qu'un \fbox{\arc local}.\index{corps réel clos!non discret@\textsl{non} discret} 
\end{definition}
%----------- fin definition -------------------------------- 

%p
%: Proposition{propCrc2}
\begin{proposition} \label{propCrc2}
Les théories \Sa{Corv}, \Sa{Crc1} et \Sa{Crc2} sont \esids (on doit définir le prédicat $>0$ que l'on ajoute à \sa{Crc2}). 
\end{proposition}
%----------- fin proposition ----------------------------- 
%
\begin{proof} 
Le \corl \ref{corCrc1Corv} donne la comparaison de \sa{Corv} et \sa{Crc1}. Le lemme \ref{lemArftr} nous dit qu'un corps ordonné \textsl{non} discret n'est autre qu'un \aftr local. En d'autres termes
la théorie \sa{Co} est \esid à la théorie \sa{Aftr} à laquelle on ajoute l'axiome \tsbf{AFRL}. Partons de \sa{Aftr}. Si l'on ajoute les \ravs puis \tsbf{AFRL} on passe à \Sa{Arc} (lemme \ref{lemArcbis} point~\textsl{1}) puis à \sa{Crc2}. Si l'on ajoute \tsbf{AFRL} puis les \ravs on passe à \sa{Co} puis à \sa{Corv}.
\end{proof}
%

%r
%: Remark{rem}
\begin{remarks} \label{remCrc}~

\noindent 1) Le corps $\RR$ est un modèle \cof de la théorie \sa{Crc2}.

\smallskip \noindent 2) La théorie \Sa{Crcd} des \crcs discrets est \esid à la théorie
obtenue en ajoutant à \sa{Crc2} l'axiome~\Edinq\  qui dit que l'\egt est décidable.

\smallskip \noindent 3) La théorie \sa{Crc2} n'est rien d'autre que la théorie des \arcs locaux. Cependant, il existe des \arcs locaux qui ne sont pas des corps au sens de Heyting. Considérons par exemple l'anneau $\gA$ des \fsagcs sur $\RRa$, et soit $\gB=S^{-1}\gA$ où $S$ est le \mo des fonctions~$f$ telles que $f(0)\neq 0$. C'est l'anneau des germes en $(0)$ des fonctions~\hbox{$f\in\gA$}. Un \elt $f\in\gA$ est $>0$ dans $\gB$ (resp. $\leq 0$ dans $\gB$) \ssi $f(0)>0$ dans $\RRa$ (resp.~\hbox{$f(x)\leq 0$} au voisinage de $0$). Ceci montre que 
\Tsbf{HOF} n'est pas satisfait dans $\gB$, car il ne suffit pas que $f(0)\leq 0$ pour que $f$ soit $\leq 0$ au voisinage de $0$. 
Notons que cet \arc local admet deux \ideps minimaux, 
avec pour localisés respectifs les germes de fonctions à droite (ou à gauche) de $0$. 

\smallskip \noindent 4) La théorie \Sa{Crc2} permet de démontrer
l'existence d'une racine carrée pour un nombre complexe de module 1.
On recouvre le cercle unité $\so{x^2+y^2=1}$ par les ouverts $\so{x>-1}$ et $\so{x<1}$, sur chacun desquels l'existence est assurée par une fonction continue. 
Cependant, cette existence ne peut pas être démontrée dans \Sa{Arc}, car dans cette théorie \peq, toute existence est certifiée par un terme, et tout terme définit une \fsagc.
\eoe
\end{remarks}
%----------- fin remark ---------------------------------- 

\Subsection{Clôture réelle d'un \afr réduit}
%: Subsection{Clôture réelle d'un \afr réduit}

Étant donné un \afr réduit $\gA$, on sait (\pstfref{Pst2}) que la théorie $\Sa{Crcdsup}(\gA)$ prouve les mêmes \ralgs que $\Sa{Afrnz}(\gA)$. Il en va de même pour toutes les théories intermédiaires, 
en particulier pour les théories \Sa{Afrrv} et \Sa{Arc}. 

Comme ces dernières sont des \tpes, l'\afr réduit $\gA$
engendre un \afrvr $\AFRRV(\gA)$ et un \arc $\ARC(\gA)$.\label{notaARC} 

Comme les théories \sa{Afrnz}, \sa{Afrrv} et \sa{Arc} prouvent les mêmes \ralgs, $\gA$ est une sous-structure (d'\afr) de $\AFRRV(\gA)$ qui est elle-même une sous-structure (d'\afrvr) de $\ARC(\gA)$. Autrement dit, l'ajout des symboles de \ravs et de fractions (avec leurs axiomes) ne change rien à $\gA$ en tant qu'\afr. 

Ces deux constructions de \gui{clôtures réelles} sont sans mystère, et uniques à \iso unique près.

On est dans la même situation que pour la construction de la clôture réelle d'un \codi (\cite{LR90,LR91}), mais ici le résultat parait complètement évident
alors qu'il nécessite un effort non négligeable dans les articles cités.
La raison principale de ce (tout petit) miracle est que nous nous appuyons ici sur une \prco du Positivstellensatz. La raison secondaire est que nous ne traitons ici que de \talgs (au lieu de \tdys).

%r
%: Remark{remcloturereelle}
\begin{remark} \label{remcloturereelle} 
Une construction de la clôture réelle d'un \codi $\gK$ peut \egmt être obtenue selon l'argument suivant. On commence par établir l'effondrement simultané de la théorie des \codis et de celle des \crcs discrets (comme dans \cite[Theorem 3.6]{CLR01}). C'est une variante du \Pst formel. Ensuite on évalue dynamiquement~$\gK$ comme un \crc discret. Cela force à introduire les zéros réels de tout \pol,
avec pour chacun d'eux un codage à la Thom (pour un \pol qui annule ce zéro). Comme aucune ambigüité n'est possible, la \sad construite est en fait une structure \agq usuelle de \crc. Cette construction est certes moins détaillée que celle expliquée dans \cite{LR91}, mais elle est essentiellement \eqve. Et l'on pourrait d'ailleurs sans doute, dans l'autre direction, déduire le Theorem 3.6 de \cite{CLR01} de la construction donnée dans \cite{LR91}. Ce en quoi \cite{Lom91} et \cite{CLR01} améliorent le résultat précédent, c'est d'une part que le \Pst formel est plus \gnl (Theorem 3.8 dans \cite{CLR01}), d'autre part, et surtout, que le \Pst concret y est démontré. 
\eoe
\end{remark}
%----------- fin remark ---------------------------------- 

\Subsection{Clôture réelle d'un corps ordonné \textsl{non} discret?}
%: Subsection{Clôture réelle d'un corps ordonné \textsl{non} discret}
 
Considérons un corps ordonné \textsl{non} discret, \textsl{i.e.} un modèle $\gK$ de la théorie \Sa{Co}.
On sait que \sa{Corv} prouve les mêmes \ralgs que \sa{Co}.

Notons $\gR$ la \sad $\sa{Corv}(\gK)$. 

Tous les termes clos de la \sad $\gR$
sont construits sur des \elts de $\gK$ au moyen des symboles de fonction donnés dans la signature (polynômes, racines virtuelles, fractions légitimes).

La structure dynamique $\gR$ est un candidat naturel pour être \und{la} structure \agq (usuelle) de type \Sa{Corv} engendrée par~$\gK$, s'il en existe une. 
Cependant, le \pb est que $\gR$ est une \sad de type \sa{Corv}, mais pas forcément un modèle de cette théorie, car cette \tdy est définie avec des axiomes non \agqs. 

On peut considérer tout d'abord la structure algébrique usuelle d'\arc $\ARC(\gK)$ qui s'identifie à la \sad $\sa{Arc}(\gK)$.
La question qui se pose est alors la suivante: l'axiome~\Tsbf{AFRL} 
est-il une règle valide dans $\ARC(\gK)$? Autrement dit, est-ce que $\ARC(\gK)$
est un modèle de~\sa{Corv}? Auquel cas on peut identifier $\gR$ (\sad) et $\ARC(\gK)$ (structure \agq usuelle).% 
\label{pagecloturevirtuelle}

La réponse n'est pas évidente (voir la question \ref{questCloturevirtuelle}). 

%:sibrouillon
\sibrouillon{
\smallskip {\bf Brouillon}. La tentative suivante n'aboutit pas. En effet on n'a pas montré que l'axiome \tsbf{IV} soit satisfait dans la structure $\gR_0$. Cela serait \ncr pour déduire que $\gR_0$ est un modèle \cof de \sa{Co--}.

Notons cependant que l'on a une réponse positive pour la variante suivante.
On considère un modèle $\gK_0$ de la théorie \Sa{Co--}.
On sait que \Sa{Co0rv} prouve les mêmes \ralgs que \sa{Co--}.
On note $\gR_0$ la structure \agq $\ASRRV(\gK_0)$. 
C'est une structure \agq usuelle et l'application naturelle $\gK_0\to\gR_0$ est injective.
Tous les \elts de $\gR_0$
sont construits sur des \elts de $\gK_0$ au moyen des symboles de fonction donnés dans la signature (polynômes, fonction $\sup$ et racines virtuelles).
La structure \agq $\gR_0$ est un candidat naturel pour être \und{la} structure \agq (usuelle) de type \sa{Co0rv} engendrée par~$\gK_0$, s'il en existe une. La question qui se pose est alors la suivante: l'axiome~\Tsbf{OTF} 
est-il une règle valide dans $\gR_0$? 

La réponse est positive pour la raison suivante. Si $x$ et $y$ sont des \elts de $\gR_0$ ce sont des \elts entiers sur $\gK_0$.
Si $x$ et $y$ sont de degré $\leq d$ sur $\gK_0$ et si $a_0,\dots,a_d$ sont des \elts distincts de $\gR_0$, au moins l'un des $x-a_i$ et au moins l'un des $y-a_j$ ont un signe strict. 
Si $x+y>0$, on considère des \elts $a_0,\dots,a_d$ dans $\gK_0$ sur l'intervalle ouvert $\,]0,\frac{x+y}2[\,$ (on peut en effet minorer $x+y$ par un \elt $>0$ de $\gK_0$). Si $x-a_i>0$, alors $x>0$. Si $y-a_j>0$, alors $y>0$. Si $x-a_i$ et $y-a_j<0$, alors $x+y<a_i+a_j$, ce qui est impossible. 

\noindent {\bf fin de Brouillon}} 
%: fin sibrouillon

%%%%%%%%%%%%%%%%%%%%%%%%%%%%%%%%%%%%%%%%%%%%%%%%%%%%%%%%%%%%%%%%%%%%
%%%%%%%%%%%%%%%%%%%%%%%%%%%%%%%%%%%%%%%%%%%%%%%%%%%%%%%%%%%%%%%%%%%%

\section{Un corps réel clos \textsl{non} discret non archimédien} \label{crcndna}

 \newcommand{\Pu}[1]{\RRa\big[\big[\vep^{1/#1}\big]\big]}

Nous décrivons dans cette section 
%:2025  ce que nous pensons être
ce que nous pensons être\footnote{La démonstration du \thref{propPuiseux} est incomplète.}
un exemple de corps réel clos de Heyting \textsl{non} discret non archimédien.

\smallskip Soit $\vep$ une \idtr. On a introduit dans la section \ref{condna} le corps ordonné \textsl{non} discret non archimédien $\gQ=\gZ[1/\vep]$ où $\gZ=\QQ[[\vep]]$ est l'anneau des séries formelles en $\vep$ à \coes rationnels avec $\vep$ un infiniment petit strictement positif.

En fait, les \coes des séries considérées auraient pu être pris dans n'importe quel \codi, en particulier dans le corps $\RRa$ des réels algébriques. Nous noterons $\gR_0=\RRa[[\vep]]$ l'analogue de $\gZ$ et $\gR=\gR_0[1/\vep]$ l'analogue de $\gQ$.

\smallskip Nous étendons maintenant ces constructions au corps $\gP$ des séries de Puiseux à \coes réels algébriques. 

On a tout d'abord les anneaux de séries $\gP_{0,d}=\Pu{d}$ pour les entiers $d\geq 1$, tous isomorphes à $\gR_0$, avec les morphismes d'inclusion 
$\gP_{0,d}\to \gP_{0,dd'}$. Cela forme un système inductif dont la limite $\gPo$ (les séries de Puiseux de valuation~\hbox{$\geq 0$}) peut être vue comme la réunion des $\gP_{0,d}$. 

Enfin les séries de Puiseux proprement dites forment l'anneau défini comme 
$\gP:=\gPo[1/\vep]$.

\smallskip On note $\gP_{j,d}=\sotq{\alpha\in\gP_{0,d}[1/\vep]}{ \alpha/\vep^{j/d}\in \gP_{0,d}}$.
On a $\gP=\bigcup_{j,d}\gP_{j,d}$.

Nous introduisons les notations qui généralisent~à~$\gP_{j,d}$ celles déjà données pour $\gR$. Ces notations sont cohérentes par rapport aux inclusions $\gP_{j,d}\subseteq \gP_{jd',dd'}$
%\begin{itemize}
%
%\item 

Soit $\alpha=\sum_{k=j}^\infty a_{k/d}\vep^{k/d}\in \gP_{j,d}\subseteq \gP_{0,d}[1/\vep]$ ($j,d\in\ZZ, d\geq 1$). On définit:
\begin{itemize}
\item $\rc_{\ell/d}(\alpha)=
\formule 
{0&\hbox{ si }\ell<j,\\ 
a_{\ell/d} &\hbox{ si } \ell\geq j.
}$
\item $\kappa_{\ell/d}(\alpha)=s_{\ell/d}\in\so{-1,0,1}$ est défini par \recu comme suit:\\
$s_{\ell/d}=\formule {\hbox{ si }\ell<j, & \hbox{alors }
0&
\\
\hbox{ si }\ell\geq j,& 
\formule{\hbox{si }s_{(\ell-1)/d}\neq 0,& \hbox{alors } s_{(\ell-1)/d},\\
\hbox{si }s_{(\ell-1)/d}= 0,& \hbox{alors }\hbox{signe de }a_{\ell/d}.} 
}$ 
\item $\alpha>0$ signifie $\exists \ell\geq j\; \kappa_{\ell/d}(\alpha)=1$. 
\item $\alpha\geq 0$ signifie $\forall \ell\geq j\; \kappa_{\ell/d}(\alpha)\geq 0$. 
\item $\rv(\alpha)>k/d$ signifie $\kappa_{k/d}(\alpha)=0$.
\item $\rv(\alpha)\leq k/d$ signifie $\kappa_{k/d}(\alpha)=\pm1$.
\item $\rv(\alpha)\geq k/d$ signifie $\kappa_{(k-1)/d}(\alpha)=0$.
\item $\rv(\alpha)=k/d$ signifie $\rv(\alpha)\geq k/d$ et $\rv(\alpha)\leq k/d$.
%
%\item 
%
\end{itemize}

%
%\item 
%%
%\item 
%\end{itemize}
 
\smallskip De l'étude précédente dans la section \ref{condna} qui a abouti à la proposition \ref{propQ[[T]]} pour l'anneau~$\gZ$ et au \thref{propQ((T))} pour l'anneau $\gQ$, on déduit les résultats analogues pour les anneaux $\gP_{0,d}$
puis pour $\gPo$, puis pour~$\gP$. 
%: Proposition{propP_0}
\begin{proposition} \label{propP_0}~
\begin{enumerate}
\item 
L'anneau $\gPo$ est un \asr réduit qui satisfait les \prts suivantes. 
\begin{itemize}
\item Il satisfait les règles \Tsbf{OTF}, \OTFx, \Tsbf{FRAC} et \Tsbf{Val2}. En particulier (lemme \ref{lemAfrnzFRAC}) la \fsagc~$\Fr$ satisfaisant les règles \Tsbf{Fr1} et \Tsbf{Fr2} est bien définie et le symbole de fonction correspondant peut être considéré comme faisant partie de la signature. 
\item C'est un \alrd hensélien. 
\item Son corps résiduel est isomorphe à $\RRa$. \\
On a $\gPo\eti=\sotq{\xi\in\gPo}{\kappa_0(\xi)=\pm1}$ et $\Rad(\gPo)=\sotq{\xi\in\gPo}{\kappa_0(\xi)=0}$.

\item Le groupe de valuation est isomorphe à $(\QQ,+,\geq)$ (la classe de $\vep$ correspond à l'\elt~$1$ de~$\QQ$).
\item Les \elts $\geq 0$ sont des carrés: l'anneau $\gPo$ est un \asr $2$-clos (théorie \Sa{Asr2c}).
\item Plus \gnlt, les \elts $\geq 0$ sont des puissances $k$-èmes d'\elts $\geq 0$. Comme il s'agit d'existence unique, on peut introduire les symboles de fonction correspondant dans la signature.
\item En outre, l'axiome de Heyting ordonné $\lnot(\xi>0)\Rightarrow\xi\leq 0$ est satisfait.
\end{itemize}
\item
L'anneau $\gP$ est un \asr réduit qui satisfait les \prts suivantes. 
\begin{itemize}
\item Un \elt est $>0$ \ssi il est $\geq 0$ et inversible.
\item Les règles \Tsbf{OTF}, \tsbf{OTF$\eti$}, \Tsbf{FRAC} et \Tsbf{IV} (à fortiori \Tsbf{Val2}) sont satisfaites. En particulier (lemme \ref{lemAfrnzFRAC}) la \fsagc~$\Fr$ satisfaisant les règles \tsbf{Fr1} et \tsbf{Fr2} est bien définie.
\item C'est un \alo avec $\Rad(\gP)=0$ (un corps de Heyting dans la terminologie de \cite{CACM} ou \cite{MRR}).
\item Les \elts $\geq 0$ sont des carrés d'\elts $\geq 0$: l'anneau $\gP$ est un \asr $2$-clos (théorie \Sa{Asr2c}).
\item Plus \gnlt, les \elts $\geq 0$ sont des puissances $k$-èmes d'\elts $\geq 0$. Comme il s'agit d'existence unique, on peut introduire les symboles de fonction correspondant dans la signature.
\item L'axiome de Heyting ordonné $\lnot(\xi>0)\Rightarrow\xi\leq 0$ est satisfait. 
\end{itemize}
\end{enumerate}

\end{proposition}
%----------- fin proposition ----------------------------- 
%
\begin{proof}
Seul le fait que les \elts $\geq 0$ sont des puissances $k$-èmes d'\elts $\geq 0$ est un point vraiment nouveau qui réclame une \demo. Elle est laissée \alec.
\end{proof}

Nous notons $\gPal$ la clôture intégrale de $\RRa(\vep)$ dans $\gP$: c'est l'anneau des séries de Puiseux qui sont entières sur le sous-\codi $\QQ(\vep)$.

Dans la suite nous utiliserons la notion de prolongement par continuité. Pour parler de prolongement par continuité, il faut définir la notion de suite convergente, et vérifier que les règles usuelles de passage à la limite fonctionnent pour cette notion. 
%d
%: Definition{deficonvinP}
\begin{definition}[suites convergentes dans $\gP$] \label{deficonvinP}
On dira que \textsl{la suite~$(\alpha_n)_{n\in\NN}$ converge vers $\alpha$ dans~$\gP$} s'il existe~$j$ et $d\in \ZZ$ avec $d\geq 1$ tels que 
\begin{itemize}
\item $\alpha$ et les $\alpha_n$ sont tous dans $\gP_{j,d}$,
\item $\lim_n\rv(\alpha_n-\alpha)=+\infty$, \cade:\\ 

\vspace{-1em} 

\centerline{$\forall k\geq j$ $\exists N$ $\forall m\geq N\, \forall \ell\in \lrb{j.. k}\; \rc_{\ell/d}(\alpha_m)=\rc_{\ell/d}(\alpha)$. } 
\end{itemize}
On écrira alors $\alpha=\lim_n\alpha_n$. 
 
\end{definition}
%----------- fin definition -------------------------------- 

On établit facilement les \prts suivantes.
%p
%: Proposition{propconvinP}
\begin{proposition} \label{propconvinP}~
\begin{enumerate}
\item $\lim_n \alpha_n=0$ \ssi $\lim_n \rv(\alpha_n)=+\infty$.
\item Si $\lim_n \alpha_n=\alpha$ alors $\alpha$ est \iv \ssi $\exists N\,\forall n>N \;\rv(\alpha_n)=\rv(\alpha_N)<\infty$. Dans ce cas $\alpha^{-1}=\lim_{n\geq N}\alpha_n^{-1}$.
\item Si $\alpha=\lim_n \alpha_n$, $\beta=\lim_n \beta_n$ et $a\in\RRa$, alors 
\begin{itemize}
\item $a\alpha=\lim_na\alpha_n$,
\item $\alpha+\beta=\lim_n(\alpha_n+ \beta_n)$,
\item $\alpha \beta=\lim_n(\alpha_n \beta_n)$,
\item $\alpha\vu \beta=\lim_n(\alpha_n\vu \beta_n)$,
\item $\Fr(\alpha,\beta)=\lim_n\Fr(\alpha_n,\beta_n)$ et
\item $({\alpha^+})^q=\lim_n(\alpha_n^+)^q$ ($q\in\QQ,q>0$).
\end{itemize}
\item Tout $\alpha\in \gP_{j,d}$ est limite de la suite des \pols de Laurent $\pi_{m}(\vep^{1/d})$ pour $m\geq j$ obtenus par troncation de la série $\alpha$, définis \prmt par 
$$\pi_{m}\in\RRa[\vep^{1/d}][1/\vep],\;\pi_{m}=\som_{k:j\leq k<m}\rc_{k/d}(\alpha)\vep^{k/d}$$ 
On note aussi que $\RRa[\vep^{1/d}][1/\vep]\subseteq \gPal$. 
\end{enumerate}
 \end{proposition}
%----------- fin proposition ----------------------------- 

%:2025 modif de la phrase ci dessous
Notre but est de démontrer le \tho suivant.

%: theorem{propPuiseux}
\begin{theorem} \label{propPuiseux}
L'anneau $\gP=\gPo[1/\vep]$ satisfait tous les axiomes de la théorie \Sa{Crc2}. C'est donc un corps de Heyting réel clos \emph{non} discret non archimédien.
\end{theorem}
%----------- fin proposition ----------------------------- 

%
\begin{proof}[Première démonstration] On va généraliser les \prts de passage à la limite décrite dans la proposition \ref{propconvinP}~à toutes les \fsagcs définies sur $\RRa$.

\noindent L'article \cite{CM2013} démontre que $\gPal$ est un \crc discret. C'est donc une clôture réelle de $\RRa(\vep)$, construite d'une manière très différente de celle proposée dans~\cite{LR91}. 
Considérons maintenant un cube $[-a,a]^r=K\subseteq \RRa^r$ et une \fsagc $f\colon K\to\RRa$. Comme $\gPal$ est un \crc discret, $f$ se prolonge de manière unique en une \fsagc $f_1\colon K_1\to\gPal$, où $K_1\subseteq \gPal^r$ est défini par le même système d'inégalités que $K$. On va démontrer que $f_1$ se prolonge par continuité en une fonction $f_2\colon K_2\to\gP$, où $K_2\subseteq \gP^r$ est défini par le même système d'inégalités que $K$. 
Cela suffira à démontrer que $\gP$ est un modèle de \sa{Crc1}\footnote{Les détails pour cette affirmation sont laissés \alec.}.
%l
%: propdef{lemconvfsagc}
\begin{propdef} \label{lemonvfsagc} On applique les notations précédentes pour $K\subseteq K_1\subseteq K_2$.\\
Soit $f\colon K\to\RRa$ une \fsagc et $f_1\colon K_1\to \gPal$ son extension à $\gPal$.
Alors pour toute suite $(\alpha_n)$ dans $\RRa[\vep,\vep^{-1}]^r\cap K_2$ qui converge vers un $\alpha\in\gP$, la suite $f_1(\alpha_n)$ converge dans $\gP$.
La limite ne dépend que de $\alpha$ et est notée $f_2(\alpha)$. 
\end{propdef}
%----------- fin propdef ----------------------------------- 
%
\begin{proof} Cela ne semble pas si simple, il faut d'abord s'assurer que, les $f_1(\alpha_n)$ sont dans un même $\gP_{j,D}$; ensuite une inégalité de \L ojasievicz assurerait la convergence. \fbox{?? À préciser}
\end{proof}
\end{proof}

\hum{Il faudra aussi démontrer que la \dfn \ref{lemonvfsagc} se comporte correctement vis à vis de la composition des \fsagcs.} 

\hum{On a sans doute un résultat précis du style suivant, que j'explicite pour une \fsagc $f\colon [-1,1]\to\RRa$ entière sur $\RRa[X]$.
Pour $\alpha\in \RRa [\vep] $ et pour chaque $k$, il existe un $m$ tel que le \coe $\rc_{k/d}(f_2(\alpha))$ est une \fsagc des $\rc_{j/d}(\alpha)$ pour $0\leq j\leq m$. En outre la dépendance de $m$ en $k$ doit être bien contrôlée. 
}

\begin{proof}[Deuxième démonstration]
Vu la proposition \ref{propP_0} il suffit de démontrer la \prt d'existence des racines virtuelles pour l'anneau $\gP$. 
Une solution est peut-être de reprendre la proposition \ref{lemonvfsagc} en se limitant à prolonger par continuité les fonctions \ravs, ce qui pourrait se faire par \recu.
\fbox{?? À préciser}
\end{proof}

%%%%%%%%%%%%%%%%%%%%%%%%%%%%%%%%%%%%%%%%%%%%%%%%%%%%%%%%%%%%%%%%%%%%
%%%%%%%%%%%%%%%%%%%%%%%%%%%%%%%%%%%%%%%%%%%%%%%%%%%%%%%%%%%%%%%%%%%%

\section{Utilisation des racines virtuelles en algèbre réelle constructive}\label{secRVARC}
Les résultats énoncés dans cette sous-section pour le corps des réels semblent aussi valides dans la \tdy \Sa{Corv}. Certains nécessitent peut-être seulement~\Sa{Co0rv} ou~\Sa{Arc}.

\hum{Il faudrait préciser tout cela dans chaque cas. En tout cas, il serait extrêmement souhaitable que les résultats s'appliquent pour les \crcs \textsl{non} discrets tels que nous les définissons. Sinon, nous aurions manqué quelque chose d'essentiel.}

%--- SubSubsection{TDY}---- 
\Subsection{Sous-ensembles \sagqs de base de la droite réelle}
%: Subsection{Sous-ensembles \sagqs de base de la droite réelle}
%-----------

Définissons un \textsl{fermé \sagq de base} de la droite réelle comme un sous-ensemble de la forme $\rF\!_f=\sotq{x\in\RR}{f(x)\geq 0}$ pour un $f\in\RRX$.

\smallskip \noindent \textsl{Premier exemple}. Prenons le cas des \pols $f(X)=X^2-b$ et $g=-f$.
%i
\begin{itemize}
\item Si $b<0$, on a $\rF\!_f=\RR$ et $\rF\!_g=\emptyset$.
\item Si $b>0$, on a $\rF\!_f=\,]-\infty,-\sqrt b]\cup[\sqrt b,+\infty[$ et $\rF\!_g=[-\sqrt b,\sqrt b]$.
\item Si $b=0$, on a $\rF\!_f=\RR$ et $\rF\!_g=\so 0$.
\end{itemize}
Pour obtenir une description sous une forme aussi précise de ces fermés \sagqs il est absolument \ncr de connaitre
le signe de $b=-f(0)$. 
\\
Si l'on note $\alpha$ et $\beta$ les racines virtuelles de $f$, on a la description alternative suivante. 
%i
\begin{itemize}
\item Si $\alpha<\beta$, \cad si $f\big(\frac {\alpha+\beta} 2\big)>0$ on a $\rF\!_f=\,]-\infty,\alpha]\cup[\beta,+\infty[$ et $\rF\!_g=[\alpha,\beta]$.
\item Si $\alpha=\beta$ et $f(\alpha)<0$, \cad si $f\big(\frac {\alpha+\beta} 2\big)<0$, on a $\rF\!_f=\RR$ et $\rF\!_g=\emptyset$.
\item Si $\alpha=\beta$ et $f(\alpha)=0$, \cad si $f\big(\frac {\alpha+\beta} 2\big)=0$, on a $\rF\!_f=\RR$ et $\rF\!_g=\so \alpha$.
\end{itemize}

\smallskip \noindent \textsl{Deuxième exemple}. \\
Le cas d'un \pol \unt $f$ de degré $\delta>2$. Notons $\mathrm{Vr}_f$ la liste de ses racines virtuelles. Le
\thref{thVirtualRoots} permet de décrire l'adhérence de $\rF\!_f\cup \mathrm{Vr}_f$ de manière exacte comme l'adhérence de la réunion des intervalles suivants 
%i
\begin{itemize}
\item $]-\infty,\rho_{\delta,1}]\;\hbox{si }\delta\equiv0 \mod 2$
\item $[\rho_{\delta,k}, \rho_{\delta,k+1}] \hbox{ pour } k\in\lrb{0..\delta-2},\,k\equiv\delta\mod 2$
\item $[\rho_{\delta,\delta}, +\infty[$
\end{itemize}
En langage imagé imprécis: \gui{on connait $\rF\!_f$ à $\mathrm{Vr}_f$ près}. 
\\
De manière \gnle, le \pb avec un \pol de degré connu provient du fait que le \thref{thVirtualRoots}
affirme quelque chose de précis concernant le signe du \pol sur un intervalle $[\rho_{\delta,j}, \rho_{\delta,j+1}]$
uniquement lorsque $\rho_{\delta,j}< \rho_{\delta,j+1}$. On a alors le résultat suivant.

%l
%: Lemma{lemrFf}
\begin{lemma} \label{lemrFf}
Soit $f\in\RRX$ un \pol de degré $\delta$ connu et $g=f/c_\delta$ le \pol \unt correspondant ($c_\delta$ est le \coe dominant, $>0$ ou $<0$).
Notons $\rho_{\delta,k}=\rho_{\delta,k}(g)$.
\begin{enumerate}
\item L'adhérence de $\rF\!_f\cup \mathrm{Vr}_f$ est égale à l'adhérence d'une réunion explicite d'intervalles fermés
avec pour bornes des $\rho_{\delta,k}$ ou $+\infty$, ou $-\infty$.
\item Lorsque l'on connait les signes des $(\rho_{\delta,k+1}-\rho_{\delta,k})$ 
et des $g(\rho_{\delta,k})$, on a une description exacte du fermé $\rF\!_f$ sous forme d'une réunion
d'intervalles fermés disjoints.
L'information \ncr est \eqve à la donnée des signes 
 des $g(\frac {\rho_{\delta,k}+\rho_{\delta,k+1}} 2)$.
\end{enumerate}

\end{lemma}
%----------- fin lemma ----------------------------------- 

\noindent \textsl{Lorsque le degré de $f$ n'est pas connu}, on perd le contrôle de la situation {en $+\infty$ et $-\infty$}. 
La situation la plus floue, dans laquelle on ne contrôle rien du tout, se présente lorsque l'on ne sait pas si le \pol est identiquement nul ou non.

\smallskip On a des résultats analogues
pour un ouvert de base $\rU_f=\sotq{x\in\RR}{f(x)> 0}$. 

%--- SubSubsection{TDY}---- 
\Subsection{Tableaux de signes et de variations}
%: Subsection{Tableaux de signes et de variations}
%-----------

Soit $\gR$ un modèle \cof de \Sa{Co--}.
On dit que deux \elts $a$ et $b$ sont \gui{distincts} si $a\neq b$,
\cad $a-b$ est \iv.

%l
%: Lemma{lemMinoration}
\begin{lemma} \label{lemMinoration}
Étant donnée une liste $L$ de $k$ \elts et une liste $L'$
de $k+\ell$ \elts distincts dans $\gR$, au moins $\ell$ \elts de $L'$ sont distincts de tous les \elts de $L$. 
\end{lemma}
%----------- fin lemma ----------------------------------- 
%
%\begin{proof}
%
%\end{proof}
%

Le \thref{thVirtualRoots}, points \textsl{\ref{ivrvariation}} et 
\textsl{\ref{ivrsigne}}, donne presque un tableau de signes et de variations complet pour le \pol unitaire $f$.

Pour le tableau des signes complet de $f$, les hésitations éventuelles concernent les signes de $f$ en les \ravs $\xi_k$ de~$f'$.
Même chose pour le tableau de variations de $f$, avec les signes de $f'$ aux racines virtuelles de $f''$.

Ceci conduit au résultat suivant.

%f
%: Proposition{factTableaucomplet}
\begin{proposition} \label{factTableaucomplet} Soit $\gR$ un \covr.
\begin{enumerate}
\item Soit $f(x)\in\Rx$ unitaire de degré $k\geq 2$ et soient $k+\ell-1$ \elts $a_i$ distincts dans~$\gR$. Pour au moins $\ell$ de ces \elts, le \pol $f(x)+a_i$
a un signe strict connu en chacune des racines virtuelles de $f'$, et son tableau complet de signes est connu de manière exacte. 
\\
Si $k=2$ le tableau complet de signes et de variations est alors connu de manière exacte.
\item Soit $f(x)\in\Rx$ unitaire de degré $k\geq 3$, soient $k+\ell-1$ \elts $a_i$ distincts dans $\gR$, \hbox{et $k+\ell-2$} \elts $b_j$ distincts dans $\gR$. Pour au moins $\ell^2$ des couples $(a_i,b_j)$, 
on a un tableau complet de signes et de variations connu de manière exacte pour le \pol $f(x)+b_jx+a_i$.
\end{enumerate} 
\end{proposition}
%----------- fin fact ----------------------------------------- 
%
%\begin{proof}
%\hum{Les \demos méritent d'être écrites en détail
%(peut-être l'énoncé n'est pas tout à fait correct)}
%\end{proof}
%

%r
%: Remark{remfactTableaucomplet}
\begin{remarks} \label{remfactTableaucomplet} ~

\noindent 1) On a probablement un résultat de perturbation du même style 
qui dit que \textsl{pour presque toutes les perturbations} d'un nombre fini de \pols unitaires $f_i$, on connait de manière sûre les \egts et les inégalités strictes entre toutes les racines virtuelles des $f_i$ et de toutes leurs dérivées, ainsi que les signes des $f_i$ en chacune de ces \ravs, ce qui donne un tableau complet de signes 
et de variations pour la famille des $f_i$ et de leurs dérivées.

\smallskip \noindent 2)
Si l'on veut un résultat analogue à la proposition \ref{factTableaucomplet} pour une \fsagc,
il faudra se placer dans la théorie \Sa{Corv} et limiter le tableau de signes et de variations recherché à un intervalle fermé borné. 
\eoe\end{remarks}
%----------- fin remark ---------------------------------- 

%--- SubSubsection{TDY}---- 
\Subsection{Une décomposition algébrique cylindrique ({CAD}) approximative?}
%: Subsection{Une décomposition algébrique cylindrique approximative?}
%-----------

On se pose le \pb de donner une CAD approximative pour une famille finie de \pols de $\RXn$ où $\gR$ est un modèle \cof de \Sa{Co0rv} (ou de \Sa{Corv}). 
Ce serait un résultat qui généraliserait de manière astucieuse
le lemme \ref{lemrFf} ou la proposition \ref{factTableaucomplet}.

En termes de déménageurs de pianos, au lieu de décider si \gui{cela passe} ou \gui{cela ne passe pas}, on obtiendrait des résultats approximatifs du genre suivant: en fonction des données du \pb et d'une précision souhaitée $\epsilon$, on calculerait \und{de manière uniforme} un $\alpha\in\gR$ tel que:
%i
\begin{itemize}
\item si $\alpha>0$, il y a moyen de passer en respectant une distance $>\epsilon$ des obstacles, et l'on vous dit comment,
\item si $\alpha<1$ il n'y a pas moyen de passer en respectant une distance $>2\epsilon$.
 \end{itemize}
Naturellement le piano doit être un compact \sagq bien défini, les obstacles itou, et l'espace dans lequel on déplace le piano itou.

\smallskip De manière \gnle, comme il est impossible de contrôler, même de manière approximative, le comportement à l'infini d'un \pol dont tous les \coes sont proches de $0$, il faut \ncrt se limiter à
calculer une CAD approximative pour une famille finie de \pols sur un compact bien défini du style $[0,1]^n$. On s'en rend compte, si l'on essaie de reproduire une CAD usuelle (pour un \crc discret) sur $\RR$, au fait que les \coes d'un \pol sous-résultant peuvent très bien a priori 
être tous très proches de $0$. Or à priori les racines virtuelles ne sont efficaces que pour les \pols unitaires.

Sur ce genre de sujet, on en est à l'époque des balbutiements.

%%%%%%%%%%%%%%%%%%%%%%%%%%%%%%%%%%%%%%%%%%%%%%%%%%%%%%%%%%%%%%%%%%%% 
\Subsection{Stratifications}
%: Subsection{Stratifications}
%-----------

Il semble que les stratifications, lorsqu'on les suppose données en hypothèse, soient un cadre reposant dans lequel beaucoup de résultats valables pour les \crcds s'étendent sans trop de difficulté au cas \textsl{non} discret. 

\hum{À vérifier. Cela pourrait remplacer avantageusement les balbutiements précédents.}
 
%%%%%%%%%%%%%%%%%%%%%%%%%%%%%%%%%%%%%%%%%%%%%%%%%%%%%%%%%%%%%%%%%%%% 
\Subsection{Le \tho fondamental de l'algèbre (\tsbf{TFA})}
%: Subsection{Le \tho fondamental de l'algèbre}
%-----------

Pour un traitement du \tsbf{TFA} sans l'axiome du choix dépendant,
voir \cite{Ric2000}.

 Puisque les racines virtuelles sont des fonctions continues, et puisqu'il est impossible de suivre par continuité les zéros d'un \pol complexe (unitaire de degré $m$ fixé et à \coes variables), on ne peut certainement pas obtenir un de ces zéros exprimé 
comme un \elt de $\SaCe_m(\RR)$. Néanmoins, on peut recouvrir les zéros d'un \pol complexe de degré $\delta$ par un nombre fini d'expressions dans $\Sace_{\delta^2}(\RR)$.

Ce que l'on aimerait ici, c'est le faire manière assez optimale.

\paragraph{1. Les racines carrées d'un nombre complexe $c=a+ib$.} \\
Les zéros du \pol $f(Z)=Z^2-c$ sont donnés sous la forme $x+iy$ par les solutions réelles du \sys \gui{$x^2-y^2=a$, $2xy=b$} et sont calculées comme suit:
%i
\begin{itemize}
\item $(x^2+y^2)^2=a^2+b^2$, donc $x=\pm u$ avec $u=\sqrt{\frac1 2 (a+\sqrt{a^2+b^2})}\in\Sace_{4}(\RR)$
\item $y=\pm v$ avec $v=\sqrt{\frac1 2 (-a+\sqrt{a^2+b^2})}$, avec la contrainte $xyb\geq 0$.
\end{itemize}

\noindent Si l'on note $z_1=u+iv$, $z_2=-z_1$, $f_1(Z)=(Z-z_1)(Z-\ov{z_1})$, $f_2(Z)=(Z-z_2)(Z-\ov{z_2})$ et $g(Z)=(Z-c)(Z-\ov c)$ on obtient l'\egt
% equation label {eqsqrC}
\vspace{-.5em}
\begin{equation} \label {eqsqrC}
g(Z^2)=f(Z)\ov f(Z)=f_1(Z)f_2(Z)
\end{equation}
% end-equation
Les \pols $f_1$, $f_2$, $g$ et $f\ov f$ sont réels, partout $\geq 0$, chacun muni d'un certificat \agq simple pour son caractère $\geq 0$ lorsque la variable est réelle. Lorsque $c\neq 0$, les zéros de $f$ se répartissent entre les zéros de $f_1$ et ceux de $f_2$.

On peut estimer que l'on a ainsi obtenu la solution optimale pour les racines carrées d'un nombre complexe dans le cadre de la \RRlg de fonctions engendrée par les fonctions \gui{racines carrées virtuelles} $\rho_{2,2}$, et plus \gnlt la solution optimale dans le cadre des \algs $\Sace_m(\RR)$.

Notons qu'à priori, $x$ et $y$ étant racines de \pols réels de degré $4$, cela faisait~$16$ choix possibles pour $x+iy$. 

\paragraph{2. Le cas général}. ~\\
On a à priori le résultat non optimal suivant.

%: Proposition{propTFA}
\vspace{-.5em}
\begin{proposition}[\tsbf{TFA} via les racines virtuelles] \label{propTFA}~ \\
Soit $f$ un \pol \unt complexe de degré $\delta$.
Il existe $\delta^4$ \pols $q_\ell$ quadratiques positifs\footnote{\Prmt: des \pols unitaires de degré $2$ partout $\geq 0$.} ayant leurs \coes construits sur des $\rho_{\delta^2,k}(\dots)$ (polyracines en les parties réelles et imaginaires des \coes de $f$)
 tels que $f\ov f$ divise le produit des $q_\ell$. \\
 Si le \crc considéré est discret, le \pol $f(z)$ se décompose en un produit de facteurs $(z-\zeta_j)$ explicites sur $\gC$, avec les $\zeta_j$ dont les parties réelles et imaginaires sont des racines de \pols réels unitaires de degré $\delta^2$, dont les \coes sont des $\QQ$-\pols en les parties réelles et imaginaires des \coes de $f$.
\end{proposition}
%--------- fin proposition ----------------------------- 
%
\begin{proof}
La partie réelle d'un zéro $\zeta_j$ de $f$ s'écrit $(\zeta_j+\ov{\zeta_j})/2$. Les $(\zeta_j+\ov{\zeta_k})/2$ sont au nombre de $\delta^2$ et ce sont les zéros d'un \pol réel $h_1$ de degré $\delta^2$ dont les \coes s'expriment comme des $\QQ$-\pols en les parties réelles et imaginaires des \coes de $f$. Parmi les
zéros réels de~$h_1$ figurent les $\frac1 2(\zeta_j+\ov{\zeta_j})$. Ce sont donc des racines virtuelles de $h_1$.
Raisonnement analogue pour la partie imaginaire, avec un \pol réel $h_2$ de degré $\delta^2$. Si $\alpha$ est une racine virtuelle de~$h_1$ et~$\beta$ une racine virtuelle de~$h_2$, on leur associe le \pol 
$$q=(z-\alpha-i\beta)(z-\alpha+i\beta)=(z-\alpha)^2+\beta^2$$ 
qui est l'un des $q_\ell$ de l'énoncé.
\end{proof}

%r
%: Remark{rem}
\begin{remark} \label{remPR2019} 
Dans l'article \cite{PR2019} les auteurs démontrent qu'un corps ordonné discret $\delta^2$-clos (\cad qui satisfait le \tho des valeurs intermédiaires pour les \pols de degré $\leq \delta^2$) satisfait le \tho fondamental de l'algèbre pour les \pols de degré $\leq \delta$. On peut déduire ce résultat
de la proposition \ref{propTFA} en utilisant le \Pst formel comme suit. On suppose que le \crc est discret. Alors, le fait que $f\ov f$ divise le produit des $q_\ell$ implique que $f$ admet au moins un zéro complexe, parmi les zéros des $q_\ell$\footnote{On a un peu mieux. Le produit des $q_\ell$ se décompose en un produit de facteurs linéaires, donc \irds, dans $\gC[Z]$. Comme $f(Z)$ divise ce produit, et comme $\gC[Z]$ est un anneau intègre à pgcd, c'est en fait un sous-produit.}. Par ailleurs les racines virtuelles de $h_1$ et de $h_2$ sont caractérisées par des \syss d'inégalités larges. Une \ralg sur le langage des corps odonnés dit que, pour un \crc discret, si l'on met en hypothèse ces \syss d'inégalités larges,
on obtient comme conséquence valide le fait que le produit des $f(\alpha\pm i\beta)$ convenables est nul. D'après le point \textsl{2} du \Pst formel \ref{Pst1}, cette \ralg est valide pour tout corps ordonné (discret ou pas)
ainsi que pour les \arcs,
car elle est valide dans la théorie \Sa{Asonz}. Et si le corps satisfait le \tsbf{TVI} pour les \pols de degré $\leq \delta^2$, les hypothèses sont satisfaites par les \ravs de $h_1$ et~$h_2$. De même, si le langage des
corps ordonnés a été enrichi en introduisant les fonctions \ravs pour les \pols de degré $\leq \delta^2$, avec les axiomes correspondants, on obtiendra aussi pour les \sads correspondantes le fait que le produit des $f(\alpha\pm i\beta)$ convenables est nul.\eoe 
\end{remark}
%----------- fin remark ---------------------------------- 

%\vspace{-1.5em}
\paragraph{3. Le cas général en termes de multisets.} \label{TFA3}
Référence: le \tsbf{TFA} dans \cite[Richman]{Ric2000}. 
\\
\textsl{À priori}, le \gui{\tsbf{TFA} version multisets} semble difficile à formuler correctement sans disposer de l'espace métrique des $n$-multisets de nombres complexes.

Pour contourner l'obstacle, on peut réduire le \gui{\tsbf{TFA} version multisets} à un ensemble de \rdys valides donnant une formulation essentiellement \eqve qui utilise le comptage
du nombre de zéros à l'intérieur de rectangles
dans le style de~\cite[Eisermann]{Eis2012}. L'article \cite{PR2019} nous semble donner toutes les précisions \ncrs.

Comme nous ne supposons pas que le corps ordonné est discret, nous devons utiliser uniquement des rectangles sur les bords desquels on est certain qu'il n'y a aucun zéro complexe du \pol considéré. 

Un test explicite montre qu'il faut éviter au plus $\delta$ lignes horizontales et au plus $\delta$ lignes verticales pour nos rectangles. Cela se formule en disant que si l'on considère $\delta+m$ lignes horizontales distinctes, on est certain que au moins $m$ d'entre elles sont bonnes (même chose pour les verticales).

Pour ces rectangles, le comptage des zéros à l'intérieur fonctionne et donne toujours un nombre entier bien défini. 

On a une règle valide qui assure qu'aucun zéro complexe
ne se trouve en dehors d'un rectangle suffisamment grand explicite. Pour ce rectangle suffisamment grand le comptage donne le nombre~$\delta$ attendu.
Et une règle valide dit que lorsqu'on coupe un rectangle en deux, la somme des deux comptages égale le comptage précédent.

Si l'on veut aussi traiter le cas non-archimédien, il faut établir des \ralgs
disant que l'on peut enfermer les $\delta$ zéros dans une réunion de rectangles de taille arbitrairement petite. 

Une étude poussée de l'article \cite{PR2019} devrait aboutir aux résultats souhaités, qui sont plus précis que le \tsbf{TFA} envisagé au point 2,
résultats que l'on peut considérer comme étant \underline{la} forme \cov satisfaisante du \tsbf{TFA},
et qui seront valables pour la théorie \sa{Corv}, formulables comme des \ralgs valides dans cette théorie. Mais ces \ralgs ne seraient pas valides dans la théorie des \arcs.

%
%L'énoncé obtenu devrait ressembler à ce qui suit pour un \pol unitaire $f$ de degré $d$ sur~$\gC=\gR[i]$: 
%
%\smallskip \noindent \textsl{Pour tout $\varepsilon>0$, on peut trouver une famille finie $(z_i,\rho_i,r_i)$ où $z_i\in\gC$, $\rho_i>0\in\gR$ et $r_i\in\NN^>$ telle que:
%\begin{itemize}
%%
%\item $\sum_ir_i=d$;
%%
%\item $\norm{\prod_{i}(X-z_i)^{r_i}-f}\leq \varepsilon$; 
%%
%\item $\abs{z_i-z_j}>2(\rho_i+\rho_j)$ pour $i\neq j$; 
%%
%\item $\forall z\in\gC$, $\big(\Vi_i\abs{z-z_i}\geq \rho_i\big)\;\Rightarrow\;\abs{f(z)}>\varepsilon$. 
%\end{itemize}
%}

% Subsection
\section{Quelques questions}

\Subsection{Fonctions \sagcs}

%: Question{Qu-continuiteRV}
\begin{question} \label{Qu-continuiteRV} Rendre plus explicite le résultat (\cof)
de continuité des fonctions racines virtuelles: point \emph{2} du \thref{thVirtualRoots}. 
\textsl{Chaque fonction $\rho_{d,j}:\gR^d\to\gR$ est uniformément continue sur toute boule $\mathrm{B}_{d,M}:=\sotQ{(a_{d-1},\dots,a_0)\,}{\,\sum_ia_i^2\leq M}$, ($M>0$).} La continuité devrait être donnée sous forme complètement explicite à la {\L}ojasiewicz.
\hum{cela doit pouvoir se faire par \recu sur le degré $d$?}
\end{question}
%--------- fin fact ---------------------------------------------- 

%: Question{questArc3}
\begin{question} \label{questArc3} 
Donner une \prco du \thref{thArc3}.
\end{question}
%----------- fin question ----------------------------- 

%: Question{Qu-polrootsafrvr}
\begin{question} \label{Qu-polrootsafrvr} Soit $\gR$ un \arc. Tout \elt de $\SaC_n(\gR)$ entier sur l'anneau des \pols $\Rxn$ est-il un \elt de $\SaCe_n(\gR)$? 
\end{question}
%--------- fin fact ---------------------------------------------- 

%Les deux questions précédentes semblent abordables. Après tout, tout vient de $\RRa$, et les résultats semblent démontrables en \clama.
%
%A priori beaucoup plus difficile est la question suivante.
%: Question{Qu-polrootsRR}
\begin{question} \label{Qu-polrootsRR} Est-ce que toute fonction continue $\RR^m\to\RR$ entière sur l'anneau des \pols est un \elt de $\SaCe_m(\RR)$? 
La réponse est positive en \clama car on peut appliquer le \thref{thPBpolyroots} à $\RR$. 
\end{question}
%--------- fin fact ---------------------------------------------- 

Les questions \ref{questArc2} ont à voir avec le caractère o-minimal
de la structure de \crc \textsl{non} discret.
Le mot \gui{compact} ci-dessous est utilisé comme signifiant \gui{fermé borné}.
%q
%: Question{questArc2}
\begin{questions} \label{questArc2} \textsl{(se rappeler la proposition \ref{propfsagccovrsup})}\\
On considère un \covr $\gR$. \hum{un \arc?}
\begin{enumerate}
\item Démontrer qu’une \fsagc qui est partout~$>0$ sur le compact~\hbox{$%\mathrm{K}_n=
[0,1]^n\subseteq \gR^n$}
est minorée (sur ce compact) par un \elt $>0$. Et que la \bif est un \elt de $\gR$.
\item Étendre le résultat à un compact \sagq \gui{bien défini} arbitraire:
on entend par là un \fsa borné $K$ pour lequel la fonction \gui{distance à $K$} est une \fsagc (un \elt de $\Sac_n(\gR)$).
\end{enumerate}
\end{questions}
%----------- fin questions ----------------------------- 

\Subsection{Clôture réelle}

%: Question{quest2cloture}
\begin{question} \label{quest2cloture} 
Si $\gK$ est un modèle de \Sa{Co} (ou de \Sa{Co--}), est-ce que sa 2-clôture $\gL$ en tant qu'\afr est encore un modèle de \sa{Co} (ou de \sa{Co--})? 
\end{question}
%----------- fin question ----------------------------- 

On reprend la question précédente (en ajoutant des précisions) pour la clôture réelle.
%: Question{questCloturevirtuelle}
\begin{question} \label{questCloturevirtuelle} Soit $\gK$ un modèle de la théorie \Sa{Co} et $\gR$ la \sad $\Sa{Corv}(\gK)$ 
 (comme \paref{pagecloturevirtuelle}). $\gR$ est-il un modèle constructif de \sa{Corv}? 
\\
En particulier, le métathéorème suivant est-il satisfait? Étant donnés deux termes clos $\alpha$ et $\beta$ de $\gR$ tels que la règle $\vd \alpha+\beta>0$ est valide, est-il vrai que l'une des deux règles $\vd \alpha>0\,$, $\vd \beta>0$ soit valide?
\\
On peut poser la même question sous la forme suivante: si $\gK$ est un modèle de \sa{Co}, est-ce que la structure \agq (usuelle) $\ARC(\gK)$ satisfait la règle \Tsbf{AFRL}? 
\end{question}
%----------- fin question ----------------------------- 

%: Question{questEis2012}
\begin{question} \label{questEis2012} 
Est-ce que l'article \cite[Eisermann]{Eis2012} peut être relu pour un \crc \textsl{non} discret, i.e. dans la théorie \Sa{Crc2}? Cette question demande de développer en détail les idées dans le point 3, \paref{TFA3} dans le paragraphe concernant le \tsbf{TFA}.
\end{question}
%----------- fin question ----------------------------- 

%q
%: Question{questTVICrc}
\begin{question} \label{questTVICrc} 
Montrer que le \tho de la valeur intermédiaire, énoncé sous la forme de la règle \RCFn\ \paref{RCFn}, n'est pas valide dans la théorie \Sa{Crc2}.
Notez qu'une forme légèrement affaiblie est valable: voir le point \textsl{\ref{i4Pst3}} 
du Positivtellensatz formel \ref{Pst3} et le point \textsl{\ref{ivrTVI}} du \thref{thVirtualRoots}.
\end{question}
%----------- fin question ----------------------------- 

%: Question{questTVICrc2}
\begin{question} \label{questTVICrc2} 
Montrer que le \tho qui affirme que tout nombre complexe a une racine carrée n'est pas une règle valide dans la théorie \Sa{Crc2}. Comparer avec la proposition \ref{propTFA} qui pourrait sembler affirmer le contraire.
\end{question}
%----------- fin question ----------------------------- 

\Subsection{Pierce-Birkhoff}
%q
%: Question{questpbring}
\begin{questions} \label{questpbring} ~

\noindent 1) la \dfn d'un anneau de Pierce-Birkhoff donnée en \ref{defiSacembis} coïncide-t-elle en \clama avec la notion définie dans \cite[Madden, 1989]{Mad89}?

\smallskip \noindent 2) Si c'est bien le cas, se pose la question de donner des \prcos pour des résultats sophistiqués, comme le fait qu'un anneau cohérent \noe régulier de dimension $\leq 2$ est un anneau de Pierce-Birkhoff \cite{LMSS2012}.

\smallskip \noindent 3) Rappelons que la conjecture de Pierce-Birkhoff usuelle est démontrée dans \cite{Mahe83} \hbox{pour $\gR[x,y]$} lorsque $\gR$ est un \crc discret mais
il n'est pas si clair qu'il y ait une \prco pour $\RR[x,y]$. 
\end{questions}
%----------- fin question ----------------------------- 

\Subsection{Le 17-ème \pb de Hilbert}
%: Subsection{Le 17-ème \pb de Hilbert}

%q
%: Question{quest17H}
\begin{question} \label{quest17H} 
Dans quelle mesure la solution \cov du 17\ieme\ \pb de Hilbert
pour~$\RR$ (voir \cite[section 6.1]{GL93}) s'applique à tout \asrvr?
Si ce n’est pas le cas, quelle théorie plus forte ferait l’affaire:
\Sa{Crc2}, \Sa{Arc}, \Sa{Crca}, \Sa{Crc3} (\paref{theorieCrc3})?
\end{question}
%----------- fin question ----------------------------- 

\Subsection{Le Graal?}
%: Subsection{Le Graal?

On se pose la question d'un \tho analogue à \ref{thTarski}, pour le cas non discret. 

Le \Pst formel \ref{Pst3} implique que la théorie \sa{Arc} est la \talg engendrée par $\gRa$, par $\RR$ ou par $\RR_{\tsbf{PR}}$ sur la signature de \sa{Arc}.

%q
%: Question{questGraal1}
\begin{question} \label{questGraal1} 
La théorie \Sa{Arc} est-elle la skolémisée de la théorie cartésienne engendrée par $\gRa$, par $\RR$ ou par $\RR_{\tsbf{PR}}$ sur la signature des anneaux commutatifs? 
\end{question}
%----------- fin question ----------------------------- 
%q
%: Question{questGraal2}
\begin{question} \label{questGraal2} 
En quel sens pourrait-on dire que la théorie \Sa{Corv} est la \tdy engendrée par $\RR$ \gui{sans axiome du choix dépendant} sur la signature de \sa{Arc}? Même question avec $\RR_{\tsbf{PR}}$.

\noindent NB: cette question semble impossible à formuler en \clama, et en \coma, il faudrait avoir une idée claire de $\RR$ \gui{sans axiome du choix dépendant}. 
\end{question}
%----------- fin question ----------------------------- 

\newpage \thispagestyle{empty}

%%%%%%%%%%%%%%%%% CHAPITRE %%%%%%%%%%%%%

\chapter{L'axiome d'archimédianité}\label{secGeomReelsArchi}
\Today
\minitoc

Dans ce chapitre, pour mieux décrire les propriétés \agqs de $\RR$, on fait une tentative qui consiste à ne pas quitter les théories dynamiques tout en conservant l'essence de la règle non dynamique~\Tsbf{HOF}.

Le langage reste cependant pour l'essentiel celui des anneaux ordonnés.

Dans la troisième partie, nous ferons une tentative beaucoup plus ambitieuse
en utilisant un langage beaucoup plus riche, qui nous fera voir essentiellement une \tgm des réels comme précurseure de la théorie des structures o-minimales. 

% section
\section{Corps réel clos archimédien \textsl{non} discret}

 La règle suivante, qui signifie que le corps est archimédien, est satisfaite sur $\RR$

\Regles{\Lab{AR1} $\dsp\vd \Vou_{n\in\NN} \;\abs x \leq n$ \qquad \qquad\qquad\qquad ({\bf Archimède 1})}

%d
%: Definition{defiCrca}
\begin{definition} \label{defiCrca}
On définit la théorie \SA{Crca} des \crcs archimédiens comme la \tgm obtenue en ajoutant l'axiome \Tsbf{AR1} à la théorie \Sa{Crc2}. 
\end{definition}
%----------- fin definition -------------------------------- 

L'exemple donné dans le point 3 de la remarque \ref{remCrc} (un anneau réel clos local avec des diviseurs de zéro, modèle de la théorie \sa{Crc2}) reste un modèle de \sa{Crca}. Les exemples \ref{exacorpsnondiscret} sont
\egmt des modèles de la théorie \sa{Crca}:
de manière \gnle les sous-anneaux de $\RR$ stables pour les fonctions racines virtuelles, la fraction $\mathrm{Fr}$ et les inverses des \elts inversibles,
sont des modèles de \sa{Crca}.

%%%%%%%%%%%%%%%%%%%%%%%%%%%%%%%%%%%%%%%%%%%%%%%%%%%%%%%%%%%%%%%%%%%%
%%%%%%%%%%%%%%%%%%%%%%%%%%%%%%%%%%%%%%%%%%%%%%%%%%%%%%%%%%%%%%%%%%%%

% Subsection
\section{Quelques questions}

\Subsection{Axiome d'archimédianité}

%: Question{questRarchi}
\begin{questions} \label{questRarchi} ~\\
On sait qu'on ne peut pas exprimer de manière fintaire le fait que $\RR$ est archimédien. On l'exprime avec la règle infinitaire \Tsbf{AR1}.
\begin{itemize}
\item Une question qui se pose est de savoir si la théorie \Sa{Crca} obtenue en ajoutant la règle \Tsbf{AR1} est une extension \cosv de \Sa{Crc2}, cela semble probable.
\item On pourrait commencer par montrer que \Sa{Crca} prouve les mêmes \ralgs que~\Sa{Crc2}. 
\item Par contre, pour la théorie formelle correspondante dans laquelle on autorise l'introduction de prédicats pour $P\Rightarrow Q$ et $\forall x P$ (avec les règles de déduction de Gentzen) il se pourrait 
qu'un énoncé comme \Tsbf{HOF} devienne prouvable. 
\end{itemize}

\end{questions}
%----------- fin question -----------------------------

%%%%%%%%%%%%%%%%%%%%%%%%%%%%%%%%%%%%%%%%%%%%%%%%%%%%%%%%%%%%%%%%%%%%
\paragraph{Le principe d'omniscience \tsbf{LPO} sans danger en \alg réelle?}

%: Question{questRLPO}
\begin{question} \label{questRLPO} ~\\ 
La règle suivante n'est pas satisfaite sur $\RR$, car elle implique la règle \Edinq\ .

\Regles{\Lab{AR2} $\dsp\vd x= 0 \; \vou \; \Vou
_{n\in\NN} \abs x > 1/2^n$ \qquad\quad ({\bf Archimède 2})
}

\noindent Mais sans doute elle est \gui{admissible}, au sens où son ajout à \Sa{Crca} fournirait une extension conservative de \Sa{Crcd}. 

\noindent Ce résultat serait une sorte de réalisation du programme de Hilbert pour \tsbf{LPO}, limité à la théorie~\Sa{Crca}. En effet, la théorie \Sa{Crcd} est elle-même inoffensive par rapport à \Sa{Crc2} car elle prouve les mêmes \ralgs.
\end{question}
%----------- fin question -----------------------------

%%%%%%%%%%%%%%%%%%%%%%%%%%%%%%%%%%%%%%%%%%%%%%%%%%%%%%%%%%%%%%%%%%%%
\paragraph{Des séries convergentes en \alg réelle?}

%: Question{questRLPO}
\begin{question}
On note $[x]^n=\frac 1{2^n}\vi(x\vu - \frac 1{2^n})$. La règle suivante n'est pas une \rdy

\UneRegle{Cauchy} { $\dsp\vd \Exists x\,\Vii_{n\in\NN}\, 
\abs {x -\som_{p=0}^n[x_p]^p} \leq 1/2^n $ 
}

Il faudrait introduire un
symbole de fonction $\sum_{n=0}^\infty [x_n]^n$ pour ces sommes infinies. 
Cela remplacerait la règle illégitime \tsbf{Cauchy} par une infinité de \ralgs légitimes. 
Mais un tel symbole de fonction est-il légitime?
\end{question}
%----------- fin question -----------------------------

%%%%%%%%%%%%%%%%%%%%%%%%%%%%%%%%%%%%%%%%%%%%%%%%%%%%%%%%%%%%%%%%%%%%
\Subsection{Le Positivstellensatz de Schmüdgen}\index{Positivstellensatz!de Schmüdgen}
%: Subsection{Le Positivstellensatz de Schmüdgen}

Références: \cite{Schm91,Schw2002,Schw2003}

%q
%: Question{questSchmudgen}
\begin{question} \label{questSchmudgen} 
La \tgm \Sa{Crca} suffit-elle pour développer les \thos
du type Schmüdgen?
 \end{question}
%----------- fin question ----------------------------- 

\newpage \thispagestyle{empty}

%%%%%%%%%%%%%%%%%%%%%%%%%%%%%%%%%%%%%%%%%%%%%%%%%%%%%%%%%%%%%%%%%%%%

\chapter*{Conclusion}
%\newpage\vspace{2cm}\noindent {\Huge{\bf Conclusion}}
\addstarredchapter{Conclusion}

Les questions les plus importantes qui restent à résoudre pour cette 2\ieme partie nous semblent les suivantes.

\begin{enumerate}
\item Question \ref{questthParamcontFsagc0}. Donner une \demo \cov du \thref{thParamcontFsagc0}.
\\
\textsl{Soit $\gR$ un \crcd et $f\colon \gR^n\to\gR$ une \fsagc. Il existe un entier $r\geq 0$, une \fsagc $g\colon \gR^{r+n}\to\gR$ définie sur~$\RRa$, et un \elt $\uy\in\gR^r$ tels que}
$$
\forall \xn\in\gR\;\; f(\xn)=g(\yr,\xn).
$$

\item Question \ref{Qu-continuiteRV}. Rendre explicite le
point \emph{2} du \thref{thVirtualRoots} affirmant la continuité uniforme des fonctions $\rho_{d,j}:\gR^d\to\gR$ sur toute boule fermée.

\item Question \ref{questArc3}. Donner une \prco du \thref{thArc3} 
\textsl{Toute \fsagc $\RRa^n\to\RRa$ peut être définie par un terme de la théorie~\sa{Arc}.} Cela permettra de clarifier définitivement la \dfn \cov
\ref{defiArc} des \arcs et son rapport en \clama avec différentes \carns \covs des \arcs. 

\item Question \ref{questCloturevirtuelle} concernant la possibilité de construire une clôture réelle d'un corps ordonné non discret. 
\textsl{Soit $\gK$ un modèle de la théorie \Sa{Co} et $\gR$ la \sad $\Sa{Corv}(\gK)$ 
 (comme \paref{pagecloturevirtuelle}). $\gR$ est-il un modèle constructif de \sa{Corv}}?

\end{enumerate}
\newpage \thispagestyle{empty}

\part[Version améliorée de la théorie des \crcs \textsl{non} discrets]{Une version améliorée de la théorie des \crcs \textsl{non} discrets et une tentative de version \cov des structures o-minimales}

\chapter*{Introduction}
\addstarredchapter{Introduction}

Nous explorons dans cette troisième partie la possibilité de mieux décrire les \prts \agqs de~$\RR$ en étendant le langage par l'introduction de sortes pour les \fsagcs sur des cubes compacts.

\smallskip Nous notons en effet que la situation \gnle s'est éclaircie quand nous avons introduit les
fonctions $\cdot\vu\cdot$, $\cdot\vi\cdot$ et les fonctions \ravs. Ces extensions naturelles du langage utilisé ont permis en bonne partie de surmonter les obstacles que la notion de \gui{corps ordonné \textsl{non} discret} semble offrir à une formalisation en \tdy finitaire. 

\smallskip Cependant, d'un point de vue \cof, il n'est pas naturel de s'intéresser aux zéros réels des seuls \pols dont le degré est fixé.
La bonne raison pour cela avec $\RR$ est que l'on ne contrôle pas les zéros au voisinage de l'infini quand le degré n'est pas clairement fixé. En remplaçant $\RR$ par l'intervalle réel $\II_\RR=[-1,1]\subseteq \RR$
cette soi-disant bonne raison disparait d'elle-même. 

\smallskip L'idée est que l'on contrôle les choses
de manière \cov uniquement dans le cadre des compacts. Il faut se désintoxiquer de $\pm\infty$ et revenir aux \maths grecques! En conséquence il faut laisser tomber~$\RR$ au profit de l'intervalle $\II_\RR$, par exemple en remplaçant $+$ par la demi-somme. Cela nécessite de reprendre l'axiomatique, mais le bénéfice sera que l'on pourra plus facilement formuler certaines \prts liées au fait que du point de vue \cof $\RR$ \textsl{n'est pas} discret.

\smallskip Notons que jusqu'à maintenant, nous sommes restés assez secs concernant certaines des \prts désirables énoncées en \ref{prptaCrcdesirable}:
en effet nous n'avons pas été capables d'énoncer correctement les principes de prolongement par continuité ou les principes de recollement avec la généralité suffisante. Nous pouvions parler de continuité uniforme uniquement depuis l'extérieur de la \tdy. En effet, la continuité uniforme réclame une alternance de \qtfs du type $\forall n \exists m \forall x,y$
ce qui nécessite à priori de sortir du cadre des \tgms. 
Cela tient aussi~à ce que nous n'avions pas de sorte des \fsagcs.

Nous essayons dans cette partie de réparer ce manque. Et nous devons nous rappeler que d'un point de vue \cof, une fonction continue sur un compact n'existe pas sans un module de continuité uniforme. Le pari que nous faisons ici est d'internaliser la question de la continuité uniforme. Cela fait que, pour le moment, nous restons dans un cadre de \tdy finitaire. 

En outre, le cadre élargi que nous proposons avec l'introduction de ces nouvelles sortes semble être un cadre correct pour aborder un traitement \cof des structures o-minimales.

\smallskip Voici maintenant une brève description du contenu de la troisième partie. 

\smallskip Le chapitre \ref{chapomin} rappelle les \prts fascinantes des structures o-minimales en \clama. Il s'agit de \prts de finitude exactement semblables à celles de la géométrie algébrique des corps réels clos discrets, et cependant dépourvues de caractère \algq par l'utilisation du test de signe sur les nombres réels dans la théorie classique. Construire une théorie \algq des structures o-minimales est un défi crucial dans le \gui{programme de Hilbert \cof} qui vise à mettre à jour les constructions cachées dans les \clama contemporaines et à reformuler les \thos purement idéaux en énoncés \cofs. Ce programme permet d'éviter le recours à la \tfo \sa{ZFC} qui décrit un univers ensembliste hypothétique ne correspondant à aucune construction \mathe avérée. 

\smallskip Le chapitre \ref{chapfnfrb} propose une première \tdy finitaire pour des sortes décrivant les fonctions réelles uniformément continues valeurs dans $\II_\RR$.

\smallskip Le chapitre \ref{chapfnbg} donne un cadre \gnl pour décrire les \prts des fonctions uniformément continues définies sur $\II_\RR^{\,m}$ à valeurs dans $\II_\RR^{\,n}$. Un aspect décisif est de prendre en compte le fait qu'un module de continuité uniforme d'une fonction $f$ peut être vu comme une autre fonction uniformément continue $g$ attachée à la fonction $f$.

\smallskip Le chapitre \ref{chapaxiomesomin} propose de nouveaux axiomes
qui sont à priori satisfaits pour la géométrie algébrique des corps réels clos et qui semblent décisifs pour aborder une théorie \cov hypothétique et hautement souhaitable des structures o-minimales. Nous sommes néanmoins très loin d'avoir formalisé en une \tdy ce qui serait une version \cov des structures o-minimales.

\newpage \thispagestyle{empty}

%%%%%%%%%%%%%%%%%%%%%%%%%%%%%%%%%%%%%%%%%

\chapter{Structures o-minimales} \label{chapomin}
\Today
\minitoc

\Subsection{Parties définissables}

Références \cite{cos99}, \cite[\hbox{1998}]{vdD}, \cite{vdD86,DD88}. 
%\cite[\hbox{2017}]{ADH2017}.

\smallskip Une \gui{structure o-minimale sur un corps réel clos $\gR$} en \clama est donnée par une collection $(S_n)_{n\in\NN}$, où chaque $S_n$ est un ensemble de parties de $\gR^n$.

Les \elts de $S_n$ sont appelés les \textsl{parties définissables} de la structure o-minimale considérée.

Comme axiomes définissant la notion de structure o-minimale, nous demandons les \prts de stabilité suivantes.
\begin{enumerate}
\item Les sous-ensembles semi-algébriques de $\gR^n$ sont dans $S_n$.
\item Chaque $S_n$ est une \agB d'ensembles (stabilité par intersection et réunion finies, et passage au complémentaire).
\item Si $A\in S_n$ et $B\in S_m$ alors $A\times B\in S_{m+n}$.
\item Si $A\in S_{n+1}$ et $p_n:\gR^{n+1}\to\gR^n$ est la projection sur le premier sous-espace $\gR^n$ de coordonnées (oubli de la dernière coordonnée), alors $p_n(A)\in S_n$.
\item Les \elts de $S_1$ sont des réunions finies d'intervalles ouverts définissables et de points définissables.
\end{enumerate}

\Subsection{Applications définissables, résultats marquants}
Une application $A\to B$ entre ensembles définissables est dite \textsl{définissable} si son graphe est définissable. 

\smallskip Rappelons quelques résultats marquants.

\smallskip 
\begin{itemize}
\labu Le domaine de \dfn et l'ensemble image d'une application définissable sont définissables.
\labu La composée de deux fonctions définissables est définissable. 
\labu Toute partie définissable est une combinaison booléenne de parties fermées définissables.
Plus \prmt on a une décomposition cylindrique définissable de $\gR^n$
adaptée à toute famille finie de parties définissables (de manière analogue à la CAD dans le cas des parties semi-\agqs pour un \crcd). Les cellules de la décomposition sont homéomorphes à des simplexes ouverts, avec des homéomorphismes définissables. 
\labu Si $A\in S_n$ est fermé (pour la distance euclidienne de $\gR^n$) non vide, alors la fonction \gui{distance à~$A$} 
$$
d_A:\gR^n\to \gR,\,x\mapsto \inf\nolimits_{y\in A}\norme{x-y}
$$ 
est (continue et) définissable.
\labu Si $f\colon \gR^n\to \gR$ est continue est définissable, les zéros de $f$ forment une partie fermée définissable. Inversement, d'après le point précédent, toute partie fermée définissable de $\gR^n$ est l'ensemble des zéros d'une fonction continue définissable. 
\labu Si $I=\;]\,a,b\,[\;\subseteq \gR$ (avec $a,b\in\gR_\infty:=\gR\cup\so{-\infty,+\infty}$) et si $f\colon I\to\gR$ est une fonction définissable, alors
\begin{itemize}
\item il y a une subdivision de $I$ 
$$a=a_0<a_1<\dots<a_k=b
$$ 
telle que sur chaque 
intervalle ouvert de la subdivision, $f$ est ou bien constante, ou bien strictement monotone et continue,
\item on a aussi une subdivision telle que sur chaque intervalle ouvert de la subdivision, $f$ est dérivable avec la dérivée définissable, continue et de signe constant ($=0$ ou $>0$ ou $<0$). 
\end{itemize}
\labu Si $A\in S_n$ est un fermé définissable
 et $f\colon A\to \gR$ est continue définissable, elle peut être prolongée en une fonction continue définissable sur $\gR^n$ tout entier.
\labu Toute fonction continue définissable $\;]\,-1,1\,[\;\to\;]\,-1,1\,[ $ se prolonge par continuité en une fonction continue définissable $[-1,1]\to [-1,1]$. 
\labu Toute fonction continue définissable $[-1,1]\to [-1,1]$ atteint ses bornes. 
\labu Si $f\colon [-1,1]^{n+1}\to \gR$ est continue et définissable, la fonction $g\colon [-1,1]^{n}\to \gR$ définie par 
\[
g(\xn):=\sup\nolimits_{y\in[-1,1]}f(\xn,y)
\] 
est continue et définissable.\\
Notons qu'en particulier si $f$ est partout $\leq 0$ et si $A$ est l'ensemble des zéros de $f$, alors $p_n(A)$ est l'ensemble des zéros de $g$. 
Si $A$ est un fermé définissable $\subseteq [-1,1]^{n+1}$, on peut prendre pour $f$ la fonction $-d_A\colon[-1,1]^{n+1}\to \gR$. 
\end{itemize}

\Subsection{Variante} 
Tout ceci implique que l'on pourrait aussi bien définir la structure o-minimale sur $\gR$ par la donnée des objets suivants.
\begin{enumerate}
\item Les fonctions continues définissables 
$[-1,1]^{n}\to [-1,1]$.
\item La bijection croissante bicontinue (définissable dans toute structure o-minimale) 
\[
]-1,1\,[\,\to \gR,\;x\mapsto x/(1-x^2)
\] 
et la bijection réciproque
\[
\gR\to \,]-1,1\,[\,,\;x\mapsto \big(\sqrt{4x^2+1}-1\big)/2x
\]
\end{enumerate}

En effet, via le codage donné au point 2, pour avoir les fonctions continues définissables de~$\gR^n$ vers $\gR$ il nous suffit de savoir
décrire les fonctions continues définissables $f\colon \,]-1,1\,[\,^{n}\to\,]-1,1\,[$.
Et pour cela il suffit de savoir
décrire les fonctions continues définissables $g\colon [-1,1]^n\to[-1,1]$.

En effet notons $\norme{x}=\sup_{i\in\lrbn}\abs{x_i}$ pour $x=(\xn)\in\gR^n$.

Dans le cas où la croissance à l'infini de toute fonction définissable $f$ de~$\gR^n$ vers $\gR$ est majorée par un \pol, pour une telle fonction $f$, on a une fonction continue définissable $g\colon [-1,1]^{n}\to [-1,1]$ écrite sous la forme $g(x)=(1-\norme{x}^2)^kf(x)$, et la fonction $f$ peut être codée par le couple $(g,k)$. La fonction $g$ tend uniformément vers $0$ lorsque $x$ tend vers le bord de
$[-1,1]^n$.

Dans le cas général, on peut remplacer $h(x):=1-\norme{x}^2$ dans $g(x)$ par une fonction $\varphi\circ h$ où $\varphi\colon [0,1]\to [0,1]$ est continue définissable et strictement positive sur $[0, 1[\,$.

\newpage \thispagestyle{empty}

%%%%%%%%%%%%%%%%% CHAPTER %%%%%%%%%%%%%

\chapter{Anneaux de fonctions réelles bornées}
\label{SubsecIcrftr}
\label{chapfnfrb}
\Today
\minitoc

%\section*{Introduction}
%\addtocontents{toc}{\vskip0.8em}
%\addcontentsline{toc}{section}{Introduction}
%\rdb
%
%\fbox{À ÉCRIRE}

\section{Quelques rappels de la deuxième partie}

La \talg \Sa{Aftr} des \aftrs est la théorie des \afrrs à laquelle on ajoute le symbole de relation $\cdot >0$ comme une abréviation de \gui{$x\geq 0\vii \exists z\,xz=1$} et le symbole de fonction $\Fr$ avec les axiomes \Tsbf{fr1} et \Tsbf{fr2} (\dfn \ref{defiAftr}).

Un \aftr est donc une \QQlg \ftm réticulée réduite dans laquelle tout \elt supérieur à un \elt positif \iv est lui même \iv, et dans laquelle
la règle~\tsbf{FRAC} est valide.

Enfin, la \tdy \Sa{Co} des corps ordonnés \textsl{non} discrets
peut être décrite comme la théorie des \aftrs locaux, ce qui revient à ajouter l'axiome \Tsbf{OTF} (lemme~\ref{lemArftr}, point \textsl{3}) à la théorie \sa{Atfr}.

%: Lemma{lemIFR}
\begin{lemma} \label{lemIFR}\label{NOTAuplus}
Soit $\gA$ un \aftr. Notons $\gI=\sotq{x\in\gA}{-1\leq x\leq 1}$. Nous définissons sur $\gI$ la loi $x\uplus y=\frac 1 2 (x+y)$. La structure obtenue sur $\gI$ pour la signature 
\Sigt{\Icro}{\cdot= \cdot,\cdot\geq \cdot, \cdot>\cdot\mathrel{;}\cdot\,\uplus\,\cdot, \cdot\times \cdot, \cdot\vu\cdot, \mathrm{Fr}(\cdot,\cdot),\,- \cdot,0} \label{NOTASigIcro}
\noindent permet de reconstruire, de manière unique, la structure de $\gA$ comme \aftr.
\end{lemma}
%----------- fin lemma ----------------------------------- 
%
%
\begin{proof}
Cela tient pour l'essentiel à ce que tout \elt $z\in\gA$ peut s'écrire sous la forme $x/y$ avec $-1<x<1$ et $0<y<1$ (par exemple $x=\frac{z}{2+\abs z}$
et $y=\frac{1}{2+\abs z}$).
\end{proof}
%

%%%%%%%%%%%%%%%%%%%%%%%%%%%%%%%%%%%%%%%%%%%%%%%%%%%%%%%%%%%%%%%%%%%%
\section[Théorie dynamique des \afrbs]{Théorie dynamique des \afrbs}
\label{secAfrb}

Nous allons utiliser une signature plus fournie
qui correspond mieux à l'intuition d'un intervalle en tant que partie convexe.
\rdb

\smallskip Nous notons $x\oplus y$ la loi de composition suivante dans un \afr: $(x,y)\mt -1\vu(1\vi(x+y))$\footnote{Il s'agit de l'addition, remise dans l'intervalle $\gI$ si elle en sort.}.
\label{oplus}

\smallskip Nous notons $\mathrm{Cb}$ l'ensemble des \textsl{systèmes de \coes barycentriques}, défini \prmt comme suit: 
\[\ndsp
\mathrm{Cb}=\sotq{(r_k)_{k\in\lrbn}}{n\geq 2,\,r_1,\dots,r_n\in\QQ,\,r_1,\dots,r_n\geq 0, \,\sum_{k=1}^n r_k=1}. 
\]
Nous notons $\II_\QQ=\sotq{r\in\QQ}{-1\leq r\leq 1}$.

Pour chaque $\rho=(r_k)_{k\in\lrbn}$ dans $\mathrm{Cb}$, $\mathrm{Brc}_{\rho}$ est un symbole de fonction d'arité $n$ correspondant à la fonction: $\gI^n\to\gI,\;(x_k)_{k\in\lrbn}\mapsto\sum_{k=1}^nr_kx_k$.\label{notacb}

Le langage de la \tdy des \textsl{\afrbs} \SA{Afrb} est défini par la signature suivante. On a une seule sorte, notée $\AfrB$
\Sigt{\AfrB}{\cdot=0,\cdot\geq 0,\cdot>0\mathrel{;}(\mathrm{Brc}_{\rho})_{\rho\in\mathrm{Cb}}, \cdot\oplus\cdot, \cdot\times\cdot,-\,\cdot,\cdot\vu\cdot,\mathrm{Fr}(\cdot,\cdot), (r)_{r\in\II_\QQ}} \label{NOTASigAfrb}

\noindent {\bf Abréviations}

\smallskip \noindent \textsl{Symboles fonctionnels}

\vspace{-1em} \DeuxCols{
\begin{itemize}
\itbu $x\vi y$ signifie $ - (-x\vu -y)$
\itbu $\abs{x}$ signifie $x \vu -x$
\end{itemize}
}
{
\begin{itemize}
\itbu ${x}^+$ signifie $ x \vu 0$
\itbu ${x}^-$ signifie $ -x \vu 0$
\end{itemize}
}

\medskip \noindent \textsl{Prédicats}

\vspace{-1em} \DeuxCols{
\begin{itemize}
\itbu $x = y $ signifie $ x - y = 0$
\itbu $x \geq y $ signifie $ x - y \geq 0$
\itbu $x > y $ signifie $ x - y > 0$
\end{itemize}
}
{\begin{itemize}
\itbu $x \perp y $ signifie $ \abs x \vi \abs y =0$
\itbu $x \leq y $ signifie $ y\geq x$
\itbu $x < y $ signifie $ y> x$
\end{itemize}
}

\smallskip \noindent {\bf Axiomes}\label{AxiomesAfrb}

\smallskip\noindent 
\textsl{Les axiomes sont toutes les \rdys énoncées 
dans le langage de \sa{Afrb} qui sont valides
pour l'intervalle $\gI=[-1,1]$
dans la théorie $\Sa{Aftr}(\QQ)$ des \QQlgs \ftm réelles.} 

%l
%: Lemma{lemAxiomesAfrb}
\begin{lemma} \label{lemAxiomesAfrb}
Les \ralgs valides dans \sa{Afrb} sont décidables.
\end{lemma}
%----------- fin lemma ----------------------------------- 

%
\begin{proof}
Conséquence du point \textsl{3} du \corl \ref{corPst2}.
\end{proof}

Notons que l'on ne sait pas si les \rdys valides sont décidables. La même question se pose d'ailleurs dans le cas local pour la \tdy \sa{Co}. 
Cette question ne semble pas très importante dans la mesure où l'on s'intéresse essentiellement au cas des théories \sa{Crc1} et \sa{Crc2}, où le \pb reste mystérieux et s'ajoute à celui de savoir si l'on a bien capturé toutes les \prts \agqs du corps $\RR$. 

\section{La \tdy d'intervalle compact réel}\label{secIcr}

La \tdy \Sa{Icr} des \textsl{intervalles compacts réels}
possède une seule sorte, notée $\Icr$. Son langage est défini par la signature suivante.
\Sigt{\Icr}{\cdot=0,\cdot>0,\cdot\geq0\mathrel{;}(\mathrm{Brc}_{\rho})_{\rho\in\mathrm{Cb}}, \cdot\oplus\cdot, \cdot\times\cdot,-\,\cdot,\cdot\vu\cdot,(\rT_n)_{n\in\NN},\mathrm{Fr}(\cdot,\cdot), (r)_{r\in\II_\QQ}}\label{NOTASigIcr}

\noindent La \tdy \SA{Icr} 
 est obtenue à partir de la théorie \Sa{Afrb} décrite 
dans la section \ref{secAfrb} en ajoutant
\begin{itemize}
\item Les symboles de fonctions $\rT_n$ pour les \pols de Chebyshev, avec les axiomes;\\
 $\vd T_0(x) =1$, $\vd T_1(x) =x$, $\vd T_{n}(x) = 2xT_{n-1}(x) - T_{n-2}(x)$ ($n\geq 2$).\\ 
Pour les principales \prts des \pols de Chebyshev nous renvoyons à l'ouvrage
\cite[Chebyshev Polynomials]{MH2003}
\item L'axiome \Tsbf{OTF} (valide pour la structure de corps ordonné \textsl{non} discret) reformulé comme suit:

\regles{\lAb{OTF$'$} $\,\, x\oplus y>0 \vdi_{x,y\Icr}\; x>0 \;\;\vou\;\;y>0$}
 
\end{itemize}

%: Question{questIcr}
\begin{remark} \label{questIcr}
%Est-ce que les théories \Sa{Afrb} et 
Les théories \Sa{Icr} et \Sa{Co} sont probablement \esids. Sinon il serait \ncr d'ajouter des axiomes à \sa{Icr} pour rendre la chose vraie.
\eoe
\end{remark}
%----------- fin question ----------------------------- 

\newpage \thispagestyle{empty}

%%%%%%%%%%%%%%%%% CHAPTER %%%%%%%%%%%%%

\stepcounter{chapter}

\chapter{Un langage renforcé et les premiers axiomes correspondants}
\label{chapfnbg}
\Today
\minitoc

\section*{Introduction}
\addtocontents{toc}{\vskip0.8em}
\addcontentsline{toc}{section}{Introduction}
\rdb

Nous introduisons maintenant les sortes des \fsagcs en vue d'obtenir une \tdy plus expressive que \Sa{Crc2} pour les \crcs \textsl{non} discrets.

Cette nouvelle \tdy, que nous baptiserons \sa{Crc3} 
 essaye ici de récapituler ce que l'on est droit d'attendre d'une structure o-minimale pour les fonctions définissables uniformément continues sur $\II_\RR$.
 
 Comme nous l'avons déjà indiqué, nous nous limitons aux fonctions uniformément continues bornées, un peu dans le même esprit que Bishop.

%%%%%%%%%%%%%%%%%%%%%%%%%%%%%%%%%%%%%%%%%%%%%%%%%%%%%%%%%%%%%%%%%%%%
% section{Les sortes du langage renforcé}
\section{Les sortes du langage renforcé}
\begin{enumerate}
\item La sorte $\Icr$, pour l'intervalle compact $\gI=[-1,1]$.
\item Pour chaque $m\geq 0,n\geq 1$, une sorte $\Df_{m,n}$ pour les
fonctions définissables uniformément continues\footnote{Fonctions \sagcs pour la théorie des corps réels clos.} $\gI^m\to \gI^n$, la sorte $\Df_{m,1}$ est notée~$\Df_m$. En particulier $\Icr=\Df_0=\Df_{0,1}$.
\item Une sorte $\Li_n$ vue comme une sous-sorte de $\Df_n$, pour certaines fonctions lisses données au départ (au moins les \pols de Chebyshev).
\item Une sorte $\Mc$ pour les modules de continuité uniformes.
On les voit comme des objets particuliers de type $\Df_1$.
\item Une sorte $\Dfmc_n$ pour des couples formés par un objet de sorte $\Df_n$ et par un module de continuité uniforme qui lui convient.
\end{enumerate}

%%%%%%%%%%%%%%%%%%%%%%%%%%%%%%%%%%%%%%%%%%%%%%%%%%%%%%%%%%%%%%%%%%%%
\section{Un principe d'abstraction}
%: Subsection Un principe d'abstraction

Pour tout terme $t(\xn)$ de type $\Icr^n\to \Icr$ de la théorie (où les $x_i$ couvrent toutes les variables libres présentes dans le terme),
terme qui fournit une fonction $\gI^n\to\gI$ dans un modèle, nous devons faire ce qu'il faut pour
qu'il existe un terme $\dot t$ dans $\Df_n$ qui \gui{s'évalue comme $t$}.
Autrement dit, on doit mettre en place ce qu'il faut pour mimer, au sein de notre \tgm, la $\lambda$-abstraction du $\lambda$-calcul.

Il faut pour cela avoir dans la signature
\begin{itemize}
\item des symboles de type $\Df_n\times\, \Icr^n\to\,\Icr$ pour l'évaluation
d'un $u:\Df^n$ en des $x_i:\Icr$;
\item des symboles pour la composition des fonctions (avec des axiomes adéquats); 
\item des symboles qui donnent un nom aux fonctions \gui{\elrs} données dans la signature (par exemple $\cdot\times \cdot$ doit avoir un nom comme objet de sorte $\Df_2$); 
\item \dots
\end{itemize}

\smallskip Cette démarche est indispensable pour pouvoir parler uniformément, et non pas seulement ponctuellement, des \prts des fonctions définissables continues\footnote{Ceci rappelle ce que fait Kleene quand il définit les fonctionelles (uniformément) primitives récursives.}. 

%: Remark{remsymbolesfonctions}
\begin{remark} \label{remsymbolesfonctions} 
 On pourrait penser que certains symboles de fonction introduits 
\textsl{à priori} pour mimer la $\lambda$-abstraction 
auraient pu être ajoutés \textsl{a posteriori} en vertu de la possibilité d'ajouter
un symbole de fonction en cas d'existence unique, ce qui fournit une \tdy \esid à la précédente. 
Mais l'existence (dans l'existence unique en question) d'une fonction bien définie de la sorte $\Icr\times \Icr$ vers la sorte $\Icr$ ne signifie pas l'existence d'un objet correspondant dans~$\Df_2$,
 ni même son unicité (car les axiomes d'extensionnalité introduits plus loin sont trop faibles). 
Ce que nous signifions en introduisant \textsl{à priori} ces fonctions comme objets de type $\Df_{m,n}$, c'est que toutes les fonctions suffisamment simples, en particulier celles décrites dans les signatures sont bien définissables continues. \eoe
\end{remark}
%----------- fin remark ---------------------------------- 

\section{Premières structures sur les sortes $\Df_{m,n}$}

\Subsection{Les sortes $\Df_{m,n}$}

La sorte $\Icr$ des intervalles compacts réels (\ftm réticulés)
est munie de la structure d'\icr telle que décrite 
dans la section \ref{secIcr}.

\smallskip Chaque sorte $\Df_{m}$ ($m\geq 0,n\geq 1$) est accompagnée des symboles de fonctions et de prédicats ainsi que des axiomes des \afrbs (\tdy \Sa{Afrb}). 

%\smallskip Chaque sorte $\Df_{m,n}$ ($m\geq 0,n\geq 2$) est accompagnée des symboles de fonctions et de prédicats ainsi que des axiomes des \afrs (\tdy \Sa{Afr}). 

\smallskip \rem L'axiome \tsbf{OTF'} \paref{secIcr} n'est pas valide pour les sortes $\Df_{m}=\Df_{m,1}$ pour $m\geq 1$.
\eoe

\paragraph{Identification de $\Df_{m,n}$ et $({\Df_{m}})^n$}~

\smallskip Pour chaque $i\in\lrbn$ on a un symbole de fonction $\pi_{m,n,i}$
 de type $\Df_{m,n}\to \Df_{m}$ correspondant à la $i$-ème \coo. On donne aussi un symbole de fonction de type $(\Df_{m})^n\to \Df_{m,n}$ pour la bijection; on le notera $(\varphi_1,\dots,\varphi_n)$ de manière certes un peu ambigüe. 
Avec les axiomes convenables cela permet d'identifier $\Df_{m,n}$ et $(\Df_{m})^n$.
 
 \Regles{ 
 \lab{ } 
$\vdi_{\varphi_1,\dots,\varphi_n:\Df_{m}}\; \varphi_i=\pi_{m,n,i}((\varphi_1,\dots,\varphi_n))$ \quad $(i\in\lrbn)$%
 \lab{ }
 $\vdi_{\varphi:\Df_{m,n}}\; \varphi=(\pi_{m,n,1}(\varphi), \dots,\pi_{m,n,n}(\varphi))$
 } 

On donne les axiomes qui disent que $\pi_{m,n,i}$ est un morphisme pour les structures d'\afrb de $\Df_{m,n}$ et $\Df_{m}$.
 
%%%%%%%%%%%%%%%%%%%%%%%%%%%%%%%%%%%%%%%%%%%%%%%%%%%%%%%%%%%%%%%%%%%%
\paragraph{Composition de fonctions}~

\smallskip On a des symboles de fonction $\rC_{m,n,p}$ de type $\Df_{n,p}\times \Df_{m,n}\to \Df_{m,p}$ correspondant à la composition des fonctions. Pour $\varphi:\Df_{m,n}$ et $\psi:\Df_{n,p}$, on notera
$\psi\circ \varphi:=\rC_{m,n,p}(\psi,\varphi)$. On a des constantes de type $\Df_{n,n}$ pour les fonctions \gui{identité} $\Id_n:\gI^n\to\gI^n$.

\smallskip On donne les axiomes pour l'associativité de la composition.

\smallskip On donne les axiomes qui disent que pour $\varphi$ fixé de type $\Df_{m,n}$, l'application $\psi\mt\psi\circ \varphi$ est un morphisme pour les structures d'\afrb
de~$\Df_{n,p}$ et~$\Df_{m,p}$.

\smallskip On note de manière abrégée $\rC_{n,m}$ le terme de type $\Df_n\times \,({\Df_m})^n\to\Df_m$, défini par $$\rC_{n,m}(\varphi,\eta_1,\dots,\eta_n)\eqdef \varphi\circ(\eta_1,\dots,\eta_n).$$

%%%%%%%%%%%%%%%%%%%%%%%%%%%%%%%%%%%%%%%%%%%%%%%%%%%%%%%%%%%%%%%%%%%%
\paragraph{Évaluation de fonctions}~

\smallskip Pour $\Icr=\Df_0$, le symbole de fonction $\rC_{n,0}$ de type $\Df_n\times \,\Icr^n\to \Icr$ définit l'évaluation d'une fonction en des variables prises dans $\gI$. L'associativité de la composition se relit alors naturellement comme suit avec des $x_i$ de type $\Icr$, $\varphi$ de type $\Df_n$ et $\psi$ de type $\Df_1$
\[
(\psi\circ \varphi)(\xn)=(\psi\circ \varphi)\circ (\xn)=\psi\circ (\varphi\circ (x_1,\dots,x_n))=\psi(\varphi(\xn)).
\]

%%%%%%%%%%%%%%%%%%%%%%%%%%%%%%%%%%%%%%%%%%%%%%%%%%%%%%%%%%%%%%%%%%%%
\paragraph{Fonctions constantes}~

\smallskip On a un symbole de fonction $\jmath_n=\jmath_{0,n}$ de type $\Icr\to\Df_n$ pour les fonctions constantes.

On donne les axiomes qui disent qu'il s'agit de morphismes pour les structures d'\afrb\footnote{Par exemple pour $r\in\II_\QQ$, un axiome dit que $\jmath_n(r)$ est égal, en tant qu'objet de type $\Df_n$, au~$r$ donné dans la structure d'\afrb.} et que l'évaluation d'une fonction constante en n'importe quels arguments est bien la constante voulue.

\smallskip Plus \gnlt, si $0\leq m<n$ on a un symbole de fonction $\jmath_{m,n}$ de type $\Df_m\to\Df_n$ pour les objets de type $\Df_n$ correspondant à des fonctions qui ne dépendent que des $m$ premières variables et qui peuvent donc 
s'exprimer à partir d'objets de type $\Df_m$.
Les axiomes sont analogues à ceux donnés pour le cas $m=0$.

%%%%%%%%%%%%%%%%%%%%%%%%%%%%%%%%%%%%%%%%%%%%%%%%%%%%%%%%%%%%%%%%%%%%
\paragraph{Réarrangement de variables}~

\smallskip Pour $m,n>0$ et une fonction $\kappa\colon \lrbm\to\lrbn$ on a un objet $\wi\kappa$ de type $\Df_{n,m}$
avec l'axiome 

\Regles{
\lab{c$_\kappa$} $\vdi_{\xn:\Icr}\; \wi\kappa(\xn)=(x_{\kappa_1},\dots,x_{\kappa_m})$
}

 On donne aussi les axiomes naturels associés: $\wi{\kappa\circ \tau}=\wi\kappa\circ \wi\tau$.
 
\smallskip On a donc pour $m<n$ l'\egt $\jmath_{m,n}(\varphi)=\varphi\circ \wi\kappa$, où $\kappa\colon \lrbm\to\lrbn$ vérifie $\kappa(i)=i$ \hbox{pour $i\in\lrbm$}. Cette \egt nous dispense d'ailleurs d'introduire le symbole $\jmath_{m,n}$.

\smallskip On a aussi, pour $\tau_{n,i}:\lst1 \to \lrbn$ défini par $\tau_{n,i}(1)=i$
et $\psi$ de sorte $\Df_{m,n}$, l'\egt \hbox{$\pi_{m,n,i}(\psi)=\wi{\tau_{n,i}}\circ \psi$}. 
%%%%%%%%%%%%%%%%%%%%%%%%%%%%%%%%%%%%%%%%%%%%%%%%%%%%%%%%%%%%%%%%%%%%

\Subsection{Recollement d'\elts de $\Df_1$ sur des intervalles consécutifs}

%%%%%%%%%%%%%%%%%%%%%%%%%%%%%%%%%%%%%%%%%%%%%%%%%%%%%%%%%%%%%%%%%%%%
\paragraph{Restriction d'un \elt de $\Df_1$ à un intervalle}~

\smallskip Si $f$ est de sorte $\Df_1$, on veut avoir un nom pour la fonction $g$
obtenue à partir de la restriction de $f$ à un intervalle $[a,b]\subseteq \gI$. 

On réalise cela au moyen d'un symbole fonction $\Rs:\Df_1 \times \Icr \times \Icr \to \Df_1$. 

Lorsque $a\leq b$, on prolonge $g$ avec $g(x)=f(b)$ si $x\geq b$ et $g(x)=f(a)$ si $x\leq a$. Lorsque $a\geq b$, on permute $a$ et $b$. On a donc les axiomes suivants

\Regles
{\labu $\vdi_{a,b:\Icr,f\colon \Df_1} \Rs(f,a,b)=\Rs(f,a\vi b,a\vu b)$
\labu $\,\,x\leq a\vi b\vdi_{a,b,x:\Icr,f\colon \Df_1} \Rs(f,a,b)(x)=f(a\vi b)$
\labu $\,\,x\geq a\vu b\vdi_{a,b,x:\Icr,f\colon \Df_1} \Rs(f,a,b)(x)=f(a\vu b)$
\labu $\,\,a\vi b \leq x\leq a\vu b\vdi_{a,b,x:\Icr,f\colon \Df_1} \Rs(f,a,b)(x))=f(x))$
}

\paragraph{Recollement}~ 

\smallskip Si $f_0,\dots,f_n$ sont de sorte $\Df_1$, et si $0\leq a_1\leq \dots\leq a_n\leq 1$ on veut avoir un nom pour la fonction qui recolle les $f_i$ restreintes à $[a_i,a_{i+1}]$, éventuellement décalées verticalement pour assurer la continuité. 

On réalise cela au moyen d'un symbole fonction $\Rc_n:(\Df_1)^{n+1} \times (\Icr)^{n} \to \Df_1$. 

On a les axiomes de base suivants (on pose $a_0=0$ et $a_{n+1}=1$)

\Regles{\lab{Rc$_{n,j}$} $\Vii_{i} (a_i\leq a_{i+1}, f_{i}(a_{i+1})=f_{i+1}(a_{i+1}))
\vdi_{a_i:\Icr,f_i:\Df_1} \Rs(\Rc(\uf,\ua),a_j,a_{j+1})=\Rs(f_j,a_j,a_{j+1})$ } 

On ajoute les axiomes adéquats pour forcer les hypothèses de \tsbf{Rc$_{n,j}$}. 

%r

\Subsection{Les axiomes d'extensionnalité faible}
%: Subsection{Les axiomes d'extensionnalité faible}

Pour chaque sorte $\Df_n$ avec $n\geq 1$ on a l'axiome d'extensionnalité faible suivant.

\UneRegle{EXT$_n$} {$\,\,a> 0 \vdi_{a:\Icr,\varphi:\Df_n}\; 
\abs {\varphi} < \jmath_n(a) \;\vou\; \Exists x\, \abs{\varphi(x)}>\frac {a} {2}$ }

Comme conséquence, une fonction qui est partout nulle est \gui{presque} nulle: elle est majorée en valeur absolue par toute constante $>0$.
Pour conclure à la nullité il faudrait invoquer \Tsbf{HOF}, axiome non \gmq dont nous ne voulons pas, ou un axiome d'archimédianité douteux tel 
que~\Tsbf{AR2} dans une \tgm infinitaire\footnote{Une solution, très bancale, à cette faiblesse de la \tdy serait de ne considérer comme modèles que ceux où les objets de type $\Df_n$ sont $\geq 0$ (resp. $>0$) exactement quand ils sont évalués $\geq 0$ (resp. $>0$) en tout point de $\gI$.}. 

\smallskip Notons que l'axiome \tsbf{EXT$_0$} dit simplement que pour $x,a:\Icr$ et $a>0$, on a $\abs x<a$ ou $\abs x>\frac a 2$, ce qui est une variante de~\Tsbf{OTF}.

\smallskip Notons enfin que la règle \tsbf{EXT$_n$} résulte de \tsbf{OTF} et des axiomes de borne supérieure dans la section \ref{AxBsup} (avec $m=0$).

\smallskip \hum{Peut-être il serait judicieux d'avoir un symbole de relation $\cdot \neq 0$ sur chaque sorte $\Df_n$ de sorte que $\varphi\neq 0$ soit une abréviation de $\Exists \ux\;\abs{\varphi(\ux)}>0$?}

\newpage \thispagestyle{empty}

%%%%%%%%%%%%%%%%% CHAPITRE %%%%%%%%%%%%%

\chapter{Des axiomes décisifs}\label{chapaxiomesomin}

\Today
\minitoc

%\section*{Introduction}
%\addtocontents{toc}{\vskip0.8em}
%\addcontentsline{toc}{section}{Introduction}
%\rdb

%%%%%%%%%%%%%%%%%%%%%%%%%%%%%%%%%%%%%%%%%%%%%%%%%%%%%%%%%%%%%%%%%%%%
%: Subsection{Les axiomes de borne supérieure}\label{AxBsup}
\section{Les axiomes de borne supérieure}\label{AxBsup}

Les axiomes de borne supérieure remplacent \textsl{à priori} l'axiome de projection pour les parties définissables dans les structures o-minimales.

\smallskip 
Pour $m>0$ on a un symbole de fonction
$\sup_{m}$ de type $\Df_{m}\to \Icr$ pour la borne supérieure.
Il satisfait les axiomes qui décrivent la borne supérieure, à savoir

\regles{
\lab{sup$^{\Df}_{m}$} $ \vdi_{\varphi:\Df_{m}}
\;\varphi\leq \jmath_m(\sup_{m}(\varphi)) $
\lab{SUP$^{\Df}_{m}$} $\,\,\epsilon>0\vdi_{\epsilon:\Icr;\varphi:\Df_{m}}
\;\Exists \uy \;\;\varphi(\uy)+\epsilon> \sup_{m}(\varphi)$
}

Plus \gnlt pour $n\geq 0$ et $m>0$ on a un symbole de fonction
$\sup_{m+n,n}$ de type $\Df_{m+n}\to \Df_{n}$ pour la borne supérieure sur les $m$ dernières variables (dans $\gI^m$) avec les axiomes suivants (donc $\sup_{m,0}$ n'est autre que $\sup_{m}$).

\regles{
\lab{sup$^{\Df}_{m+n,n}$} $\vdi_{\varphi:\Df_{n+m}}
\;\varphi\leq \jmath_{n,m+n}(\sup_{m+n,n}(\varphi))$
\lab{SUP$^{\Df}_{m+n,n}$} $\,\,\epsilon>0\vdi_{\epsilon,\xn:\Icr;\varphi:\Df_{n+m}}
\;\Exists \uy \;\;\varphi(\ux,\uy)+\epsilon> \sup_{m+n,n}(\varphi)(\ux)$
%\lab{SUp$_{m+n,n}$}
%$\,\,\varphi<0\vdi_{\varphi:\Df_{m+n}}\;\sup_{m+n,n}(\varphi)< 0$
}

%\rem \eoe 

\smallskip \hum{On ne doit pas mettre d'axiome disant que la borne supérieure est atteinte (c'est faux \cot même pour les \pols),
mais on aura des axiomes de finitude qui impliqueront qu'elle est presque toujours atteinte sous forme d'un sup fini de valeurs en des points que l'on calcule. Ceci règlera peut-être la question de la borne supérieure $<0$
si la fonction est $<0$?
}

%%%%%%%%%%%%%%%%%%%%%%%%%%%%%%%%%%%%%%%%%%%%%%%%%%%%%%%%%%%%%%%%%%%%
\section{Les axiomes de continuité uniforme}
%: Subsection{Les axiomes de continuité uniforme}

On explique maintenant comment un système d'axiomes convenable peut traduire le fait que toute fonction définissable continue admet un module de continuité uniforme, tout en restant dans le cadre d'une \tgm. Cela est possible parce qu'une fonction définissable continue admet un module de continuité uniforme qui est lui-mêmes une fonctions fonction définissable particulière.
Les sortes $\Mc$ et $\Dfmc_n$ avec leurs axiomes sont ici cruciales.

%\footnote{On utilise ici la terminologie usuelle des \coma (Bishop). Ce qui est souvent appelé \gui{module de continuité uniforme} dans les textes de \clama mériterait mieux le nom de \gui{module d'oscillation uniforme}.}

\smallskip On commence par donner un symbole de fonction $\jmath_\Mc$ pour une injection de type $\Mc\to\Df_1$. Un axiome
précise que $\jmath_\Mc$ est injective.

\smallskip On a un prédicat $\mathrm{Mcu}_{n}$ sur $\Df_n\times\, \Mc$ qui exprime que $\mu$ est un module pour $\varphi$ au moyen de l'abréviation suivante.

\Regles
{
\labu $\,\,\mathrm{Mcu}_n(\varphi,\mu)$ est une abréviation
pour: $\mu \big(\abs{\varphi(\ux)-\varphi(\ux')}\big)\leq \norme{(\ux)-(\ux') }$
}

%\vspace{-.2em}
\noindent où $\norme{(\uz)}=\sup_i\abs{z_i}$. Ici l'inégalité semble écrite, sous forme de fonctions évaluées, avec des $x_i$ et~$x'_i$ de sorte~$\Icr$.
Mais en fait, \textsl{on doit lire cette inégalité comme liant deux objets de sorte $\Df_{2n}$}. Cela permet d'éviter l'utilisation du quantificateur universel sur les $x_i$ et $x'_j$
dans la \dfn de la continuité uniforme! Les \tdys ne permettent pas la création de nouvelles formules utilisant des quantificateurs universels, nous contournons la difficulté en mimant la $\lambda$-abstraction!

\medskip Les axiomes suivants précisent les contraintes sur les objets $\mu$ de sorte $\Mc$.

\DeuxRegles
{
\lab{Mc$_1$} $\,\, 0<b<c\vdi_{\mu:\Mc;b,c:\Icr}
\;0<\mu(b)< \mu(c)$
\lab{mc$_1$} $\vdi_{\mu:\Mc} 
\;\mathrm{Mcu}_1(\jmath_\Mc(\mu),\mu)$}
{\lab{Mc$_2$} $\,\, a\leq 0\vdi_{\mu:\Mc;a:\Icr}
\;\mu(a)=0$
}

\smallskip 
\rem Dans le cas où l'on considère uniquement les \fsagcs, \L ojasievicz nous assure que tout module de continuité uniforme peut être pris parmi les seules fonctions $\epsilon>0 \mt c \,\epsilon^n$ (avec $c>0$)
\eoe

\smallskip 
\rem L'axiome peu intuitif $\tsbf{mc}_1$ sera une règle valide si l'on demande dans un autre axiome que tout objet
de type $\Mc$ corresponde à une fonction convexe. 
\eoe

\smallskip La sorte $\Dfmc_n$ est définie comme sous-sorte de la sorte produit $\Df_n\times \Mc$. Elle est accompagnée de deux symboles de fonctions $\mathrm{df}_n$ et $\mathrm{mc}_n$, avec les axiomes adéquats, qui font qu'un objet de sorte 
$\Dfmc_n$ peut être considéré comme un couple d'objets $(\varphi,\mu)$ de sortes respectives\footnote{Il n'est pas \ncr de créer la sorte produit en tant que telle. L'axiome suivant sera suffisant: 
$\mathrm{df}_n(\theta)=\mathrm{df}_n(\theta'), \mathrm{mc}_n(\theta)=\mathrm{mc}_n(\theta')\vdi_{\theta,\theta':\Dfmc_n} \theta=\theta'$.}
$\Df_n$ et $\Mc$. 
L'axiome $\tsbf{mcu}_n$ dit comment la sous-sorte est définie: si $(\varphi,\mu)$ est de sorte
$\Dfmc_n$, alors le prédicat $\mathrm{Mcu}_n(\varphi,\mu)$
est satisfait. 

\regles
{\lab{mcu$_n$} $
\vdi_{\psi:\Dfmc_n}\;
\mathrm{Mcu}_n(\mathrm{df}_n(\psi), \mathrm{mc}_n(\psi))$
}

\smallskip Enfin, on a l'axiome $\tsbf{DFMC}_n$ qui dit que tout objet $\varphi$ de sorte $\Df_n$ est l'image par $\mathrm{df}_n$ d'un objet $\theta$ de sorte
$\Dfmc_n$.

\regles{\lab{DFMC$_n$} $\vdi_{\varphi:\Df_n}\; \Exists \theta\;\; \mathrm{df}_n(\theta)=\varphi$} 

Avec toute cette machinerie on garantit de manière explicite la continuité uniforme des fonctions représentées par des objets de type $\Df_n$.
 
%r
%: Remark{remMcupartout}
\begin{remark} \label{remMcupartout} 
Chaque fois que nous avons introduit une constante de type $\Df_n$,
il fallait en fait introduire une constante \gui{au dessus d'elle} de type $\Dfmc_n$. Ceci n'offre pas de difficulté car dans chaque cas un module de continuité uniforme est évident.\eoe
\end{remark}
%----------- fin remark ---------------------------------- 

%%%%%%%%%%%%%%%%%%%%%%%%%%%%%%%%%%%%%%%%%%%%%%%%%%%%%%%%%%%%%%%%%%%%
\section{Les axiomes pour les fonctions lisses}
%: Subsection{Les axiomes de lissité}

Les objets de type $\Li_n$ sont vus comme des objets de type $\Df_n$
qui définissent certaines fonctions lisses (i.e.~$\Cin$).
La signature comporte un symbole de fonction $\jmath_{\Li_n}$ pour l'injection correspondante, avec les axiomes qui disent qu'il s'agit d'un morphisme injectif pour les lois convenables (celles qui conservent les fonctions lisses). 

\hum{Donner quelques précisions.}

On donne les fonctions constantes et les fonctions \coos comme objets de type $\Li_n$.

La sorte $\Li_1$ contient les \pols de Chebyshev. 

On peut envisager d'introduire d'autres fonctions Nash dans $\Li_n$, cela ne devrait pas changer la \tdy mais pourrait faciliter certaines \demos. 

Dans le cas où l'on vise une structure o-minimale particulière (autre que celle fournie par les \fsagcs), d'autres fonctions peuvent être données qui serviront de base pour la \dfn de la structure. 

\Subsection{Axiome de densité} 

Les zéros d'une fonction lisse non nulle (dans une structure o-minimale) forment un fermé $F$ d'intérieur vide. On peut exprimer (au moins partiellement) cette \prt de densité (pour le complémentaire de $F$) au moyen de l'axiome suivant

\Regles{
\lab{Dens$_n$}
$\,\,\abs{\varphi(\an)}>0,\; \varphi\times \psi=0 \vdi_{a_1,\dots,a_n:\Icr;\varphi:\Li_n;\psi:\Df_n}\;\psi=0
$
}

\Subsection{La dérivation}

Il nous semble confortable d'introduire la dérivée (ou la dérivée partielle) en suivant la \dfn de Bridger-Stolzenberg (voir \cite{AD1994} et \cite{BS1999}). Une fonction $\varphi\colon \II\to\RR$ est continument dérivable
si la fonction \gui{taux d'accroissement} peut être prolongée par continuité, i.e. s'il existe une fonction uniformément continue $\psi\colon \II^2\to\RR$ satisfaisant l'identité \fbox{$\varphi(x_1)-\varphi(x_2)=\psi(x_1,x_2)\times (x_1-x_2)$}.
La dérivée de $\varphi$ est alors donnée par $\varphi'(x)=\psi(x,x)$.

Comme nous ne voulons que des fonctions de $\gI^n$ dans $\gI$, nous
devons utiliser un codage implicite $(x,p)$ avec $x\in \gI$ et $p\in \NN$
pour le réel $px$. 

La fonction $\psi$ est uniquement déterminée par $\varphi$ (voir ci-après la règle valide \Tsbf{Der}) de sorte que
dans notre \tdy nous pouvons l'introduire au moyen d'un symbole 
 de fonction $\Delta=\Delta_{1,1}$ de type $\Li_1\to\Li_2$
qui satisfait l'axiome

\Regles{
\Lab{der}
$\vdi_{\varphi:\Li_1}\;\varphi(x_1)-\varphi(x_2)=\Delta(\varphi)(x_1,x_2)\times (x_1-x_2)$
}

\noindent \rem Cette \egt semble écrite sous forme de fonctions évaluées comme liant deux objets de sorte~$\rR$, mais en fait \textsl{elle doit être lue comme liant deux objets de sorte $\Li_2$}, qui s'évaluent en $(x_1,x_2)$
sous la forme indiquée dans l'axiome tel qu'il semble écrit.

\smallskip \rem En fait nous devons utiliser le codage implicite auquel
nous avons fait allusion plus haut et la règle \tsbf{der} doit en fait être écrite
sous forme

\Regles{
\lab{der}
$\vdi_{\varphi:\Li_1}\;\frac 1 {2p} (\varphi(x_1)-\varphi(x_2))=\Delta_p(\varphi)(x_1,x_2)\times (x_1-x_2)$
}

Pour ne pas compliquer l'exposé nous faisons dans la suite comme si $\Delta(\varphi)$ était la vraie fonction \gui{taux d'accroissement}. 

\smallskip La règle d'unicité suivante découle de l'axiome \tsbf{Dens$_2$}: dans le premier membre on doit lire une \egt entre objets de type $\Df_2$ et la fonction lisse est $x_1-x_2$ vue comme \elt de $\Li_2$. 

\Regles{
\Lab{Der}
$\,\,\varphi(x_1)-\varphi(x_2)=\psi(x_1,x_2)\times (x_1-x_2) \vdi_{\varphi:\Li_1;\psi:\Df_2}\;\Delta(\varphi)=\psi$
}

\smallskip 
De même en plusieurs variables on demande les axiomes analogues pour chaque dérivée partielle. En particulier, pour $n\geq 2$ et $i\in\lrbn$ on a un symbole de fonction $\Delta_{n,i}$ de type $\Li_n\to\Li_{n+1}$
qui satisfait l'axiome

\Regles{
\lab{der$_{n,i}$} {\mathrigid 1.7mu
$\vdi_{\varphi:\Li_n}\;\varphi(x_1,\dots,x_i,\dots,x_n)-\varphi(x_1,\dots,x'_i,\dots,x_n)=\Delta_{n,i}(\varphi)(x_1,\dots,x_i,x'_i,\dots,x_n)\times (x_i-x_i')$}
}

\noindent \textsl{Cette \egt doit être lue comme liant deux objets de sorte $\Li_{n+1}$}.

\smallskip 
\rem On obtient, en utilisant l'axiome de borne supérieure, que les fonctions lisses sont lipschitziennes,
 ce qui donne un module de continuité uniforme particulièrement simple.
\eoe

\Subsection{Quels autres axiomes pour la dérivation?} 

Il faut examiner ici quels axiomes il est \ncr d'introduire correspondant aux \prts usuelles de la dérivation. La plupart de ces \prts devraient résulter de la \dfn (axiomes $\tsbf{der}_{n,i}$) et des axiomes $\tsbf{Dens}_{n}$.

%%%%%%%%%%%%%%%%%%%%%%%%%%%%%%%%%%%%%%%%%%%%%%%%%%%%%%%%%%%%%%%%%%%%
\Subsection{Les axiomes de racines virtuelles}

À priori on peut définir les \ravs pour toute fonction lisse dont la dérivée d'ordre~$k$ est $>0$ (sur $\II$), ceci en vertu du lemme \ref{lemBasicVirtualRoots} et de la version \cov uniforme du \tho des accroissements finis. On obtient alors l'essentiel des \prts décrites dans la \dfn \ref{prdfVirtualRoots} et le \thref{thVirtualRoots}. Le \pol $f(X) = X^{d} - ( a_{d-1} X^{d-1} + \cdots +a_1X+ a_0)$ qui dépend de $d+1$ variables peut être remplacé par n'importe quelle fonction $\varphi$ lisse de $d+1$ variables $X,a_1,\dots,a_{d}$ dont la dérivée partielle $k$-ème par rapport à $X$ est $>0$ en tant qu'objet de type $\Li_{d+1}$.

Si $\inf(\varphi^{(k)})=\phi(a_1,\dots,a_{d})$, on peut traiter la fonction 
$\psi_k=\varphi+(c-\phi)^+X^k/k!$, pour une constante $c>0$.
Sa dérivée $k$-ème par rapport à $X$ est $\geq c$, et elle est égale à $\varphi$ si $\phi\geq c$. On peut alors introduire les~$k$ \ravs de $\varphi$ sur $\gI$ en tant qu'objets de sorte $\Df_d$ comme dans la \dfn \ref{prdfVirtualRoots} et le \thref{thVirtualRoots}, mais en utilisant notre $\lambda$-abstraction. Plus \prmt, on a des symboles de fonction \gui{racines virtuelles} $\mathrm{Rv}_{d,k,j}$ de type $\Li_{d+1}\to\Df_d$. Et on a les axiomes correspondants traductions directes de la \dfn \ref{prdfVirtualRoots} et du \thref{thVirtualRoots} (en remplaçant $-\infty$ et $+\infty$ par $-1$ et $+1$). 

\smallskip \hum{Écrire les détails.}

%%%%%%%%%%%%%%%%%%%%%%%%%%%%%%%%%%%%%%%%%%%%%%%%%%%%%%%%%%%%%%%%%%%%
\section{Les axiomes de clôture réelle ou o-minimale}
%: Subsection{Les axiomes de clôture réelle}

Nous traitons à partir de maintenant des axiomes qui correspondent à l'idée \gnle de
clôture réelle et de structure o-minimale.
%%%%%%%%%%%%%%%%%%%%%%%%%%%%%%%%%%%%%%%%%%%%%%%%%%%%%%%%%%%%%%%%%%%%

%%%%%%%%%%%%%%%%%%%%%%%%%%%%%%%%%%%%%%%%%%%%%%%%%%%%%%%%%%%%%%%%%%%%
\Subsection{Les axiomes de finitude}
Les axiomes de racines virtuelles sont déjà des axiomes de finitude, mais indépendants de toute structure o-minimale. 

On devrait avoir un analogue à la proposition \ref{factTableaucomplet} (tableau de signes et de variations) pour les \fsagcs $\II\to\II$, et cela devrait aussi fonctionner pour les structures o-minimales. En \clama les tableaux de signes et de variations existent pour les fonctions définissables $\II\to\II$ d'une structure o-minimale, et la proposition \ref{factTableaucomplet} montre comment transformer l'énoncé classique en un énoncé \cof. 
Ici aussi le \pb est celui de formuler des axiomes dynamiques qui capturent ce type de résultat. 
Une solution serait d'avoir une \tdy infinitaire avec des axiomes qui disent grosso modo qu'une fonction définissable continue est \gui{lisse monotone par morceaux} dans un énoncé à préciser, semblable au point \textsl{2} 
de la proposition~\ref{factTableaucomplet}. 

%%%%%%%%%%%%%%%%%%%%%%%%%%%%%%%%%%%%%%%%%%%%%%%%%%%%%%%%%%%%%%%%%%%%
\Subsection{Recollement de fonctions définies sur un recouvrement ouvert}

Un recouvrement fini de $\gI^n$ par des ouverts définissables est donné ici sous la forme 
$$V_i=\sotq{\ux\in\gI^n}{g_i(\ux)>0}\quad i\in\lrbp$$ 
où les $g_i$ sont de sorte $\Df_n$ et satisfont \fbox{$\sum_{i}g_i^+>0$ (1)}. Des fonctions $h_i$ de sorte $\Df_n$ sont considérées, dont à priori seules les restrictions $h_i\frt{V_i}$ sont pertinentes. Le fait que $h_i$ et $h_j$ coïncident sur $V_i\cap V_j$
se traduit par l'\egt \fbox{$h_ig_i^+g_j^+=h_jg_j^+g_i^+$ (2)}. Sous les hypothèses (1) et (2) on demande l'existence et l'unicité d'une $f$ de sorte $\Df_n$
vérifiant $fg_i^+=h_ig_i^+$ pour chaque $i$ (ce qui signifie que $f\frt{V_i}=h_i\frt{V_i}$). À priori on doit avoir $f=(\sum_{i}h_ig_i^+)(\sum_{i}g_i^+)^{-1}$ (d'où l'unicité). Et l'on obtient bien
$$fg_k^+=\frac{\sum_{i}h_ig_i^+g_k^+}{\sum_{i}g_i^+}=
\frac{\sum_{i}h_kg_i^+g_k^+}{\sum_{i}g_i^+}=\frac{(h_kg_k^+)\,\sum_{i}g_i^+}{\sum_{i}g_i^+}=h_kg_k^+.
$$

%%%%%%%%%%%%%%%%%%%%%%%%%%%%%%%%%%%%%%%%%%%%%%%%%%%%%%%%%%%%%%%%%%%%
\Subsection{Recollement de fonctions définies sur un recouvrement fermé}

Un recouvrement fini de $\gI^n$ par des fermés définissables est donné ici sous la forme 
$$F_i=\sotq{\ux\in\gI^n}{g_i(\ux)\geq 0}\quad i\in\lrbp$$ 
où les $g_i$ sont de sorte $\Df_n$ et satisfont \fbox{$\sup_{i}g_i\geq 0$}. Des fonctions $h_i$ de sorte $\Df_n$ sont considérées, dont à priori seules les restrictions $h_i\frt{F_i}$ sont pertinentes. Le fait que $h_i$ et $h_j$ coïncident sur $F_i\cap F_j$
se traduit par la validité des règles ($i,j\in\lrbp$)

\Regles{
\lab{} $\,\,g_i(\ux)\geq 0\vet g_j(\ux)\geq 0\vdi_{\xn:\Icr;g_i,h_i,g_j,h_j:\Df_n)}\;h_i(\ux)=h_j(\ux) $
}

Une version uniforme \agq de cette validité peut être énoncé comme suit

\Regles{
\labu $\vdi_{g_i,h_i,q_{ij},q_{ji}:\Df_n} \;({h_i-h_j})^2+g_iq_{ij}^++g_jq_{ji}^+=0$
}

\noindent où les $q_{k\ell}$ sont de sorte $\Df_n$. 
Abrégeons le second membre sous la forme $E_{ij}$. Sous l'hypothèse des \egts $E_{ij}$, on veut avoir une fonction $f$ (un objet $f$ de sorte $\Df_n$) satisfaisant une \idt qui signifie que $f\frt{F_i}=h_i\frt{F_i}$.
Ceci peut s'exprimer sous la forme de la règle suivante

\Regles{
\lab{RCVF}$\,\,\sup_ig_i\geq 0\vet E_{1,2}\vet\dots\vet E_{p-1,p}\vdi_{g_i,h_i,q_{ij},q_{ji}:\Df_n}\;\Exists f,q_1,\dots,q_p\;\;\Vi_i(f-h_i)^2+g_iq_i^+=0$
}

\noindent Toutes les variables (libres ou muettes) dans cette règle sont de sorte $\Df_n$.

En \clama ce type de règle est valide dans le cadre des structures o-minimales.
Il se peut cependant que d'un point de vue \cof il faille se limiter aux recouvrements par des \textsl{fermés situés}\footnote{Un fermé est dit \textsl{situé} lorsque la fonction distance au fermé est bien définie d'un point de vue \cof. Il semble nécessaire d'ajouter les axiomes qui disent que la fonction distance aux zéros d'une fonction continue définissable est elle-même définissable.}. Ce qui compliquera l'écriture des axiomes.

Notons que l'objet $f$ dont on postule l'existence est prouvablement unique en vertu d'un calcul classique pour les Positivstellensätze: on utilise l'\idt $(a+b)^2+(a-b)^2=2(a^2+b^2)$.
%%%%%%%%%%%%%%%%%%%%%%%%%%%%%%%%%%%%%%%%%%%%%%%%%%%%%%%%%%%%%%%%%%%%
\Subsection{Les axiomes de prolongement par continuité}

 Typiquement, les règles \FRACn\ sont des axiomes de 
prolongement par continuité particuliers.
On cherche ici à énoncer des règles différentes qui s'appliquent de manière plus \gnle (sans que le prolongement par continuité donne $0$ aux valeurs litigieuses) mais avec un dénominateur lisse.

Par exemple une fonction définissable en dehors des zéros d'une fonction lisse (non nulle) et continue sur son domaine de \dfn se prolonge par continuité de manière unique si elle est uniformément continue.

Tout le \pb est de formuler ceci dans le cadre de notre \tdy.

\smallskip Ce sera déjà bien de réussir à le formuler pour un quotient $f/g$ (bien défini en dehors des zéros de $g$) avec $g$ lisse.

Le fait que $f$ s'annule aux zéros de $g$ peut être mis comme hypothèse sous la forme forte suivante: il existe un $\alpha$ de sorte $\Mc$ tel que \fbox{$\alpha\circ \abs f\leq \abs g$}.

La continuité uniforme de $f/g$ en dehors des zéros de $g$ semble pouvoir s'énoncer en utilisant la bijection réciproque d'un objet de sorte $\Mc$ sur l'intervalle $[0,1]$. En effet on désire écrire quelque chose comme
$$
\mu\left(\abs{\frac{f(\ux)g(\ux')-f(\ux')g(\ux)}{g(\ux)g(\ux')}}\right)\leq \norme
{(\ux)-(\ux')}
$$
pour $\norme{(\ux)-(\ux')}>0$, qui pourrait se réécrire 
sans l'hypothèse $\norme{(\ux)-(\ux')}>0$ dans le cadre de la \tgm sous forme
$$
\abs{f(\ux)g(\ux')-f(\ux')g(\ux)}\leq \abs{g(\ux)g(\ux')}\,\nu(\norme{(\ux)-(\ux')})
$$
avec la bijection réciproque $\nu$ (sur $[0,\mu(1)]$) de $\mu\frt{[0,1]}$.

\smallskip \hum{écrire les détails}

%%%%%%%%%%%%%%%%%%%%%%%%%%%%%%%%%%%%%%%%%%%%%%%%%%%%%%%%%%%%%%%%%%%%
\Subsection{Conclusion: la structure améliorée de \crc}
%: Subsection{Conclusion: la structure améliorée de \crc}

La théorie \SA{Crc3} sera obtenue une fois qu'on aura mis au point tous les axiomes. On a vu que la théorie \Sa{Icr} pouvait être considérée comme une variante de la théorie \Sa{Co}. La théorie \sa{Crc3}, que l'on pourrait aussi appeler théorie des \icrcs, est une variante améliorée de \Sa{Crc2}, dans laquelle les structures o-minimales (qui sont des structures enrichies de \crcs) pourraient avoir droit de cité en tant que \sads particulières.

%t
%: Theorem{thicrc1}
\begin{theorem} \label{thicrc1}
En \coma, l'intervalle réel $\II=\II_\RR$ et les \fsagcs ${\II}^n\to\II$, fournit un modèle de la \tdy \sa{Crc3}.
\end{theorem}
%----------- fin theorem ----------------------------- 
%
\begin{proof}
{\tt Il semble que la \dfn ad hoc des \fsagcs adoptée en \ref{defiFSAGC2+}
ramène ce \tho à un \tho concernant pour l'essentiel $\RRa$. Mais il faut vérifier tous les détails et cela nous amènera peut-être à modifier la formulation de certains axiomes.}
\end{proof}

Cette théorie \sa{Crc3}
devrait permettre de démontrer des résultats \cofs qui échappent à la théorie plus \elr \Sa{Crc2} pour la simple raison qu'ils ne correspondent pas à des \rdys de \sa{Crc2}. Par ailleurs la question se pose pour les \rdys de \sa{Crc2} elles-mêmes.

%%%%%%%%%%%%%%%%%%%%%%%%%%%%%%%%%%%%%%%%%%%%%%%%%%%%%%%%%%%%%%%%%%%%
\section{Structures o-minimales}

Il semble que les axiomes proposés ici pour la structure d'\icrc soient 
presque corrects pour décrire \cot certaines structures o-minimales définies en \clama: celles engendrées par les restrictions au cube compact $\II^n$ de certaines fonctions lisses au voisinage de $\II^n$. 

Le point le plus faible semble la stabilité par projection. à priori le système d'axiomes actuel ne garantit cette stabilité que pour les parties définissables fermées bornées. 

La structure obtenue dépend des fonctions lisses données au départ dans les sortes $\Li_n$.

Nous sommes prioritairement intéressés par la structure obtenue en prenant pour fonctions lisses de départ les fonctions analytiques réelles au voisinage du cube. En dimension 1, cela se passe sans doute agréablement avec les séries de Chebyshev. 

C'est un véritable enjeu que de donner une version constructive de la théorie classique, par exemple en partant des exposés donnés dans \cite{vdD86} et \cite{DD88}. Il faudrait au moins démontrer \cot que les fonctions analytiques réelles au voisinage du cube engendrent une structure qui est un modèle de la \tdy \Sa{Crc3}.

Notons aussi que du point de vue proprement calculatoire, nous sommes à priori plus intéressés par le corps énumérable $\RR_{\tsbf{PR}}$ des réels calculables en temps primitif récursif, ou par le corps énumérable $\RR_{\tsbf{Ptime}}$ des réels calculables en temps polynomial (voir l'exemple \ref{exacorpsnondiscret}). Concernant
les fonctions continues définissables correspondant à ces corps (pour une structure o-minimale fixée), elles peuvent sans doute être elles aussi énumérées en utilisant les séries de Chebyshev. 

Enfin, soulignons qu'en l'état actuel le système d'axiomes envisagé ne semble pas suffire à décrire vraiment les structures o-minimales, car il ne garantit la stabilité par projection que pour les parties définissables fermées bornées. 

%En l'absence des axiomes de finitude, la théorie obtenue peut êre considérée comme une première approche pour une version \gui{\tdy}

%%%%%%%%%%%%%%%%%%%%%%%%%%%%%%%%%%%%%%%%%%%%%%%%%%%%%%%%%%%%%%%%%%%%
\section{Quelques questions}

%q
%: Question{questRRmodele}
\begin{question} \label{questRRmodele}~

\noindent Est-ce que la théorie \Sa{Crc3} prouve plus de \rdys que la théorie de l'intervalle $[-1,1]$ pour un \crc décrit par \Sa{Crc1}, ou par \Sa{Crc2}?
\end{question}
%----------- fin question ----------------------------- 
%q
%: Question{questRROmin}
\begin{questions} \label{questRROmin}~

\noindent 1) 
Est-ce que $\RR$ (pour la sorte $\rR$), avec $\II=\II_\RR$ (pour la sorte $\Icr$) et pour $\Li_n$ les fonctions ${\II}^n\to\II$ qui sont analytiques dans un voisinage de ${\II}^n$, \gui{engendre} un modèle de la \tdy \Sa{Crc3}? 

\smallskip \noindent 2)
Si oui, est-ce que les objets de type $\Df_n$ correspondent exactement aux \elts de l'\aftr engendré par les
fonctions associées aux objets de type $\Li_n$.
\end{questions}
%----------- fin question ----------------------------- 
%%%%%%%%%%%%%%%%%%%%%%%%%%%%%%%%%%%%%%%%%%%%%%%%%%%%%%%%%%%%%%%%%%%%
%%%%%%%%%%%%%%%%%%%%%%%%%%%%%%%%%%%%%%%%%%%%%%%%%%%%%%%%%%%%%%%%%%%%
%%%%%%%%%%%%%%%%%%%%%%%%%%%%%%%%%%%%%%%%%%%%%%%%%%%%%%%%%%%%%%%%%%%%

\newpage \thispagestyle{empty}

%%%%%%%%%%%%%%%%% CHAPITRE %%%%%%%%%%%%%

\chapter*{Conclusion générale}
\addstarredchapter{Conclusion générale}

Ce mémoire, et les questions sans réponse qu'il contient, permet de mesurer notre ignorance de l'\alg réelle. 
\rdb

%%%%%%%%%%%%%%%%%%%%%%%%%%%%%%%%%%%%%%%%%%%%%%%%%%%%%%%%%%%%%%%%%%%%%

\part*
{Références et index}
\addstarredpart{Références et index}
%:biber

\rdb
\addcontentsline{toc}{section}{Références. Livres}

\setlength\labelnumberwidth{4em}
\printbibliography[title=Références. Livres,keyword=book,resetnumbers=true]

\rdb
\addcontentsline{toc}{section}{Références. Articles}
\setlength\labelnumberwidth{2em}
\printbibliography[title=Références. Articles,keyword=paper,resetnumbers=true]

%:fin biblio

\newpage \thispagestyle{empty}

\normalsize

\catcode`\@=11
%\@makesindexhead{Index des notations}

\newbox\toto
\newbox\tata
\newlength\largeurtoto
\newlength\largeurtata

\newcommand \ttt[3]{\setbox\tata=\hbox{#1}%
\ifdim\wd\tata<.17\textwidth\relax%
\setlength{\largeurtoto}{.80\textwidth}%
\setlength{\largeurtata}{.15\textwidth}%
\else%
\setlength{\largeurtoto}{.95\textwidth}%
\addtolength{\largeurtoto}{-\wd\tata}%
\setlength{\largeurtata}{\wd\tata}%
\fi%
\setbox\toto=\hbox{\parbox[b]{\largeurtoto}{\leftskip10pt\parindent-\leftskip\strut#2\dotfill\par}}%
\smallskip\noindent\mbox{%
\parbox[b][\ht\toto][t]{\largeurtata}{\strut#1}%
\parbox[b]{\largeurtoto}{\leftskip10pt\parindent-\leftskip\strut#2\dotfill\par}%
\hspace{.01\textwidth}%
\parbox[b]{.04\textwidth}{\hfill\strut#3}%
}%
\par}

\newcommand\NOTA[3]{\ttt{\smash{#1}}{\smash{#2}}{\pageref{#3}}}
\newcommand\NOTAx[2]{\ttt{\smash{\tsbf{#1}}}{\smash{#2}}{\pageref{Ax#1}}}
\newcommand\NOTAt[2]{\ttt{\smash{\sa{#1}}}{\smash{#2}}{\pageref{theorie#1}}}

\newbox\Toto
\newbox\Tata
\newlength\largeurToto
\newlength\largeurTata

\newcommand \ttT[3]{\setbox\Tata=\hbox{#1}%
\setlength{\largeurToto}{.85\textwidth}%
\setlength{\largeurTata}{.10\textwidth}%
\setbox\Toto=\hbox{\parbox[b]{\largeurToto}{\leftskip10pt\parindent-\leftskip\strut#2\dotfill\par}}%
\smallskip\noindent\mbox{%
\parbox[b][\ht\Toto][t]{\largeurTata}{\strut#1}%
\parbox[b]{\largeurToto}{\leftskip10pt\parindent-\leftskip\strut#2\dotfill\par}%
\hspace{.01\textwidth}%
\parbox[b]{.04\textwidth}{\hfill\strut#3}%
}%
\par}

\newcommand\NOTAA[3]{\ttT{\smash{#1}}{\smash{#2}}{\pageref{#3}}}

\newpage
\chapter*{Index des notations}
\addcontentsline{toc}{section}{Index des notations}
\markboth{Index des notations}{Index des notations}

%\begin{flushright}{page}\end{flushright}%
\vspace{.5em}
\subsection*{Logique}

\NOTA{$\vd$}{règle de déduction}{NOTAvou}
\NOTA{$\vou$}{ouvrir des branches de calcul}{NOTAvou}
\NOTA{$\Exists$}{introduire une variable fraiche}{NOTAExists}
\NOTA{$\Bot$}{symbole de collapsus}{NOTABot}
\NOTA{$\vii$}{et logique}{NOTAvii}
\NOTA{$\vuu$}{ou logique}{NOTAvuu}
\NOTA{$\exists$}{Il existe logique}{NOTAexists}

\subsection*{Symboles fonctionnels}

\NOTA{$a\vu b$}{$\sup(a,b)$}{deficodisup}
\NOTA{$\Fr(a,b)$}{$a/b$ (supposé bien défini)}{lemCo0FRAC}
\NOTA{$\fsac_F$}{\fsagc de graphe $F$}{Notafsa}
\NOTA{$a\vi b$}{$\inf(a,b)$}{secgrl}
\NOTA{$\Sqr$}{$\Sqr(x)=\sqrt{ x^+}$}{subsecclotrlRR}
\NOTA{$\rho_{d,j}$}{\rav}{prdfVirtualRoots}
\NOTA{$\uplus$}{\({\frac 1 2} \,(x+y)\) dans l'intervalle \([-1,+1]\)}{NOTAuplus}
\NOTA{$\oplus$}{addition forcée dans l'intervalle \([-1,+1]\)}{oplus}
\NOTA{$\mathrm{Cb}$}{\coes barycentriques}{notacb}
\NOTA{$\mathrm{Brc}$}{barycentres}{notacb}
\NOTA{$\rT_n$}{\pol de Chebyshev}{NOTASigIcr}

\subsection*{Théories}
%\NOTAt{Al}{anneaux locaux}
%\NOTAt{Anz}{anneaux réduits}
\NOTAt{Cd}{corps discrets}
\NOTAt{Al}{anneaux locaux}
\NOTAt{Ac}{anneaux commutatifs}
\NOTAt{Asdz}{anneaux sans diviseurs de zéro}
\NOTAt{Ai}{anneaux intègres}
\NOTAt{Al1}{anneaux locaux avec unités}
\NOTAt{Alrd}{anneaux locaux résiduellement discrets}
\NOTAt{Cod}{corps ordonnés discrets}
\NOTAt{Crcd}{corps réels clos discrets}
\NOTAt{Apo}{anneaux préordonnés}
\NOTAt{Ao}{anneaux ordonnés}
\NOTAt{Aonz}{anneaux ordonnés \stm réduits}
\NOTAt{Ato}{anneaux totalement ordonnés}
\NOTAt{Atonz}{anneaux totalement ordonnés réduits}
\NOTAt{Apro}{anneaux proto-ordonnés}
\NOTAt{Aso}{anneaux strictement ordonnés}
\NOTAt{Asto}{anneaux strictement totalement ordonnés}
\NOTAt{Asonz}{anneaux \stm ordonnés réduits}
\NOTAt{Aito}{anneaux intègres totalement ordonnés}
\NOTAt{Codsup}{\codis avec sup}
\NOTAt{Atosup}{anneaux totalement ordonnés avec sup}
\NOTAt{Astosup}{anneaux \stm totalement ordonnés avec sup.}
\NOTAt{Crcdsup}{\crcds avec sup}
\NOTAt{Co--}{théorie minimale des corps ordonnés \emph{non} discrets}
\NOTAt{Co}{corps ordonnés \emph{non} discrets}
\NOTAt{Crc1}{corps réels clos (\emph{non} discrets)}
\NOTAt{Tr0}{treillis (bornés)}
\NOTAt{Tr}{treillis non triviaux}
\NOTAt{Trdi}{\trdis}
\NOTAt{Grl}{groupes réticulés}
\NOTAt{Gtosup}{groupes totalement ordonnés avec sup}
\NOTAt{Afr}{\afrs}
\NOTAt{Afrnz}{\afrs réduits}
\NOTAt{Afrsdz}{\afrs \sdz}
\NOTAt{Asr}{\asrs}
\NOTAt{Asrnz}{\asrs réduits}
\NOTAt{Aftr}{anneaux (réticulés) \ftm réels}
\NOTAt{Afr2c}{\afrs $2$-clos}
\NOTAt{Asr2c}{\asrs $2$-clos}
\NOTAt{Co2c}{corps ordonnés (\emph{non} discrets) $2$-clos}
%\NOTAt{Co--rv}{}
\NOTAt{Afrrv}{\afrs avec racines virtuelles}
\NOTAt{Asrrv}{\asrs avec racines virtuelles}
\NOTAt{Aitorv}{anneaux intègres totalement ordonnés avec racines virtuelles}
\NOTAt{Arc}{anneaux réels clos}
\NOTAt{Corv}{corps ordonnés (\emph{non} discrets) avec racines virtuelles}
\NOTAt{Co--rv}{}
\NOTAt{Crc2}{corps réels clos \emph{non} discrets, 2; essentiellement identique à \sa{Corv} et à \sa{Crc1}}
\NOTAt{Crca}{corps réels clos \emph{non} discrets archimédiens}
\NOTAt{Afrb}{\afrbs}
\NOTAt{Icr}{intervalles compacts réels}
\NOTAt{Crc3}{corps réels clos \emph{non} discrets, 3}

\subsection*{Signatures}

\NOTA{$\Sigma_\Ac$}{$(\cdot=0\mathrel{;}\cdot+\cdot,\cdot\times \cdot,-\,\cdot,0,1)$\quad Anneaux commutatifs }{NOTASigAc}

\NOTA{$\Sigma_\Ai$}{$(\cdot=0,\cdot\neq0\mathrel{;}\cdot+\cdot,\cdot\times \cdot,-\,\cdot,0,1)$\quad Anneaux intègres }{NOTASigAi}

\NOTA{$\Sigma_\Alu$}{$(\cdot=0,\U(\cdot)\mathrel{;}\cdot+\cdot,\cdot\times \cdot,-\,\cdot,0,1)$\quad Anneaux locaux avec unités }{NOTASigAlu}

\NOTA{$\Sigma_\Alrd$}{$(\cdot=0,\U(\cdot),\Rn(\cdot)\mathrel{;}\cdot+\cdot,\cdot\times \cdot,-\,\cdot,0,1)$ Anneaux locaux résiduellement discrets}{NOTASigAlrd}

\NOTA{$\Sigma_\Aso$}{$(\cdot=0,\cdot\geq 0,\cdot>0\mathrel{;}\cdot+\cdot, \cdot\times\cdot,-\cdot, 0,1)$\quad Anneaux strictement ordonnés }{NOTASigAso}

\NOTA{$\Sigma_\Ao$}{$(\cdot=0,\cdot\geq 0\mathrel{;}\cdot+\cdot, \cdot\times\cdot,-\,\cdot,0,1)$\quad Anneaux ordonnés }{NOTASigAo}

\NOTA{$\Sigma_{\CoO}$}{$(\cdot=0,\cdot\geq 0,\cdot>0\mathrel{;}\cdot+\cdot, \cdot\times\cdot,\cdot\vu\cdot,-\,\cdot,0,1)$\quad Corps ordonné \nds, 0 }{NOTASigCoO}

\NOTA{$\Sigma_\Co$}{$(\cdot=0,\cdot\geq 0,\cdot>0\mathrel{;}\cdot+\cdot, \cdot\times\cdot,\cdot\vu\cdot,-\,\cdot,\mathrm{Fr}(\cdot,\cdot),0,1)$\;Corps ordonné \nds }{NOTASigCo}

\NOTA{$\Sigma_\TR$}{$(\cdot=\cdot\mathrel{;}\cdot \vi \cdot, \cdot \vu \cdot, 0,1)$\quad Treillis bornés }{NOTASigTr}

\NOTA{$\Sigma_\Gao$}{$(\cdot=0,\cdot\geq 0\mathrel{;}\cdot+\cdot,-\,\cdot,0)$\quad Groupes abéliens ordonnés }{NOTASigGao}

\NOTA{$\Sigma_\GRL$}{$(\cdot=0\mathrel{;}\cdot+\cdot,-\,\cdot,\cdot\vu\cdot,0)$\quad Groupes réticulés }{NOTASigGrl}

\NOTA{$\Sigma_\AfR$}{$(\cdot=0\mathrel{;}\cdot+\cdot, \cdot\times\cdot,\cdot\vu\cdot,-\,\cdot,0,1)$\quad Anneaux fortement réticulés }{NOTASigAfr}

\NOTA{$\Sigma_{\AfR'}$}{$(\cdot=0,\cdot\geq 0\mathrel{;}\cdot+\cdot, \cdot\times\cdot,\cdot\vu\cdot,-\,\cdot,0,1)$\; Anneaux fortement réticulés, bis}{NOTASigAfr'}

\NOTA{$\Sigma_\AsR$}{$(\cdot=0,\cdot\geq 0,\cdot>0\mathrel{;}\cdot+\cdot, \cdot\times\cdot,\cdot\vu\cdot,-\,\cdot,0,1)$\quad Anneaux strictement réticulés }{NOTASigAsr}

\NOTA{$\Sigma_\AftR$}{$(\cdot=0,\cdot\geq 0,\cdot>0\mathrel{;}\cdot+\cdot, \cdot\times\cdot,\cdot\vu\cdot,-\,\cdot,\Fr(\cdot),0,1)$\quad Anneaux fortement réels }{NOTASigAftr}

\NOTA{$\Sigma_\AfRdc$}{$(\cdot=0,\mathrel{;}\cdot+\cdot, \cdot\times\cdot,\cdot\vu\cdot,-\,\cdot,\Sqr(\cdot),0,1)$\quad Anneaux fortement réticulés $2$-clos }{NOTASigAfr2c}

\NOTA{$\Sigma_\AsRdc$}{$(\cdot=0,\cdot\geq 0,\cdot>0\mathrel{;}\cdot+\cdot, \cdot\times\cdot,\cdot\vu\cdot,-\,\cdot,\Sqr(\cdot),0,1)$ }{NOTASigAsr2c}

\NOTA{$\Sigma_\AftRdc$}{$\cdot=0,\cdot\geq 0,\cdot>0\mathrel{;}\cdot+\cdot, \cdot\times\cdot,\cdot\vu\cdot,-\,\cdot,\Fr(\cdot),\Sqr(\cdot),0,1)$ }{NOTASigAftr2c}

\NOTA{$\Sigma_\AfRrv$}{$(\cdot=0 \mathrel{;}\cdot+\cdot, \cdot\times\cdot,\cdot\vu\cdot,-\,\cdot,(\rho_{d,j})_{1\leq j\leq d},0,1)$\quad Anneaux \dots  \ravs }{NOTASigAfrrv}

\NOTA{$\Sigma_\ArC$}{$(\cdot=0\mathrel{;}\cdot+\cdot, \cdot\times\cdot,\cdot\vu\cdot,-\,\cdot,(\rho_{d,j})_{1\leq j\leq d},\mathrm{Fr}(\cdot,\cdot),0,1)$\quad Anneaux réels clos }{NOTASigArc}

\NOTA{$\Sigma_\CoRv$}{$(\cdot=0,\cdot\geq 0,\cdot>0\mathrel{;}\cdot+\cdot, \cdot\times\cdot,\cdot\vu\cdot,-\,\cdot,(\rho_{d,j})_{1\leq j\leq d)},\mathrm{Fr}(\cdot,\cdot),0,1)$  }{NOTASigCorv}

\NOTA{$\Sigma_{\Icro}$}{$(\cdot= \cdot,\cdot\geq \cdot, \cdot>\cdot\mathrel{;}\cdot\,\uplus\,\cdot, \cdot\times \cdot, \cdot\vu\cdot, \mathrm{Fr}(\cdot,\cdot),\,- \cdot,0)$\quad Intervalle compact réel, 0}{NOTASigIcro}

\NOTA{$\Sigma_\AfrB$}{$(\cdot=0,\cdot\geq 0,\cdot>0\mathrel{;}(\mathrm{Brc}_{\rho})_{\rho\in\mathrm{Cb}}, \cdot\oplus\cdot, \cdot\times\cdot,-\,\cdot,\cdot\vu\cdot,\mathrm{Fr}(\cdot,\cdot), (r)_{r\in\II_\QQ})$ }{NOTASigAfrb}

\NOTA{$\Sigma_\Icr$}{$(\cdot=0,\cdot>0,\cdot\geq0\mathrel{;}(\mathrm{Brc}_{\rho})_{\rho\in\mathrm{Cb}}, \cdot\oplus\cdot, \cdot\times\cdot,-\,\cdot,\cdot\vu\cdot,(\rT_n)_{n\in\NN},\mathrm{Fr}(\cdot,\cdot), (r)_{r\in\II_\QQ})$ }{NOTASigIcr}

%\NOTA{$\Sigma_\$}{$()$\quad }{NOTASig}
%
%\NOTA{$\Sigma_\$}{$()$\quad }{NOTASig}

%: axiomes
%\NOTA{}{}{}
%\NOTAx{}{}{}

\subsection*{Quelques axiomes et \rdys}
\NOTAx{AL}{$\,\, (x+y)z=1 \vd \Exists u \; xu=1 \;\vou\;\Exists v\;yv=1$ }
\NOTAx{Anz}{$\,\, x^2= 0 \vd x = 0$ }
\NOTAx{ASDZ}{$\,\,xy=0 \vd x=0 \;\;\vou\;\;y=0$ }
\NOTAx{CD}{$\vd x=0\;\;\vou\;\;\Exists {y}\;\;xy=1$ }
\NOTA{$\tsbf{CL$_{\Ac}$}$}{$\,\,1=_\Ac0\vd \Bot$ }{AxCLnqAc}
\NOTA{$\tsbf{ED\inq}$}{$\vd x=0 \vou x\neq 0$ }{AxEdnq}
\NOTA{$\tsbf{col\inq}$}{$\,\,0\neq 0 \vd 1=0$ } {Axcolnq}
\NOTAx{Uv}{$\,\, xy = 1 \vd \U(x) $ }
\NOTAx{UV}{$\,\, \U(x) \vd \Exists y\, xy = 1$ }
\NOTAx{AL1}{$\,\, \U(x+y) \Vd \U(x) \vou \U(y)$ }
\NOTAx{NIL}{$\,\,\mathrm{Z}(x)\vd\Vou_{n\in\NN^+} x^n=0$ }
\NOTA{\tsbf{aso1} à \tsbf{aso4}}{\,\, }{Axaso1}
\NOTAx{Iv}{$\,\, xy = 1 \vd x\neq 0 $ }
\NOTAx{IV}{$\,\, x\neq 0 \vd \Exists y\, xy = 1$ }
\NOTA{$\tsbf{col\igt}$}{$\,\,0> 0 \vd 1=0$ }{Axcoligt}
\NOTAx{OT}{$\vd x \geq 0 \;\vou\; x\leq 0$ }
\NOTAx{Aonz}{$\,\, c\geq 0\vet x(x^2+c)\geq 0 \vd x\geq 0$ }
\NOTA{\tsbf{RCF$_n$}}{$\,\, a< b\vet P(a)P(b)<0 \vd \Exists x\, 
\big(P(x)=0,\,a<x<b\big)$ \quad ($P(x)=\sum_{k=0}^na_kx^k$) }{RCFn}\NOTAx{Aso1}{$\,\, x> 0\vet xy\geq 0 \vd y\geq 0 $ }
\NOTAx{Aso2}{$\,\, x\geq 0\vet xy> 0 \vd y> 0$ } 
\NOTAx{OTF}{$\,\, x+y> 0 \vd x > 0\;\vou\;y>0$ }
\NOTA{\tsbf{OTF}$\eti$}{$\,\, xy< 0 \vd x <0 \;\vou\; y<0 $ }{AxOTFx}
\NOTAx{Ato1}{$\,\, x\geq 0 \vet xy=1\vd y\geq 0$ }
\NOTAx{Ato2}{$\,\, c\geq 0\vet x(x^2+c)\geq 0 \vd x^3\geq 0$ }
\NOTAx{Aonz3}{$\,\, a\geq 0\vet b\geq 0\vet a^2=b^2 \vd a=b$ }
\NOTAx{CVX}{$\,\,0\leq a\leq b\vd \Exists z\; zb=a^2$ }\NOTAx{FRAC}{$\,\,0\leq a\leq b\vd \Exists z\; (zb=a^2\vet 0\leq z\leq a)$ }\NOTAx{fr1}{$\vd \mathrm{Fr}(a,b)\, \abs b=(\abs a\!\vi\! \abs b)^2$ }\NOTAx{fr2}{$\vd 0\leq {\mathrm{Fr}(a,b)}\leq \abs a\!\vi\! \abs b$ }
\NOTA{\tsbf{FRAC$_n$}}{\smash{$\,\,\abs u^n\leq \abs v^{n+1}\vd \Exists z\; (zv=u\vet \abs z^{n}\leq \abs v)$ \quad ($n\geq 1$)}}{AxFRACn}
\NOTAx{Gao}{$\,\, x\geq 0\vet x\leq 0 \vd x = 0$ }
\NOTAx{grl}{$\vd x+(y\vu z)=(x+y)\vu(x+z)$ }
\NOTAx{sup1}{$\vd x\vu y\,\geq x$ }
\NOTAx{sup2}{$\vd x\vu y\,\geq y$ }
\NOTAx{Sup}{$\,\, z\geq x\vet z\geq y\vd z\geq x\vu y$ }
\NOTAx{afr}{$\vd x^+\, (y\vu z)=(x^+\, y)\vu(x^+\, z)$ }
\NOTAx{sup}{$\vd ((x\vu y)- x)\,((x\vu y)- y)=0$ }
\NOTA{\tsbf{afr1} à \tsbf{afr7} }{}{Axafr1}
\NOTA{\tsbf{Afr1} à \tsbf{Ato2} }{}{AxAfr1}
\NOTA{\tsbf{Afrnz1} - \tsbf{Afrnz2} }{}{AxAfrnz1}
\NOTAx{AFRL}{$\,\, z(x+y)=1\vet x+y\geq 0\vd \Exists u\;(ux=1\vet x\geq 0 )\;\vou\; \Exists v\;(vy=1\vet y\geq 0 )$ }
\NOTAx{sqr0}{$\vd \Sqr(0)=0$ }
\NOTAx{sqr1}{$\vd \Sqr(x)=\Sqr(x)^+ $ }
\NOTAx{sqr2}{$\vd \Sqr(x)=\Sqr(x^+)$ }
\NOTAx{sqr3}{$\vd \Sqr(x)^2=x^+ $ }
\NOTAx{sqr4}{$\vd \Sqr(x^+y^+)=\Sqr(x)\Sqr(y)$ }
\NOTA{\tsbf{vr$_{i,j,k}$}}{\,\, axiomes pour les racines virtuelles, exemples}{exavr} 
%\NOTAx{ }{$ $}

\subsection*{Structures algébriques, modèles \gnqs}

\NOTA{$\RRa$}{le corps des réels algébriques}{NOTARa}
\NOTA{$\RR_{\tsbf{PR}}$}{le corps des réels primitifs récursifs}{exacorpsnondiscret}
\NOTA{$\RR_{\tsbf{Ptime}}$}{le corps des réels calculables en temps polynomial}{exacorpsnondiscret}
\NOTA{$\RR_{\tsbf{Rec}}$}{le corps des réels récursifs}{exacorpsnondiscret}
\NOTA{$\Sac_n(\RR)$}{\ref{defiFSAGC2}: fonctions \sagcs en $n$ variables, voir aussi \ref{defiFSAGC2+}}{defiFSAGC2}

\medskip \noindent  Ci-après $\gA$ désigne un anneau commutatif ou une structure \agq sur une \tdy convenable dans le contexte considéré, $\gR$ désigne un \afrvr.
 
\NOTA{$\AFR(\gA)$}{\afr librement engendré par $\gA$}{notaAFR}
\NOTA{$\SIPD_n(\gA)$}{anneau des \spos en~$n$ variables sur $\gA$}{defiSIPD}
\NOTA{$\SIPD_n(\gA,\gB)$}{anneau des $\gA$-\spos en~$n$ variables sur $\gB$}{defiSIPD}
\NOTA{$\AFRdC(\gA)$}{2-clôture d'un \afr}{notaAFR2C}
\NOTA{$\Sace_n(\gR)$}{anneau des fonctions \sagcs entières}{defiSacem}
\NOTA{$\Sac_n(\gR)$}{anneau des fonctions \sagcs en $n$ variables}{defiFSAGC2+}
\NOTA{$\AFRNZ(\gA)$}{\afr réduit engendré par $\gA$}{defiSacembis}
\NOTA{$\AFRRV(\gA)$}{\afr avec \ravs engendré par~$\gA$}{defiSacembis}
\NOTA{$\PPM(\gA)$}{anneau des fonctions \polles par morceaux sur $\gA$}{defiSacembis}
\NOTA{$\ARC(\gA)$}{\arc engendré par $\gA$}{notaARC}
%\NOTA{$$}{}{}

%\newpage
\rdb
\addcontentsline{toc}{section}{Index des termes}
\printindex

\end{document}